%% file: hdinferrev3.tex
\documentclass[aos]{imsart}

\RequirePackage{amsthm,amsmath,amsfonts,amssymb}
\RequirePackage[authoryear]{natbib}
\RequirePackage{graphicx}

\startlocaldefs
\theoremstyle{plain}

\theoremstyle{remark}

\usepackage{epsfig,enumitem,bm,multirow,float,color}
\usepackage{amsbsy}
\usepackage{algorithm,mathabx}
\usepackage{lineno}
\usepackage{url}
\usepackage{booktabs}
\usepackage{epstopdf}
\usepackage{verbatim}
\usepackage{xr}
\newtheorem{Th}{\underline{\bf Theorem}}

\newtheorem{Def}{Definition}

\newtheorem{Lem}{\underline{\bf Lemma}}
\newtheorem{Cor}{\underline{\bf Corollary}}

\def\be{\begin{eqnarray*}}
	\def\ese{\end{eqnarray*}}
\def\be{\begin{eqnarray}}
\def\ee{\end{eqnarray}}
\def\bsq{\begin{equation*}}
\def\esq{\end{equation*}}
\def\bq{\begin{equation}}
\def\eq{\end{equation}}

\def\supp{\hbox{supp}}

\def\cov{\hbox{cov}}

\def\wh{\widehat}
\def\wt{\widetilde}
\def\wc{\widecheck}

\def\mL{\mathcal L}
\def\mM{\mathcal M}
\def\mA{\mathcal A}
\def\mR{\mathbb{R}}
\def\n{\nonumber}

\def\cov{\mbox{cov}}

\def\argmin{\mbox{argmin}}

\def\sumi{\sum_{i=1}^n}

\def\trans{^{\rm T}}

\def\btheta{\boldsymbol\theta}

\def\bmu{\boldsymbol\mu}
\def\bnu{\boldsymbol\nu}

\def\bb{\boldsymbol\beta}
\def\bGamma{\boldsymbol\Gamma}
\def\bOmega{\boldsymbol\Omega}
\def\bomega{\boldsymbol\omega}

\def\bPsi{\boldsymbol\Psi}
\def\bmu{\boldsymbol\mu}
\def\bg{{\boldsymbol\gamma}}

\def\0{{\bf 0}}
\def\A{{\bf A}}
\def\Q{{\bf Q}}
\def\U{{\bf U}}

\def\e{{\bf e}}
\def\R{{\bf R}}

\def\a{{\bf a}}
\def\B{{\bf B}}

\def\f{{\bf f}}
\def\h{{\bf h}}
\def\b{{\bf b}}
\def\I{{\bf I}}
\def\mU{\mathcal{U}}

\def\M{{\bf M}}

\def\t{{\bf t}}

\def\bQ{{\bf Q}}

\def\U{{\bf U}}
\def\S{{\bf S}}
\def\s{{\bf s}}
\def\u{{\bf u}}

\def\v{{\bf v}}
\def\W{{\bf W}}

\def\w{{\bf w}}
\def\X{{\bf X}}

\def\x{{\bf x}}
\def\I{{\bf I}}

\def\C{{\bf C}}

\def\Z{{\bf Z}}
\def\z{{\bf z}}

\def\Ubar{{\overline{\U}}}

\def\bSig{{\bf \Sigma}}

\def\bq{\begin{equation}}
\def\eq{\end{equation}}
\def\pr{\hbox{pr}}
\def\wh{\widehat}
\def\wt{\widetilde}
\def\trans{^{\rm T}}

\def\log{\hbox{log}}

\def\squarebox#1{\hbox to #1{\hfill\vbox to #1{\vfill}}}
\def\btheta{{\boldsymbol \theta}}

\def\mA{\mathcal{A}}

\def\mF{\mathcal{F}}

\def\mM{\mathcal{M}}

\def\cov{\hbox{cov}}

\def\spacingset#1{\renewcommand{\baselinestretch}%
	{#1}\small\normalsize} \spacingset{1}
\def\bse{\begin{eqnarray*}}
	\def\ese{\end{eqnarray*}}
\def\be{\begin{eqnarray}}
\def\ee{\end{eqnarray}}
\def\bsq{\begin{equation*}}
\def\esq{\end{equation*}}
\def\bq{\begin{equation}}
\def\eq{\end{equation}}
\def\pr{\hbox{pr}}
\def\wh{\widehat}
\def\wt{\widetilde}

\def\trans{^{\rm T}}

\def\boxit#1{\vbox{\hrule\hbox{\vrule\kern6pt\vbox{\kern6pt#1\kern6pt}\kern6pt\vrule}\hrule}}

\endlocaldefs

\begin{document}

\begin{frontmatter}
\title{On High dimensional Poisson models with measurement error:
	hypothesis testing for nonlinear nonconvex optimization}
\runtitle{On High dimensional Poisson measurement error models}

\begin{aug}


\author[A]{\fnms{Fei} \snm{Jiang}\ead[label=e1]{fei.jiang@ucsf.edu}},
\author[B]{\fnms{Yeqing}
  \snm{Zhou}\ead[label=e2,mark]{zhouyeqing@tongji.edu.cn}},
\author[C]{\fnms{Jianxuan}
  \snm{Liu}\ead[label=e3,mark]{jliu193@syr.edu}}
and
\author[D]{\fnms{Yanyuan}
  \snm{Ma}\ead[label=e4,mark]{yzm63@psu.edu}}
\address[A]{Department of Epidemiology and Biostatistics,
       The University of California, San Francisco,
\printead{e1}}
\address[B]{School of Mathematical Sciences, Tongji University, 
\printead{e2}}
\address[C]{Department of Mathematics, Syracuse University,
\printead{e3}}
\address[D]{Department of Statistics, Pennsylvania State University,
\printead{e4}}
\end{aug}

\begin{abstract}
We study  estimation and testing in the
  Poisson regression model with noisy high 
  dimensional covariates, which has wide applications in analyzing
  noisy big
  data. Correcting for the estimation bias due to
  the covariate noise leads to a non-convex target function to
  minimize. Treating the high dimensional issue further leads us to
augment an amenable penalty term to the target function. We propose to
estimate the regression parameter through minimizing the penalized
target function. 
We derive the $L_1 $ and $L_2$ convergence rates of the
estimator and prove the variable selection consistency. We further
establish the asymptotic normality of any subset of the parameters,
where the subset can have infinitely many components as long as
its cardinality grows sufficiently slow.
 We develop Wald and
score tests based on the asymptotic normality of the estimator, which
permits testing of linear functions of the members if the subset.
We examine the finite sample
performance of the proposed tests by extensive simulation. Finally,
the proposed method is successfully applied to the
Alzheimer's Disease Neuroimaging Initiative study, which motivated
this work initially.
\end{abstract}



\end{frontmatter}

\section{Introduction}\label{sec:intro}
Count data are routinely encountered in practice.
For example, 
cognitive scores in a neuroscience study, the number of
deaths in an infectious disease study, and the number of clicks on a
particular product on an e-commerce platform, are all count data.  Because most of the
count data are concentrated on a few
small discrete values rather than expanded on the entire real line and because
the distribution of  count variables is often skewed, the familiar
linear model becomes less ideal to capture these features. In the
literature, Poisson regression \citep{mccullagh2019} 
is arguably the most popular model to describe count
outcomes, because it naturally models the skewed distribution for
positive outcomes. On the other hand, together with the count data, a
large number of covariates are often collected thanks to the ever
advancing capability of modern 
technologies. However, these covariates are often contaminated with
errors due to imperfect data acquisition and processing procedures. Ignoring
these errors can produce biased results, which can finally lead to misleading
statistical inference on the model parameters
\citep{carroll2006} that explain the association between 
covariates and outcomes.  Our goal is to develop rigorous statistical inference
procedures to 
test linear hypotheses in the high dimensional Poisson model with noisy
covariates. Such inference tools will enable explaining the
association between the count outcome and the individual covariate or
combination of covariate,  quantifying the 
uncertainties of the estimated association,
 and controlling the
false discovery rate when testing scientifically
important hypotheses.

Let $Y$ be the count outcome and $\X$ be its associated covariate
vector. In the Poisson model, $Y$ is related to
$\X$ as
\be\label{eq:pois}
\pr(Y\mid\X) =
e^{-\exp(\bb\trans \X)}\{\exp(\bb\trans \X)\}^Y/Y!. 
\ee
We study the testing problem in (\ref{eq:pois}) under 
the situation that the covariate vector $\X$ is both high dimensional
and  contaminated with noise.
When $\X$ is accurately observed, the testing problem has been
extensively discussed in the 
literature \citep{ning2017, zhang2017simultaneous,
  van2014asymptotically, shi2019linear}.
However, when  $\X$  is not accurately observed,
it is unclear that any of the existing proposed tests are
applicable, and
testing in the high dimensional noisy Poisson regression
model has not been explored. The major obstacles
in constructing valid hypothesis testing procedures are as
follows.  1) The existing lasso-type penalized Poisson estimator \citep{jiang2021}
does not enjoy the variable selection consistency when the
number of parameters is much larger than
the sample size.
2) The
asymptotic normality of the estimator
has not been established. We develop Wald and score tests targeting at
linear hypothesis on the parameters of 
interest in  (\ref{eq:pois}). To overcome obstacle
1), we improve the penalized Poisson estimator proposed in
\cite{jiang2021} by using a
class of ``amenable'' penalty functions first defined in
\cite{loh2015, loh2017} in combination with a modified
log-likelihood function
to construct estimators.  We establish the
estimation consistency and variable selection consistency of the resulting
estimators.  To bypass obstacle 2), we derive the asymptotic linear form
of the estimators, and establish the asymptotic normality.
The asymptotic normal estimator has a wider range of applications than
the lasso type estimator does, because it
facilitates subsequent inference procedures such as constructing
hypothesis testing procedures.

Even after establishing the asymptotic normal properties, it is still
challenging to generalize Wald and
score tests to the high dimensional setting for
Poisson regression with noisy data. This is because 
under the amenable
penalties \citep{loh2015, loh2017},
the asymptotic normality of the
estimators is built on  a  minimal
signal condition, which requires the nonzero elements in
$\bb$ to be at least of order $\lambda$. Here
$\lambda$ is the penalty parameter which goes to zero when
sample size increases. Now consider testing the
null hypothesis $\beta_1 = 0$ versus the alternative $\beta_1 =
h_n$, where $\beta_1$ is the first element of $\bb$.  The minimal
signal condition implies that the  test will have no power
in testing the local alternative when $|h_n|_2<<\lambda$.
To resolve this issue, we remove the penalties on the
subvector of the parameters
involved in the test.  However, it is still unclear how
fast the  dimension of the
subvector can grow while still ensuring sufficient power. To this
end, we derive
the convergence property of the estimators, which provides the explicit
rate at which
the dimension of the subvector is allowed to grow with the sample size
in order to
achieve consistency, asymptotic normality, and sufficient
  power in testing.
Furthermore, to implement the score test, we need to estimate the
regression parameters under the null hypothesis, which involves
optimization under linear equality constraints. This type of
constrained parameter estimation for noisy Poisson
model has not yet been developed. To fill this gap, we develop a
general procedure for parameter
estimation under linear constraints. The constraints include inequality
constraints for the parameter estimation under general Poisson
model and an additional equality constraint imposed by the
null hypothesis, which leads to great challenge in
establishing the convexity.  Incorporating inequality
  constraints is 
practically important because it allows to
incorporate additional parameter
information, which will reduce the
estimation variation and in turn the sample size needed to achieve
satisfactory estimation accuracy.

We briefly summarize our contributions as follows.  First, we develop
a new estimation procedure of the Poisson
model with amenable regularization for noisy data. Second, we
show the variable selection consistency and the
consistency of the resulting estimator. We provide
explicit convergence rate of the estimator. Third, we derive the
asymptotic normality of the estimator for the nonzero
parameters and the parameters to test. Fourth, we propose the Wald
and score test procedures by constructing  the corresponding
test statistics. Fifth, we derive the asymptotic
distributions of the Wald and score test statistics. These five essential
elements combined together finally
allow us to perform hypothesis testing for 
Poisson model with high dimensional noisy covariates,  which
allows us to answer important questions such as ``if the left
inferior temporal gyrus has a significant impact on the development of
Alzheimer's disease''. These estimation and inference tasks
are not straightforward to achieve, they require building up a series of
theoretical properties first, which involves techniques related to analyzing
conditional sub-Gaussian distribution tails, utilizing and modifying
various concentration 
inequalities, constructing the prime-dual equivalence, 
carefully bounding various quantities, linking 
different vector and matrix norms, and establishing
 a Lyapunov-type bound  \citep{bentkus2005} on the probability distribution to derive the asymptotic distribution of 
proposed test statistics. All these analyses are performed under the
unusual constraints involving both linear equality constraints and
parameter restrictions.
We also modify the alternating direction method of multipliers (ADMM)
algorithm to solve a regularized optimization problem under linear
constraint in constrained parameter space. 
Although each individual technique in its basic form has been used in the literatures of 
mathematical analysis, statistics, combinatorics, operations research and computer
science, a seamless combination of all these into a general tool to
solve the problem under study is very challenging and difficult.

\begin{figure}
  \centering
  \includegraphics[scale = 0.25]{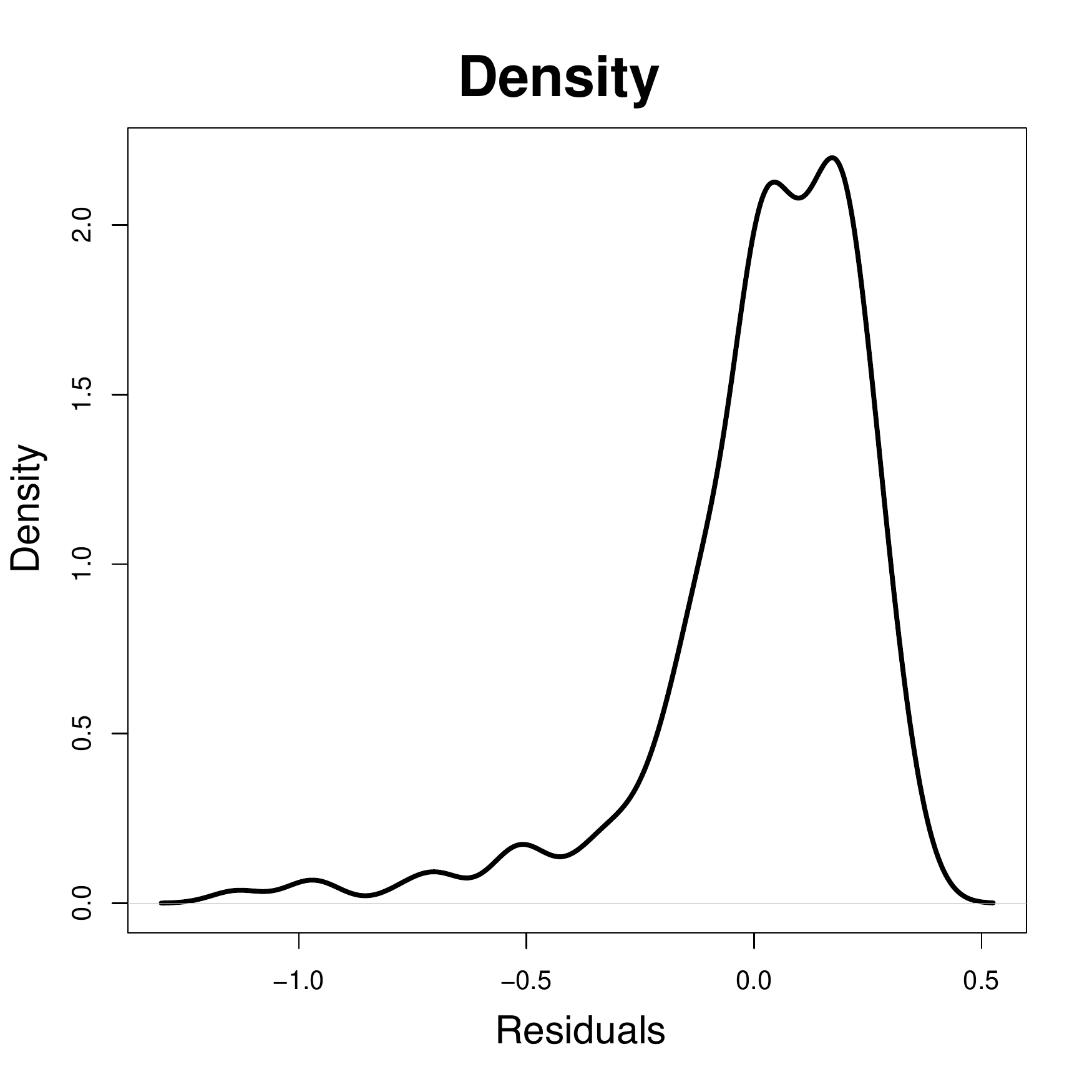}
  \includegraphics[scale = 0.25]{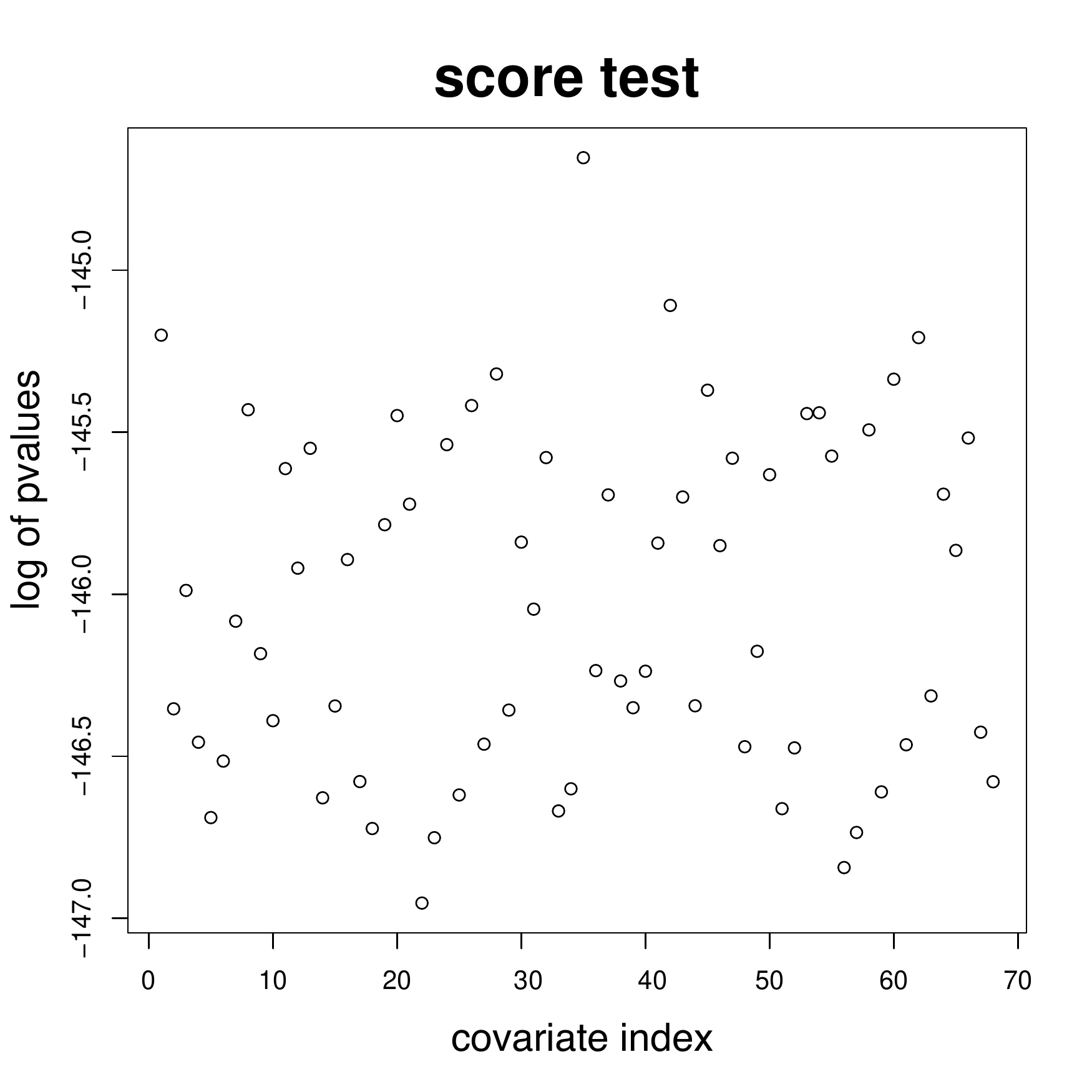}
  \includegraphics[scale = 0.25]{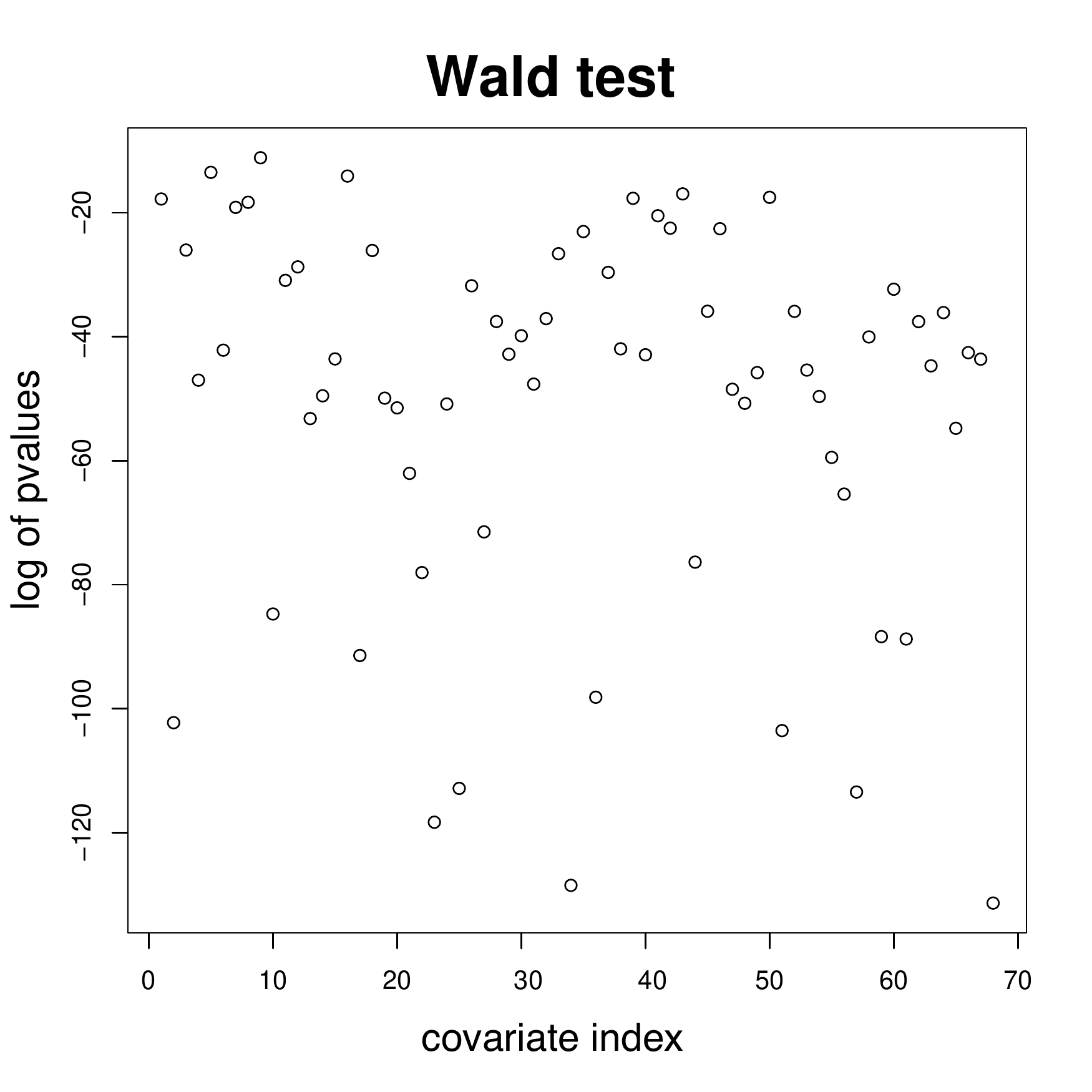}
  \caption{Left: The density of the residuals from lasso regression. The lasso regression does not provide a satisfactory fit
    for the data. Middle and Right: The logarithm of the p values from the Wald and score tests
    proposed by \cite{shi2019linear} for testing whether the SUVR
    from each cortical regions is significant predictor for the
    cognitive score. The Wald and score tests suggest  that the SUVRs at
    all the cortical regions have significant association with the
    cognitive score.  }\label{fig:linear}
\end{figure}

Count data occur frequently in practice, and it is a
  rule rather than exception that the covariates can be
  contaminated. In modern data collection mechanism, 
  covariates are almost always high dimensional. Hence, estimation and
  inference in Poisson
  regression with high dimensional noisy covariates is a general
  problem with wide
  applications. A direct motivation of this work is the Alzheimer's
  Disease Neuroimaging Initiative 
(ADNI) study,  which is a multi-site longitudinal study investigating
early detection of Alzheimer's disease (AD) and tracking disease progression 
biomarkers \citep{weiner2017}. Recently,
the advent of
tau-targeted positron emission tomography (PET) tracers such
as flortaucipir ($^{18}$F-AV-1451) has made it possible to
investigate the relative (to patient's body weight) tissue
radioactivity concentration of the tracers, 
quantified as standardized uptake value ratio (SUVR), 
in
relationship to the cognitive function.   Therefore, we aim to study
the association between cognitive scores and SUVRs from PEG image
data. We extract Montreal
Cognitive Assessment (MoCa)
scores  ($Y$) and SUVRs ($\X$) from the PET
image in the ADNI study taken within 14 days of the cognitive tests from 196
subjects in the ADNI phase 3 study. We first perform a linear lasso
regression between the  logarithm of MoCa score and the 218 covariates
including age, gender, SUVRs, and
volumes of whole brain ROIs.  Figure \ref{fig:linear} shows the
density  of the 
residuals from the lasso regression, which suggests that the residuals
are skewed and hence the  linear lasso regression does not
provide a satisfactory fit for the data.  This motivates us to consider Poisson
regression. 
 We utilize the Poisson high
dimensional hypothesis testing procedure developed in
\cite{shi2019linear} to examine which SUVRs are significantly
associated with the MoCa scores. For each covariate of interest, we
test the hypothesis that the corresponding coefficient is
greater than zero. We plot the logarithm of the p-values from the score and Wald
tests proposed in  \cite{shi2019linear}  for the coefficients of the
SUVRs at cortical  
ROIs.  Figure \ref{fig:linear} shows that if using 
0.05/218 as a cut off for the p-value, 
both the Wald and score test identify the SUVRs at all cortical 
ROIs  as
significant predictors, which contradicts the fact that 
the cognitive functions are controlled by a subset of brain ROIs
\citep{leisman2016}. This unsatisfactory result likely attributes to
the fact the \cite{shi2019linear}'s method relies on the assumption that
the expectation of the exponential of the distance between outcome and
regression function is bounded (Condition (A3) in \citep{shi2019linear})
while neuroimage data are often subject to
data acquisition and processing errors, which likely leads to violation
of this assumption. This motivating example demonstrates the
necessity of developing novel statistical inference procedure  to  test  linear hypothesis in the high dimensional Poisson
model with noisy data. 



The rest of the paper is organized as follows. Section
\ref{sec:related} discusses related work. In Section
\ref{sec:model}, we describe our model assumptions and
the overall estimation strategy.  We further detail the
estimation with
and without the null constraint, and the construction of the test
statistics.
The fundamental theoretical developments are
provided in Section \ref{sec:theory}, where we establish convergence
rates,
the asymptotic normal results, and the properties of
the test procedures.
We study the practical implementation and the numerical
performances in Section \ref{sec:num}, where a detailed algorithm is
provided, extensive
simulations are carried out, and a
ADNI data set is analyzed. We conclude the paper in
Section \ref{sec:con}. The main mathematical proofs are provided in
an Appendix given in a Supplementary Document.

\section{Related Works and Notations}\label{sec:related}
Nonlinear models with high dimensional noisy data are in
general hard problems to work with, partly because existing
treatments usually lead to non-convex optimization,
which violates standard requirements in the high
dimensional data analysis literature.
Thus, only linear models, which are the simplest
in all noisy data problems, have received relatively thorough
investigation \citep{loh2012, belloni2016,
	datta2017cocolasso,  belloni2017linear,  belloni2017inference, lilima2021}.
Expanding the research framework to the Poisson regression context is
difficult because the link function in the Poisson model is
nonlinear. Subsequently, it is
not easy to construct noise adjusted quantities such as a noise
adjusted Hessian matrix like in the linear case. In addition, the
Hessian matrix involves heavy tailed random variables due to the
exponential link, even if all the covariates are
sub-Gaussian in their original scale. These difficulties require
additional restrictions on the moments
of the covariate distribution as well as on the parameter searching
space, which complicates all the subsequent computation and analysis.
Indeed, the only works we are aware of in the high dimensional Poisson
model with noisy data are  \cite{jiang2021, sorensen2015, sorensen2018covariate,
  brown2019meboost}, while only the estimator in
\cite{jiang2021} has been shown to be consistent. 
However, because all these methods use lasso-type $L_1$
penalty in the estimation, the resulting estimators do not enjoy variable
selection consistency and their asymptotic
distribution results are not established.

There is extensive literature on the linear hypothesis testing under high
dimensional noise free setting. \citep{ning2017}  introduced a decorrelated score
function to construct confidence regions for low dimensional
components in high dimensional models.
\cite{zhang2017simultaneous} used the
desparsifying lasso estimator \citep{van2014asymptotically} to
propose a maximal-type statistic allowing the number of
parameters that are involved in the test to grow with the sample
size.  Moreover, \cite{shi2019linear}
proposed a partial penalized likelihood ratio test, a score
test, and a Wald test for testing the linear hypothesis of the
parameters in high dimensional generalized linear
models.

\subsection*{Notations}\label{sec:notation}
We introduce some general notation that will be used throughout the text.
For a matrix $\M$, let $\|\M\|_{\max}$ be the matrix maximum norm,
$\|\M\|_\infty$ be the $L_\infty$ norm and
$\|\M\|_p$ be the $L_p$ norm. Let $\mF(\bb)$ be the $\sigma$-field
generated by $\X_i, \bb\trans \W_i, i = 1, \ldots, n$. Further, let
$\mF_x$ be the sigma-field generated by $\X_i, i = 1, \ldots,
n$.   For a general vector $\a$, let
$\|\a\|_\infty$ be the vector sup-norm, $\|\a\|_p$ be the vector
$l_p$-norm. Let $\e_j$ be the unit vector with 1 on its $j$th
entry. For a vector ${\bf v} = (v_1,
\ldots, v_p)\trans$, let
$\supp({\bf v})$ be the set of indices with $v_i \neq 0$ and $\|{\bf v}\|_0 =
|\supp({\bf v})|$, where $|\mU|$ stands for the cardinality of the set $\mU$.
For a vector ${\bf v}\in \mathbb{R}^p$ and a subset $S \subseteq (1,
\ldots, p)$, we use ${\bf v}_S \in \mathbb{R}^S$ to  denote the vector
obtained by restricting
${\bf v}$  on the set $S$.

Following \cite{fletcher1980},  for an arbitrary norm
$\|\cdot\|_A$ and its dual normal $\|\cdot\|_D$, we define $\partial
\|\x\|_A$ as the set $({\bf v}: \|\x\|_A = {\bf v}\trans \x, \|{\bf v}\|_D \leq
1)$. Thus, for an arbitrary  vector
$\x = (x_1, \ldots, x_p)\trans$,
$\partial \|\x\|_1 = \{{\bf v} =
(v_1, \ldots, v_p)\trans: v_j =
{\rm sign}(x_j) \text{ if } x_j\ne0,  \text{ and } |v_j| \leq 1 \text{ if } x_j=
0  \}$, and
$\partial \|\x\|_2 = \{{\bf v} =
(v_1, \ldots, v_p)\trans: v_j = x_j/\|\x\|_2\}$.

\section{Model,  Estimation and Test Statistics}\label{sec:model}

\subsection{Problem Formulation: High dimensional Poisson model with noisy data}\label{sec:submodel}

Let $\X_i$ be a $p$-dimensional covariate, for example the  image 
features,  and let $Y_i$ be a count
random response variable, for example the MoCa score from the
ADNI data.  We model the relationship between $Y_i$ and
$\X_i$ ($i=1, \ldots, n$) through a Poisson model $\pr(Y_i=y\mid\X_i=\x) =
e^{-\exp(\bb_t\trans \x)}\{\exp(\bb_t\trans \x)\}^y/y!$. Here,
$\bb_t$ is a $p$-dimensional sparse parameter
vector. We allow the number of nonzero
  entries in $\bb_t$ to grow with the sample size. We consider Poisson
model here because our response is a count, and Poisson model is
arguably the most standard model for count data. 
Indeed, Poisson model  has  been widely used to model the
distribution of cognitive scores \citep{katz2021t, fallah2011applying,
mitnitski2014multi}.
We use
$e^{\bb_t\trans\x}$ to model the conditional mean of the Poisson model to
 ensure the positiveness of the mean, and to allow possible
skewness in the distribution
\citep{mccullagh2019}. We assume $\bb_t$ to be sparse
because it often happens that only a few covariates  have effect on the outcome.
For example, in the ADNI data, because the cognitive functions are
controlled by a subset of brain ROIs 
\citep{leisman2016}, only a subset of brain features 
contributes to the cognitive function.

Furthermore,  we assume the covariate $\X_i$ is not precisely observed and
instead, a contaminated version of $\X_i$, denoted  $\W_i$, is
observed, where $\W_i =
\X_i + \U_i$, and $\U_i$ is the
noise that is independent of both $\X_i$ and $Y_i$. For example, in 
the ADNI data, $\X_i$ can be
the true image features, while $\W_i$ represent the observed image
features which can deviate from the truth due to imperfect data
collection and processing procedure. 
Without loss of generality,
assume that $E(\X_i)=\0$, which can always be achieved by centering
the observed covariates in practice. Furthermore, we assume $\U_i$ is a normally
distributed random noise vector with mean zero and
known covariance matrix $\bOmega$. The normal assumption for $\U_i$
is the common assumption at the state of the art in the
 Poisson measurement error literature and
allows to derive analytic form of the loss function, which
is the only setting that we can directly examine the convexity of the
loss function. The known $\bOmega$ assumption is only for
convenience of presentation. In practice, it is often
  replaced by an estimated version based on 
multiple observations, validation data or other standard instruments
under both low and high dimensional settings \citep{carroll2006,
  loh2012}, and the corresponding analysis is routine. 
Let $(\X_i, \W_i, Y_i, \U_i)$ be
independent and identically distributed (iid)
and assume $(\W_i, Y_i), i=1, \dots, n$ are the iid observations. In
this work, we devise estimation procedures for $\bb$ and
  establishing theoretical properties of the estimator, 
we further aim at performing inference, such as conduct
hypothesis testing. Throughout, we allow the covariate dimension to be
much higher than the number of observations, i.e. $p>>n$.
We assume $\bb_t$ is in the
feasible set: $\{\bb: \|\bb\|_0 \leq k$, $\|\bb\|_2 \leq b_0\}$, which
is practically sensible.
A vector $\bb$ in the feasible set automatically satisfies
$\|\bb\|_1\leq b_0\sqrt{k}$.

\subsection{General Estimation Strategy}\label{sec:estimation}
If the true covariates $\X_i$ can be observed and  the dimension $p$
is fixed,  this is a standard regression
model and we routinely estimate $\bb$ by minimizing the  negative
loglikelihood, which is
proportional to
\bse
-n^{-1}\sumi \{Y_i \X_i\trans \bb- \exp (\bb\trans \X_i)\}.
\ese
Here we use $\exp(\bb\trans\X_i)$ to model
the mean of $Y_i$ because it preserves the positiveness of the
mean estimate, and it is a standard choice in the generalized linear
model \citep{mccullagh2019}.
It is useful to note that for normal noise $\U_i$,
we have the relation
\be
E\{\exp(\bb_t\trans \W_i - \bb_t\trans \bOmega\bb_t/2 )\mid\X_i\} &=&
\exp(\bb_t\trans\X_i),\label{eq:merr1}\\
E\{\exp(\bb_t\trans \W_i - \bb_t\trans \bOmega\bb_t/2 )(\W_i-\bOmega\bb_t)\mid\X_i\} &=&
\exp(\bb_t\trans\X_i)\X_i,\label{eq:merr2}\\
E[\exp(\bb_t\trans \W_i - \bb_t\trans \bOmega\bb_t/2
) \{(\W_i-\bOmega\bb_t)^{\otimes2} - \bOmega\} \mid\X_i] &=&
\exp(\bb_t\trans\X_i)\X_i^{\otimes2}. \label{eq:merr3}
\ee
Due to the conditional independence of $\W_i$ and $Y_i$  given $\X_i$,
(\ref{eq:merr1}) leads to
\bse
E \{Y_i \W_i\trans \bb_t - \exp (\bb_t\trans \W_i - \bb_t\trans \bOmega
\bb_t/2)\mid\X_i, Y_i\} = Y_i\X_i\trans \bb_t - \exp(\bb_t\trans\X_i).
\ese
Consequently, it is a reasonable practice to estimate $\bb$ by
minimizing the loss function
\be\label{eq:likelihood}
\mL (\bb)= - n^{-1}\sumi
\{Y_i \W_i\trans \bb - \exp (\bb\trans \W_i - \bb\trans \bOmega
\bb/2)\},
\ee
which has the same mean as the negative log-likelihood function
when $\X_i$ is accurately observed.  When $n >p$, the estimator for
$\bb$ can be obtained by minimizing $\mL (\bb)$ using the standard 
gradient descent method. However,
when $n <p$, without addition regularization,
optimizing (\ref{eq:likelihood}) is an ill-posed mathematical problem
because it does not have a unique solution.
To take into account the ultra-high dimension nature of
the model, using the fact that $\bb$ is sparse,
we propose to estimate
$\bb$ through solving the following constrained minimization problem
\be\label{eq:lap}
\min_{\|\bb\|_1 \leq R_1,  \|\bb\|_2 \leq R_2}\left\{\mL(\bb) + \rho_{\lambda}(\bb)\right\}
\ee
at suitable $R_1, R_2$, where $\rho_{\lambda}(\bb)$ is a suitable
penalty function. 
For convenience, define the
set $\{\bb: \|\bb\|_1 \leq
R_1, \|\bb\|_2 \leq R_2\}$ as the feasible set \citep{fletcher1980}. Here $R_1, R_2$ can be any
  constants that are greater than the true $\|\bb\|_1$ and $\|\bb\|_2$,
  respectively. The condition
  $\|\bb\|_1 \leq R_1$ is imposed to guarantee that the objective
  function satisfies the restricted eigenvalue condition discussed in
  \cite{loh2012} and therefore the objective function is convex in
  the feasible set, while the condition $\|\bb\|_2 \leq R_2$ is
  imposed to avoid the explosion of the mean
  function  $\exp (\bb\trans \W_i - \bb\trans \bOmega
\bb/2)$. In
  practice, we often set $R_1, R_2$ to be a constant times
the $L_1$, $L_2$ norms of the initial estimators of $\bb$.  
Here, with a slight abuse of notation, we
use the same symbol $\rho_\lambda$ to denote both multivariate and
univariate penalty functions and
let  $\rho_{\lambda}(\bb) =
\sum_{j=1}^{\|\bb\|_0}\rho_{\lambda}(\beta_j)$, where $\beta_j$ is the
$j$th element of $\bb$ and $\|\bb\|_0$ is the number of  nonzero elements in $\bb$. 

\subsection{Estimation under Hypotheses}\label{sec:hypo}

Consider testing the  hypothesis that $\C\bb_{t\mM} = \t +
\h_n$ for some $\h_n  \in \mathbb{R}^r$, where $\C$ is a
$r\times m$ matrix with  $r \leq m$, $\bb_{t\mM}$ is a $m$-dimensional
sub-vector of $\bb$ with index set $\mM$. The null
hypothesis holds  when $\h_n = \0$,  while the alternative hypothesis
holds
when $\h_n\ne\0$. For example if $\t = 0$, $\h_n = 1$,  $\C = (1, 0)$, $\mM$ contains the
index of the first element in $\bb_t$, then
testing $\C\bb_{t\mM} = \t +
\h_n$  is  testing the null hypothesis that
$\beta_{t1} = 0$ versus the alterative that $\beta_{t1} = 1$. Similarly, we
can test $\beta_{t1} - \beta_{t2} = 0$  versus $\beta_{t1} - \beta_{t2} \neq 0$ by
choosing $\C=(1,  -1)$, $\t=0$, $\h_n=0$ or nonzero, and $\mM=\{1,2\}$. In summary, by varying 
$\C$, $\t$, $\h_n$, and $\mM$, we can generate different linear
hypotheses to test. 
Without loss of generality, we assume
$\bb_\mM$ contains the first $m$ elements of $\bb$.  Further, let
$\bb_\mM^c$ be the
vector containing elements that are not in $\mM$, i.e. the last
$p-m$ components of $\bb$. Let $S \subseteq
\mM^c$ be the index
set of the nonzero  elements of
$\bb_{t\mM^c}$.
We assume  $\bb_{t\mM^c}$
to be $k$ sparse, i.e. $|S| = k$. Note that $k$ is allowed to diverge
with $n$. Without loss of generality, we assume
the first $k$ elements in $\bb_{t\mM^c}$ are none zero.

Suppose we are interested in testing whether $\C\bb_{t\mM}  = \t$ or
not.  Under the null hypothesis that $H_0: \C\bb_{t\mM} =
\t, $  we modify the general estimation
strategy slightly and  consider the estimator resulting
from the equality and inequality constrained minimization:
\be\label{eq:lapHT}
\wh {\bb} =\argmin_{\|\bb\|_1 \leq R_1,  \|\bb\|_2 \leq R_2}\left\{\mL
(\bb) + \rho_{\lambda}(\bb_{\mM^c})\right\}, \mbox{ s.t. } \C\bb_{\mM} = \t
\ee
for  suitable $R_1, R_2$.
Without assuming the null hypothesis,
we consider a similar estimator resulting
from the inequality constrained minimization:
\be\label{eq:lapHT1}
\wh {\bb}_a =\argmin_{\|\bb\|_1 \leq R_1,  \|\bb\|_2 \leq R_2}\left\{\mL
(\bb) + \rho_{\lambda}(\bb_{\mM^c})\right\}.
\ee
Note that here, both (\ref{eq:lapHT}) and (\ref{eq:lapHT1}) are
slightly different from the general strategy in (\ref{eq:lap}), in that
we do not place the penalty $\rho_{\lambda}$ on the parameters in
$\mM$, which are to be tested for the linear relation $\C\bb_{t\mM} = \t$. This special
treatment is to avoid the situation that the penalty forces some
components in $\bb_{\mM}$ to be zero,  and therefore
the null hypothesis $\C\bb_{t\mM} = \t$ is affected not only by the data
but also by our penalization.
\subsection{Test statistics}\label{sec:testing}
We define
\be \label{eq:Qtruedefine} {\Q}(\bb) \equiv E\{\exp(\bb\trans\X)\X\X\trans\},
\ee
define the covariance of the residuals
\bse
\bSig(\bb)&\equiv& E [\{Y_i \W_i - \exp(\bb\trans\W_i - \bb\trans
\bOmega \bb/2)  (\W_i -
\bOmega \bb)\}^{\otimes2}],
\ese
and define
\bse
{\bPsi} (\bSig, \Q, \bb)\equiv (\C [\I_{m\times m}, {\bf 0}_{m \times k}]  {\Q}_{\mM \cup S,
	\mM\cup S}^{-1}({\bb}) \bSig_{\mM \cup S,
	\mM\cup S} ({\bb}) {\Q}_{\mM \cup S,
	\mM\cup S}^{-1} ({\bb}) [\I_{m\times m}, {\bf 0}_{m \times k}]\trans
\C\trans).
\ese
Furthermore, let $\wh{\bSig}(\bb)$ and $\wh{\Q}(\bb)$ be a sample
estimator of $\bSig(\bb)$ and $\Q(\bb)$, respectively.
To test $\C\bb_\mM= \t$, we introduce two statistics, the  Wald statistic
\be
T_W =n (\C \wh{\bb}_{a\mM} - \t)\trans  {\bPsi} (\wh\bSig, \wh\Q, \wh\bb_a) ^{-1}(\C
\wh{\bb}_{a\mM} - \t),\label{eq:waldtest}
\ee
and the score statistic
\be
T_S &=& n \left\{\frac{\partial \mL(\wh{\bb})}{\partial \bb\trans
}\right\}_{\mM \cup S} (\C [\I_{m\times m}, {\bf 0}_{m \times k}] \wh{\Q}_{ \mM \cup
	S,  \mM \cup S}^{-1} (\wh\bb))\trans \nonumber \\
&& \times \bPsi^{-1}(\wh\bSig, \wh\Q, \wh\bb) \C [\I_{m\times m}, {\bf 0}_{m \times k}] \wh{\Q}_{ \mM \cup
	S,  \mM \cup S}^{-1} (\wh\bb) \left\{\frac{\partial
	\mL(\wh{\bb})}{\partial \bb }\right\}_{\mM \cup S}. \label{eq:scoretest}
\ee
As we will show later in Section \ref{sec:asytest} that $T_W$ and
$T_S$ are both asymptotically chi-square distributed with $r$
degrees of freedom under the null hypothesis. Therefore, to control
the false discovery rate at level $\alpha$, we reject the null
hypothesis if $T_W >\chi^2_{1 - \alpha}(r)$ when we perform Wald test,
or if $T_S >\chi^2_{1 - \alpha}(r)$ when we perform score test.
Here  $\chi^2_{1 - \alpha}(r)$ is the
$1 - \alpha$ quantile of the chi-square distribution.

\section{Theoretical Properties}\label{sec:theory}
Define
\bse
\wc{\bb}_\mM \equiv \bb_{t\mM} - \C\trans (\C\C\trans)^{-1}\h_n,
\ese
and let $\wc\bb=(\wc\bb_{t\mM}\trans,\bb_{t\mM^c}\trans)\trans$.
Thus,  the last $p-m$ components of $\wc\bb$,
i.e. $\wc{\bb}_{\mM^c}$,  and the last $p-m$ components of $\bb_0$, i.e.
${\bb}_{0\mM^c}$, are identical.
However, the first $m$ components of
$\wc\bb$ and $\bb$ are different, in that $\C\wc{\bb}_\mM = \t$ under
both null and alternative, while
$\C\bb_{t\mM} = \t$ under the null alone.
Under some conditions,
we first show that the inequality and equality constrained estimator
$\wh{\bb}$  is a consistent estimator of
$\wc\bb$ regardless the null or the alternative holds, and when $\|\h_n\|_2 $ vanishes,  $\wh{\bb}$  is also
consistent as an estimator of the true parameter $\bb_t$. Furthermore,
we show that $\wh{\bb}_a$  is a consistent estimator of
$\bb_t$ regardless the null or the alternative holds. We then
establish the asymptotic linear form of the estimators of a subvector
$\wh\bb$ and a subvector of $\wh\bb_a$, which are formed by
components of $\bb_t$ that are either to be tested or nonzero.
Finally, using the asymptotic linear forms, we construct test statistics and prove
the convergence properties of these test statistics under both null
and alternative.

\subsection{Conditions}\label{sec:condition}
Before we proceed with the specific results,
we first list a set of assumptions on the univariate penalty function
$\rho_{\lambda}$ which are similar to those in \cite{loh2015} and
\cite{loh2017} .
\begin{enumerate}[label=(A\arabic*)]
	\item The function $\rho_{\lambda}(t)$ satisfies $\rho_{\lambda}(0) =0$
	and is symmetric around zero.
	\item On the nonnegative real line $t\ge 0$, the function $\rho_{\lambda}(t)$ is
	nondecreasing. Furthermore, $\rho_{\lambda}(t)$ is subadditive, i.e.
	$\rho_{\lambda}(t_1 + t_2)\leq \rho_{\lambda}(t_1) +
	\rho_{\lambda}(t_2)$ for all $t_1, t_2 \geq 0$.
	\item For $t>0$, the function $\rho_{\lambda}(t)/t$ is non-increasing
	in t.
	\item The function $\rho_{\lambda}(t)$ is differentiable at all $t\neq0$
	and sub-differentiable at $t=0$, with $\lim_{t \to 0+}
	\rho'_{\lambda}(t) = \lambda $, where $\rho'(t)$ denotes the derivative of
	$\rho(t)$. Together with the symmetric Condition in  (A1), this
	leads to  $\lim_{t \to 0-}  \rho'_{\lambda}(t) = -\lambda $.
	\item There exists $\mu >0$ so that $\rho_{\lambda}(t) + \mu t^2/2$ is
	convex.
	\item There exists a $\gamma \in (0, +\infty)$ such that
	$\rho_{\lambda}'(t) = 0$ for all $t \geq \gamma\lambda$.
\end{enumerate}
Conditions (A1)--(A3) are some general requirements as
discussed in \cite{zhang2012}. Condition (A4)
restricts the class of penalties by
excluding regularizers that  are not differentiable at 0,
for example, the lasso penalty is excluded.
Condition (A5) is known as weak convexity \citep{vial1982, chen2014} and is a type of curvature constraint that controls the level
of nonconvexity of $\rho_{\lambda}$. Condition (A6) is imposed to
allow penalty to be zero if the estimator is
$\gamma \lambda$ away from zero, which
removes the estimation bias for the nonzero parameters.
We say $\rho_{\lambda}$ is  $\mu$-amenable if Conditions
(A1)--(A5) hold, and
we name $\rho_{\lambda}$ $(\mu, \gamma)$-amenable if Conditions
(A1)--(A6) hold. The $(\mu, \gamma)$-amenable penalty includes
the smoothly clipped
absolute deviation (SCAD) and the minimax concave penalty
\citep{loh2017}.

We need some additional regularity conditions to support
the theoretical development. These conditions impose upper and lower
bounds on various quantities to ensure that the upper bounds are finite and
the lower bounds are positive. They also restrict the relation between the
sample size and parameter number so that $\log(p)/n \to 0$ in a
slow rate of $1/\{\log(n)\}^2$. 
To save space, we only provide a discussion of these
conditions here, while provide the details in the supplementary
material. Specifically,
Condition \ref{con:boundX} (a) is a standard assumption 
  used in noisy data problem such as that used in
  \cite{csenturk2005} and is 
  usually satisfied in practice. Condition \ref{con:boundX} (b)
  guarantees the boundedness and the invertibility of the Hessian matrix (\ref{eq:merr3}),
  i.e. the second derivative of the noise free log likelihood. 
Conditions \ref{con:Wx} and \ref{con:Yx} bound the total
  variability of both the response $Y$ and the noise $\U$
  marginally and conditionally on the covariates $\X$. Similar
  requirement is also assumed in \cite{loh2012}. Condition
  \ref{con:pn} shows that the dimension of the covariate can grow
  exponentially faster than the sample size. Finally, \cite{jiang2021}
  have discussed the Conditions 
\ref{con:ww1}--\ref{con:ww3} and provided examples showing that the conditions
are 
usually satisfied in practice.

\subsection{Consistency}\label{sec:consistency}
We first show that the equality and inequality constrained estimator
$\wh\bb$ is a consistent estimator of $\wc\bb$ in Theorems
\ref{from:th1} and \ref{from:lem:th1}, which is the same as
the true parameter $\bb_t$, except that the first $m$ components are
adjusted to ensure that $H_0$ holds for $\wc\bb$.

\begin{Th}\label{from:th1}
	Define
	\bse
	\alpha_1 \equiv \min_{\|\bb\|_1\le R_1, \|\bb\|_2\le R_2}
	\alpha_{\min}
	[E\{\exp(\bb\trans\X_i)\X_i\X_i\trans\}]  /2.
	\ese
	Assume $\|\C_r^{-1
	}\C_{m -r}\|_2  = O(1)$, $\rho_\lambda$ satisfies
	Conditions (A1) -- (A6) and Conditions \ref{con:boundX} --
        \ref{con:ww2} in the supplementary material hold. Assume $\alpha_1 > 3/4\mu$, and $\bb$ is in the
	feasible set.
	Let $\lambda$ satisfy
	\bse
	4 \max\left\{\|\partial \mL (\wc \bb)/\partial \bb\|_\infty,
	\alpha_1 (\log(p)/n)^{1/4}\right\}\leq \lambda \leq \frac{\alpha_1}{6 R_1}
	\ese
	and $n \geq
	\log (p) \max(16R_1^4\tau_1^4/\alpha_1^4, 64 R_1^4
	\tau_1^2/\alpha_1^2)$.
	Write $t_1\equiv\sqrt{r}\|\C_r^{-1}\C_{m-r}\|_2 + \sqrt{m
		-r}$ and
	$
	t\equiv (6 \lambda \sqrt{k} + 2  \lambda t_1 )(4 \alpha_1 - 3\mu)^{-1}.
	$
	Then the local minimum of
	(\ref{eq:lapHT}) satisfies the error bounds
	\bse
	\|\wh {\bb} -\wc \bb\|_2 \leq t.
	\ese
	and
	\bse
	\|\wh {\bb} -\wc  \bb\|_1 \leq (4 \sqrt{k} +  t_1) t.
	\ese
\end{Th}
Following the similar argument, we also show that the inequality constrained estimator
$\wh{\bb}_a$ is a consistent estimator of
the true parameter $\bb$.
\begin{Th}\label{from:lem:th1} Let
	\bse
	\alpha_1 = \min_{\|\bb\|_1\le R_1, \|\bb\|_2\le R_2}
	\alpha_{\min}
	[E\{\exp(\bb\trans\X_i)\X_i\X_i\trans\}]  /2
	\ese
	and let $\rho_\lambda$ satisfy
	Conditions (A1) -- (A6) and Conditions \ref{con:boundX} --
        \ref{con:ww2} in the supplementary material hold. Assume $\alpha_1 > 3/4\mu$, and $\bb$ is in the
	feasible set. Let $\lambda$ satisfy
	\bse
	4 \max\left\{\|\partial \mL (\bb_t)/\partial \bb\|_\infty,
	\alpha_1 (\log(p)/n)^{1/4}\right\}\leq \lambda \leq \frac{\alpha_1}{6 R_1}
	\ese
	and $n \geq
	\log (p) \max(16R_1^4\tau_1^4/\alpha_1^4, 64 R_1^4 \tau_1^2/\alpha_1^2)$. Then the local minimum of
	(\ref{eq:lapHT1}) satisfies the error bounds
	\bse
	\|\wh {\bb}_a - \bb_t\|_2 \leq  \frac{6 \lambda \sqrt{k} + 2  \lambda \sqrt{m} }{4 \alpha_1 - 3\mu}.
	\ese
	and
	\bse
	\|\wh {\bb}_a -  \bb_t\|_1 &\leq&  (4 \sqrt{k} + \sqrt{m}) \frac{6 \lambda \sqrt{k} + 2  \lambda   \sqrt{m} }{4 \alpha_1 - 3\mu}.
	\ese
\end{Th}
Theorems \ref{from:th1} and \ref{from:lem:th1} suggest that
when $\log(p)/n \to 0$,  and when $\lambda$ is suitably
chosen, for example,  $\lambda$ is at least no smaller than
$O[\{\log(p)/n\}^{1/4}]$,
both $\wh{\bb}$ and $\wh{\bb}_a$ converge
to their corresponding true values in terms of both
$l_1$ and $l_2$ norms, as long as
$k$ and $m$ grow slower
than $\{n/\log(p)\}^{1/2}$.
These theoretical results suggest that
the dimension of $\bb_{t\mM}$, i.e.,  the number of parameters involved
in the tests, and the number of
nonzero entries in $\bb_t$ can grow at a slower rate of
$\{n/\log(p)\}^{1/2}$ under noisy Poisson model. These
results also assist us to find reasonable ranges for $\lambda$ in practice to
obtain consistent estimators.

\subsection{Asymptotic linear forms}\label{sec:normality}
We denote $\wt{\bb}$ as a stationary point
of (\ref{eq:lapHT}), which satisfies
the first order
condition that
\be\label{eq:firstorderHT}
\{\partial \mL(\wt \bb)/\partial \bb\trans+ \partial \rho_\lambda(\wt
\bb_{\mM^c})/\partial \bb_{\mM^c}\trans \A \}(\bb - \wt {\bb})\geq 0,
\ee
for all $\bb \in \mathbb R^p$ in the feasible set and satisfies $\C\bb_{\mM} = \t$.
Here $\A =
(\0_{p-m, m}, \I_{p-m, p-m})$ is a matrix that
satisfies  $\|\A\|_\infty = \|\A\|_1 = 1$.
Likewise, we denote $\wt{\bb}_a$ as a stationary point of
(\ref{eq:lapHT1}), which satisfies
the first order
condition that
\be\label{eq:firstorderHT1}
\{\partial \mL(\wt \bb_a)/\partial \bb\trans+ \partial \rho_\lambda(\wt
\bb_{a \mM^c})/\partial \bb_{a \mM^c}\trans \A \}(\bb_a - \wt
{\bb}_a)\geq 0,
\ee
for all $\bb_a \in \mathbb R^p$ in the feasible set.

To show the asymptotic normality of $\wh\bb$ and $\wh\bb_a$, our first
step is to establish that the local minimizers
$\wt{\bb}$ and  $\wt{\bb}_a$ achieve variable selection consistency.
To do this, we follow the prime-dual construction
introduced in \cite{wainwright2009}. We first show that both
\be\label{eq:lapHTsub}
\min_{\|\bb\|_1 \leq R_1,  \|\bb\|_2 \leq R_2, \bb \in
	\mathbb{R}^{\mM \cup S}}\left\{\mL
(\bb) + \rho_{\lambda}(\bb_{\mM^c})\right\}, \text{ such that }
\C\bb_{\mM} = \t
\ee
and
\be\label{eq:lapHTsub1}
\min_{\|\bb\|_1 \leq R_1,  \|\bb\|_2 \leq R_2, \bb \in
	\mathbb{R}^{\mM \cup S}}\left\{\mL
(\bb) + \rho_{\lambda}(\bb_{\mM^c})\right\}
\ee
have unique local minimizer in the interior of the feasible
set.
Then we show that all
stationary points of (\ref{eq:lapHT}) and  (\ref{eq:lapHT1})
must have support in
$\mM\cup S$. Since the local minimizers of (\ref{eq:lapHT}) and
(\ref{eq:lapHT1}) are automatically stationary points of
(\ref{eq:lapHT}) and  (\ref{eq:lapHT1}) respectively,  the
local minimizers of (\ref{eq:lapHT}) and
(\ref{eq:lapHT1}) must also have support in $\mM\cup S$.
Therefore, the local minimizers of (\ref{eq:lapHT}) and
(\ref{eq:lapHT1}) are actually the local minimizers of (\ref{eq:lapHTsub})
and (\ref{eq:lapHTsub1}) respectively, so are also unique.  In other
words, $\wh{\bb}$ and $\wh{\bb}_a$ are respectively the unique solution of
(\ref{eq:lapHTsub}) and  (\ref{eq:lapHTsub1}) hence
achieve the variable
selection consistency.
The details of the above analysis are presented in
Theorem \ref{th:1} and Theorem \ref{from:th:1} in the Appendix A in
the supplementary material.

In our second step to  establish the asymptotic distribution properties of
$\wh\bb$ and $\wh\bb_a$, we
define
\bse
\wh {\Q}(\bb) = \frac{\partial^2  \mL({\bb})}{\partial \bb \partial
	\bb\trans},
\ese
and define
\bse
\A_2 &=&  [\I_{m\times m}, {\bf 0}_{m
	\times k}]\trans \C\trans (\C [\I_{m\times m}, {\bf 0}_{m \times k}]\{\wh{\Q}_{\mM \cup S,
	\mM\cup S}(\bb^*)\}^{-1}\\
&&\times [\I_{m\times m}, {\bf 0}_{m \times k}]\trans
\C\trans)^{-1} \C [\I_{m\times m}, {\bf 0}_{m \times  k}],
\ese
where $\bb^*$ is the point in between $\wh{\bb}$ and $\bb_t$ and
\bse
\A_2^*  &=&  [\I_{m\times m}, {\bf 0}_{m
	\times k}]\trans \C\trans (\C [\I_{m\times m}, {\bf 0}_{m \times k}]
\{{\Q}_{\mM \cup S,
	\mM\cup S}(\bb)\}^{-1} \\
&&\times [\I_{m\times m}, {\bf 0}_{m \times k}]\trans
\C\trans)^{-1} \C [\I_{m\times m}, {\bf 0}_{m \times  k}],
\ese
where ${\Q}(\bb) =E\{\exp(\bb\trans\X)\X\X\trans\} $ is
defined in (\ref{eq:Qtruedefine}).
Based on the variable selection consistency established in
the first step,  we derive the
asymptotic  linear form of $\wh{\bb}_{\mM\cup S}$ and
$\wh{\bb}_{a \mM\cup S}$  under null and alternative
hypothesis in Theorems \ref{th:2} and  \ref{from:th:2}, respectively.
\begin{Th}\label{th:2}
	Assume $\rho_\lambda$ satisfies
	Conditions (A1) -- (A6) and Conditions \ref{con:boundX} -- \ref{con:ww3} in the
	supplementary material hold, $\lambda =
	O_p[\{\log(p)/n\}^{1/4}]$,
	$\|\C_r^{-1}\C_{m-r}\|_2=O(1)$,
	and $\lambda \leq
	\alpha_1/(8R_1)$. Further we assume the boundedness $\| \{\Q_{(\mM \cup S), \mM \cup
		S}^{-1}(\bb_t)\}\|_\infty \leq c_\infty$,  and $\|\{ \Q_{(\mM \cup S),
		\mM \cup S}(\bb_t)\}^{-1}\A_2  \Q_{(\mM \cup S), \mM \cup
		S}^{-1}(\bb_t)\|_\infty \leq c_\infty$. In addition assume $\|\h_n\| _2
	= O\{\sqrt{\max (m + k - r,
		r)/n}\}$, $\min(|\beta_{j}|)\geq \lambda (\gamma + 5c_\infty )$ for $j
	\in S$ and  $n\geq
	c_\infty (m+ k)^4  \log(p)$.  Then we have
	\bse
	&&\wh{\bb}_{\mM \cup S}- \bb_{t\mM \cup S}\\
	&=&  - (\{{\Q}_{\mM \cup S,
		<  \mM\cup S}(\bb_t)\}^{-1} - \{{\Q}_{\mM \cup S, \mM\cup
		S}(\bb_t)\}^{-1} \A_2^* \\
	&&\times \{{\Q}_{\mM
		\cup S, \mM\cup S}(\bb_t)\}^{-1} )
	\left\{\frac{\partial
		\mL( \bb_t)}{\partial \bb}\right\}_{\mM
		\cup S}\{1 + o_p(1)\}\\
	&& +\{{\Q}_{\mM \cup S, \mM\cup S}(\bb_t)\}^{-1} \A_2^*  [\{(\C\C\trans)^{-1}\C, {\bf 0}_{ r \times
		k}\}\trans \h_n ]\{1 + o_p(1)\}
	\ese
	and $\wh{\bb}_{(\mM \cup S)^c} = \bf 0$.
\end{Th}

\begin{Th}\label{from:th:2}
	Assume $\rho_\lambda$ satisfies
	Conditions (A1) -- (A6) and Conditions \ref{con:boundX} -- \ref{con:ww3} in the
	supplementary material hold, $\lambda =
	O_p[\{\log(p)/n\}^{1/4}]$, and $\lambda \leq
	\alpha_1/(8R_1)$. Further we assume $\| \{\Q_{(\mM \cup S), \mM \cup
		S}(\bb_t)\}^{-1}\|_\infty \leq c_\infty$, $\min(|\beta_{j}|)\geq \lambda (\gamma + 5c_\infty )$ for $j
	\in S$ and $n\geq
	c_\infty (m+ k)^4  \log(p)$. Then we have
	\bse
	\wh{\bb}_{a \mM \cup S}- \bb_{t\mM \cup S}
	&=&  -\{ {\Q}_{\mM \cup S,
		\mM\cup S}(\bb_t)\}^{-1} \left\{\frac{\partial
		\mL( \bb_t)}{\partial \bb}\right\}_{\mM
		\cup S}\{1 + o_p(1)\}
	\ese
	and $\wh{\bb}_{(\mM \cup S)^c} = \bf 0$.
\end{Th}

Theorems \ref{th:2} and  \ref{from:th:2} suggest that the
asymptotic linear forms of $\wh{\bb}_{\mM\cup  S}$ and $\wh{\bb}_{a \mM\cup  S}$ are
the usual product of the inverse of Hessian matrix  and the
score function. Furthermore, only the first  $(m+ k)\times (m+k)$
block in the Hessian matrix and the first $m+k$ elements in the score
function contribute  to the asymptotic distribution. Therefore,  when
$m+k$ grows slower than $\{n/\log(p)\}^{1/4}$ and $\|\h_n\| \to 0$, it
is easy to see that the asymptotic
linear forms converge in distribution to Gaussian random vectors.  It
is  worth mentioning that the minimal signal condition
$\min(|\beta_{j}|)\geq \lambda (\gamma + 5c_\infty )$ for $j 
	\in S$ is a standard requirement for the 
        optimization using  nonconvex penalty such as SCAD
        \citep{fan2001}. This condition is also very weak because $\lambda
        \to 0$, which allows the minimal signal vanishing to zero.
 
\subsection{Asymptotic distribution of the test statistics}\label{sec:asytest}

To study the asymptotic behavior of  $T_S$ and $T_W$, we first
investigate the distribution of their asymptotic form $T_0$
defined by \bse
T_0 \equiv (\bomega_n + \sqrt{n}\h_n) \trans  \bPsi^{-1}(\bSig, \Q, \bb_t)  (\bomega_n +
\sqrt{n}\h_n),
\ese
where
\bse
\bomega_n = - \sqrt{n} \C [\I_{m\times m}, {\bf 0}_{m \times k}] {\Q}_{ \mM \cup
	S,  \mM \cup S}^{-1} (\bb_t)\left\{\frac{\partial \mL(\bb_t)}{\partial
	\bb}\right\}_{\mM \cup S}.
\ese
As shown in Lemma \ref{th:chisq}, $T_0$ is asymptotically  noncentral
chi-square distributed with the noncentral parameter approaches
$n \h_n\trans \bPsi^{-1}(\bSig, \Q, \bb_t) \h_n $.

\begin{Lem}\label{th:chisq}
	Assume $\rho_\lambda$ satisfies
	Conditions (A1) -- (A6)  and Conditions \ref{con:boundX} and \ref{con:3ord} in the
	supplementary material hold and $n\geq
	c_\infty (m+ k)^4  \log(p)$, then
	\bse
	\lim_{n\to \infty}\sup_{\mathcal{C}} |\Pr(T_0 \leq x) -
	\Pr\{\chi^2(r, n \h_n\trans \bPsi^{-1}(\bSig, \Q, \bb_t) \h_n ) \leq x\}| = 0,
	\ese
	where $\chi^2(r, \gamma)$ is a non-central chi-square
	random variable, with non-centrality parameter $\gamma$.
      \end{Lem}
      Here Condition \ref{con:3ord} provides upper bound of the third
      moment of each
      summand in $\bomega$ (note that $\partial \mL(\bb_t)/\partial
	\bb$ is the summation of the derivatives of the negative
        log-likelihood from  $n$ samples), which is a necessary  condition to establish
      convergence in distribution. See Theorem 3.1 in \cite{shi2019linear} for
      example. 
To establish the asymptotic distribution of $T_W$ and $T_S$, in Theorems \ref{th:TW} and
\ref{th:TS} respectively, we show
that $T_W$ and $T_S$ are close to $T_0$, hence has the same testing
property asymptotically when $r$ is finite.
\begin{Th}\label{th:TW}
	Assume the conditions in Theorem \ref{from:th:2} and Conditions
\ref{con:3ord}	and \ref{con:psieigen} in the Section \ref{sec:assD}
in the supplementary material hold, we have
	$T_W - T_0 = o_p(r)$. Therefore,
	\bse
	\lim_{n\to \infty}\sup_{\mathcal{C}} |\Pr(T_W \leq x) -
	\Pr\{\chi^2(r, n \h_n\trans \bPsi^{-1}(\bSig, \Q, \bb_t)\h_n ) \leq x\}| = 0,
	\ese
	where $\chi^2(r, \gamma)$ is a non-central chi-square
	random variable, with non-centrality parameter $\gamma$.
\end{Th}

\begin{Th}\label{th:TS}
	Assume the conditions in Theorem \ref{th:2}, Conditions
\ref{con:3ord}	and \ref{con:psieigen} in the Section \ref{sec:assD}
in the supplementary material hold, we have
	$T_S - T_0 = o_p(r)$.  Therefore,
	\bse
	\lim_{n\to \infty}\sup_{\mathcal{C}} |\Pr(T_S \leq x) -
	\Pr\{\chi^2(r, n \h_n\trans \bPsi^{-1}(\bSig, \Q, \bb_t)\h_n ) \leq x\}| = 0,
	\ese
	where $\chi^2(r, \gamma)$ is a non-central chi-square
	random variable, with non-centrality parameter $\gamma$.
\end{Th}
Here 
Condition \ref{con:psieigen}  in the Section \ref{sec:assD}  is a regularity condition ensures 
${\bPsi} (\bSig, \Q, \bb_t)$ to be positive definite.  Theorems \ref{th:TW} and \ref{th:TS} show that the two test statistics
$T_W$ and $T_S$ indeed have the same $\chi^2(r,\gamma)$ distribution
as $T_0$ in large samples, hence can be used to perform the standard
chi-square test. A curious question is whether or not a
  likelihood ratio type of test can also be constructed. We feel it is
  hard in this context because it is almost impossible
  to obtain a likelihood function in the functional measurement error
  context. Much work is
  needed to overcome this obstacles.

\section{Numerical Implementation}\label{sec:num}
\subsection{Computational algorithms}\label{sec:comp}

We compute the estimators $\wh {\bb}$ and $\wh {\bb}_a$ using
the popular ADMM. In what follows,
we only detail the algorithm to estimate $\wh {\bb}$. The estimator $\wh
{\bb}_a$  can be computed in a similar way. For a given $\lambda$, we
consider
\bse
\wh {\bb}=\argmin_{\|\bb\|_1 \leq R_1,  \|\bb\|_2 \leq R_2}\left\{\mL
(\bb) + \rho_{\lambda}(\bb_{\mM^c})\right\}, \mbox{ s.t. } \C\bb_{\mM} = \t
\ese
for constants  $R_1, R_2$. Similar to \cite{shi2019linear}, this optimization problem is equivalent to
\bse
(\wh {\bb},\wh \btheta)=\argmin_{\|\bb\|_1 \leq R_1,  \|\bb\|_2 \leq R_2}\left\{\mL
(\bb) + \rho_{\lambda}(\bb_{\mM^c})\right\}, \mbox{ s.t. } \C\bb_{\mM} = \t, \bb_{\mM^c}=\btheta.
\ese
By the augmented Lagrangian method, the estimators can be obtained by minimizing
\bse
L(\bb,\btheta,\bm v)=\mL
(\bb) + \rho_{\lambda}(\bb_{\mM^c})+\bm v^{\trans}
\begin{pmatrix}
	
	\C\bb_{\mM}-\t  \\
	
	\bb_{\mM^c} -\bm\theta
	
\end{pmatrix}+\frac{\rho}{2} \begin{Vmatrix}
	
	\C\bb_{\mM}-\t  \\
	
	\bb_{\mM^c} -\bm\theta
	
\end{Vmatrix}^2_2,
\ese
with $\|\bb\|_1 \leq R_1$, $ \|\bb\|_2 \leq R_2$, where the dual variables $\bm v$
are Lagrange multipliers and $\rho>0$ is a given penalty parameter. We
compute the estimators of $(\bb,\btheta,\bm v)$ through iterations.
Let the sup-script (t) indicate the $t$-th iteration, we describe the
main steps of ADMM methods in Algorithm \ref{alg1}.

In the implementation, the initial value ${\bb}^{(0)}$ can be
computed by a penalized Poisson regression following
  \cite{jiang2021}. For the radii $R_1$ and 
$R_2$, we consider $R_1=\sqrt{2}R_2$ and $R_2=1.5\|\bb\|_2^{(0)}$. 
In the implementation, if the algorithm converges to the
  boundary, we can increase the corresponding norm $R_1$ or $R_2$ slightly. In
  contrast, if multiple minimum problems are encountered, we can
  decrease $R_1$ and when the estimation procedure leads to
a very large $\exp(\bb\trans\X)$,  we can decrease $R_2$, gradually.
The
tuning parameter $\lambda$ is selected by minimizing
\be\label{eq:BIC}
\textrm{BIC}(\lambda) =n\mL(\wh{\bb})+c_n\|\wh {\bb}\|_0
\ee
with respect to $\lambda$, where $c_n$ is a positive number that may
depend on $n$.  In our analysis,  we follow \cite{shi2019linear} to
adopt $c_n=\max\{\log n, \log (\log(n))\log p\}$.
For simplicity, we set $\rho=1$.
\begin{algorithm}[!t]
	\caption{ADMM Algorithm for estimating $\wh
          {\bb}$. }\label{alg1}
        \flushleft
	For $t= 0, 1, \ldots, t_{\text{max}}$, perform:
	
	\hspace{0.5cm} Step 1. Use the Newton-Raphson algorithm to solve \eqref{step1} to obtain  $\widetilde{\bb}^{(t+1)}$.
	\be\label{step1}
	\widetilde{\bb}^{(t+1)}&=&\argmin_{\bb}\left\{\mL
	(\bb)+\bm v^{(t)\trans}
	\begin{pmatrix}
		
		\C\bb_{\mM}-\t  \\
		
		\bb_{\mM^c} -\bm\theta^{(t)}
		
	\end{pmatrix}+\frac{\rho}{2} \begin{Vmatrix}
		
		\C\bb_{\mM}-\t  \\
		
		\bb_{\mM^c} -\bm\theta^{(t)}
		
	\end{Vmatrix}^2_2\right\}.
	\ee
	\hspace{0.5cm} Step 2.  Project $\widetilde{\bb}^{(t+1)}$ to a $L_1$ ball with radius $R_1$ to obtain $\breve{\bb}^{(t+1)}$ by the simplex projection method \citep{duchi2008}. If $||\breve{\bb}^{(t+1)}||_2>R_2$, we shrink it to get ${\bb}^{(t+1)}=\breve{\bb}^{(t+1)}R_2/||\breve{\bb}^{(t+1)}||_2$. Otherwise, ${\bb}^{(t+1)}=\breve{\bb}^{(t+1)}$.
	
	\hspace{0.5cm} Step 3. Obtain  $\bm\theta^{(t+1)}$ by solving \eqref{step1}, where the penalty term we use is SCAD with $a=3.7$.
	\be\label{step3}
	\bm\theta^{(t+1)}&=&\argmin_{\bm\theta}\left\{\rho_{\lambda}(\bm\theta)+\frac{\rho}{2} \begin{Vmatrix}
		
		\bb^{(t+1)}_{\mM^c} -\bm\theta
		
	\end{Vmatrix}^2_2+\bm v^{(t)\trans}
	\begin{pmatrix}
		
		\C\bb^{(t+1)}_{\mM}-\t  \\
		
		\bb^{(t+1)}_{\mM^c} -\bm\theta
		
	\end{pmatrix}\right\}.
	\ee
	\hspace{0.5cm} Step 4. Update the dual variables $\bm v$ by
	\bse
	\bm v^{(t+1)}&=&\bm v^{(t)}+\rho
	\begin{pmatrix}
		
		\C\bb^{(t+1)} _{\mM}-\t  \\
		
		\bb^{(t+1)} _{\mM^c} -\bm\theta^{(t+1)}
		
	\end{pmatrix}.
	\ese
	\hspace{0.5cm} Step 5. If stopping rule $\|{\bb}^{(t+1)}-{\bb}^{(t)}\|_2 \le \delta_{\textrm{tol}}$ or $\|{\btheta}^{(t+1)}-{\btheta}^{(t)}\|_2 \le \delta_{\textrm{tol}}$ is satisfied, where $\delta_{\textrm{tol}}$ denotes the tolerance of error, then terminate
	the algorithm.
	
	End of the main loop.
\end{algorithm}

\subsection{Simulation Experiments}\label{sec:simulation}

We generate the outcome $Y_i$ from the Poisson model
\bse
\Pr(Y_i=y \mid \X_i)=\exp\{-\exp(\bb\trans\X_i) \}\exp(y
\bb\trans\X_i) /y!,
\ese
where the covariates $\X_i=(X_{i,1},\ldots, X_{i,p})\trans$ are generated
from two distributions: (I) the  multivariate normal distribution with mean zero
and covariance matrix $\Sigma$. (II) the uniform distribution in the
interval $(-\sqrt{6}/2, \sqrt{6}/2)$.
To generate correlated uniform distribution, we first draw covariates independently from $\cal{U}$$(-\sqrt{6}/2, \sqrt{6}/2)$, and then  transform these  covariates by multiplying the
Choleski factorization of covariance $\Sigma$.
We consider two forms of the
covariance matrix: uncorrelated structure $\Sigma=0.5\I_p$ and  correlated with auto-regressive AR(1) structure $\Sigma=(0.5^{|i-j|+1})_{p\times p}$ for $i,j=1,\ldots, p$.
Furthermore,
the noise $\U_i$ is drawn from the multivariate normal
distribution with mean zero and covariance matrix $\Omega=0.1\Sigma$.
The true coefficient
$\bb =(\beta_1,\ldots, \beta_p)\trans=(0.75,
-0.75+h_2,h_3, 0, \ldots, 0, h_p)\trans$. Here $h_j, j = 2, 3, p$ are
assigned various values to check the empirical powers of the
tests. We set  $h_j=0$ when
$j\ne 2, 3$ or $p$.
For simplicity, the initial ${\bb}^{(0)}$ is set to be a $p$-dimensional zero vector. We
select parameter $\lambda$ as described in Section \ref{sec:comp}. The
candidate list for $\lambda$ is
$\{e^{-2.5},e^{-2.245},\ldots,e^{0.5}\}$ of length 41.
We consider sample size $n = 300, 500$
  and covariate dimension $p=50, 350, 600$. The
tolerance of error $\delta_{\textrm{tol}}=10^{-4}$.
We repeat each setting 500 times, and report the size and power of the
proposed tests under different hypotheses. We perform the tests at type I error
$\alpha = 0.05$ in the following scenarios.

\subsubsection{Univariate parameter testing}
We first consider the following three hypotheses on a single element in $\bb$.
\bse
H_{0,1}:\beta_{2}=-0.75, &\ v.s.& \ H_{a,1}:  \beta_{2} \neq-0.75.\\
H_{0,2}:\beta_{3}=0,  &\ v.s. & \ H_{a,2}:  \beta_{3} \neq0.\\
H_{0,3}:\beta_{p}=0,  &\ v.s. & \  H_{a,3}:\beta_{p} \neq0.
\ese
To test a hypothesis set regarding $\beta_j$,
we simulate data with $h_j=0, 0.1, 0.2, 0.4$, while
set $h_k=0$ for $k\neq j$. For example, to test $H_{0, 1}$ and $H_{a, 1}$, we simulate data with  $h_3=0$, $h_p=0$, and $h_2=0, 0.1,
0.2, 0.4$. When $h_2=0$, the null hypothesis $H_{0,1}$ holds, we study
the type I error of the test. On the other hand,
when $h_2=0.1$ to $0.4$, the alternative hypothesis is true, which allows us to examine
the power of the test.
Tables \ref{S1} and \ref{S1_500} summarize  the empirical type I error and powers of
the Wald and score tests. It
is clear that  the empirical type I errors are controlled at the
nominal  level 0.05 in all scenarios,  indicating that the proposed
tests are consistent. The powers of the Wald
and score tests increase gradually when the  magnitude of  $|h_j|$'s
increases, and have satisfactory powers in general.
The Wald and score tests yield similar performances in all
scenarios. This finding is in accordance with theoretical analysis.
\begin{table}
  \footnotesize
  \centering
	\caption{ The empirical sizes and powers of Wald and score
          tests for univariate parameter testing with $n = 300$. }\label{S1}
	\begin{tabular}{c|c|cc|cc|cc|cc}\hline
		&&\multicolumn{4}{c|}{$\X \sim$ Normal}&\multicolumn{4}{c}{$\X \sim$ Uniform}\\\hline
		&&\multicolumn{2}{c|}{$\Sigma=0.5\I_p$}&\multicolumn{2}{c|}{$\Sigma=0.5^{|i-j|+1}$} &\multicolumn{2}{c|}{$\Sigma=0.5\I_p$}&\multicolumn{2}{c}{$\Sigma=0.5^{|i-j|+1}$} \\\hline
		&&$T_W$&$T_S$&$T_W$&$T_S$&$T_W$&$T_S$&$T_W$&$T_S$\\\hline
	\multirow{12}{*}{$p=50$}&	$\beta_{2}$&\multicolumn{8}{c}{$H_{0,1}:\beta_{2}=-0.75,\ v.s. \ H_{a,1}:  \beta_{2} \neq-0.75$}\\\cline{2-10}

		&-0.75&0.068&0.054&0.056&0.050&0.054&0.044&0.066&0.062   \\                                       
		&-0.65&0.352&0.288&0.222&0.176&0.292&0.276&0.232&0.208    \\                                                             			
		&-0.55&0.826&0.778&0.680&0.592&0.736&0.726&0.632&0.598     \\                                                            			
		&-0.35&1.000&0.996&0.996&0.976&1.000&1.000&0.990&0.972    \\\cline{2-10}	                                                       			
		&$\beta_{3}$&\multicolumn{8}{c}{$H_{0,2}:\beta_{3}=0, \ v.s. \ H_{a,2}:  \beta_{3} \neq0$ } \\\cline{2-10}			
		&0.0&0.056&0.046&0.058&0.038&0.056&0.046&0.068&0.060      \\
		&0.1&0.302&0.272&0.204&0.172&0.250&0.240&0.214&0.182       \\                                                                           			
		&0.2&0.752&0.724&0.554&0.524&0.692&0.682&0.530&0.508        \\                                                                          			
		&0.4&0.996&0.996&0.984&0.960&1.000&1.000&0.976&0.942       \\\cline{2-10}	                                                            			
		&$\beta_{p}$&\multicolumn{8}{c}{$H_{0,3}:\beta_{p}=0, \ v.s.  \  H_{a,3}:\beta_{p} \neq0$ }    \\\cline{2-10}	            			
		&0.0&0.060&0.044&0.056&0.042&0.062&0.054&0.052&0.050      \\
		&0.1&0.246&0.222&0.234&0.210&0.276&0.258&0.240&0.226       \\                                              			
		&0.2&0.708&0.672&0.666&0.634&0.714&0.698&0.708&0.682        \\                                             			
		&0.4&0.998&0.998&0.998&0.996&0.998&0.998&0.998&0.994       \\\hline     
\multirow{12}{*}{$p=350$}&	$\beta_{2}$&\multicolumn{8}{c}{$H_{0,1}:\beta_{2}=-0.75,\ v.s. \ H_{a,1}:  \beta_{2} \neq-0.75$}\\\cline{2-10}	
       &-0.75 &0.050&0.036&0.054&0.034&0.064&0.068&0.066&0.056\\
       &-0.65 &0.312&0.322&0.260&0.242&0.268&0.266&0.230&0.202\\
       &-0.55 &0.750&0.766&0.650&0.644&0.752&0.750&0.612&0.598\\
       &-0.35 &0.998&0.998&0.980&0.878&0.998&0.998&0.978&0.892\\\cline{2-10}	                                                          			
		&$\beta_{3}$&\multicolumn{8}{c}{$H_{0,2}:\beta_{3}=0, \ v.s. \ H_{a,2}:  \beta_{3} \neq0$ } \\\cline{2-10}		
&0.0&0.064&0.048&0.066&0.066&0.066&0.058&0.068&0.054\\
&0.1&0.328&0.330&0.224&0.200&0.270&0.262&0.198&0.164\\
&0.2&0.770&0.770&0.590&0.546&0.708&0.706&0.568&0.504\\
&0.4&1.000&1.000&0.950&0.846&1.000&1.000&0.942&0.830\\\cline{2-10}                                                      			
		&$\beta_{p}$&\multicolumn{8}{c}{$H_{0,3}:\beta_{p}=0, \ v.s.  \  H_{a,3}:\beta_{p} \neq0$ }    \\\cline{2-10}	            			
		&0.0&0.072&0.066&0.050&0.046&0.058&0.050&0.066&0.056\\
  &0.1      &0.346&0.342&0.208&0.198&0.250&0.250&0.220&0.206\\
  &0.2      &0.736&0.742&0.662&0.646&0.782&0.768&0.654&0.638\\
  &0.4      &1.000&1.000&0.996&0.994&1.000&1.000&0.994&0.992\\\hline                                            			
	\end{tabular}                                                                                                         			
      \end{table}

      \begin{table}
	\footnotesize
	\centering
	\caption{ The empirical sizes and powers of Wald and score tests for univariate parameter testing with $n=500$. }\label{S1_500}
	\begin{tabular}{c|c|cc|cc|cc|cc}\hline
		&&\multicolumn{4}{c|}{$\X \sim$ Normal}&\multicolumn{4}{c}{$\X \sim$ Uniform}\\\hline
		&&\multicolumn{2}{c|}{$\Sigma=0.5\I_p$}&\multicolumn{2}{c|}{$\Sigma=0.5^{|i-j|+1}$} &\multicolumn{2}{c|}{$\Sigma=0.5\I_p$}&\multicolumn{2}{c}{$\Sigma=0.5^{|i-j|+1}$} \\\hline
		&&$T_W$&$T_S$&$T_W$&$T_S$&$T_W$&$T_S$&$T_W$&$T_S$\\\hline
		\multirow{12}{*}{$p=50$}&	$\beta_{2}$&\multicolumn{8}{c}{$H_{0,1}:\beta_{2}=-0.75,\ v.s. \ H_{a,1}:  \beta_{2} \neq-0.75$}\\\cline{2-10}	
		&-0.75&0.066 & 0.054&0.044 & 0.040 &0.064&0.060&0.060&0.056 \\                                       
		&-0.65& 0.488 & 0.450& 0.346 & 0.316  &0.442&0.422&0.324 & 0.304\\                                                             			
		&-0.55& 0.954 & 0.950&  0.864 & 0.838 &0.910&0.906&0.800 & 0.792\\                                                            			
		&-0.35& 1.000 & 1.000& 1.000 & 1.000 &1.000&1.000&1.000 & 1.000\\\cline{2-10}	                                                       			
		&$\beta_{3}$&\multicolumn{8}{c}{$H_{0,2}:\beta_{3}=0, \ v.s. \ H_{a,2}:  \beta_{3} \neq0$ } \\\cline{2-10}			
		&0.0&0.048 & 0.046  &0.066 & 0.058 &0.052&0.050&0.056&0.056    \\
		&0.1&  0.402 & 0.382 &0.324 & 0.300&0.390&0.386&0.302&0.296     \\                                                                           			
		&0.2&   0.890 & 0.888&0.780 & 0.770 &0.892&0.884&0.778&0.766     \\                                                                          			
		&0.4& 1.000 & 1.000 &1.000 & 0.998 &1.000&1.000&1.000&1.000    \\\cline{2-10}	                                                            			
		&$\beta_{p}$&\multicolumn{8}{c}{$H_{0,3}:\beta_{p}=0, \ v.s.  \  H_{a,3}:\beta_{p} \neq0$ }    \\\cline{2-10}	      			
		&0.0& 0.066 &0.060&  0.064 & 0.058&0.050&0.048&0.050&0.048   \\
		&0.1& 0.400 &0.368&  0.350 & 0.338&0.390&0.380&0.336 & 0.320   \\                                              			
		&0.2&  0.922 &0.914&  0.896 & 0.878 &0.892&0.890&0.884 & 0.878  \\                                             			
		&0.4&  1.000 &1.000&  1.000 & 1.000 &1.000&1.000 &1.000 & 1.000\\\hline     
		\multirow{12}{*}{$p=600$}&	$\beta_{2}$&\multicolumn{8}{c}{$H_{0,1}:\beta_{2}=-0.75,\ v.s. \ H_{a,1}:  \beta_{2} \neq-0.75$}\\\cline{2-10}	
		&-0.75 &0.052&0.056&0.046&0.046&0.062&0.056&0.040&0.038\\
		&-0.65 &0.458&0.478&0.328&0.330&0.390&0.392&0.396&0.402\\
		&-0.55 &0.920&0.930&0.864&0.868&0.926&0.928&0.902&0.902\\
		&-0.35 &0.842&0.918&0.988&0.990&1.000&1.000&1.000&1.000\\\cline{2-10}	                                                          			
		&$\beta_{3}$&\multicolumn{8}{c}{$H_{0,2}:\beta_{3}=0, \ v.s. \ H_{a,2}:  \beta_{3} \neq0$ } \\\cline{2-10}		
		&0.0&0.076&0.066&0.062&0.058&0.070&0.066&0.066&0.056\\
		&0.1&0.454&0.452&0.344&0.342&0.392&0.392&0.454&0.454\\
		&0.2&0.904&0.904&0.832&0.822&0.902&0.894&0.870&0.862\\
		&0.4&0.974&0.974&0.998&0.980&1.000&1.000&0.998&0.894\\\cline{2-10}                                                      			
		&$\beta_{p}$&\multicolumn{8}{c}{$H_{0,3}:\beta_{p}=0, \ v.s.  \  H_{a,3}:\beta_{p} \neq0$ }    \\\cline{2-10}	            			
		&0.0&0.048&0.046&0.050&0.046&0.060&0.058&0.058&0.054\\
		&0.1&0.444&0.444&0.404&0.390&0.414&0.416&0.450&0.448\\
		&0.2&0.930&0.928&0.870&0.860&0.912&0.912&0.876&0.874\\
		&0.4&0.980&0.980&1.000&1.000&1.000&1.000&1.000&1.000\\\hline                                            			
	\end{tabular}                                                                                                         			
\end{table}

\subsubsection{Linear hypothesis testing}
We also consider the hypotheses that contain the linear combinations of two coefficient parameters:
\bse
H_{0,4}:\beta_{1}+\beta_{2}=0, &\ v.s.& \ H_{a,4}:  \beta_{1}+\beta_{2}\neq0.\\
H_{0,5}:\beta_{3}+\beta_{4}=0,  &\ v.s. & \ H_{a,5}: \beta_{3}+\beta_{4}\neq 0.\\
H_{0,6}:\beta_{1}+\beta_{p}=0.75,  &\ v.s. & \ H_{a,6}: \beta_{1}+\beta_{p}\neq 0.75.\\
H_{0,7}:\beta_{2}+\beta_{3}=-0.75,  &\ v.s. & \ H_{a,7}: \beta_{2}+\beta_{3}\neq -0.75.
\ese
For the first three sets of hypotheses, we still set $h_j=0, 0.1, 0.2, 0.4$ if
the hypothesis involves $\beta_{j}$ for  $j=2, 3, p$, and set $h_k=0$ if
the corresponding $\beta_{k}$ is not involved in the hypotheses. For
the last hypothesis $H_{0,7}$, we set $h_2=0$,  $h_p=0$  and vary
$h_3$ from $0$ to $0.4$.  Tables \ref{S2} and
  \ref{S2_500} show that the Wald and score 
tests control the type I error at nominal level, and their powers
improve when  $h_j$ increases.
\begin{table}
  \footnotesize
  \centering
	\caption{ The empirical size and power of Wald and score tests
          for linear hypothesis testing with $n = 300$. }\label{S2}
	\begin{tabular}{c|c|cc|cc|cc|cc}\hline
		&&\multicolumn{4}{c|}{$\X \sim$ Normal}&\multicolumn{4}{c}{$\X \sim$ Uniform}\\\hline
		&&\multicolumn{2}{c|}{$\Sigma=0.5\I_p$}&\multicolumn{2}{c|}{$\Sigma=0.5^{|i-j|+1}$} &\multicolumn{2}{c|}{$\Sigma=0.5\I_p$}&\multicolumn{2}{c}{$\Sigma=0.5^{|i-j|+1}$} \\\hline
		&&$T_W$&$T_S$&$T_W$&$T_S$&$T_W$&$T_S$&$T_W$&$T_S$\\\hline
	\multirow{16}{*}{$p=50$}&	$\beta_{1} + \beta_{2}$&\multicolumn{8}{c}{$H_{0,4}:\beta_{1}+\beta_{2}=0, \ v.s. \ H_{a,4}:  \beta_{1}+\beta_{2}\neq0$}\\\cline{2-10}
&0.0&0.056&0.048&0.040&0.042&0.056&0.050&0.052&0.042\\
&0.1&0.154&0.132&0.174&0.162&0.138&0.136&0.208&0.184\\
&0.2&0.420&0.386&0.574&0.532&0.374&0.354&0.574&0.530\\
&0.4&0.930&0.906&0.990&0.984&0.908&0.902&0.994&0.990  \\\cline{2-10}                                                         	   
&$\beta_{3} +\beta_{4}$&\multicolumn{8}{c}{$H_{0,5}:\beta_{3}+\beta_{4}=0,  \ v.s. \ H_{a,5}: \beta_{3}+\beta_{4}\neq 0.$}\\\cline{2-10}
&0.0&0.058&0.050&0.066&0.056&0.056&0.050&0.068&0.066  \\
&0.1&0.130&0.116&0.154&0.138&0.152&0.136&0.190&0.174\\
&0.2&0.450&0.428&0.470&0.450&0.412&0.388&0.450&0.424\\
&0.4&0.930&0.920&0.946&0.930&0.932&0.920&0.958&0.932 \\\cline{2-10}	
&$\beta_{1} + \beta_{p}$&\multicolumn{8}{c}{$H_{0,6}:\beta_{1}+\beta_{p}=0.75, \ v.s.  \ H_{a,6}: \beta_{1}+\beta_{p}\neq0.75$}\\\cline{2-10}	   
&0.75	&0.050&0.044&0.060&0.046&0.058&0.048&0.044&0.036  \\
&0.85	&0.172&0.146&0.138&0.116&0.162&0.154&0.154&0.130  \\
&0.95	&0.490&0.444&0.408&0.354&0.462&0.436&0.418&0.376   \\
&1.15	&0.970&0.950&0.936&0.916&0.970&0.958&0.942&0.914  \\\cline{2-10}	
&$\beta_{2} + \beta_{3}$&\multicolumn{8}{c}{$H_{0,7}: \beta_{2}+ \beta_{3}=-0.75,  \ v.s.  \ H_{a,7}: \beta_{2} + \beta_{3} \neq -0.75$}\\ \cline{2-10}	   
&-0.75&0.056&0.050&0.054&0.044&0.060&0.044&0.062&0.062  \\
&-0.65&0.200&0.176&0.182&0.172&0.148&0.138&0.180&0.174  \\
&-0.55&0.484&0.444&0.516&0.486&0.422&0.408&0.468&0.466  \\
&-0.35&0.922&0.910&0.966&0.962&0.920&0.910&0.960&0.962 \\\hline
\multirow{16}{*}{$p=350$}&	$\beta_{1} + \beta_{2}$&\multicolumn{8}{c}{$H_{0,4}:\beta_{1}+\beta_{2}=0, \ v.s. \ H_{a,4}:  \beta_{1}+\beta_{2}\neq0$}\\\cline{2-10}	
&0.0&0.062&0.056&0.062&0.056&0.050&0.046&0.048&0.046\\
&0.1&0.164&0.160&0.216&0.202&0.106&0.096&0.206&0.184\\
&0.2&0.472&0.438&0.612&0.572&0.402&0.378&0.536&0.510\\
&0.4&0.940&0.934&0.988&0.988&0.910&0.900&0.982&0.980\\\cline{2-10}
&$\beta_{3} +\beta_{4}$&\multicolumn{8}{c}{$H_{0,5}:\beta_{3}+\beta_{4}=0,  \ v.s. \ H_{a,5}: \beta_{3}+\beta_{4}\neq 0.$}\\ \cline{2-10}	   
&0.0&0.058&0.046&0.070&0.040&0.038&0.040&0.068&0.048\\
&0.1&0.192&0.188&0.174&0.138&0.126&0.124&0.172&0.136\\
&0.2&0.454&0.442&0.462&0.410&0.392&0.378&0.404&0.356\\
&0.4&0.952&0.952&0.912&0.814&0.944&0.942&0.916&0.828\\\cline{2-10}
&$\beta_{1} + \beta_{p}$&\multicolumn{8}{c}{$H_{0,6}:\beta_{1}+\beta_{p}=0.75, \ v.s.  \ H_{a,6}: \beta_{1}+\beta_{p}\neq0.75$}\\\cline{2-10}	   
&0.75 &0.046&0.044&0.058&0.038&0.056&0.042&0.060&0.054\\
&0.85 &0.148&0.142&0.280&0.292&0.110&0.110&0.240&0.242\\
&0.95 &0.466&0.472&0.562&0.566&0.394&0.390&0.496&0.512\\
&1.15 &0.942&0.944&0.960&0.940&0.938&0.940&0.954&0.932\\\cline{2-10}
&$\beta_{2} + \beta_{3}$&\multicolumn{8}{c}{$H_{0,7}: \beta_{2}+ \beta_{3}=-0.75,  \ v.s.  \ H_{a,7}: \beta_{2} + \beta_{3} \neq -0.75$}\\ \cline{2-10}	   
&-0.75 &0.052&0.038&0.052&0.046&0.052&0.038&0.062&0.038\\
&-0.65 &0.154&0.138&0.174&0.138&0.134&0.120&0.142&0.130\\
&-0.55 &0.450&0.420&0.478&0.442&0.400&0.374&0.488&0.456\\
&-0.35 &0.932&0.916&0.968&0.966&0.914&0.912&0.960&0.960\\\hline
                                           			
	\end{tabular}                                                                                                         			
      \end{table}

      \begin{table}
	\footnotesize
	\centering
	\caption{ The empirical size and power of Wald and score tests for linear hypothesis testing with $n=500$. }\label{S2_500}
	\begin{tabular}{c|c|cc|cc|cc|cc}\hline
		&&\multicolumn{4}{c|}{$\X \sim$ Normal}&\multicolumn{4}{c}{$\X \sim$ Uniform}\\\hline
		&&\multicolumn{2}{c|}{$\Sigma=0.5\I_p$}&\multicolumn{2}{c|}{$\Sigma=0.5^{|i-j|+1}$} &\multicolumn{2}{c|}{$\Sigma=0.5\I_p$}&\multicolumn{2}{c}{$\Sigma=0.5^{|i-j|+1}$} \\\hline
		&&$T_W$&$T_S$&$T_W$&$T_S$&$T_W$&$T_S$&$T_W$&$T_S$\\\hline
		\multirow{16}{*}{$p=50$}&	$\beta_{1} + \beta_{2}$&\multicolumn{8}{c}{$H_{0,4}:\beta_{1}+\beta_{2}=0, \ v.s. \ H_{a,4}:  \beta_{1}+\beta_{2}\neq0$}\\\cline{2-10}
		&0.0&0.058 & 0.052&0.048 & 0.046 &0.062&0.056&0.046&0.044\\
		&0.1&0.240 & 0.220&0.286 & 0.280&0.202&0.196&0.296&0.274\\
		&0.2& 0.670 & 0.656&0.810 & 0.794&0.586&0.580&0.758&0.754\\
		&0.4&0.996 & 0.996 &1.000 & 1.000&0.988&0.988&0.998&0.998  \\\cline{2-10}                                                         	   
		&$\beta_{3} +\beta_{4}$&\multicolumn{8}{c}{$H_{0,5}:\beta_{3}+\beta_{4}=0,  \ v.s. \ H_{a,5}: \beta_{3}+\beta_{4}\neq 0.$}\\\cline{2-10}
		&0.0&0.046 & 0.040&0.064 & 0.058&0.060&0.048&0.060&0.062  \\
		&0.1& 0.244 & 0.216&0.254 & 0.244&0.242&0.238& 0.260 & 0.244\\
		&0.2&0.646 & 0.616&0.696 & 0.692&0.666&0.648&0.686 & 0.678\\
		&0.4&0.996 & 0.996 &0.998 & 0.996&0.996&0.998 &0.998 & 0.998\\\cline{2-10}	
		&$\beta_{1} + \beta_{p}$&\multicolumn{8}{c}{$H_{0,6}:\beta_{1}+\beta_{p}=0.75, \ v.s.  \ H_{a,6}: \beta_{1}+\beta_{p}\neq0.75$}\\\cline{2-10}	   

		&0.75	&0.068 & 0.062 & 0.058 & 0.056&0.050&0.046&0.054&0.048 \\   
		&0.85	&0.242 & 0.222 & 0.190 & 0.162&0.226&0.212&0.176&0.160 \\   
		&0.95	& 0.676 & 0.640 & 0.620 & 0.580&0.646&0.626&0.646&0.614 \\   
		&1.15	&0.996 & 0.994 & 0.986 & 0.980 &0.996&0.994&0.998&0.990\\\cline{2-10}	
		&$\beta_{2} + \beta_{3}$&\multicolumn{8}{c}{$H_{0,7}: \beta_{2}+ \beta_{3}=-0.75,  \ v.s.  \ H_{a,7}: \beta_{2} + \beta_{3} \neq -0.75$}\\ \cline{2-10}	   
		&-0.75& 0.060 & 0.058 & 0.056 & 0.044&0.056&0.052&0.050&0.048   \\
		&-0.65&0.238 & 0.234 & 0.302 & 0.302 &0.214&0.206&0.226&0.218\\
		&-0.55& 0.640 & 0.622& 0.730 & 0.720 &0.608&0.598&0.664&0.660\\
		&-0.35& 0.992 & 0.992 &1.000 & 1.000&0.994&0.992&0.998&0.998\\\hline
		\multirow{16}{*}{$p=600$}&	$\beta_{1} + \beta_{2}$&\multicolumn{8}{c}{$H_{0,4}:\beta_{1}+\beta_{2}=0, \ v.s. \ H_{a,4}:  \beta_{1}+\beta_{2}\neq0$}\\\cline{2-10}	
		&0.0&0.054&0.044&0.056&0.050&0.046&0.042&0.046&0.042\\
		&0.1&0.192&0.180&0.286&0.268&0.190&0.182&0.292&0.288\\
		&0.2&0.602&0.594&0.790&0.786&0.578&0.558&0.824&0.816\\
		&0.4&0.688&0.700&0.998&1.000&0.986&0.984&1.000&1.000\\\cline{2-10}
		&$\beta_{3} +\beta_{4}$&\multicolumn{8}{c}{$H_{0,5}:\beta_{3}+\beta_{4}=0,  \ v.s. \ H_{a,5}: \beta_{3}+\beta_{4}\neq 0.$}\\ \cline{2-10}	   
		&0.0&0.072&0.064&0.066&0.060&0.044&0.046&0.060&0.056\\
		&0.1&0.222&0.218&0.244&0.238&0.206&0.194&0.264&0.264\\
		&0.2&0.654&0.650&0.678&0.670&0.614&0.608&0.712&0.714\\
		&0.4&0.974&0.974&0.994&0.980&0.992&0.990&0.994&0.978\\\cline{2-10}
		&$\beta_{1} + \beta_{p}$&\multicolumn{8}{c}{$H_{0,6}:\beta_{1}+\beta_{p}=0.75, \ v.s.  \ H_{a,6}: \beta_{1}+\beta_{p}\neq0.75$}\\\cline{2-10}	   
		&0.75 &0.044&0.040&0.056&0.048&0.058&0.056&0.052&0.050\\
		&0.85 &0.238&0.242&0.316&0.316&0.196&0.200&0.386&0.390\\
		&0.95 &0.670&0.672&0.716&0.722&0.604&0.616&0.790&0.788\\
		&1.15 &0.986&0.994&1.000&1.000&0.994&0.994&0.9996&0.994\\\cline{2-10}
		&$\beta_{2} + \beta_{3}$&\multicolumn{8}{c}{$H_{0,7}: \beta_{2}+ \beta_{3}=-0.75,  \ v.s.  \ H_{a,7}: \beta_{2} + \beta_{3} \neq -0.75$}\\ \cline{2-10}	   
		&-0.75 &0.058&0.052&0.052&0.044&0.056&0.052&0.044&0.046\\
		&-0.65&0.222&0.200&0.230&0.208&0.162&0.146&0.264&0.246\\
		&-0.55&0.630&0.608&0.692&0.680&0.544&0.516&0.686&0.678\\
		&-0.35&0.974&0.990&1.000&1.000&0.994&0.992&1.000&1.000\\\hline
		
	\end{tabular}                                                                                                         			
\end{table}

\subsubsection{Performance regarding $m$}
We further investigate how the testing performance changes as
$m$ changes. We consider three sets of hypotheses:
\bse
H_{0,8}:\sum_{j=1}^4\beta_{j}=0, &\ v.s.& \ H_{a,8}:  \sum_{j=1}^4\beta_{j} \neq0.\\
H_{0,9}:\sum_{j=1}^8\beta_{j}=0,  &\ v.s. & \ H_{a,9}: \sum_{j=1}^8\beta_{j} \neq0.\\
H_{0,10}:\sum_{j=1}^{12}\beta_{j}=0,  &\ v.s. & \ H_{a,10}: \sum_{j=1}^{12}\beta_{j} \neq 0,
\ese
corresponding to $m=4, 8$ and $12$.
We set $h_2=0$, $h_p=0$, and $h_3=0, 0.2, 0.4,
0.8$. The empirical sizes and powers are displayed in Table
\ref{S3}. These results suggest that under different $m$, the empirical sizes
remain close to the nominal significance level for both the Wald and
score tests. On the other hand, the empirical power decreases in general when
$m$ increases.
For instance, as shown in Table \ref{S3}, when $\X$ follows the
multivariate normal distribution  with mean zero and covariance
$\Sigma=0.5\I_p$, $p=350$ and $h_3=0.8$, the powers of the Wald test
are 1.000, 0.950 and 0.854 for $m=4, 8$ and 12, respectively. This is
intuitively sensible,
and suggests that larger sample size is needed  to reach a desired
power when the hypothesis concerns more parameters.

\begin{table}
  \footnotesize
  \centering
	\caption{ The empirical size and power of Wald and
	score tests under different $m$. }\label{S3}     
	\begin{tabular}{c|c|cc|cc|cc|cc}\hline
		&&\multicolumn{4}{c|}{$\X \sim$ Normal}&\multicolumn{4}{c}{$\X \sim$ Uniform}\\\hline
		&&\multicolumn{2}{c|}{$\Sigma=0.5\I_p$}&\multicolumn{2}{c|}{$\Sigma=0.5^{|i-j|+1}$} &\multicolumn{2}{c|}{$\Sigma=0.5\I_p$}&\multicolumn{2}{c}{$\Sigma=0.5^{|i-j|+1}$} \\\hline
		&&$T_W$&$T_S$&$T_W$&$T_S$&$T_W$&$T_S$&$T_W$&$T_S$\\\hline
	\multirow{12}{*}{$p=50$}&	$\sum_{j=1}^4\beta_{j}$&\multicolumn{8}{c}{$H_{0,8}:\sum_{j=1}^4\beta_{j}=0, \ v.s. \ H_{a,8}:  \sum_{j=1}^4\beta_{j}\neq0$}\\\cline{2-10} 
&0.0&0.062&0.052&0.046&0.038&0.062&0.060&0.060&0.052\\
&0.2&0.238&0.220&0.434&0.406&0.244&0.228&0.404&0.386\\       
&0.4&0.754&0.724&0.914&0.904&0.660&0.646&0.914&0.900\\       
&0.8&1.000&1.000&1.000&1.000&0.998&0.998&1.000&1.000\\\cline{2-10} 
	&$\sum_{j=1}^8\beta_{j}$&\multicolumn{8}{c}{$H_{0,9}:\sum_{j=1}^8\beta_{j}=0,  \ v.s.  \ H_{a,9}: \sum_{j=1}^8\beta_{j}\neq0$}\\\cline{2-10}
&0.0&0.054&0.052&0.062&0.058&0.052&0.054&0.062&0.052 \\
&0.2&0.162&0.142&0.316&0.308&0.178&0.172&0.296&0.282 \\       
&0.4&0.410&0.388&0.764&0.732&0.424&0.406&0.764&0.742 \\       
&0.8&0.966&0.952&1.000&1.000&0.920&0.906&1.000&1.000 \\\cline{2-10} 
&$\sum_{j=1}^{12}\beta_{j}$&\multicolumn{8}{c}{$H_{0,10}:\sum_{j=1}^{12}\beta_{j}=0,  \ v.s.  \ H_{a,10}: \sum_{j=1}^{12}\beta_{j} $}  \\\cline{2-10}
&0.0&0.046&0.044&0.062&0.062&0.046&0.052&0.052&0.050\\
&0.2&0.116&0.124&0.250&0.240&0.084&0.096&0.210&0.210\\       
&0.4&0.288&0.298&0.610&0.604&0.318&0.330&0.642&0.642\\       
&0.8&0.854&0.796&0.996&0.994&0.802&0.756&0.994&0.992\\\hline
\multirow{12}{*}{$p=350$}&	$\sum_{j=1}^4\beta_{j}$&\multicolumn{8}{c}{$H_{0,8}:\sum_{j=1}^4\beta_{j}=0, \ v.s. \ H_{a,8}:  \sum_{j=1}^4\beta_{j}\neq0$}\\\cline{2-10}
&0.0&0.062&0.052&0.062&0.052&0.060&0.054&0.036&0.036\\
&0.2&0.260&0.244&0.420&0.400&0.230&0.222&0.402&0.382\\
&0.4&0.718&0.684&0.918&0.916&0.710&0.672&0.924&0.912\\
&0.8&1.000&0.998&1.000&1.000&1.000&1.000&1.000&1.000\\\cline{2-10}
&$\sum_{j=1}^8\beta_{j}$&\multicolumn{8}{c}{$H_{0,9}:\sum_{j=1}^8\beta_{j}=0,  \ v.s.  \ H_{a,9}: \sum_{j=1}^8\beta_{j}\neq0$}\\\cline{2-10}
&0.0&0.062&0.060&0.068&0.060&0.066&0.058&0.042&0.040\\
&0.2&0.132&0.130&0.266&0.258&0.180&0.166&0.262&0.238\\
&0.4&0.452&0.420&0.770&0.746&0.380&0.360&0.764&0.752\\
&0.8&0.950&0.936&1.000&1.000&0.936&0.926&1.000&1.000\\\cline{2-10}
&$\sum_{j=1}^{12}\beta_{j}$&\multicolumn{8}{c}{$H_{0,10}:\sum_{j=1}^{12}\beta_{j}=0,  \ v.s.  \ H_{a,10}: \sum_{j=1}^{12}\beta_{j} $}  \\\cline{2-10}
&0.0&0.056&0.056&0.054&0.054&0.052&0.050&0.038&0.038\\
&0.2&0.098&0.106&0.170&0.174&0.100&0.100&0.178&0.172\\
&0.4&0.296&0.296&0.616&0.604&0.276&0.274&0.584&0.574\\
&0.8&0.854&0.792&0.998&0.996&0.838&0.818&0.992&0.990\\\hline
		\end{tabular}

\end{table}

\subsubsection{Comparison with naive test}

We further compare the performances of our proposed tests with the
naive Wald
and score
tests developed under the noise free framework. We
consider the covariates $\X_i=(X_{i,1},\ldots, X_{i,p})\trans$
generated from the  multivariate normal distribution with mean zero
and covariance matrix $0.7\I_p$. The noise $\U_i$ follows
the multivariate normal
distribution with mean zero and covariance matrix $0.3\I_p$. Other
settings remain unchanged. We consider the hypotheses  on a single
element in $\bb$: $H_{0,2}$, and the linear
combinations of two coefficient parameters: $H_{0,5}$ and $H_{0,7}$ as described previously.
We report the empirical sizes and powers of the Wald and score tests
with/without noises for $p=50$ in Table \ref{S4}. It is
clear that while the
proposed tests achieve  Type I errors reasonably close to the nominal
level under different null
hypotheses, the naive tests lead to
precarious performance.
For instance, the Type I errors of Wald and Score tests for
$H_{0,5}$ are as large as  0.474 and 0.554, respectively. These Type
I errors are far beyond the significance level.
Because they cannot control the significance level, we do not recommend
consider using them in practice.
\begin{table}
  \footnotesize
  \centering
		\caption{The empirical size and power of Wald and score tests with/without noise considered. }\label{S4}
		\begin{tabular}{c|cc|cc|cc|cc|cc|cc}\hline
			&\multicolumn{2}{c|}{With noise}&\multicolumn{2}{c|}{Without noise}&\multicolumn{2}{c|}{With noise}&\multicolumn{2}{c|}{Without noise} &\multicolumn{2}{c|}{With noise}&\multicolumn{2}{c}{Without noise}\\\cline{2-13}
			&$T_W$&$T_S$&$T_W$&$T_S$&$T_W$&$T_S$&$T_W$&$T_S$&$T_W$&$T_S$&$T_W$&$T_S$\\\hline
			$\beta_{3}$&\multicolumn{4}{c|}{$H_{0,2}: \beta_{3}=0$}&\multicolumn{4}{c|}{$H_{0,5}: \beta_{3}+\beta_{4}=0$}&\multicolumn{4}{c}{$H_{0,7}: \beta_{2}+\beta_{3}=-0.75$}\\\hline
			0.0&  0.084 & 0.078 &0.774 & 0.824 &0.064 & 0.066 &0.474 & 0.554&0.056 & 0.054 &0.538 & 0.604 \\
			0.1  & 0.340 & 0.316 & 0.332 & 0.410 &0.106 & 0.094& 0.784 & 0.838 & 0.166 & 0.188& 0.234 & 0.288  \\
			0.2  & 0.692 & 0.646 & 0.104 & 0.092 & 0.270 & 0.194& 0.930 & 0.956& 0.388 & 0.378& 0.096 & 0.090   \\
			0.4  & 0.914 & 0.954& 0.690 & 0.362 & 0.700 & 0.540& 0.996 & 1.000& 0.882 & 0.880& 0.406 & 0.196\\ \hline
		\end{tabular}

\end{table}

\subsection{Neuroimage application}\label{sec:real}
We apply our proposed testing procedures to study how the SUVRs from
PET image data affect the MoCa score. We download the preprocessed
$^{18}$F-AV-1451 PET 
image features, and demographic and cognitive  assessments 
from the ADNI database. The image features include $^{18}$F-AV-1451
SUVRs and volumes of the cortical, sub-cortical regions, brainstem,
ventricles and sub-divisions of corpus callosum. Furthermore,  the demographic
variables include gender  and standardized age (divided by the
standard deviation) at the  image 
examining time. For each subject, we obtain his/her MoCa
score within 14 days of his/her image examining time as the outcome, which
ranges from 9 to 30.  Furthermore, we remove the covariates with more than 100
missing values. We standardize the volumes of ROIs by
subtracting the means and dividing by the standard
deviations. We use the
SUVR from
inferior cerebellum as a reference and divide the rest of SUVRs by
this reference as suggested in \citep{Landau2016}.  Finally, we have
$n = 196$ complete samples with $p=218$ 
covariates in the analysis.

Since the neuroimage data are longitudinally collected, we estimate
the covariance matrix of $\U$ using repeatedly
measured image features, while assuming that age and gender are recorded
precisely. More specifically, let
$\wt{\W}_{ij}$ denote the observed image features at the $j$th
examining time. We 
first perform the regression between
$\wt{\W}_{ij} $ and age of the $i$th patient at the $j$th examining
time, and obtain $\wt{\U}_{ij}$ as the residual of the
regression. Then we obtain the estimator for the covariance matrix
\bse
\wt{\bOmega} = \frac{\sumi \sum_{j = 1}^{n_i}(\wt{\U}_{ij} - \Ubar_{i}) (\wt{\U}_{ij}
  - \Ubar_{i})\trans}{\sumi (n_i -1)},
\ese
where $n_i$ is the number of repeated measurements of
$\wt{\W}_{ij}$, and $\Ubar_{i} = \sum_{j = 1}^{n_i}\wt{\U}_{ij}
/n_i$. Finally, because the first two covariates, age and gender, are
measured precisely, the first two columns and rows of the estimated
$\bOmega$, denoted by $\wh{\bOmega}$, are zeros. We set the rest $(p-
2) \times 
(p-2)$ sub-matrix of $\wh{\bOmega}$ to be $\wt{\bOmega}$.

\begin{figure}[!h]
  \centering
\includegraphics[scale = 0.5]{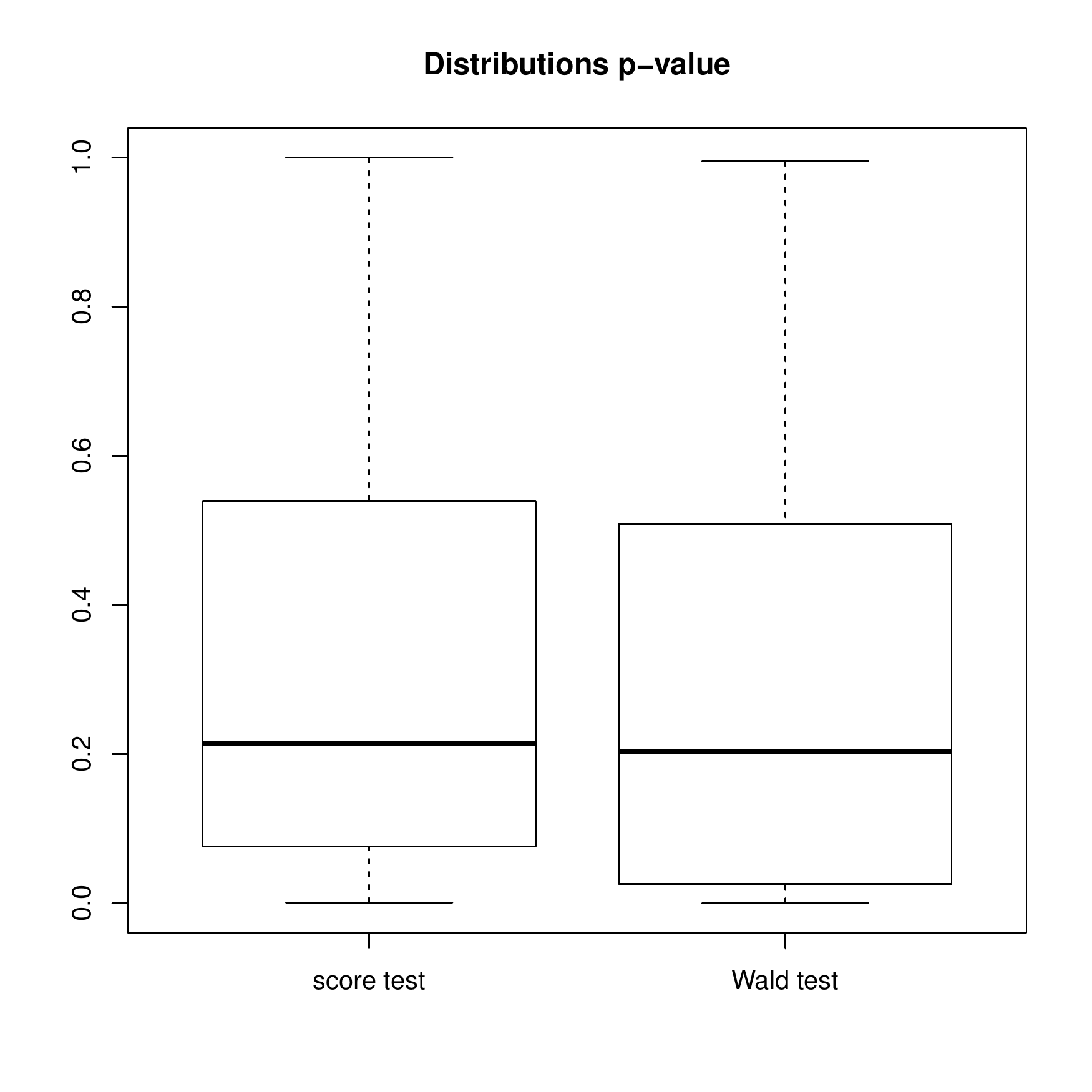}
	\caption{\label{fig:box} The boxplot of the $p$-values based on
          the score and Wald tests. The distributions of the
          $p$-values are similar from both methods. }
      \end{figure}

\begin{figure}[!h]
   \centering
   \includegraphics[scale = 0.2]{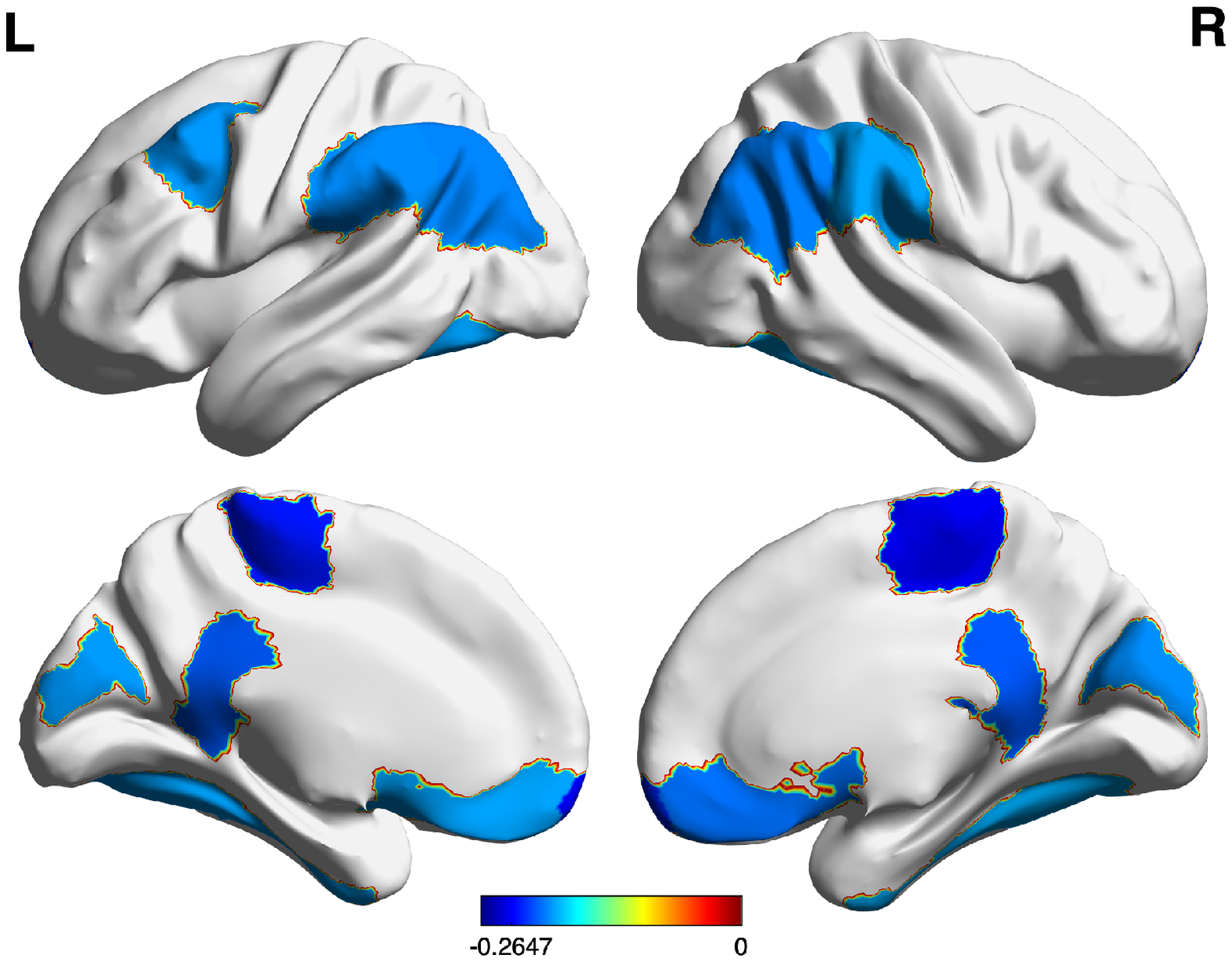}
   \includegraphics[scale = 0.2]{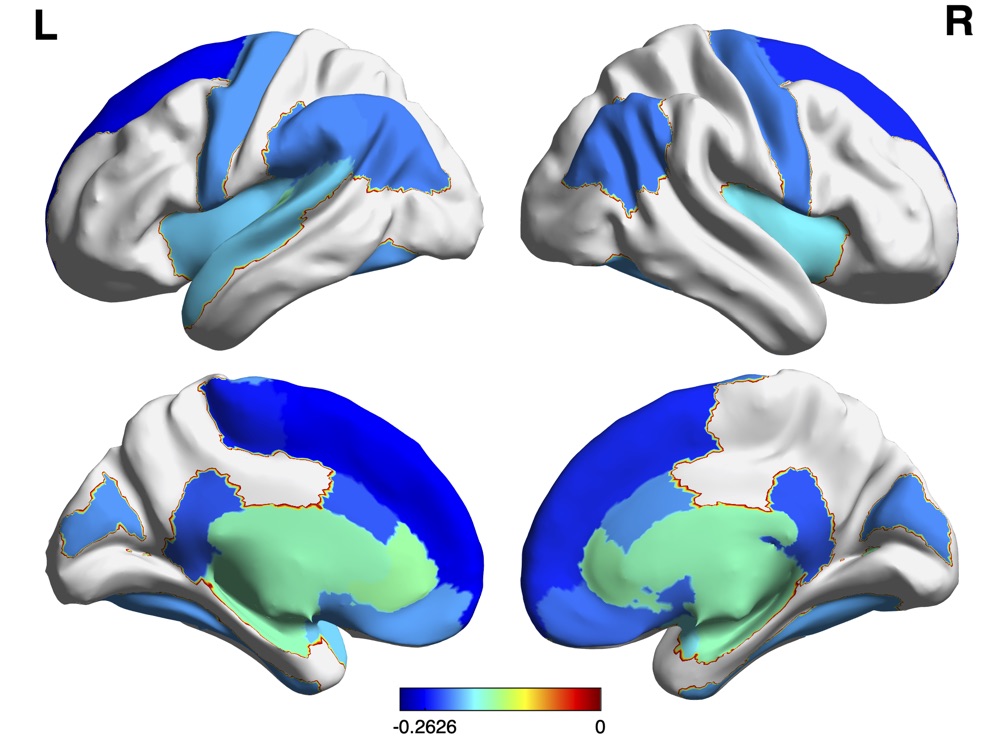}
	\caption{\label{fig:brain} The effects of SUVRs at the cortical
          regions. The colors represent the values (indicated by the
          color bars) of the estimated
          coefficients of the SUVRs. We only  plot the coefficient values
         corresponding to the significant brain regions with p-value less than
         0.05 from score test (left) and Wald test (right).  The white
         areas are the  non-significant brain regions. The $L$
         and $R$ letters in the plot represent the left and right
         hemispheres.  }
      \end{figure}

We test $p$ hypotheses, each of the form
\be\label{eq:test}
H_0: \beta_j=0 \quad \textrm{v.s.} \quad H_{a}:  \beta_j \neq 0, \quad
\ee
for $j=1, \dots, p$ at 0.05 nominal level.
To implement the hypothesis testing procedure, in each
  test, we first fit a standard 
penalized Poisson regression 
model to obtain the initial values of the 
coefficient estimators.  Then we construct the score test  and Wald test
statistics based on (\ref{eq:scoretest}) and (\ref{eq:waldtest}),
respectively.  The
tuning parameter $\lambda$ is selected by minimizing
(\ref{eq:BIC}). We obtain the p-value as the probability of a 
$\chi^2(1)$ random variable that is greater than the resulting score
and Wald test
statistics. 
There are  33 and 69 covariate coefficients with significant p-values
at 0.05 nominal level
based on the score and Wald tests, respectively. Furthermore, we plot
the boxplot of the resulting $p$-values in
Figure \ref{fig:box}. 
It is clear that the distribution of the $p$-values
are similar for the score and the Wald test. 
For each covariate $j$, we obtain the estimated $j$th
 coefficient  based on
 (\ref{eq:lapHT1}) under the corresponding alternative
 hypothesis, and plot the 
 estimated coefficients of the SUVRs at the 
 cortical regions on a template brain in Figure \ref{fig:brain}. The
 results show that  
 the SUVRs have negative effects on the cognitive score, suggesting
 that the higher the SUVR values, the
 lower the MoCa score and in turn the worse the cognitive function,
 which is consistent with the scientific evidences \citep{braak1991,
   scholl2016, baker2017}. 
Furthermore, the score test is more
 stringent and gives less number of  significant SUVRs. Among 33
significant predictors from the score test, 27 of them are also
significant in the  Wald test. Based on this 
 high agreement between the score and Wald tests,  we believe the difference
 between the two tests is a small sample phenomenon.

\begin{figure}
        \centering
        \includegraphics[scale = 0.5]{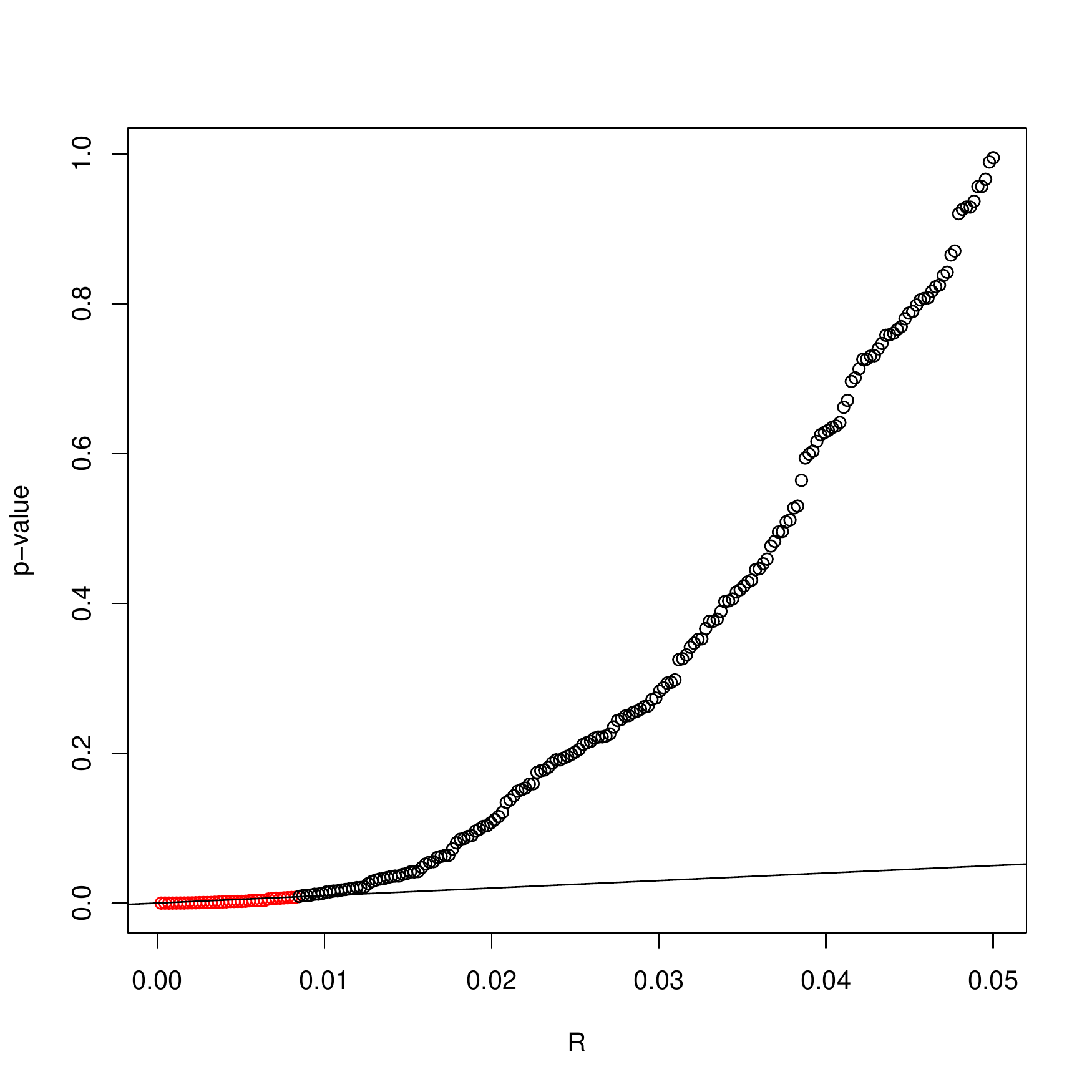}
        \caption{Sorted p-value versus $R = 0.05 j/218, j = 1, \ldots,
          218$. There are 36 important predictors corresponding to
        the p-values (in red) that below the line. }\label{fig:FDR}
        \end{figure}

To adjust for the multiple testing, we further performed
an analysis to control false discovery rate (FDR)
\citep{benjamini1995} within 0.05 by treating the p-values 
as independent. Since the score test is too stringent, no significant
covariate has been identified at 0.05 FDR by using the score
test. Therefore, we only 
present the results from the Wald test. We plot the 
p-values versus $0.05j/218$ in Figure \ref{fig:FDR} in an
increasing order, which
suggests 36 covariates are selected as the important predictors. There are 13
cortical SUVRs among the 36 important predictors that are
significant. We present their estimated coefficient, p-values 
from the Wald test in Table \ref{tab:suvrs}. The results show
that the majority of the significant cortical SUVRs are in the temporal
lobe, which consists of structures that are vital for declarative or
long-term memory \citep{smith2008}.

      \begin{table}
        \centering
        \footnotesize
        \begin{tabular}{l|c|c|c}
          \hline
          Cortical regions& Brain lobes  & Estimated coefficient & Wald
                                                           test p-value
                                                           \\
          \hline
               left middle temporal gyrus  & Temporal lobe &       -0.214     &
                                                                   0.0002\\
           left inferior parietal cortex   &     Parietal lobe      &      -0.214   &
                                                               0.0003\\
          left inferior temporal    gyrus  & Temporal lobe  &   -0.223&0.0007\\
          right inferior parietal cortex  &   Temporal lobe    &      -0.211&0.0007\\
          left  BANKSSTS& Temporal lobe      &    -0.174      &    0.0015\\
          left fusiform gyrus  & Temporal lobe&            -0.262
                                                                 &      0.0016\\
          right middle temporal gyrus  & Temporal lobe&
                                                        -0.229& 
          0.0024\\
          left caudal  middle frontal gyrus & Frontal lobe    &    -0.236
                                                                 &      0.0030\\
          left precuneus cortex & Parietal lobe&           -0.215     &     0.0034\\
          left entorhinal  cortex     & Temporal lobe&          -0.217
                                                                 &     0.0036\\
           right inferior temporal &      Temporal lobe&
                                                         -0.225& 
          0.0059\\
         right left entorhinal  cortex &     Temporal lobe&
                                                            -0.221
                                                                 &     0.0065\\
          right BANKSSTS & Temporal lobe&        -0.168    &      0.0076\\
                                                                                   \hline
        \end{tabular}
        \caption{The estimated coefficients, p-values from score and
          Wald tests for the significant SUVRs at 27 cortical
          regions. We also include the specific brain lobe that contains each cortical
          region. BANKSSTS  stands for banks of the superior temporal sulcus.}\label{tab:suvrs}
        \end{table}

\begin{figure}
  \centering
  \includegraphics[scale = 0.5]{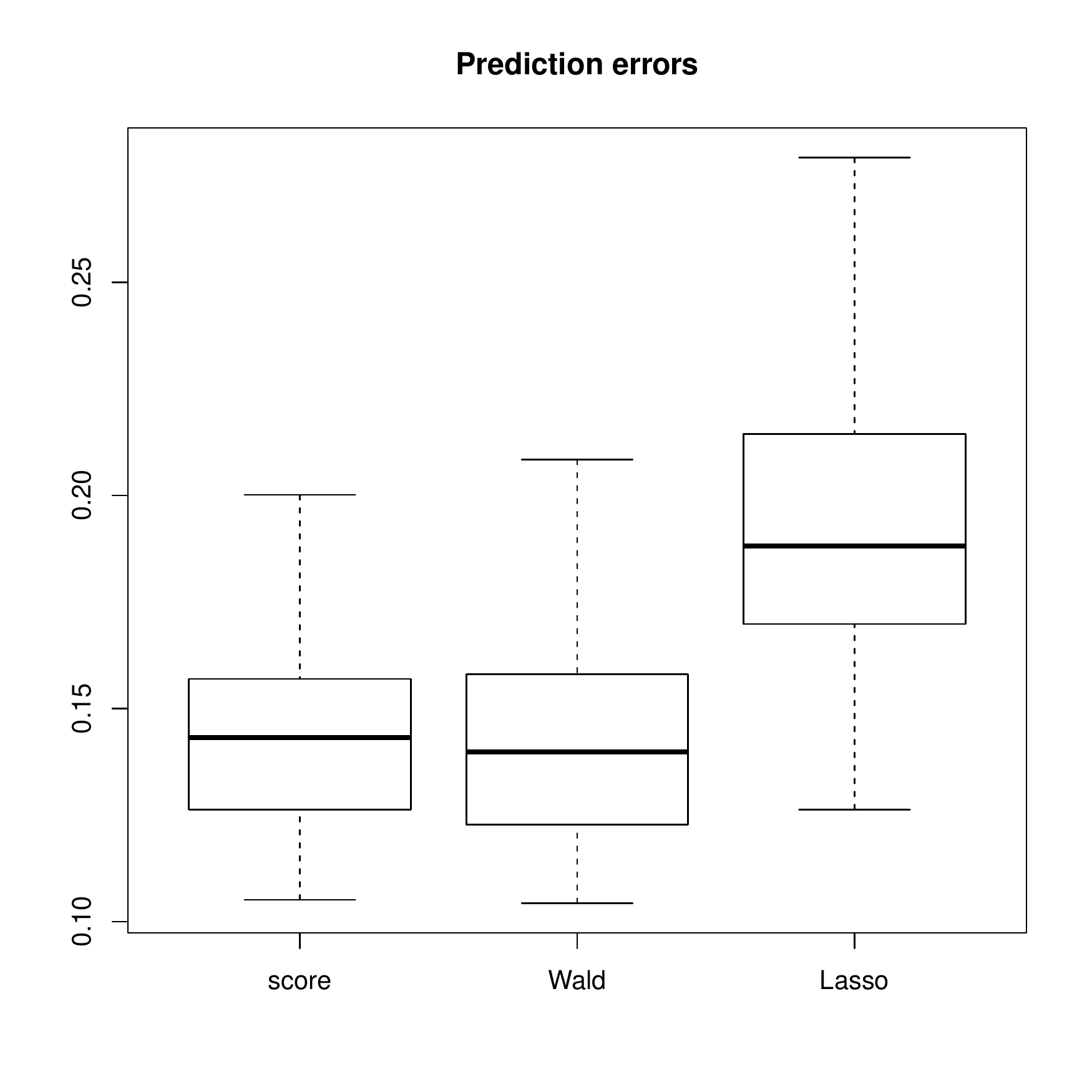}
  \caption{The distribution of the prediction errors from 100 runs of
the 
5-fold cross-validation based on Methods
    (a), (b) and (c). Method (d) breaks down because the algorithm does not
converge. }\label{fig:prederrors}
  \end{figure}

In addition, we perform a 5-fold cross
validation and compare the prediction errors among the four methods:
(a)  We select the important predictors 
as those with p-value less than 0.05 in the test (\ref{eq:test}) based
on the score statistics and then use  formula $\exp(\wh{\bb}_S\trans\W_{Si} - \wh{\bb}_S\trans
\wh{\bOmega}_S\wh{\bb}_S)$ to predict the outcome in the test sample,
where $\W_{Si}$ is  the selected 
covariates, $\wh{\bb}_S$ is estimator from (\ref{eq:lap}) using selected
covariates, $\wh{\bOmega}_S$ is the subset of $\wh{\bOmega}$
corresponding to the selected covariates.  
(b) We select the important predictors 
as those with p-value less than 0.05 in the test (\ref{eq:test}) based
on the Wald statistics and then use  formula $\exp(\wh{\bb}_W\trans\W_{Wi} - \wh{\bb}_W\trans
\wh{\bOmega}_W\wh{\bb}_W)$ to predict the outcome in the test sample,
where $\W_{Wi}$ is  the selected 
covariates, $\wh{\bb}_W$ is estimator from (\ref{eq:lap}) using selected
covariates, $\wh{\bOmega}_W$ is the subset of $\wh{\bOmega}$
corresponding to the selected covariates.   
(c) We select the
important predictors using  the 
standard lasso regression between the logarithm of the MoCa score and
all covariates and then use formula $\exp(\wh{\bb}\trans\W_i)$ to predict the
outcome, where $\wh{\bb}$ is the estimator from the lasso regression.    
(d) We select the
important predictors using 
the penalized Poisson regression between the logarithm of the MoCa score and
all covariates and then use formula $\exp(\wh{\bb}\trans\W_i)$ to predict the
outcome, where $\wh{\bb}$ is the estimator from the penalized Poisson
regression. The penalty parameters in the lasso and penalized Poisson
regression are selected using a sub-routine of 10-folder
cross-validation. Method (d) breaks down because the algorithm does not
converge for any selections of the penalty parameters. Therefore, in
Figure \ref{fig:prederrors}, we
show the distributions of the prediction 
errors, defined as $\sumi |Y_i - \wh{Y_i}|/|Y_i|$,  only
  for the methods (a), (b) and (c) after 100 runs of
the 
5-fold cross-validation. The results
shows that Method (a) and (b) have similar performance and both
outperform Method (c) with much smaller prediction errors on average. 

Finally,  we perform the score and the Wald tests to test whether any
SUVRs from any
composite regions may have significant association with the MoCa
score, where the composite regions, namely BRAAK12,  BRAAK34, BRAAK56,
are defined in \citep{braak1991}
and used in \cite{Landau2016} and \cite{scholl2016}. 
We provide
the list of ROIs in each composite regions in Appendix
\ref{sec:composite}. 
 Let $\bb_{\S_k}$ be the coefficients of the SUVRs
from the ROIs  that belong
to the composite region $k$. We test the null hypothesis that $\bb_{\S_k}
= \bf 0$. The results in Table \ref{tab:composit} show that all the
tests are significant, 
suggesting that at least one ROI  in each of the composite region has
significant association with the cognitive function. This result partially agrees
with results in \citep{scholl2016} that the SUVRs from the composite
regions are significantly different in healthy subjects and patients
with a diagnosis of probable Alzheimer's disease. 

\begin{table}[!h]
  \centering
  \small
  \begin{tabular}{c|c|c|c|c|c}

    Composite regions &  TS& score test p-value& TW &Wald test p-value &
                                                                 DF\\
        \hline
BRAAK12& 10.10 &  0.039 &  15.41 &   0.0039& 4\\
BRAAK34 &42.51 & 0.0113& 89.45 &  1.774e-09& 24\\
BRAAK56 &  66.29 &  0.0165 & 71.44 &  0.0055&44\\
  \end{tabular}
 \caption{The score and Wald test results of hypothesis that $\bb_{\S_k}
= \bf 0$. TS and TW are the score and Wald test statistics, DF is the degree of
freedom in the asymptotic distribution of the score and
Wald statistics,
which equals the number of the ROIs in the composite regions. }\label{tab:composit}
\end{table}

\section{Conclusion and discussion}\label{sec:con}

We have proposed an amenably penalized noise corrected  Poisson model
to study the relationship between the cognitive score and high
dimensional noisy neuroimage data. Under the sparsity assumption, we
established the parameter convergence rates in both $l_1$ and $l_2$
norms, the variable selection consistency property and the asymptotic
normality of a subvector with  possibly infinitely many
components. Inference
tools are subsequently developed.  The neuroimage application shows that the
inference tools generate scientifically meaningful results, which have
potential to be used to study the cognitive function and cognitive changes for 
neurodegenerative diseases. Further research along this line is ongoing
in our group. The
neuroimage dataset
and computational code are available at \cite{github}.

Thanks to an anonymous referee, we would like to point out one
important extension. Instead of a constant matrix $\Omega$, we can
further allow $\Omega$ to depend on both the covariate $\X$ and the
response $Y$, hence $\Omega(Y,\X)$, and assume $E(\U|Y, \X) = \0$. This would include heteroscedastic 
measurement error and to allow dependent relation between $\W$ and $Y$
given $\X$. All the estimation and inference results will still hold
and the regularity conditions and proofs in the Suppement also 
do not need to be further modified to accomodate this extension.

Establishing similar results in generalized linear models  beyond
Poisson or general regression models with non-Gaussian 
noise turns out to be surprisingly
difficult due to various technical obstacles. 
The main difficult lies in being unable to construct a
  loss function that is positive-definite at the true parameter
  value. In the case when an estimating equation is available,
  although one may be tempted to treat the $l_2$ norm square of the
  estimating equation as a  loss function, we find other technical
  issues arise partially because the Hessian of the loss function may
  involve the response, hence some of the techniques used here cannot
  be directly applied.
Likewise, extending the Poisson model to allow
  overdispersion also turns out challenging, regardless if we use a 
  negative binomial model, or incorporate random effects, or use extra
 observed covariates. All these will lead to models different from Poisson.
The biggest hurdle of
considering general regression model and/or non-Gaussian noise is to rigorously establish that the loss
function is locally convex. More investigation and 
dedicated effort are needed in this aspect.

The assumption that the covariance of the measurement is known is
widely adopted 
in the low and dimensional noisy data literature
\citep{stefanski1989, cook1995, loh2012, sorensen2015}, because the
parameter estimation in the
noisy model with unknown noise
covariance is a challenging, especially in high dimensional
setting where the covariance is a high dimensional unknown parameter
to be estimated. Thresholding techniques as those proposed in
\cite{bickel2008, cai2011, fan2011}
can be used for the covariance estimation, but the theoretical
properties of the resulting estimators are involved, requiring careful
treatment of the additional error from the covariance
estimator. In a relatively simple situation when the error
  variance can be estimated through estimating a parameter $\bg$ via
  solving $\f_\bg(\bg)=\0$, then writing $\mL(\bb)$ as $\mL(\bb, \bg)$, 
 we can accomodate the additional parameter by concatenating $\bb$
 with $\bg$ and carrying out the subsequent analysis. 
 For example, in this case the result in Theorem
 \ref{from:th:2} will be updated to
\bse
\left(
\begin{array}{c}
\wh{\bb}_{a \mM \cup S}- \bb_{t\mM \cup S}\\
\wh\bg-\bg\end{array}\right)&=&-
\left\{\begin{array}{cc}
\bQ_{\mM\cup S, \mM\cup S}(\bb_t,\bg)&\partial^2 {\cal
  L}(\bb,\bg)/\partial\bb_{\mM\cup S}\partial\bg\trans\\
\partial \f_\bg(\bg)/\partial\bb_{\mM\cup S}\trans&
\partial \f_\bg(\bg)/\partial\bg\trans
\end{array}\right\}^{-1}\\
&&\times\left[\begin{array}{c}
\left\{\frac{\partial
		\mL( \bb_t,\bg)}{\partial \bb}\right\}_{\mM
		\cup S}\\
\f_\bg(\bg)\end{array}\right]\{1 + o_p(1)\}.
\ese
Letting $\M\equiv {\bQ}_{\mM \cup S,
  \mM\cup S}(\bb_t)
-\{\partial^2 {\cal
  L}(\bb_t,\bg)/\partial\bb_{\mM\cup S}\partial\bg\trans\}
\{\partial \f_\bg(\bg)/\partial\bg\trans\}^{-1}
\{\partial \f_\bg(\bg)/\partial\bb_{\mM\cup S}\trans\}$, then this
 leads to
\bse
\wh{\bb}_{a \mM \cup S}- \bb_{t\mM \cup S}
&=&  -\M^{-1} \left[\left\{\frac{\partial
    \mL( \bb_t,\bg)}{\partial \bb}\right\}_{\mM
  \cup S}
+\{\partial^2 {\cal
  L}(\bb_t,\bg)/\partial\bb_{\mM\cup S}\partial\bg\trans\}\{\partial \f_\bg(\bg)/\partial\bg\trans\}^{-1}\right.\\
&&\left.\times
\f_\bg(\bg)\right]
\{1 + o_p(1)\}.
\ese
We further conduct simulations to evaluate the proposed
  adjustment in Section 
\ref{sec:addsimulation} of the supplementary document. The results
suggest that the proposed adjustment controls type I error rate when
$\bOmega$ contain small number of unknown parameters. Estimation and
testing when $\bOmega$  has a large number of unknown parameters
are challenging problems and deserve much more extensive
investigation.

\bibliographystyle{agsm}
\bibliography{hdmerrorinfer}

\input{hdinfersupprev3}
\end{document}

%% file: hdinfersupprev3.tex
\newcommand{\mystrut}{\vphantom{\int_0^1}}
\allowdisplaybreaks
\thispagestyle{empty}
\baselineskip 17pt
\renewcommand {\thepage}{}
\include{titre}
\pagenumbering{arabic}
\begin{center}
{\Large\bf Supplementary Materials to ``On High dimensional Poisson models with measurement error:
	hypothesis testing for nonlinear nonconvex optimization''}
\end{center}

\renewcommand{\theequation}{A.\arabic{equation}}
\renewcommand{\thesubsection}{A.\arabic{subsection}}
\renewcommand{\thesubsubsection}{A.\arabic{subsection}.\arabic{subsubsection}}
\renewcommand{\thesection}{A.\arabic{section}}
\setcounter{Lem}{0}\renewcommand{\theLem}{A.\arabic{Lem}}
\setcounter{Th}{0}\renewcommand{\theTh}{A.\arabic{Th}}
\setcounter{Cor}{0}\renewcommand{\theCor}{A.\arabic{Cor}}

\section*{Appendix A}
We define an auxiliary function $q_{\lambda}(t) = \lambda |t| - \rho_{\lambda}
(t)$ to facilitate the theoretical derivation, where $q_{\lambda}(t)  - \mu/2t^2$ is concave and everywhere
differentiable as shown in Lemma \ref{lem:fromlemm5support} in the
supplementary material.

\subsection{Conditions for the estimation consistency}

Define $\alpha_{\min}(\M)$ and  $\alpha_{\max}(\M)$ as the minimal and
maximal eigenvalues of the matrix $\M$, respectively. Further, we
define  the sub-exponential norm and sub-Gaussian norm as 
\bse
\|X\|_{\psi_1} = \sup_{k \geq 1} 1/{k} E(|X|^k)^{1/k}, 
\ese
and 
\bse
\|X\|_{\psi_2} = \sup_{k \geq 1} 1/ \sqrt{k} E(|X|^k)^{1/k}. 
\ese
For notational convenience, let  
\bse A(\bb\trans \W_i) := \exp(\bb\trans \W_i - \bb\trans
\bOmega\bb/2), 
\ese
\bse
g(\W_i, \bb, \v, \w) :=  \v \trans\{(\W_i -
\bOmega\bb) ^{\otimes2}- \bOmega\} \w,
\ese 
and 
\bse
g_1(\W_i, \bb, \v, \w) :=  \v \trans\{(\W_i -
\bOmega\bb) ^{\otimes2}\} \w. 
\ese 

\begin{Def}\label{def:def1}\cite{loh2012}
(Lower-RE condition). The matrix $\bGamma$ satisfies a lower
restricted eigenvalue condition with curvature $\alpha_1 >0$ and
tolerance $\tau(n, p) >0$ if 
\bse
\bb\trans \bGamma \bb \geq \alpha_1\|\bb\|_2^2 - \tau(n, p)
\|\bb\|_1^2, \forall \bb\in {\mathbb R}^p. 
\ese
\end{Def}

\begin{Def}\label{def:def2}\cite{loh2012}
(Upper-RE condition). The matrix $\bGamma$ satisfies a upper
restricted eigenvalue condition with smoothness $a_2 >0$ and
tolerance $\tau(n, p) >0$ if 
\bse
\bb\trans \bGamma \bb \leq a_2 \|\bb\|_2^2 +  \tau(n, p)
\|\bb\|_1^2, \forall \bb\in {\mathbb R}^p. 
\ese
\end{Def}
We first state the regularity conditions as follows:

\begin{enumerate}[label=(C\arabic*)]

\item \label{con:boundX}  (a) $\sup_{i =1, \ldots, n, \|\v\|_2\leq 1} |\W_{i}\trans \v| \leq 
  M_W \sqrt{\|\v\|_0}$ for a positive constant $M_W$. $\|\bOmega\|_2 =
  O(1)$. \\
(b)\bse
D_1 \leq \alpha_{\min} [E\{\exp(\bb\trans \X) \X\X\trans \}] \leq \alpha_{\max}
[E\{\exp(\bb\trans \X) \X\X\trans \}]\leq D_2, 
\ese
\bse
&&D_{W1} \leq \alpha_{\min} [E\{\exp(2 \bb\trans\W - \bb\trans \bOmega \bb) (\W -
 \bOmega\bb)^{\otimes2}\}] \\
&\leq& \alpha_{\max}
[E\{\exp(2 \bb\trans\W - \bb\trans \bOmega \bb) (\W -
\bOmega\bb)^{\otimes2}\}]\leq D_{W2},
\ese
\bse
\alpha_{\max}
[E\{\exp(\bb\trans\W - \bb\trans \bOmega \bb/2) (\W -
\bOmega\bb)^{\otimes2}\}]\leq D_{W3},
\ese
and $E\{\exp(2\bb\trans\X)\} = O(1)$ for any $\bb$ with $\|\bb\|_2 \leq 2 R_2$.\\
(c) $E(\|\W_i - \bOmega\bb\|_2^2) \leq D_{\Omega}$, 
for any $\bb$ with $\|\bb\|_2 \leq 2 R_2$. Here $D_1, D_2, D_{W1},
D_{W2}, D_{W3}, D_{\Omega}$ are positive
constants. \\
(d) $\|\C\|_2 = O(1)$ and $\|(\C\C\trans)^{-1}\|_2 = O(1)$. \\
(e) The $L_2$ norm of the true parameter $\bb$ is bounded, that is
$\|\bb\|_2 \leq b_0$ for some $0<b_0 <\infty$. \\

\item\label{con:Wx} 
 For $j = 1, \ldots, p$, define 
\bse
K_{j} :=\|U_{ij}\|_{\psi_2}
=(2\Omega_{jj})^{1/2}\sup_{k \geq 1} k^{-1/2} \frac{\Gamma ^{1/k}\{(k+1)/2\}}{\pi
  ^{1/(2 k)}},
\ese
where $\Gamma$ is the Gamma function, 
then there exist constants $m_0, M_0$ so that
\bse
m_{0} <K_{j}^2 \sumi  Y_i ^2 / n < M_{0}, 
\ese
uniformly for all $j$ almost surely.

\item\label{con:Yx} Define $$K_{Y}(\X_i)= \sup_{k \geq 1} k^{-1}
  E[| Y_i- \exp(\bb_t\trans
\X_i) |^k|\X_i]^{1/k}. $$ 
There exist constants $m_1, m_2, M_1,
M_2$ so that
\bse
&&m_{1} <\sumi X_{ij}^2 K_{Y}(\X_i) ^2/ n < M_{1}, \\
&&\max_i |X_{ij} |K_{Y}(\X_i) /\sqrt{\log{n}} <
M_{2}, 
\ese
uniformly for all $j = 1, \ldots, p$ almost surely.

\item \label{con:pn} Sample size $n$ and  dimension
  of covariates $p$ satisfy 
\bse
\log(n)\sqrt{\log (p) / n}\leq C
\ese 
for an absolute constant $C$.

\item \label{con:ww1}  For $\e_j$, $j = 1, \ldots, p$, define 
\bse
  K_{wij} (\bb)&=&\sup_{k \geq 1} k^{-1/2 } E[|(\W_i -
\bOmega\bb) \trans\e_j -E\{(\W_i -
\bOmega\bb) \trans\e_j|\bb\trans\W_i, \X_i\}
|^k|\bb\trans\W_i, \X_i]^{1/k}, 
\ese
we assume $E\{K_{wij} (\bb_t)^4\}< Q_0$. 
Then there exist constants $m_3, M_3, Q_1$ so that
\begin{itemize}
\item[(i)]
\bse
m_{3} < \sumi K_{wi j}  (\bb_t)^2 \exp(2\bb_t\trans\W_i - \bb_t\trans
\bOmega \bb_t) /n <M_{3}
\ese 
and
\item[(ii)]
\bse
|\sum_{i = 1}^n \frac{ E\{ \exp(\bb_t\trans\W_i - \bb_t\trans
\bOmega \bb_t/2) (\W_i - 
\bOmega\bb_t) \trans\e_j  |\bb_t\trans\W_i,\X_i \}  -  \exp(\bb_t\trans \X_i)
\X_i \trans\e_j } {  \sqrt{n \log (p)} }|< Q_1
\ese
\end{itemize}
uniformly for all $j=1, \dots, p$ in probability.

\item \label{con:ww2}  For vectors $\v, \w \in {\mathbb R}^p$, $\|\v\|_2\leq 1,
  \|\w\|_2\leq 1$, for $\bb$ with
 $\|\bb\|_2 \leq 2 R_2$, 
$$K_{gvwi}(\bb):=  \sup_{k > 1} 1/k E (|[g(\W_i, \bb, \v, \w) -
E\{g(\W_i, \bb, \v, \w)
|\bb\trans \W_i, \X_i\}]|^k| \bb\trans\W_i, \X_i )^{1/k}. $$
We assume $E\{K_{gvwi}(\bb) ^4\}< Q_{01}$, and $E[\exp \{A
  ^2 (\bb\trans\W_i) 
K_{gvwi}^2(\bb)\} ]<Q_{02}$.
We also assume that for all $\v, \w$, 
\be
&&m_4 < \sumi  | A (\bb\trans\W_i)|^2K_{gvwi}(\bb)^2/ n < M_4,\label{eq:M4} \\
&&m_5 < \max_i |A (\bb\trans\W_i) | K_{gvwi}(\bb)/\sqrt{\log{n}} < M_5,\label{eq:M5}
\ee 
and
\be\label{eq:Q2}
&&n^{-1/2} |\sup_{\v, \w}\sumi \v\trans (A(\bb\trans\W_i) E[ \{(\W_i-\bOmega\bb)^{\otimes2} - \bOmega\} |\bb\trans\W_i, \X_i]\nonumber\\
&&-
E\{\exp(\bb\trans
  \X_i) \X_i\X_i\trans\} ) \w |_2 <Q_2, 
\ee
in probability.

\item \label{con:ww3}   For vectors $\v, \w \in {\mathbb R}^p$, $\|\v\|_2\leq 1,
  \|\w\|_2\leq 1$, for $\bb$ with
 $\|\bb\|_2 \leq 2 R_2$, 
$$K_{g_1 vwi}(\bb):=  \sup_{k > 1} 1/k E (|[g_1(\W_i, \bb, \v, \w) -
E\{g_1(\W_i, \bb, \v, \w)
|\bb\trans \W_i, \X_i\}]|^k| \bb\trans\W_i, \X_i )^{1/k}. $$
We assume $E\{K_{g_1 vwi}(\bb) ^4\}< Q_{11}$, and $E[\exp \{A
  ^4 (\bb\trans\W_i) 
K_{g_1 vwi}^2(\bb)\} ]<Q_{12}$.
We also assume that for all $\v$, 
\bse
&&m_6 < \sumi  | A^2 (\bb\trans\W_i)|^2K_{g_1 vwi}(\bb)^2/ n < M_6,\label{eq:M6} \\
&&m_7 < \max_i |A^2 (\bb\trans\W_i) | K_{g_1vwi}(\bb)/\sqrt{\log{n}} < M_7,\label{eq:M7}
\ese 
and 
\bse
&&m_{61} < \sumi  | A^2 (\bb\trans\W_i)|^2K_{g_1 vwi}(\bb)^2/ n < M_{61},\label{eq:M61} \\
&&m_{71} < \max_i |A^2 (\bb\trans\W_i) |
K_{g_1vwi}(\bb)/\sqrt{\log{n}} < M_{71},\label{eq:M71}. 
\ese 
Further 
\bse\label{eq:Q3}
&&n^{-1/2} |\sup_{\v, \w} \sumi \v\trans (A^2(\bb\trans\W_i) E[ \{(\W_i-\bOmega\bb)^{\otimes2}\} |\bb\trans\W_i, \X_i]\nonumber\\
&&-
E\{\exp(2 \bb\trans \W_i - \bb\trans \bOmega\bb)
(\W_i-\bOmega\bb)^{\otimes2} \} )\w |  <Q_3, 
\ese
in probability and 
\bse\label{eq:Q31}
&&n^{-1/2} |\sup_{\v, \w} \sumi \v\trans (A(\bb\trans\W_i) E[ \{(\W_i-\bOmega\bb)^{\otimes2}\} |\bb\trans\W_i, \X_i]\nonumber\\
&&-
E\{\exp(\bb\trans \W_i - \bb\trans \bOmega\bb/2)
(\W_i-\bOmega\bb)^{\otimes2} \} )\w |  <Q_{31}, 
\ese
in probability. 
\end{enumerate}

\subsection{Proof of Theorems in Section \ref{sec:consistency}}
\subsubsection{Proof of Theorem \ref{from:th1}}
First denote $\wh{\v} = \wh{\bb} - \wc \bb$, by the Taylor expansion of the first order derivative
\bse
\{\partial \mL(\wh\bb)/\partial \bb\trans - \partial \mL(\wc \bb)/\partial \bb\trans\}\wh{\v}
=  \wh \v \trans \partial ^2 \mL(\bb^*)/\partial \bb\partial \bb\trans
\wh \v,
\ese
where $\bb^* $ is a point on the line connecting $\wc \bb$ and
$\wh{\bb}$ and hence is in the feasible set. Therefore, by Lemma
\ref{lem:RE1} we have
\be\label{eq:lower2}
\{\partial \mL(\wh \bb)/\partial \bb\trans - \partial
\mL(\wc \bb)/\partial \bb\trans\}\wh{\v} \geq \alpha_1\|\wh \v\|_2^2 -
\tau_1\sqrt{\log(p)/n} \|\wh \v\|_1^2.
\ee
We first show that $\|\wh \v\|_2 \leq 1$. If not, we have
\bse
\{\partial \mL(\wh \bb)/\partial \bb\trans - \partial
\mL(\wc \bb)/\partial \bb\trans\}\wh{\v} \geq \alpha_1\|\wh \v\|_2 -2
\tau_1\sqrt{\log(p)/n} R_1 \|\wh \v\|_1.
\ese
Together with (\ref{eq:firstorderHT}), we obtain
\be\label{eq:rholower}
\{-\partial \rho_{\lambda}(\wh \bb_{\mM^c})/\partial \bb_{\mM^c}\trans\A- \partial
\mL(\wc \bb)/\partial \bb\trans\}\wh {\v}\geq \alpha_1\|\wh \v\|_2 -
2\tau_1\sqrt{\log(p)/n} R_1 \|\wh \v\|_1.
\ee
Further,
\be\label{eq:rhoupp}
&&\{-\partial \rho_{\lambda}(\wh \bb_{\mM^c})/\partial \bb_{\mM^c}\trans\A- \partial
\mL(\wc \bb)/\partial \bb\trans\}\wh {\v} \nonumber\\
&\leq& \{\|-\partial \rho_{\lambda}(\wh \bb_{\mM^c})/\partial \bb_{\mM^c}\trans\|_\infty\|\A\|_\infty+  \|\partial
\mL(\wc \bb)/\partial \bb\trans\|_\infty \}\|\wh {\v}\|_1\nonumber\\
&\leq& (\lambda +  \|\partial
\mL(\wc \bb)/\partial \bb\trans\|_\infty ) \|\wh {\v}\|_1\nonumber\\
&\leq&  3\lambda /2 \|\wh {\v}\|_1.
\ee
The second inequality holds because the maximum row sum of $\A$ is 1 and $\|\partial \rho_{\lambda}(\wh \bb_{\mM^c})/\partial \bb_{\mM^c}\|_\infty\leq \lambda $
by Condition (A1)--(A6) and Lemma \ref{lem:fromlemma4}. The last equality holds because $\|\partial
\mL(\wc \bb)/\partial \bb\trans\|_\infty \leq \lambda /2$ by the
statement assumption.
Now combine (\ref{eq:rhoupp}) and (\ref{eq:rholower}), we have
\bse
\|\wh{\v}\|_2 &\leq& \alpha_1^{-1} (2 \tau_1\sqrt{\log(p)/n} R_1 + 3\lambda
/2) \|\wh{\v}\|_1\\
&\leq&\alpha_1^{-1} (2 \tau_1\sqrt{\log(p)/n} R_1 + 3\lambda/2) 2 R_1.
\ese
By the assumption that  $\lambda \leq \alpha_1/(6 R_1)$ and $n \geq
\log(p) (64 \tau_1^2 R_1^4 )/\alpha_1^2$ in the statement, we conclude that
the right hand side is
at most one, which contradict to the hypothesis that $\|\wh\v\|_2>1$. Therefore $\|\wh{\v}
\|_2\leq 1$.

Further, since function $\rho_{\lambda}(\bb_{\mM^c})
+ \mu/2 \|\bb_{\mM^c}\|_2^2$ is convex function of $\bb_{\mM^c}$ by Condition (A5), we
have
\bse
&&\rho_{\lambda}(\wc\bb_{\mM^c})
+ \mu/2 \|\wc \bb_{\mM^c}\|_2^2 - \rho_{\lambda}(\wh {\bb}_{\mM^c})
-  \mu/2 \|\wh {\bb}_{\mM^c}\|_2^2 \\
&\geq&
 \{\partial \rho_{\lambda}(\wh
{\bb}_{\mM^c})/\partial \bb_{\mM^c}\trans + \mu {\wh
  \bb}_{\mM^c}\trans \}( \wc \bb_{\mM^c} - \wh{\bb}_{\mM^c})
\ese
which implies
\bse
&&\rho_{\lambda}( \wc \bb_{ \mM^c})
- \rho_{\lambda}(\wh {\bb}_{\mM^c})+ \mu/2 \|\wh {\bb}_{\mM^c} -
\bb_{\mM^c}\|^2\\
& \geq& \{\partial \rho_{\lambda}(\wh
{\bb}_{\mM^c})/\partial \bb_{\mM^c}\trans \}( \wc \bb_{\mM^c} -
\wh{\bb}_{\mM^c}) \\
&=& \{\partial \rho_{\lambda}(\wh
{\bb}_{\mM^c})/\partial \bb_{\mM^c}\trans \}\A ( \wc \bb -
\wh{\bb}).
\ese
Combine with (\ref{eq:firstorderHT}), we have
\bse
\partial \rho_\lambda(\wh
\bb_{\mM^c})/\partial \bb_{\mM^c}\trans \A(\wc \bb - \wh {\bb})\geq -
\partial \mL(\wh \bb)/\partial \bb\trans(\wc \bb - \wh {\bb}),
\ese
and hence
\bse
\rho_{\lambda}( \wc \bb_{ \mM^c})
- \rho_{\lambda}(\wh {\bb}_{\mM^c})+ \mu/2 \|\wh {\bb}_{\mM^c} -
\wc \bb_{\mM^c}\|^2 \geq
\partial \mL(\wh \bb)/\partial \bb\trans(\wh {\bb} - \wc \bb ).
\ese
Now combine with (\ref{eq:lower2}), we have
\bse
&&\alpha_1 \|\wh {\v}\|_2^2 - \tau_1 \sqrt{\frac{\log(p)}{n}} \|\wh
\v\|^2_1\\
&\leq& - \partial \mL(\wc \bb)/\partial \bb\trans \wh{\v} +
\rho_{\lambda}( \wc \bb_{ \mM^c})
- \rho_{\lambda}(\wh {\bb}_{\mM^c}) + \mu/2
\|\A\wh{\v}\|_2^2\\
&\leq& - \partial \mL(\wc \bb)/\partial \bb\trans \wh{\v} +
\rho_{\lambda}(\wc  \bb_{ \mM^c})
- \rho_{\lambda}(\wh {\bb}_{\mM^c}) + \mu/2
\|\A\|_1\|\A\|_\infty\|\wh{\v}\|_2^2\\
&=& - \partial \mL(\wc \bb)/\partial \bb\trans \wh{\v} +
\rho_{\lambda}( \wc \bb_{ \mM^c})
- \rho_{\lambda}(\wh {\bb}_{\mM^c}) + \mu/2\|\wh{\v}\|_2^2.
\ese
This implies
\be\label{eq:a1}
&&(\alpha_1  - \mu/2) \|\wh{\v}\|_2^2 \nonumber\\
&\leq& \tau_1\sqrt{\frac{\log
    (p)}{n}}\|\wh \v\|_1^2 + \|\partial \mL(\wc \bb)/\partial
\bb\|_\infty \|\wh{\v}\|_1 + \rho_{\lambda}( \wc \bb_{ \mM^c})
- \rho_{\lambda}(\wh {\bb}_{\mM^c})\nonumber\\
&\leq& \left\{2 R_1  \tau_1\sqrt{\frac{\log
    (p)}{n}} + \|\partial \mL(\wc \bb)/\partial
\bb\|_\infty\right\} \{\|\wh{\v}_{\mM ^c}\|_1+ \|\wh{\v}_{\mM}\|_1\} +
\rho_{\lambda}( \wc \bb_{ \mM^c})
- \rho_{\lambda}(\wh {\bb}_{\mM^c}).
\ee
Note that by the assumption that $n \geq
\log (p) (16 R_1^4\tau_1^4) /\alpha_1^4$, we have
\bse
2 R_1 \tau_1 \left\{\frac{\log
    (p)}{n}\right\}^{1/2 } = 2 R_1 \tau_1 \left\{\frac{\log
    (p)}{n}\right\}^{1/4 }\left\{\frac{\log
    (p)}{n}\right\}^{1/4 }
\leq \alpha_1 \left\{\frac{\log
    (p)}{n}\right\}^{1/4 }.
\ese
Further by the assumption that $4 \|\partial \mL(\wc \bb)/\partial
\bb\|_\infty  \leq \lambda$ and  $4 \alpha_1 \{\log
(p)/n\}^{1/4} \leq \lambda$ in the lemma statement we obtain
\bse
 2R_1  \tau_1\sqrt{\frac{\log
    (p)}{n}} + \|\partial \mL(\wc \bb)/\partial
\bb\|_\infty\leq \lambda /4 +\lambda/4\leq \lambda /2.
\ese
Combine with (\ref{eq:a1}) and Lemma \ref{lem:fromlemma4}  and
subadditivity of $\rho_{\lambda}$ in Condition (A2), we have
\be\label{eq:v2lev1}
&& (\alpha_1  - \mu/2) \|\wh{\v}\|_2^2\nonumber\\
 &\leq&\rho_{\lambda}( \wc \bb_{ \mM^c})
- \rho_{\lambda}(\wh {\bb}_{\mM^c})+ \lambda /2
\left\{\frac{\rho_{\lambda}(\wh {\v}_{\mM ^c})}{\lambda } +
  \|\wh{\v}_{\mM}\|_1 + \frac{\mu}{2
    \lambda } \|\wh \v _{\mM ^c}\|_2^2 \right\}\nonumber\\
&\leq& \rho_{\lambda}( \wc \bb_{ \mM^c})
- \rho_{\lambda}(\wh {\bb}_{\mM^c})+ \frac{\rho_{\lambda}( \wc \bb_{ \mM^c})
+\rho_{\lambda}(\wh {\bb}_{\mM^c})}{2} +
\lambda \|\wh{\v}_{\mM}\|_1/2+ \mu/4 \|\wh {\v}\|_2^2.
\ee
Further,  let $\b^r$ be the vector containing the first $r$ element of vector
$\b$, and $\b^{-r}$ the vector without the first $r$
element.
By the condition that $\C\bb_{\mM} = \t$, we have
$\wh{\bb}_{\mM}^r = \C_r^{-1}(\t - \C_{m-r}  \wh{\bb}_{\mM}^{-r})$, and
$\wc{\bb}_{\mM}^r = \C_r^{-1}(\t - \C_{m-r}
\wc{\bb}_{\mM}^{-r})$.
We have
\be\label{eq:wtv1}
\|\wh{\v}_{\mM} \|_1 &=&  \|\C_r^{-1}\C_{m-r} \wh{\v}_{\mM}^{-r}\|_1 +
\|\wh{\v}_{\mM}^{-r}\|_1\nonumber\\
&\leq& \sqrt{r}\|\C_r^{-1}\C_{m-r} \wh{\v}_{\mM}^{-r}\|_2 + \sqrt{m
  -r}\|\wh{\v}_{\mM}^{-r}\|_2\nonumber\\
&\leq& (\sqrt{r}\|\C_r^{-1}\C_{m-r}\|_2 + \sqrt{m
  -r}) \|\wh {\v}\|_2.
\ee
Now (\ref{eq:v2lev1}) becomes
\be\label{eq:from24}
0\leq \left(\alpha_1 - \frac{3 \mu}{4}\right) \|\wh{\v}\|_2^2 \leq 3/2 \rho_{\lambda}(\wc \bb_{ \mM^c})
- 1/2 \rho_{\lambda}(\wh {\bb}_{\mM^c}) + \lambda \|\wh{\v}_{\mM}\|_1/2.
\ee
We consider two cases, $3/2 \rho_{\lambda}( \wc \bb_{ \mM^c})
- 1/2 \rho_{\lambda}(\wh {\bb}_{\mM^c}) > 0 $ and $3/2
\rho_{\lambda}( \wc \bb_{ \mM^c})
- 1/2 \rho_{\lambda}(\wh {\bb}_{\mM^c}) \le0$.
When $3 \rho_{\lambda}( \wc \bb_{ \mM^c})
- \rho_{\lambda}(\wh {\bb}_{\mM^c}) \geq 0$,
by Lemma \ref{lem:fromlemma5}, we have
\be\label{eq:from25}
0 \leq 3 \rho_{\lambda}( \wc \bb_{ \mM^c})
- \rho_{\lambda}(\wh {\bb}_{\mM^c})\leq 3
\lambda  \|\wh {\v}_{\mM^c\mathcal{A}}\|_1 - \lambda  \|\wh
{\v}_{\mM^c \mathcal{A}^c}\|_1.
\ee
Now from (\ref{eq:from25}) we further have
\bse
\|\wh {\v}_{\mM^c \mA^c}\|_1 \leq 3\|\wh {\v}_{\mM^c \mA}\|_1.
\ese
Substitue (\ref{eq:wtv1}) and (\ref{eq:from25}) into (\ref{eq:from24}), we have
\bse
(2 \alpha_1 - \frac{3\mu}{2})\|\wh {\v}\|_2^2 &\leq& 3 \lambda \|\wh
{\v}_{\mM^c \mathcal{A}}\|_1 - \lambda \|\wh {\v}_{\mM^c\mA^c}\|_1 +
\lambda \|\wh {\v}_{\mM}\|_1\\
& \leq&  3
\lambda \|\wh{\v}_{\mM^c\mA}\|_1 + \lambda   (\sqrt{r}\|\C_r^{-1}\C_{m-r}\|_2 + \sqrt{m
  -r}) \|\wh\v\|_2\\
&\leq&  3
\lambda  \sqrt{k}\|\wh{\v}_{\mM^c\mA}\|_2 + \lambda   (\sqrt{r}\|\C_r^{-1}\C_{m-r}\|_2 + \sqrt{m
  -r}) \|\wh\v\|_2\\
&\leq&   \{3
\lambda  \sqrt{k} + \lambda  (\sqrt{r}\|\C_r^{-1}\C_{m-r}\|_2 + \sqrt{m
  -r}) \}\|\wh{\v}\|_2.
\ese
Hence we have that  when  $3/2 \rho_{\lambda}( \wc \bb_{ \mM^c})
- 1/2 \rho_{\lambda}(\wh {\bb}_{\mM^c}) >0 $
\be\label{eq:wtv2}
\|\wh {\v}\|_2 \leq \frac{6 \lambda \sqrt{k} + 2  \lambda   (\sqrt{r}\|\C_r^{-1}\C_{m-r}\|_2 + \sqrt{m
  -r}) }{4 \alpha_1 - 3\mu}.
\ee
When $3/2 \rho_{\lambda}(\wc  \bb_{ \mM^c})
- 1/2 \rho_{\lambda}(\wh {\bb}_{\mM^c})  \le 0$, by (\ref{eq:from24}) and (\ref{eq:wtv1}) we have
\bse
\left(\alpha_1 - \frac{3 \mu}{4}\right) \|\wh{\v}\|_2^2 &\leq& \lambda
\|\wh{\v}_{\mM}\|_1/2\\
&\leq&\lambda  (\sqrt{r}\|\C_r^{-1}\C_{m-r}\|_2 + \sqrt{m
  -r}) \|\wh {\v}\|_2/2,
\ese
which implies that when $3/2 \rho_{\lambda}( \wc \bb_{ \mM^c})
- 1/2 \rho_{\lambda}(\wh {\bb}_{\mM^c})  <0$,
\bse
\|\wh{\v}\|_2 \leq \frac{2 \lambda  (\sqrt{r}\|\C_r^{-1}\C_{m-r}\|_2 + \sqrt{m
  -r}) }{4\alpha_1 - 3\mu}.
\ese
Together with (\ref{eq:wtv2}), we always have
\bse
\|\wh {\v}\|_2 \leq \frac{6 \lambda \sqrt{k} + 2  \lambda  (\sqrt{r}\|\C_r^{-1}\C_{m-r}\|_2 + \sqrt{m
  -r}) }{4 \alpha_1 - 3\mu}.
\ese
Further, the $L_1$ distance is
\bse
\|\wh \v\|_1 &\leq&  \|\wh {\v}_{\mM^c \mA}\|_1 + \|\wh
{\v}_{\mM^c \mA^c}\|_1  + \|\wh {\v}_{\mM}\|_1 \\
&\leq& 4
\|\wh {\v}_{\mM^c \mA}\|_1 + \|\wh {\v}_{\mM}\|_1\\
&\leq&  4 \sqrt{k}\|\wh{\v}_{\mM^c}\|_2 +
(\sqrt{r}\|\C_r^{-1}\C_{m-r}\|_2 + \sqrt{m
  -r}) \|\wh {\v}\|_2\\
& \leq&  4 \sqrt{k}\|\wh{\v}\|_2 +  (\sqrt{r}\|\C_r^{-1}\C_{m-r}\|_2 + \sqrt{m
  -r}) \|\wh {\v}\|_2\\
&\leq&  (4 \sqrt{k} +  \sqrt{r}\|\C_r^{-1}\C_{m-r}\|_2 + \sqrt{m
  -r})\frac{6 \lambda \sqrt{k} + 2  \lambda  (\sqrt{r}\|\C_r^{-1}\C_{m-r}\|_2 + \sqrt{m
  -r}) }{4 \alpha_1 - 3\mu}.
\ese
This proves the result. \qed

\subsubsection{Proof of Theorem \ref{from:lem:th1}}
This proof is very similar to the proof of Theorem \ref{from:th1} and
is simpler, hence is omitted.

\subsection{Proof of Theorems in Section \ref{sec:normality}}

\subsubsection{Supporting Theorem of Theorem \ref{th:2} and its Proof}
\begin{Th}\label{th:1} Let $\alpha_1 =\min_{\|\bb\|_1\le R_1, \|\bb\|_2\le R_2}
\alpha_{\min}
  [E\{\exp(\bb\trans\X_i)\X_i\X_i\trans\}]  /2
$.
Assume $\rho_{\lambda}$ is $\mu$-amenable, for $\mu <
\alpha_1$, and Conditions \ref{con:boundX} -- \ref{con:ww2}
hold. Define $\A = \{\0_{(p-m)\times m}, \I_{(p-m)\times (p-m)}\}$ and
$\A_1 = \{\I_{m\times m}, \0_{m\times (p-m)}\}$.
Write $t_1\equiv\sqrt{r}\|\C_r^{-1}\C_{m-r}\|_2 + \sqrt{m
  -r}$ and
$
t\equiv (6 \lambda \sqrt{k} + 2  \lambda t_1 )(4 \alpha_1 - 3\mu)^{-1}.
$
Further we state
the following two conditions.

(a) The parameters $\lambda, R_1, R_2$ satisfy
\bse
4 \max\left\{\|\partial \mL (\wc \bb)/\partial \bb\|_\infty,
  \alpha_1 (\log(p)/n)^{1/4}\right\}\leq \lambda \leq \left[ \alpha_1
  \left\{\mystrut 16 (4 \sqrt{k} +  t_1)  t/\lambda\right\}^{-1}\right]^{1/2},
\ese
\bse
\max\left\{\mystrut 2 \|\wc \bb\|_1, 2 (4 \sqrt{k} +
t_1) t\right\} \leq R_1
\leq
\min\left(\frac{\alpha_1}{8\lambda}, [n \alpha_1^2/\{64 \tau_1^2 \log
(p)\}]^{1/4}, [n \alpha_1^4/\{16 \tau_1^4 \log
(p)\}]^{1/4}\right).
\ese
with
$
 \|\wc \bb\|_1\neq  (4 \sqrt{k} + t_1) t
$
and
$
\max\left\{2 \|\wc \bb\|_2, 2 t\right\} \leq R_2,
$
with
$
 \|\wc \bb\|_2\neq t.
$

(b)
Let $\wh{\bb}_{\mM\cup S}$ be the minimizer for
(\ref{eq:lapHTsub}),  $\wh{\bb} = (\wh{\bb}_{\mM\cup S}\trans, \0_{p -
m - k}\trans)\trans$, $\wh{\z}$ and $\bmu_4$ satisfy
\be
\frac{\partial \mL(\wh{\bb} )}{\partial \bb} - \A\trans \frac{\partial
  q_{\lambda}(\wh{\bb}_{\mM^c})}{\partial \wh{\bb}_{\mM^c}} +
\lambda (\A\trans \wh{\z})  + \A_1 \trans \C\trans \bmu_4=\bf 0. \label{eq:construct}
\ee
There exist $\delta \in [4 R_1 \tau_1 \sqrt{\log (p)/n} /\lambda,
1]$ so that
$\|(\A\trans \wh{\z}) _{(\mM \cup S)^c} \|_\infty\leq 1-\delta$, where
we name $(\A\trans \wh{\z})$ as extended sub-gradient. In addition,
\bse
n \geq \max[\log (p) \tau_1^2 (m + k)^2/(\alpha_1-\mu)^2,
 4 \log (p) \tau_1^2 \{t_1+  (2/ \delta + 1/2) \sqrt{k} \}^4/(\alpha_1-\mu)^2], ~~\alpha_1 > \mu
\ese  and
$\bb_{\mM^c}$ is $k$ sparse.

Under the above conditions and the conditions (a) and (b),
(\ref{eq:lapHT}) has a unique local
minimizer $\wh{\bb}$.
\end{Th}

\noindent Proof: We follow the primal-dual witness
construction introduced in \cite{wainwright2009}.

\noindent {\bf Step i:}
Following Lemma \ref{lem:fromlemma1} in the supplementary material,
let
$\wh{\bb}_{\mM \cup S}$ be the minimizer
for (\ref{eq:lapHTsub}).
We can easily check that when replacing $\bb$ by $\bb_{\mM\cup S}$, $\X$ by
$\X_{\mM\cup S}$, and $p$ by $m+k$,  all the conditions of Theorem \ref{from:th1} are still
satisfied. Here to check the condition regarding $\alpha_1$,  we note that for any $\bb$,
\bse
&&\min_{\|\bb\|_2\le R_2, \|\bb\|_1\le R_1}
\alpha_{\min} [E\{\exp(\bb\trans\X_i)\X_i\X_i\trans\}]  \\
&=&\min_{\|\bb\|_2\le R_2, \|\bb\|_1\le R_1} \inf_{\|\v\|_2=1,\v\in\mR^p}\v\trans E\{\exp(\bb\trans\X_i)\X_i\X_i\trans\} \v\\
&\le&\min_{\|\bb\|_2\le R_2, \|\bb\|_1\le R_1} \inf_{\|\v\|_2=1,\v\in\mR^{m+k}}
(\v\trans,\0\trans) E\{\exp(\bb\trans\X_i)\X_i\X_i\trans\} (\v\trans,
\0\trans)\trans\\
&=&\min_{\|\bb\|_2\le R_2, \|\bb\|_1\le R_1} \alpha_{\min}
 [E\{\exp(\bb\trans\X_i)\X_{i, \mM\cup S}\X_{i, \mM\cup S}\trans\}] \\
&=&\min_{\|\bb_{\mM\cup S}\|_2\le R_2, \|\bb_{\mM\cup S}\|_1\le R_1} \alpha_{\min}
 [E\{\exp(\bb_{\mM\cup S}\trans\X_{i, \mM\cup S})\X_{i, \mM\cup S}\X_{i, \mM\cup S}\trans\}].
\ese
Therefore, Theorem \ref{from:th1} applied to the $m+k$ dimensional case
leads to
\bse
\|\wh {\bb}_{\mM\cup S} - \wc \bb_{ \mM\cup S}\|_1 \le (4\sqrt{k}+t_1)t
\leq R_1/2,
\ese
and
\bse
\|\wh {\bb}_{\mM\cup S} - \wc \bb_{{\mM\cup S}}\|_2 \le t
\leq  R_2/2.
\ese
Therefore
\bse
\|\wh \bb_{\mM\cup S} \|_1  \leq \|\wh {\bb}_{\mM\cup S} -
\wc \bb_{{\mM\cup S}}\|_1  + \|\wc \bb\|_1 < R_1
\ese
and
\bse
\|\wh \bb_{\mM\cup S} \|_2  \leq \|\wh {\bb}_{\mM\cup S} -
\wc \bb_{{\mM\cup S}}\|_2  + \|\wc \bb\|_2
< R_2.
\ese
Hence $\wh \bb_{\mM\cup S}$ and $\wh \bb$ must be in the interior of the feasible
region.

\noindent {\bf  Step ii:} We show that $\wh{\bb}$ is a local  minimum
for  (\ref{eq:lapHT}) by verifying the conditions in Lemma
\ref{lem:fromlemma10} in the supplementary material. Because
\bse
\mL (\bb) + \rho_{\lambda}(\bb_{\mM^c})
&=& \mL(\bb) - q_{\lambda}(\bb_{\mM^c}) + \lambda \|\bb_{\mM^c}\|_1,
\ese
we can write $f = \mL$, $g = q_{\lambda}$, and $(\x^*, \v^*,
\w_1^*, \w_2^*, \mu_1^*, \mu_2^*, \bmu_3^*) = (\wh{\bb},\wh{\z}, \A\trans
\wh{\z}, \wh \z_1, 0, 0, \bmu_4 )$, where $\z_1 \in \partial
\|\wh{\bb}\|_2$. Lemma \ref{lem:fromlemm5support} in the supplementary
material ensures the
concavity and differentiability of $g(\x) - \mu/2
\|\x\|_2^2$. Further, since $\mu_1^* = \mu_2^* = 0$, (\ref{eq:x1}) and
(\ref{eq:x2}) are satisfied.   (\ref{eq:x1x2}) is satisfied by our
construction in (\ref{eq:construct}). Therefore, it remains to verify
(\ref{eq:G2c}). We first show that $G^* \subseteq {\mathbb R}^{\mM \cup
  S}$ so that (\ref{eq:G2c}) only needs to be satisfies for the
vectors belong to  ${\mathbb R}^{\mM \cup
  S}$. Suppose this does not  hold, let $\bnu \in G^* $ such that
$\supp (\bnu) \not\subseteq \mM\cup S$. This implies there is an index $j
\in {(\mM \cup
  S)^c}$ such that $\nu_j \neq 0$.

Now we define $\A \trans {\z}'$ such that $(\A\trans {\z}')_k = (\A\trans\wh{\z})_k $ for $k\neq
j$, and $(\A\trans {\z}')_j = {\rm sign}(\nu_j)$, where $a_k$ is the
$k$th element in vector $\a$. Clearly ${\z}'  \in \partial  \|\wh\bb_{\mM^c}\|_1$ and
\bse
\lambda  \bnu \trans\A\trans {\z}' > \lambda  \bnu \trans\A\trans \wh{\z}
\ese
because $\|(\A\trans\wh{\z})_{(\mM \cup S)^c}\|_\infty <1$.
Therefore,
\bse
&& \bnu \trans \left[\frac{\partial \mL(\wh{\bb} )}{\partial \bb} -  \left\{\A\trans \frac{\partial
  q_{\lambda}(\wh{\bb}_{\mM^c})}{\partial \wh{\bb}_{\mM^c}}\right\}
\right]  + \lambda  \bnu \trans\A\trans {\z}' + \bnu\trans  \A_1\trans\C\trans \bmu_4\\
&> & \bnu \trans \left[\frac{\partial \mL(\wh{\bb} )}{\partial \bb} -  \left\{\A\trans \frac{\partial
  q_{\lambda}(\wh{\bb}_{\mM^c})}{\partial \wh{\bb}_{\mM^c}}\right\}
\right]  + \lambda  \bnu \trans\A\trans\wh{\z}  + \bnu  \A_1\trans \C\trans \bmu_4= 0.
\ese
Now because $\bnu \in G^*$, we have
\bse
\bnu  \A_1\trans \C\trans \bmu_4 = 0.
\ese
This implies
\bse
\bnu \trans \left[\frac{\partial \mL(\wh{\bb} )}{\partial \bb} -  \left\{\A\trans \frac{\partial
  q_{\lambda}(\wh{\bb}_{\mM^c})}{\partial \wh{\bb}_{\mM^c}}\right\}
\right]  + \lambda  \bnu \trans\A\trans {\z}' >0
\ese
which contradicts with the requirement  of $G^*$ that
\bse
\sup_{\v \in \partial
  \|\wh\bb_{\mM^c}\|_1}  \bnu \trans \left[\frac{\partial \mL(\wh{\bb} )}{\partial \bb} -  \left\{\A\trans \frac{\partial
  q_{\lambda}(\wh{\bb}_{\mM^c})}{\partial \wh{\bb}_{\mM^c}}\right\}
\right]  + \lambda  \bnu \trans\A\trans\v = {0 }.
\ese
Therefore, $G^* \subseteq \mathbb{R}^{\mM \cup S}$. Now by
construction $\supp(\wh\bb) \subset \mM \cup S$,  using Lemma
\ref{lem:fromlemma1}  in the supplementary
material, we conclude that (\ref{eq:G2c}) holds with $\kappa =
\mu$. Hence, all conditions of Lemma \ref{lem:fromlemma10}  in the supplementary
material  are
satisfied, so we conclude
$\wh{\bb}$ is an isolated  local minimum of
(\ref{eq:lapHT}).

Now by Lemma \ref{lem:fromlem3}, because $\supp (\wh{\bb}) \subseteq
\mM \cup S$ and  $\wh{\bb}$ is an interior minimizer, the support of any stationary point of
(\ref{eq:lapHT}) is a subset of $\mM \cup
S$. Hence, we can write any stationary point in the form of  $\wt{\bb}
= \{\wt{\bb}_{\mM \cup S}\trans,
\0_{p - m -k}\trans\}\trans$, where $\wt{\bb}_{\mM \cup S}$ is a
stationary point for  (\ref{eq:lapHTsub}).
Further, note  that
(\ref{eq:lapHTsub}) is strictly convex by Lemma
\ref{lem:fromlemma1}, and hence $\wt{\bb}_{\mM \cup S}$ is unique in
the feasible set and
therefore $\wh{\bb}_{\mM\cup S}$ and
$\wt{\bb}$ are unique. Hence $\wh{\bb} = \wt{\bb}$ is the unique
local minimum. This proves the result. \qed

\subsubsection{Supporting Theorem of Theorem \ref{from:th:2} and its Proof}
\begin{Th}\label{from:th:1}
Let $\alpha_1 =\min_{\|\bb\|_1\le R_1, \|\bb\|_2\le R_2}
\alpha_{\min}
  [E\{\exp(\bb\trans\X_i)\X_i\X_i\trans\}]  /2
$.
Assume $\rho_{\lambda}$ is $\mu$-amenable, for $\mu <
\alpha_1$, and Conditions \ref{con:boundX} -- \ref{con:ww2} in the
supplementary material
hold. Define $\A = \{\0_{(p-m)\times m}, \I_{(p-m)\times (p-m)}\}$ and
$\A_1 = \{\I_{m\times m}, \0_{m\times (p-m)}\}$.    Further assume
that

(a) The parameters $\lambda, R_1, R_2$ satisfy
\bse
4 \max\left\{\|\partial \mL (\bb_t)/\partial \bb\|_\infty,
\alpha_1 (\log(p)/n)^{1/4}\right\}\leq \lambda \leq \left[ \alpha_1
  \left\{\mystrut 16 (4 \sqrt{k} + \sqrt{m})
\frac{6 \sqrt{k} + \sqrt{m}) }{4 \alpha_1- 3\mu} \right\}^{-1}\right]^{1/2},
\ese
\bse
&&\max\left\{\mystrut 2 \|\bb_t\|_1, 2 (4 \sqrt{k} + \sqrt{m})\lambda \frac{6 \sqrt{k} + 2\sqrt{m} }{4 \alpha_1 - 3\mu} \right\} \leq R_1 \\
&& \leq
\min\left(\frac{\alpha_1}{8\lambda}, [n \alpha_1^2/\{64 \tau_1^2 \log
(p)\}]^{1/4}, [n \alpha_1^4/\{16 \tau_1^4 \log
(p)\}]^{1/4}\right).
\ese
with
\bse
 \|\bb_t\|_1\neq  (4 \sqrt{k} + \sqrt{m})
\frac{6 \lambda \sqrt{k} + 2  \lambda
   \sqrt{m} }{4 \alpha_1 - 3\mu}
\ese
and
\bse
\max\left\{2 \| \bb_t\|_2, 2 \frac{6 \lambda \sqrt{k} + 2  \lambda  \sqrt{m} }{4 \alpha_1- 3\mu} \right\} \leq R_2,
\ese
with
\bse
 \| \bb_t\|_2\neq \frac{6 \lambda \sqrt{k} + 2\lambda\sqrt{m} }{4 \alpha_1 - 3\mu}.
\ese
(b) Let $\wh{\bb}_{a \mM\cup S}$ be the minimizer for
(\ref{eq:lapHTsub1}),  $\wh{\bb}_a = (\wh{\bb}_{a \mM\cup S}\trans, \0_{p -
m - k}\trans)\trans$, $\wh{\z}$ satisfy
\be \label{eq:construct1}\left\{\frac{\partial \mL(\wh{\bb}_a )}{\partial \bb_a}\right\} - \left\{\A\trans \frac{\partial
  q_{\lambda}(\wh{\bb}_{a \mM^c})}{\partial \wh{\bb}_{a \mM^c}}\right\} +
\lambda (\A\trans \wh{\z})  = \bf 0.
\ee
There exist $\delta \in [4 R_1 \tau_1 \sqrt{\log (p)/n} /\lambda,
1]$ so that
$\|(\A\trans \wh{\z}) _{(\mM \cup S)^c} \|_\infty\leq 1-\delta$, where
we name $(\A\trans \wh{\z})$ as extended sub-gradient.
In addition,
\bse n &\geq& \max[\log (p) \tau_1^2 (m + k)^2/(\alpha_1-\mu)^2, \\
&& 4 \log (p) \tau_1^2 \{ \sqrt{m} +  (2/ \delta + 1/2) \sqrt{k}
\}^4/(\alpha_1-\mu)^2], ~~ \alpha_1 >\mu
\ese  and
$\bb_{\mM^c}$ is $k$ spar
se. Under the above conditions and the conditions (a) and (b),
(\ref{eq:lapHT1}) has a unique local
minimizer $\wh{\bb}_a$.
\end{Th}
\noindent Proof: The proof follows the same argument as those lead
to Theorem \ref{th:1} hence is omitted.

\subsubsection{Definition needed to prove Theorems \ref{th:2} and \ref{from:th:2}}

Let $\bb^*$ is the point on the line connecting $\wh{\bb}$ and
$\wc \bb$. We write
\be\label{eq:first0}
\wh{\Q} (\bb^*)= \left\{\begin{array}{cc}\wh{\Q}_{(\mM \cup S)(\mM \cup S)
                   }(\bb^*)&\wh{\Q}_{(\mM \cup S)  (\mM \cup S) ^c}(\bb^*)\\
\wh{\Q}_{(\mM \cup S)^c  (\mM \cup S)}(\bb^*)& \wh{\Q}_{(\mM \cup S)^c  (\mM \cup S)^c} (\bb^*)\end{array}\right\}
\ee
Then by the construction in (\ref{eq:construct}), let $\wh{\bb} =
(\wh{\bb}_{\mM \cup S}\trans,{\bf 0}_{(\mM \cup S)^c}\trans)\trans $, $\wh{\bb}_{\mM \cup S}$ is the minimizer for
(\ref{eq:lapHTsub}), then  we have
\bse
&& \wh{\Q} (\bb^*)(\wh{\bb} - \wc \bb) + \left[\begin{array}{l}\left\{\frac{\partial
                               \mL(\wc \bb)}{\partial \bb}\right\}_{\mM \cup S} -
                               \left\{\A\trans \frac{\partial
                               q_{\lambda}(\wh{\bb}_{\mM^c})}{\partial
                               \bb_{\mM^c}}\right\}_{\mM \cup S}\\
\left\{\frac{\partial
                               \mL(\wc \bb)}{\partial \bb}\right\}_{(\mM \cup S)^c} -
                               \left\{\A\trans \frac{\partial
                               q_{\lambda}(\wh{\bb}_{\mM^c})}{\partial
                               \bb_{\mM^c}}\right\}_{(\mM \cup
                                     S)^c}\end{array}\right]\\
&& +\lambda \left\{
                               \begin{array}{l}(\A\trans \wh{\z})_{\mM
                                 \cup S }\\(\A\trans \wh{\z})_{(\mM
                                 \cup S)^c} \end{array}\right\} +
                           \left\{\begin{array}{c} (\A_1\trans \C\trans
                                   \bmu_3^*)_{\mM \cup S}\\\bf 0_{(\mM \cup S)^c}\end{array}\right\}= {\bf
                             0}.
\ese
Taking the upper $m+k$ non-zero component, we get
\be\label{eq:diffbb}
\wh{\bb}_{\mM \cup S} - \wc \bb_{ \mM \cup S}&=& \{\wh{\Q}_{\mM \cup S,
  \mM\cup S}(\bb^*)\}^{-1}\left[- \left\{\frac{\partial
                               \mL(\wc \bb)}{\partial \bb}\right\}_{\mM
                             \cup S} + \left\{\A\trans \frac{\partial
                               q_{\lambda}(\wh\bb_{\mM^c})}{\partial
                               \bb_{\mM^c}}\right\}_{\mM \cup S}\right.\nonumber\\
&&\left.
                           -\lambda(\A\trans \wh{\z})_{\mM\cup
                             S} - (\A_1\trans \C\trans \bmu_3^*)_{\mM
                             \cup S}\mystrut\right],
\ee
while taking the lower $p-m-k$ components, this leads to
\bse
(\A\trans \wh{\z})_{(\mM \cup S)^c} &=&  \lambda^{-1}  \left[ \left\{\A\trans \frac{\partial
                               q_{\lambda}(\wh\bb_{\mM^c})}{\partial
                               \bb_{\mM^c}}\right\}_{(\mM \cup
                                     S)^c} - \left\{\frac{\partial
                               \mL(\wc \bb)}{\partial \bb}\right\}_{(\mM
                             \cup S)^c} \right] \\
&&- \lambda^{-1}\wh{\Q}_{(\mM \cup
                           S)^c  (\mM \cup S)}(\bb^*) (\wh{\bb}_{\mM \cup S} -
                         \bb_{ \mM\cup S}) \\
&=& \lambda^{-1}  \left[ \left\{\A\trans \frac{\partial
                               q_{\lambda}(\wh{\bb}_{\mM^c})}{\partial
                               \bb_{\mM^c}}\right\}_{(\mM \cup
                                     S)^c} - \left\{\frac{\partial
                               \mL(\wc \bb)}{\partial \bb}\right\}_{(\mM
                             \cup S)^c} \right] \\
&&+ \lambda^{-1} \wh{\Q}_{(\mM \cup
                           S)^c  (\mM \cup S)}(\bb^*)\{\wh{\Q}_{(\mM \cup
                           S)(\mM \cup S) }(\bb^*)\}^{-1}\\
&&\times \left(\left[\left\{\frac{\partial
                               \mL(\wc \bb)}{\partial \bb}\right\}_{\mM \cup S} -
                              \left\{\A\trans \frac{\partial
                               q_{\lambda}(\wh{\bb}_{\mM^c})}{\partial
                               \bb_{\mM^c}}\right\}_{\mM \cup
                             S}\right] +
                           \lambda (\A\trans \wh{\z})_{\mM \cup
                             S} + (\A_1\trans \C\trans
                                   \bmu_3^*)_{\mM \cup S}\right).
\ese
Further by Condition (A4), we have
\bse
\left\{\A\trans \frac{\partial
                               q_{\lambda}(\wh \bb_{\mM^c})}{\partial
                               \bb_{\mM^c}}\right\}_{(\mM \cup
                                     S)^c}
&=& \left[\left\{ \frac{\partial \lambda |\wh\beta_{j}|}{\partial
                               \beta_{ j} }\right\} -  \left\{\frac{\partial
                               \rho_{\lambda}(\wh\beta_{j})}{\partial
                               \beta_{ j} }\right\}, j \in (\mM \cup
                           S)^c \right ]\trans\\
&=& \bf 0.
\ese
 Therefore, we have
\bse
(\A\trans \wh{\z})_{(\mM \cup S)^c} &=& \lambda^{-1}  \left[ - \left\{\frac{\partial
                               \mL(\wc \bb)}{\partial \bb}\right\}_{(\mM
                             \cup S)^c} \right] + \lambda^{-1} \wh{\Q}_{(\mM \cup
                           S)^c  (\mM \cup S)}(\bb^*)\{\wh{\Q}_{(\mM \cup
                           S)(\mM \cup S) }(\bb^*)\}^{-1}\\
&&\times \left(\left[\left\{\frac{\partial
                               \mL(\wc \bb)}{\partial \bb}\right\}_{\mM \cup S} -
                             \left\{\A\trans \frac{\partial
                               q_{\lambda}(\wh\bb_{\mM^c})}{\partial
                               \bb_{\mM^c}}\right\}_{\mM \cup
                             S}\right] +
                           \lambda (\A\trans \wh{\z})_{\mM \cup S} + (\A_1\trans \C\trans
                                   \bmu_3^*)_{\mM \cup S}\right).
\ese
By the condition that $\C\wh{\bb}_{\mM} = \C\wc {\bb} = \t$, from
(\ref{eq:diffbb}), we
have
\bse
\bf 0&=& \C [\I_{m\times m}, {\bf 0}_{m \times  k}]\{\wh{\Q}_{\mM \cup S, \mM\cup S}(\bb^*)\}^{-1}\left[- \left\{\frac{\partial
                               \mL(\wc{\bb} )}{\partial \bb}\right\}_{\mM
                             \cup S} + \left\{\A\trans \frac{\partial
                               q_{\lambda}(\wh\bb_{\mM^c})}{\partial
                             \bb_{\mM^c}}\right\}_{\mM \cup S}
                           -\lambda(\A\trans \wh{\z})_{\mM\cup S}
                         \right]\\
&&+  \C [\I_{m\times m}, {\bf 0}_{m \times k}]\{\wh{\Q}_{\mM \cup S,
  \mM\cup S}(\bb^*)\}^{-1} [\I_{m\times m}, {\bf 0}_{m \times k}]\trans
\C\trans \bmu_3^*,
\ese
which leads to
\bse
\bmu_3^* &=& (\C [\I_{m\times m}, {\bf 0}_{m \times k}]\{\wh{\Q}_{\mM \cup S,
  \mM\cup S}(\bb^*)\}^{-1} [\I_{m\times m}, {\bf 0}_{m \times k}]\trans
\C\trans)^{-1} \C [\I_{m\times m}, {\bf 0}_{m \times  k}]\{\wh{\Q}_{\mM
  \cup S, \mM\cup S}(\bb^*)\}^{-1}\\
&&\times \left[- \left\{\frac{\partial
                               \mL(\wc \bb)}{\partial \bb}\right\}_{\mM
                             \cup S} + \left\{\A\trans \frac{\partial
                               q_{\lambda}(\wh\bb_{\mM^c})}{\partial
                               \wc \bb_{\mM^c}}\right\}_{\mM \cup S}
                           -\lambda(\A\trans \wh{\z})_{\mM\cup S}
                         \right]
\ese
and
\bse
&&(\A_1\C\trans\bmu_3^*)_{\mM \cup S}\\
 &=&\A_2 \{\wh{\Q}_{\mM
  \cup S, \mM\cup S}(\bb^*)\}^{-1}\left[- \left\{\frac{\partial
                               \mL(\wc \bb)}{\partial \bb}\right\}_{\mM
                             \cup S} + \left\{\A\trans \frac{\partial
                               q_{\lambda}(\wh\bb_{\mM^c})}{\partial
                               \bb_{\mM^c}}\right\}_{\mM \cup S}
                           -\lambda(\A\trans \wh{\z})_{\mM\cup S}
                         \right],
\ese
where
\bse
\A_2 =  [\I_{m\times m}, {\bf 0}_{m
  \times k}]\trans \C\trans (\C [\I_{m\times m}, {\bf 0}_{m \times k}]\{\wh{\Q}_{\mM \cup S,
  \mM\cup S}(\bb^*)\}^{-1} [\I_{m\times m}, {\bf 0}_{m \times k}]\trans
\C\trans)^{-1} \C [\I_{m\times m}, {\bf 0}_{m \times  k}].
\ese
Hence, to use Theorem \ref{th:1}, we must show that $\|(\A\trans
\wh{\z})_{(\mM \cup S)^c}\|_\infty <1$.
Define ${\Q}(\bb) =
E\{\exp(\bb\trans\X)\X\X\trans\}$, and
\bse
\A_2^*  =  [\I_{m\times m}, {\bf 0}_{m
  \times k}]\trans \C\trans (\C [\I_{m\times m}, {\bf 0}_{m \times k}]
\{{\Q}_{\mM \cup S,
  \mM\cup S}(\bb)\}^{-1} [\I_{m\times m}, {\bf 0}_{m \times k}]\trans
\C\trans)^{-1} \C [\I_{m\times m}, {\bf 0}_{m \times  k}].
\ese

\subsubsection{Proof of Theorem \ref{th:2}}

First of all, $\wh{\bb}_{(\mM \cup S)^c} = \bf 0$ by construction.
For any unit vector, recall that
\bse
\v\trans \frac{\partial^2 \mL(\bb)}{\partial \bb\partial \bb\trans }
\w
= \exp(\bb\trans \W_i - \bb\trans
\bOmega\bb/2) \v \trans\{(\W_i -
\bOmega\bb) ^{\otimes2}- \bOmega\} \w
\ese
Denoting $\nabla ^3 \mL(\bb)$ to be the third order
gradient of $\mL$, we have
\bse
&& \v\trans \nabla ^3 \mL(\bb) \w\\
&=&n^{-1}\sumi \exp(\bb\trans \W_i - \bb\trans
\bOmega\bb/2) \v\trans\{(\W_i -
\bOmega\bb) ^{\otimes2}- \bOmega\} \w \{\W_i - \bOmega\bb\}\trans \\
&&-n^{-1}\sumi \exp(\bb\trans \W_i - \bb\trans
\bOmega\bb/2) \{\v \trans (\W_i -
\bOmega\bb) \w\trans\bOmega + \w \trans (\W_i -
\bOmega\bb) \v\trans\bOmega \}.
\ese
Hence  define vectors $\v, \w$ such that their $j$th element $|v_j|>0,
|w_j|>0$ for $j \in \mM\cup S$, and $|v_j|=|w_j|=0$ for $j\notin \mM\cup
S$, and $\|\v\|_2 = \|\w\|_2 = 1$. Firstly by Theorem \ref{from:th1} and
the condition that $\|\h_n\| _2
= O\{\sqrt{\max (m + k - r,
    r)/n}\}$, we have
\bse
\|\wh{\bb} - \bb_t\|_2 &=& \|\wh{\bb} - \wc \bb\|_2  + \|\wc{\bb} -
\bb_t\|_2 \\
&\leq &  C_1 \lambda \max(\sqrt{r}, \sqrt{m - r},\sqrt{k}) + \sqrt{(m+k)/n}\\
&\leq & C_2\left\{\frac{ \log(p)}{n}\right\}^{1/4}  \sqrt{m+ k}
\ese
for some constants $C_1, C_2 > 0$.  Further recall that ${\mathbb K}
\equiv \{\v\in {\mathbb R}^{\mM\cup S}: \|\v\|_2\leq 1\}$.
\be\label{eq:QL}
&&\sup_{\v, \w \in {\mathbb K}}\bigg| \v \trans \left\{\wh{\Q}(\bb^*) - \frac{\partial^2 \mL(\bb_t)}{\partial \bb\partial
  \bb\trans }\right\}\w\bigg|\nonumber\\
&\le& \sup_{\v, \w \in {\mathbb K}}\left[\bigg|n^{-1}\sumi \exp(\bb^{*\rm T}  \W_i - \bb^{*\rm T}
\bOmega\bb^*/2) \v\trans\{(\W_i -
\bOmega\bb^*) ^{\otimes2}- \bOmega\} \w \{\W_i - \bOmega\bb^*\}\trans
(\wh{\bb} - \bb_t)\bigg| \right. \nonumber\\
&&\left.+\bigg| n^{-1}\sumi \exp(\bb^{*\rm T} \W_i - \bb^{*\rm T}
\bOmega\bb^*/2) \{\v \trans (\W_i -
\bOmega\bb^*) \w\trans\bOmega + \w \trans (\W_i -
\bOmega\bb^*) \v\trans\bOmega \} (\wh{\bb} -
\bb_t)\bigg|\right]\nonumber\\
&\leq& \sup_{\v, \w \in {\mathbb K}}\bigg|n^{-1}\sumi \exp(\bb^{*\rm T}  \W_i - \bb^{*\rm T}
\bOmega\bb^*/2)\{ \v\trans(\W_i -
\bOmega\bb^*) ^{\otimes2}\w (\W_i - \bOmega\bb^*)\trans(\wh{\bb}
- \bb_t) \nonumber\\
&&- \v\trans\bOmega \w  (\W_i - \bOmega\bb^*)\trans(\wh{\bb} -
\bb_t)\}\bigg| \nonumber\\
&& + \sup_{\v, \w \in {\mathbb K}} \bigg| n^{-1}\sumi \exp(\bb^{*\rm T} \W_i - \bb^{*\rm T}
\bOmega\bb^*/2) \{\v \trans (\W_i -
\bOmega\bb^*) \w\trans\bOmega + \w \trans (\W_i -
\bOmega\bb^*) \v\trans\bOmega \} (\wh{\bb} -
\bb_t)\bigg|\nonumber\\
&\leq& \sup_{\v, \w \in {\mathbb K}}\bigg|\frac{n^{-1}}{2}\sumi \exp(\bb^{*\rm T}  \W_i - \bb^{*\rm T}
\bOmega\bb^*/2) \frac{(\v + \w) \trans}{\sqrt{2}}(\W_i -
\bOmega\bb^*) ^{\otimes2}\frac{\v + \w}{\sqrt{2}} (\W_i - \bOmega\bb^*)\trans(\wh{\bb}
- \bb_t) \bigg|\nonumber\\
&& + \sup_{\v, \w \in {\mathbb K}}\bigg|\frac{n^{-1}}{2}\sumi \exp(\bb^{*\rm T}  \W_i - \bb^{*\rm T}
\bOmega\bb^*/2) \frac{(\v- \w)\trans}{\sqrt{2}}(\W_i -
\bOmega\bb^*) ^{\otimes2}\frac{\v- \w}{\sqrt{2}}  (\W_i - \bOmega\bb^*)\trans(\wh{\bb}
- \bb_t) \bigg|\nonumber\\
&&+ \sup_{\v, \w \in {\mathbb K}} \bigg|\frac{(\v + \w)\trans}{\sqrt{2}} \bOmega \frac{\v + \w}{\sqrt{2}} \frac{n^{-1}}{2}\sumi \exp(\bb^{*\rm T}  \W_i - \bb^{*\rm T}
\bOmega\bb^*/2) (\W_i - \bOmega\bb^*)\trans(\wh{\bb} -
\bb_t)\bigg|  \nonumber\\
&&+  \sup_{\v, \w \in {\mathbb K}} \bigg|\frac{(\v - \w)\trans}{\sqrt{2}} \bOmega \frac{\v - \w}{\sqrt{2}} \frac{n^{-1}}{2}\sumi \exp(\bb^{*\rm T}  \W_i - \bb^{*\rm T}
\bOmega\bb^*/2)(\W_i - \bOmega\bb^*)\trans(\wh{\bb} -
\bb_t)\bigg|\nonumber\\
&& + \sup_{\v, \w \in {\mathbb K}} \bigg| n^{-1}\sumi \exp(\bb^{*\rm T} \W_i - \bb^{*\rm T}
\bOmega\bb^*/2) \{\v \trans (\W_i -
\bOmega\bb^*) \w\trans\bOmega + \w \trans (\W_i -
\bOmega\bb^*) \v\trans\bOmega \} (\wh{\bb} -
\bb_t)\bigg|\nonumber\\
&\leq&  \sup_{\v\in {\mathbb K}}\bigg|n^{-1}\sumi \exp(\bb^{*\rm T}  \W_i - \bb^{*\rm T}
\bOmega\bb^*/2) \v\trans (\W_i -
\bOmega\bb^*) ^{\otimes2}\v (\W_i - \bOmega\bb^*)\trans(\wh{\bb}
- \bb_t) \bigg|\nonumber\\
&&+  \sup_{\v\in {\mathbb K}} \bigg|\v\trans \bOmega \v n^{-1}\sumi \exp(\bb^{*\rm T}  \W_i - \bb^{*\rm T}
\bOmega\bb^*/2)(\W_i - \bOmega\bb^*)\trans(\wh{\bb} -
\bb_t)\bigg|  \nonumber\\
&& + \sup_{\v, \w \in {\mathbb K}} \bigg| n^{-1}\sumi \exp(\bb^{*\rm T} \W_i - \bb^{*\rm T}
\bOmega\bb^*/2) \{\v \trans (\W_i -
\bOmega\bb^*) \w\trans\bOmega + \w \trans (\W_i -
\bOmega\bb^*) \v\trans\bOmega \} (\wh{\bb} -
\bb_t)\bigg|\nonumber\\
&\leq& \sup_{\v\in {\mathbb K}}n^{-1}\sumi \exp(\bb^{*\rm T}  \W_i - \bb^{*\rm T}
\bOmega\bb^*/2) \v\trans (\W_i -
\bOmega\bb^*) ^{\otimes2}\v  |(\W_i - \bOmega\bb^*)\trans (\wh{\bb}
- \bb_t) |\label{eq:t1}\\
&&+ \sup_{\v\in {\mathbb K}} \v\trans \bOmega \v \bigg\| n^{-1}\sumi \exp(\bb^{*\rm T}  \W_i - \bb^{*\rm T}
\bOmega\bb^*/2) (\W_i - \bOmega\bb^*)\bigg\|_2\|(\wh{\bb} -
\bb_t)\|_2 \label{eq:t2}\\
&& + 2\sup_{\v, \w \in {\mathbb K}} \bigg| n^{-1}\sumi \exp(\bb^{*\rm T} \W_i - \bb^{*\rm T}
\bOmega\bb^*/2) \v \trans (\W_i -
\bOmega\bb^*) \{\w\trans\bOmega(\wh{\bb} -
\bb_t)\}\bigg|. \label{eq:t3}
\ee
Now for (\ref{eq:t1}),  because by Condition
\ref{con:boundX} in the
supplementary material, we have \bse
&&  |(\W_i - \bOmega\bb^*)\trans (\wh{\bb}
- \bb_t) | \\
&=&    |(\W_i -\bOmega\bb^*)\trans (\wh{\bb}
- \bb_t) /\| \wh{\bb}
- \bb_t\|_2 |  \| \wh{\bb}
- \bb_t\|_2\\
&\leq&  M_W \sqrt{m + k}\| \wh{\bb}
- \bb_t\|_2 + \|\bOmega\|_2\|\bb^*\|_2  \| \wh{\bb}
- \bb_t\|_2\\
&\leq &2 M_W \sqrt{m + k}\| \wh{\bb}
- \bb_t\|_2.
\ese
The last equality holds because $\|\bOmega\|_2 = O(1)$, $\|\bb^*\|_2
\leq \|\bb_t\|_2 + \|\wh{\bb} - \bb_t\|_2 = O_p(1) +O_p\left\{{ \log(p)}/{n}\right\}^{1/4}  \sqrt{m+ k}
= o_p(\sqrt{m + k})$.
From Corollary
\ref{cor:fromlemma9jiangext1}  in the supplementary
material, we have
\bse
&& \sup_{\v \in {\mathbb K}}\bigg|n^{-1}\sumi \exp(\bb^{*\rm T}  \W_i - \bb^{*\rm T}
\bOmega\bb^*/2) \v\trans (\W_i -
\bOmega\bb^*) ^{\otimes2} \v \sup_i| (\W_i - \bOmega\bb^*)\trans  (\wh{\bb}
- \bb)| \nonumber\\
&& \mystrut -  E\{\exp(\bb^{*\rm T}  \W_i - \bb^{*\rm T}
\bOmega\bb^*/2) \v\trans (\W_i -
\bOmega\bb^*) ^{\otimes2}  \v\} \sup_i| (\W_i - \bOmega\bb^*)\trans (\wh{\bb}
- \bb)|\bigg|\\
& \leq&  2 M_W \sqrt{(m + k)/n} \sqrt{m + k}\| \wh{\bb}
- \bb\|_2
\ese
with probability $1 - O[\exp\{- (m+k) \}]$.
Further since $ \sup_{\v\in {\mathbb K}}  E\{\exp(\bb^{*\rm T}  \W_i - \bb^{*\rm T}
\bOmega\bb^*/2) \v\trans (\W_i -
\bOmega\bb^*) ^{\otimes2}\v\}  = O(1)$ due to Condition C1(b), hence
we get
\bse
&&\sup_{\v \in {\mathbb K}}\bigg|n^{-1}\sumi \exp(\bb^{*\rm T}  \W_i - \bb^{*\rm T}
\bOmega\bb^*/2) \v\trans (\W_i -
\bOmega\bb^*) ^{\otimes2} \v (\W_i - \bOmega\bb^*)\trans  (\wh{\bb}
- \bb)\bigg|\\
&\le&\sup_{\v \in {\mathbb K}}\left\{n^{-1}\sumi \exp(\bb^{*\rm T}  \W_i - \bb^{*\rm T}
\bOmega\bb^*/2) \v\trans (\W_i -
\bOmega\bb^*) ^{\otimes2} \v \right\}\sup_i| (\W_i - \bOmega\bb^*)\trans  (\wh{\bb}
- \bb)|\\
&\le& \{2 M_W \sqrt{(m + k)/n} +O_p(1)\}\sqrt{m + k}\| \wh{\bb}
- \bb\|_2\\
&=&O_p (  \sqrt{m + k}\| \wh{\bb}
- \bb\|_2).
\ese
 For (\ref{eq:t2}), we first have
\bse
&&\| n^{-1} \sumi \exp(2\bb^{*\rm T}  \W_i - \bb^{*\rm T}
\bOmega\bb^*)   (\W_i - \bOmega\bb^*)^{\otimes2}_{\mM\cup S} \\
&&- E
\{\exp(2 \bb^{*\rm T}  \W_i - \bb^{*\rm T}
\bOmega\bb^*)  (\W_i -
\bOmega\bb^*) ^{\otimes2}_{\mM\cup S}\}\|_2\\
& =& O_p\{\sqrt{
    (m + k)/n}\}
\ese
by Corollary \ref{cor:fromlemma9jiangext}  in the supplementary
material.
Further $\| E \{\exp(2\bb^{*\rm T}  \W_i - \bb^{*\rm T}
\bOmega\bb^*)  (\W_i -
\bOmega\bb^*) ^{\otimes2}_{\mM\cup S}\}\|_2 = O(1)$.
Hence (\ref{eq:t2}) is of order $O_p(\|\wh{\bb} - \bb\|_2)$.
Now for (\ref{eq:t3}), because
\bse
&&\bigg|n^{-1}\sumi \exp(\bb^{*\rm T} \W_i - \bb^{*\rm T}
\bOmega\bb^*/2) \v \trans (\W_i -
\bOmega\bb^*) \bigg |^2\\
& \leq& n^{-1}\sumi \{\exp(2 \bb^{*\rm T} \W_i - \bb^{*\rm T}
\bOmega\bb^*) \v \trans (\W_i -
\bOmega\bb^*)^{\otimes2}\v,
\ese
by Corollary \ref{cor:fromlemma9jiangext}  in the supplementary
material, we have
\bse
&& n^{-1}\sumi\exp(2 \bb^{*\rm T} \W_i - \bb^{*\rm T}
\bOmega\bb^*) \v \trans (\W_i -
\bOmega\bb^*)^{\otimes2}\v \\
&&- E\{\exp(2 \bb^{*\rm T} \W_i - \bb^{*\rm T}
\bOmega\bb^*) \v \trans (\W_i -
\bOmega\bb^*)^{\otimes2}\v \} = O_p(\sqrt{(m+k) /n})
\ese
with probability of the order $1  - O\{\exp(m + k)\}$. Further because $E\{\exp(2 \bb^{*\rm T} \W_i - \bb^{*\rm T}
\bOmega\bb^*) \v \trans (\W_i -
\bOmega\bb^*)^{\otimes2}\v \} = O(1)$, the third term of the last line
is of the order $O_p(\|\wh{\bb} - \bb_t\|_2)$. Therefore, we have
\bse
\sup_{\v, \w \in {\mathbb K}}\bigg| \v \trans \left\{\wh{\Q}(\bb^*) - \frac{\partial^2 \mL(\bb_t)}{\partial \bb\partial
  \bb_t\trans }\right\}\w\bigg| = O_p(\sqrt{m + k} \|\wh{\bb} - \bb_t\|_2)
\ese
Hence for some positive constants $C_3$,
\be
 \sup_{\v, \w \in {\mathbb K}} \bigg|\v \trans \left\{\wh{\Q}(\bb^*) - \frac{\partial^2 \mL(\bb_t)}{\partial \bb\partial
  \bb\trans }\right\}\w\bigg|
\leq  C_{3}\left\{\frac{\log(p)}{n}\right\}^{1/4}  (m+k)\label{eq:2normbound1}
\ee
with probability $1 - O[\exp\left\{-  (m +
  k)\right\}]$.
Further, by Corollary \ref{cor:fromlemma9jiang}  in the supplementary
material, we also have
\be
 \sup_{\v, \w \in {\mathbb K}} \bigg|\v \trans \left[\frac{\partial^2 \mL(\bb_t)}{\partial \bb\partial
  \bb\trans } - \Q(\bb_t)\right]\w\bigg|
&=& \sup_{\v, \w \in {\mathbb K}}\bigg|\v_{\mM \cup S} \trans \left[\frac{\partial^2 \mL(\bb_t)}{\partial \bb\partial
  \bb\trans } -\Q(\bb_t)\right]_{\mM \cup S, \mM
  \cup S}\w_{\mM \cup S}\bigg| \nonumber\\
& = &
O_p\{\sqrt{(m + k)/n}\} \nonumber\\
&=& o_p\left[\left\{\frac{\log(p)}{n}\right\}^{1/4}  (m+ k)\right]
\label{eq:2normbound2}
\ee
 with probability $1 - 2 \exp\left\{-  (m +
  k)\right\}$.
Further, recall that
${\mathbb K}_1
\equiv \{\v\in {\mathbb R}^{(\mM\cup S)^{c}}: \|\v\|_2= 1, \|\v\|_0
= 1\}$.
\be\label{eq:QL1}
&&\sup_{\v_1 \in {\mathbb K}_1, \w \in {\mathbb K}}\bigg| \v_1 \trans \left\{\wh{\Q} (\bb^*)- \frac{\partial^2 \mL(\bb_t)}{\partial \bb\partial
  \bb\trans }\right\}\w\bigg|\nonumber\\
&\le& \sup_{\v_1 \in {\mathbb K}_1, \w \in {\mathbb K}}\left[\bigg|n^{-1}\sumi \exp(\bb^{*\rm T}  \W_i - \bb^{*\rm T}
\bOmega\bb^*/2) \v_1\trans\{(\W_i -
\bOmega\bb^*) ^{\otimes2}- \bOmega\} \w \{\W_i - \bOmega\bb^*\}\trans
(\wh{\bb} - \bb_t)\bigg| \right. \nonumber\\
&&\left.+\bigg| n^{-1}\sumi \exp(\bb^{*\rm T} \W_i - \bb^{*\rm T}
\bOmega\bb^*/2) \{\v_1 \trans (\W_i -
\bOmega\bb^*) \w\trans\bOmega + \w \trans (\W_i -
\bOmega\bb^*) \v\trans\bOmega \} (\wh{\bb} -
\bb_t)\bigg|\right]\nonumber\\
&\leq& \sup_{\v_1 \in  {\mathbb K}_1 , \w \in {\mathbb K}}\bigg|n^{-1}\sumi \exp(\bb^{*\rm T}  \W_i - \bb^{*\rm T}
\bOmega\bb^*/2)\{ \v_1\trans(\W_i -
\bOmega\bb^*) ^{\otimes2}\w (\W_i - \bOmega\bb^*)\trans(\wh{\bb}
- \bb_t) \nonumber\\
&&- \v_1\trans\bOmega \w  (\W_i - \bOmega\bb^*)\trans(\wh{\bb} -
\bb_t)\}\bigg| \nonumber\\
&& + \sup_{\v_1 \in  {\mathbb K}_1 , \w \in {\mathbb K}}\bigg| n^{-1}\sumi \exp(\bb^{*\rm T} \W_i - \bb^{*\rm T}
\bOmega\bb^*/2) \{\v_1 \trans (\W_i -
\bOmega\bb^*) \w\trans\bOmega \nonumber\\
&&+ \w \trans (\W_i -
\bOmega\bb^*) \v_1\trans\bOmega \} (\wh{\bb} -
\bb_t)\bigg|\nonumber\\
&\leq& \sup_{\v_1 \in  {\mathbb K}_1 , \w \in {\mathbb K}}\bigg|n^{-1}\sumi \exp(\bb^{*\rm T}  \W_i - \bb^{*\rm T}
\bOmega\bb^*/2) \v_1\trans(\W_i -
\bOmega\bb^*) ^{\otimes2}\w (\W_i - \bOmega\bb^*)\trans(\wh{\bb}
- \bb_t) \bigg|  \nonumber\\
&&+ \sup_{\v_1, \w \in {\mathbb K}}\bigg|n^{-1}\sumi  \exp(\bb^{*\rm T}  \W_i - \bb^{*\rm T}
\bOmega\bb^*/2)  \v_1\trans\bOmega \w  (\W_i - \bOmega\bb^*)\trans(\wh{\bb} -
\bb_t) \bigg|\nonumber\\
&& +\sup_{\v_1 \in  {\mathbb K}_1 , \w \in {\mathbb K}}\bigg| n^{-1}\sumi \exp(\bb^{*\rm T} \W_i - \bb^{*\rm T}
\bOmega\bb^*/2) \{\v_1 \trans (\W_i -
\bOmega\bb^*) \w\trans\bOmega \nonumber\\
&&+ \w \trans (\W_i -
\bOmega\bb^*) \v_1\trans\bOmega \} (\wh{\bb} -
\bb_t)\bigg|\nonumber\\
&\leq&  \sup_{\v_1 \in  {\mathbb K}_1 , \w \in {\mathbb K}} n^{-1}\sumi \bigg| \exp(\bb^{*\rm T}  \W_i - \bb^{*\rm T}
\bOmega\bb^*/2) \v_1\trans(\W_i -
\bOmega\bb^*) ^{\otimes2}\w \bigg| \bigg| (\W_i - \bOmega\bb^*)\trans(\wh{\bb}
- \bb_t) \bigg| \nonumber\\
&&+ \sup_{\v_1 \in  {\mathbb K}_1 , \w \in {\mathbb K}}\bigg|n^{-1}\sumi  \exp(\bb^{*\rm T}  \W_i - \bb^{*\rm T}
\bOmega\bb^*/2)  \v_1\trans\bOmega \w  (\W_i - \bOmega\bb^*)\trans(\wh{\bb} -
\bb_t) \bigg|\nonumber\\
&& + \sup_{\v_1 \in  {\mathbb K}_1 , \w \in {\mathbb K}} \bigg| n^{-1}\sumi \exp(\bb^{*\rm T} \W_i - \bb^{*\rm T}
\bOmega\bb^*/2) \{\v_1 \trans (\W_i -
\bOmega\bb^*) \w\trans\bOmega \nonumber\\
&&+ \w \trans (\W_i -
\bOmega\bb^*) \v_1\trans\bOmega \} (\wh{\bb} -
\bb_t)\bigg| \nonumber\\
&\leq& \sup_{\v_1 \in  {\mathbb K}_1}  \frac{1}{2} n^{-1}\sumi \exp(\bb^{*\rm T}  \W_i - \bb^{*\rm T}
\bOmega\bb^*/2) \v_1 \trans(\W_i -
\bOmega\bb^*) ^{\otimes2} \v_1  \sup_{i}\bigg| (\W_i - \bOmega\bb^*)\trans(\wh{\bb}
- \bb_t) \bigg| \label{eq:t11}\\
&&+ \sup_{\w  \in {\mathbb K}} \frac{1}{2}  n^{-1}\sumi \exp(\bb^{*\rm T}  \W_i - \bb^{*\rm T}
\bOmega\bb^*/2) \w\trans(\W_i -
\bOmega\bb^*) ^{\otimes2}\w  \sup_{i}\bigg| (\W_i - \bOmega\bb^*)\trans(\wh{\bb}
- \bb_t) \bigg| \label{eq:t12}\\
&& + \sup_{\v_1 \in  {\mathbb K}_1 , \w \in {\mathbb K}}\bigg|n^{-1}\sumi  \exp(\bb^{*\rm T}  \W_i - \bb^{*\rm T}
\bOmega\bb^*/2)  \v_1\trans\bOmega \w  (\W_i - \bOmega\bb^*)\trans(\wh{\bb} -
\bb_t) \bigg|\label{eq:t13}\\
&& +\sup_{\v_1 \in  {\mathbb K}_1 , \w \in {\mathbb K}} \bigg| n^{-1}\sumi \exp(\bb^{*\rm T} \W_i - \bb^{*\rm T}
\bOmega\bb^*/2) \{\v_1 \trans (\W_i -
\bOmega\bb^*) \w\trans\bOmega \nonumber\\
&&+ \w \trans (\W_i -
\bOmega\bb^*) \v_1\trans\bOmega \} (\wh{\bb} -
\bb_t)\bigg|.\label{eq:t14}
\ee
For (\ref{eq:t11}), from Corollary
\ref{cor:fromlemma9jiangext1}, we have
\bse
&& \sup_{\v_1 \in {\mathbb K}_1}\bigg|n^{-1}\sumi \exp(\bb^{*\rm T}  \W_i - \bb^{*\rm T}
\bOmega\bb^*/2) \v_1\trans (\W_i -
\bOmega\bb^*) ^{\otimes2} \v_1 \sup_i| (\W_i - \bOmega\bb^*)\trans  (\wh{\bb}
- \bb_t)| \nonumber\\
&& \mystrut -  E\{\exp(\bb^{*\rm T}  \W_i - \bb^{*\rm T}
\bOmega\bb^*/2) \v_1\trans (\W_i -
\bOmega\bb^*) ^{\otimes2}  \v_1\} \sup_i| (\W_i - \bOmega\bb^*)\trans (\wh{\bb}
- \bb_t)|\bigg|\\
& \leq&  2 M_W \sqrt{\log(p)/n} \sqrt{m + k}\| \wh{\bb}
- \bb_t\|_2
\ese
with probability $1 - O[\exp\{-\log(p) \}]$.
Further since $ \sup_{\v_1\in {\mathbb K}_1}  E\{\exp(\bb^{*\rm T}  \W_i - \bb^{*\rm T}
\bOmega\bb^*/2) \v_1\trans (\W_i -
\bOmega\bb^*) ^{\otimes2}\v_1\}  = O(1)$ due to Condition C1(b), hence
we get
\bse
&&\sup_{\v_1 \in {\mathbb K}_1}\bigg|n^{-1}\sumi \exp(\bb^{*\rm T}  \W_i - \bb^{*\rm T}
\bOmega\bb^*/2) \v_1\trans (\W_i -
\bOmega\bb^*) ^{\otimes2} \v_1 (\W_i - \bOmega\bb^*)\trans  (\wh{\bb}
- \bb_t)\bigg|\\
&\le&\sup_{\v_1 \in {\mathbb K}_1}\left\{n^{-1}\sumi \exp(\bb^{*\rm T}  \W_i - \bb^{*\rm T}
\bOmega\bb^*/2) \v_1\trans (\W_i -
\bOmega\bb^*) ^{\otimes2} \v_1 \right\}\sup_i| (\W_i - \bOmega\bb^*)\trans  (\wh{\bb}
- \bb_t)|\\
&\le& \{2 M_W \sqrt{\log(p)/n} +O_p(1)\}\sqrt{m + k}\| \wh{\bb}
- \bb_t\|_2\\
&=&O_p (  \sqrt{m + k}\| \wh{\bb}
- \bb\|_2).
\ese
For (\ref{eq:t12}), we use the same argument as those lead to the
order of (\ref{eq:t1}), we have (\ref{eq:t12}) is of  order $O_p (  \sqrt{m + k}\| \wh{\bb}
- \bb_t\|_2)$.
For (\ref{eq:t13}), we use the same argument as those lead to the
order of (\ref{eq:t2}), we have  (\ref{eq:t13}) is of order $O_p ( \| \wh{\bb}
- \bb_t\|_2)$.

Now for (\ref{eq:t14}), because
\bse
&&\bigg|n^{-1}\sumi \exp(\bb^{*\rm T} \W_i - \bb^{*\rm T}
\bOmega\bb^*/2) \v_1 \trans (\W_i -
\bOmega\bb^*) \bigg |^2\\
& \leq& n^{-1}\sumi \{\exp(2 \bb^{*\rm T} \W_i - \bb^{*\rm T}
\bOmega\bb^*) \v_1 \trans (\W_i -
\bOmega\bb^*)^{\otimes2}\v_1,
\ese
by Corollary \ref{cor:fromlemma9jiangext}  in the supplementary
material, we have
\bse
&& n^{-1}\sumi\exp(2 \bb^{*\rm T} \W_i - \bb^{*\rm T}
\bOmega\bb^*) \v_1 \trans (\W_i -
\bOmega\bb^*)^{\otimes2}\v_1 \\
&&- E\{\exp(2 \bb^{*\rm T} \W_i - \bb^{*\rm T}
\bOmega\bb^*) \v_1 \trans (\W_i -
\bOmega\bb^*)^{\otimes2}\v_1 \} = O_p(\sqrt{\log(p) /n})
\ese
with probability of the order $1  - O[\exp\{-\log(p)\}]$. Further because $E\{\exp(2 \bb^{*\rm T} \W_i - \bb^{*\rm T}
\bOmega\bb^*) \v_1 \trans (\W_i -
\bOmega\bb^*)^{\otimes2}\v_1 \} = O(1)$, and because of the same
argument as those lead to (\ref{eq:t13}),  we conclude that
(\ref{eq:t14})
is of the order $O_p(\|\wh{\bb} - \bb_t\|_2)$.
Hence follow (\ref{eq:QL}), and by Conditions
\ref{con:boundX},  Corollaries
\ref{cor:fromlemma9jiang}--\ref{cor:fromlemma15loh2012}  in the supplementary
material, for
positive constants $C_4$  we have
\be
\sup_{\v_1 \in {\mathbb K}_1, \w \in {\mathbb K}}\bigg|\v_1 \trans \left\{\wh{\Q}(\bb^*) - \frac{\partial^2 \mL(\bb)}{\partial \bb\partial
  \bb\trans }\right\}\w\bigg|
&=& O_p(\sqrt{m + k} \|\wh{\bb} - \bb_t\|_2)\nonumber\\
 &\leq& C_4
\left\{\frac{\log(p)}{n}\right\}^{1/4}  (m+ k)
\label{eq:supnormbound1}
\ee
and
\be
 && \sup_{\v_1 \in {\mathbb K}_1, \w \in {\mathbb K}}\bigg|\v_{1}\trans \left[\frac{\partial^2 \mL(\bb_t)}{\partial \bb\partial
  \bb\trans } - \Q(\bb_t)\right]\w\bigg| \nonumber\\
&=& \sup_{\v_1 \in {\mathbb K}_1, \w \in {\mathbb K}} \bigg|\v_{1 (\mM \cup S)^c} \trans \left[\frac{\partial^2 \mL(\bb_t)}{\partial \bb\partial
  \bb\trans } - \Q(\bb_t)\right]_{(\mM \cup S)^c, \mM
  \cup S}\w_{\mM \cup S}\bigg|\nonumber \\
&= &
O_p[\sqrt{\max\{(m + k ), \log(p)\}/n}] \nonumber\\
&=& o_p\left[\left\{\frac{\log(p)}{n}\right\}^{1/4} (m+ k)\right]
\label{eq:supnormbound2},
\ee
with probabilities to the order of $1 - 2\exp[- \max\{(m + k )
, \log(p)\}]$.

Combine (\ref{eq:2normbound1}) and (\ref{eq:2normbound2}) with  Lemma
\ref{lem:fromlemma11loh2017}  in the supplementary
material,
\be\label{eq:2normdiff}
&&\|(\{\wh \Q_{(\mM \cup S), \mM \cup S})(\bb^*)\}^{-1} - \{\Q_{ (\mM \cup S), \mM \cup S}(\bb_t)\}^{-1}\|_2\nonumber \\
&\leq& \frac{\|\{\Q_{(\mM \cup S), \mM \cup S}(\bb_t)\}^{-1}\|_2^2 \|\wh \Q_{(\mM \cup S), \mM \cup S} (\bb^*)- \Q_{(\mM \cup S), \mM \cup S}(\bb_t)\|_2} {\{1 -\|\{\Q_{ (\mM \cup S), \mM \cup S}(\bb_t)\}^{-1} \|_2\|\wh \Q_{(\mM \cup S), \mM \cup S}(\bb^*)- \Q_{(\mM \cup S), \mM \cup S}(\bb_t)\|_2\}}\nonumber\\
&=& O_p(\| \wh \Q_{(\mM \cup S), \mM \cup S} (\bb^*)-\Q_{0 (\mM \cup S), \mM \cup S}(\bb_t)\|_2)\nonumber\\
&=& C_9
\left[\left\{\frac{\log(p)}{n}\right\}^{1/4}  (m+ k)
\right],
\ee
and hence
\be\label{eq:qdiffsup}
&&\|(\{\wh \Q_{(\mM \cup S), \mM \cup S})(\bb^*)\}^{-1} - \{\Q_{(\mM \cup S), \mM
  \cup S}(\bb_t)\}^{-1}\|_\infty\nonumber\\
&=& C_9  (m + k)^{3/2 }\left\{\frac{\log(p)}{n}\right\}^{1/4}.
\ee
Further,
\be
&& \bigg\|\wh \Q_{(\mM \cup S)^c, \mM \cup S} (\bb^*) \{\wh \Q_{(\mM \cup S), \mM \cup
  S}(\bb^*)\}^{-1}\left\{\frac{\partial \mL(\wc \bb)}{\partial \bb}\right\}_{\mM \cup S}\bigg\|_\infty\nonumber\\
&\le& \bigg\|[\wh\Q_{(\mM \cup S)^c, \mM \cup S} (\bb^*)- \Q_{(\mM \cup S)^c, \mM \cup S}(\bb_t)]\{\wh \Q_{(\mM \cup S), \mM \cup
  S}(\bb^*)\}^{-1} \left\{\frac{\partial \mL(\wc\bb)}{\partial
  \bb}\right\}_{\mM \cup S}\bigg\|_\infty \nonumber\\
&& + \bigg\|\Q_{(\mM \cup S)^c,
    \mM \cup S}(\bb_t)(\{\wh \Q_{(\mM \cup S), \mM \cup
  S}(\bb^*)\}^{-1} \nonumber\\
&&- \{\Q_{(\mM \cup S), \mM \cup S}(\bb_t)\}^{-1})\left\{\frac{\partial \mL(\wc\bb)}{\partial
  \bb}\right\}_{\mM \cup S}\bigg\|_\infty \nonumber\\
&& + \bigg\|\Q_{(\mM \cup S)^c, \mM \cup S} (\bb_t)\{\Q_{(\mM \cup S), \mM \cup S}(\bb_t)\}^{-1} \left\{\frac{\partial \mL(\wc\bb)}{\partial
  \bb}\right\}_{\mM \cup S}\bigg\|_\infty\nonumber\\
&\leq& \sup_{\v_1 \in {\mathbb K}_1, \w \in {\mathbb K}}
\bigg |\v_{1 (\mM \cup S)^c} \trans [\wh\Q_{(\mM \cup S)^c, \mM \cup S} (\bb^*)-
\Q_{(\mM \cup S)^c, \mM \cup
  S}(\bb_t)]\w_{\mM \cup S}\bigg |\nonumber\\
&&\times  \|\{\wh \Q_{(\mM \cup S), \mM \cup
  S}(\bb^*)\}^{-1}\|_2 \bigg\|\left\{\frac{\partial \mL(\wc \bb)}{\partial
  \bb}\right\}_{\mM \cup S}\bigg\|_2\nonumber\\
&& + \|\Q_{ (\mM \cup S)^c, \mM \cup S} (\bb_t)\|_2 \bigg\|\{\wh \Q_{(\mM \cup S), \mM \cup
  S}(\bb^*)\}^{-1} - \{\Q_{(\mM \cup S), \mM \cup
  S}(\bb_t)\}^{-1}\bigg\|_2 \bigg\|\left\{\frac{\partial \mL(\wc \bb)}{\partial
  \bb}\right\}_{\mM \cup S}\bigg\|_2\nonumber\\
&& + \bigg\|\Q_{(\mM \cup S)^c, \mM \cup S}(\bb_t) \{\Q_{0(\mM \cup S), \mM
  \cup S}(\bb_t)\}^{-1} \left\{\frac{\partial \mL(\wc \bb)}{\partial
  \bb}\right\}_{\mM \cup S}\bigg\|_\infty\nonumber\\
&\leq&  C_{10 }
\left[\left\{\frac{\log(p)}{n}\right\}^{1/4}  (m+ k) \right] \sqrt{m + k}
\sqrt{\log (p)/n}\label{eq:ineq1}\\
&& + C_{11}\left[\left\{\frac{\log(p)}{n}\right\}^{1/4}  (m+ k)
\right]\sqrt{m + k}
\sqrt{\log (p)/n}
\label{eq:ineq2}\\
&& + \bigg\|\Q_{(\mM \cup S)^c, \mM \cup S} (\bb_t)\{\Q_{(\mM \cup S), \mM
  \cup S}(\bb_t)\}^{-1} \left\{\frac{\partial \mL(\wc \bb)}{\partial
  \bb}\right\}_{\mM \cup S}\bigg\|_\infty. \nonumber
\ee
 (\ref{eq:ineq1}) is obtained by using (\ref{eq:supnormbound1}) and
(\ref{eq:supnormbound2}) and the fact that
\bse
&&\|\{\wh \Q_{(\mM \cup S), \mM \cup
  S}(\bb^*)\}^{-1} \|_2\\
&\leq& \|\{\wh \Q_{(\mM \cup S), \mM \cup
  S}(\bb^*)\}^{-1} -[\Q_{ (\mM \cup S), \mM \cup S}(\bb_t)]^{-1}\|_2 +\| \{\Q_{(\mM \cup S), \mM \cup S}(\bb_t)\}^{-1}\|_2\\
&\leq& O_p(1).
\ese
(\ref{eq:ineq2})  is obtained by using (\ref{eq:2normbound1}),
(\ref{eq:2normbound2})  and the fact that
$\|\Q_{(\mM \cup S)^c, \mM \cup S}(\bb)
\|_2 = O(1)$.
Therefore, now note that $n\geq c_\infty (m+ k)^4  \log(p)$ and
$\lambda = O[\{\log(p)/n\}^{1/4}]$
by the statement
assumption, together with Condition
\ref{con:pn}, we have
\be\label{eq:QQL}
&& \bigg\|\wh \Q_{(\mM \cup S)^c, \mM \cup S} (\bb^*)\{\wh \Q_{(\mM \cup S), \mM \cup
  S}(\bb^*)\}^{-1}\left\{\frac{\partial \mL(\wc \bb)}{\partial
    \bb}\right\}_{\mM \cup S}\bigg\|_\infty \nonumber\\
&\leq&\bigg\|\Q_{(\mM \cup S)^c, \mM \cup S}(\bb_t) \{\Q_{(\mM \cup S), \mM
  \cup S}(\bb_t)\}^{-1} \left\{\frac{\partial \mL(\wc \bb)}{\partial
  \bb}\right\}_{\mM \cup S}\bigg\|_\infty +  o_p(\lambda)\nonumber\\
&\leq& \bigg\|\Q_{(\mM \cup S)^c, \mM \cup S} (\bb)\{\Q_{(\mM \cup S), \mM
  \cup S}(\bb_t)\}^{-1}\bigg\|_2 \bigg\|\left\{\frac{\partial \mL(\wc \bb)}{\partial
  \bb}\right\}_{\mM \cup S}\bigg\|_2+  o_p(\lambda)\nonumber\\
&=& \sqrt{m + k}
\sqrt{\log (p)/n}
+  o_p(\lambda)\nonumber\\
&=& o_p(\lambda ).
\ee
Therefore
\be\label{eq:AZ}
&& \|(\A\trans \wh{\z})_{(\mM \cup S)^c}\|_\infty \nonumber\\
&=& \|\lambda^{-1}  \left[ - \left\{\frac{\partial
                               \mL(\wc{\bb})}{\partial \bb}\right\}_{(\mM
                             \cup S)^c} \right] + \lambda^{-1} \wh{\Q}_{(\mM \cup
                           S)^c  (\mM \cup S)}(\bb^*)\{\wh{\Q}_{(\mM \cup
                           S)(\mM \cup S) }(\bb^*)\}^{-1}\nonumber\\
&&\times \left(\left[\left\{\frac{\partial
                               \mL(\wc{\bb})}{\partial \bb}\right\}_{\mM \cup S} -
                             \left\{\A\trans \frac{\partial
                               q_{\lambda}(\wh\bb_{\mM^c})}{\partial
                               \bb_{\mM^c}}\right\}_{\mM \cup
                             S}\right] +
                           \lambda (\A\trans \wh{\z})_{\mM \cup
                             S}+   (\A_1\trans \C\trans \bmu_3^*)_{(\mM
                             \cup S)} \right) \|_\infty \nonumber\\
&\le& \lambda^{-1}\bigg\| \left[ - \left\{\frac{\partial
                               \mL(\wc{\bb})}{\partial \bb}\right\}_{(\mM
                             \cup S)^c} \right] \bigg\|_\infty
                         \nonumber\\
&& + \lambda^{-1}\bigg \|\wh{\Q}_{(\mM \cup
                           S)^c  (\mM \cup S)}(\bb^*)\{\wh{\Q}_{(\mM \cup
                           S)(\mM \cup S) }(\bb^*)\}^{-1} \left\{\frac{\partial
                               \mL(\wc{\bb})}{\partial \bb}\right\}_{\mM
                             \cup S} \bigg\|_\infty \nonumber\\
&&+ \lambda^{-1} \bigg\|\wh{\Q}_{(\mM \cup
                           S)^c  (\mM \cup S)}(\bb^*)\{\wh{\Q}_{(\mM \cup
                           S)(\mM \cup S) }(\bb^*)\}^{-1}  \left(\left[ -
                             \left\{\A\trans \frac{\partial
                               q_{\lambda}(\wh\bb_{\mM^c})}{\partial
                               \bb_{\mM^c}}\right\}_{\mM \cup
                             S}\right] +
                           \lambda (\A\trans \wh{\z})_{\mM \cup
                             S}\right)\bigg\|_\infty \nonumber\\
&& + \lambda^{-1} \|\wh{\Q}_{(\mM \cup
                           S)^c  (\mM \cup S)}(\bb^*)\{\wh{\Q}_{(\mM \cup
                           S)(\mM \cup S) }(\bb^*)\}^{-1}  (\A_1\trans \C\trans \bmu_3^*)_{\mM
                             \cup S}\|_\infty\nonumber\\
&=& \lambda^{-1} \bigg\|\wh{\Q}_{(\mM \cup
                           S)^c  (\mM \cup S)}(\bb^*)\{\wh{\Q}_{(\mM \cup
                           S)(\mM \cup S) }(\bb^*)\}^{-1}  \left(\left[ -
                             \left\{\A\trans \frac{\partial
                               q_{\lambda}(\wh\bb_{\mM^c})}{\partial
                               \bb_{\mM^c}}\right\}_{\mM \cup
                             S}\right] +
                           \lambda (\A\trans \wh{\z})_{\mM \cup
                             S}\right)\bigg\|_\infty \nonumber\\
&&+\lambda^{-1}\bigg\|\wh{\Q}_{(\mM \cup
                           S)^c  (\mM \cup S)}(\bb^*)\{\wh{\Q}_{(\mM \cup
                           S)(\mM \cup S) }(\bb^*)\}^{-1} \A_2\{ \wh{\Q}_{\mM
  \cup S, \mM\cup S}(\bb^*)\}^{-1}\nonumber\\
&&\times \left[ \left\{\A\trans \frac{\partial
                               q_{\lambda}(\wh\bb_{\mM^c})}{\partial
                               \bb_{\mM^c}}\right\}_{\mM \cup S}
                           -\lambda(\A\trans \wh{\z})_{\mM\cup S}
                         \right]\bigg\|_\infty+ o(1).
\ee
Now recall (\ref{eq:diffbb})
and
\be\label{eq:a3}
&&(\A_1\C\trans\bmu_3^*)_{\mM \cup S} \\
&=&\A_2 \{\wh{\Q}_{\mM
  \cup S, \mM\cup S}(\bb^*)\}^{-1}\left[- \left\{\frac{\partial
                               \mL(\wc \bb)}{\partial \bb}\right\}_{\mM
                             \cup S} + \left\{\A\trans \frac{\partial
                               q_{\lambda}(\wh\bb_{\mM^c})}{\partial
                               \bb_{\mM^c}}\right\}_{\mM \cup S}
                           -\lambda(\A\trans \wh{\z})_{\mM\cup S}
                         \right]. \nonumber
\ee
Hence
\be\label{eq:difbb1}
&&\|\wh{\bb}_{\mM \cup S} - \wc \bb_{ \mM \cup S}\|_\infty \n\\
&\le& \|\{\wh{\Q}_{\mM
  \cup S, \mM\cup S}(\bb^*)\}^{-1} \left\{\frac{\partial
                               \mL(\wc \bb)}{\partial \bb}\right\}_{\mM
                             \cup S} \|_\infty\n\\
&&+ \bigg\|\{\wh{\Q}_{\mM \cup S,
                             \mM\cup S}(\bb^*)\}^{-1} \left[\left\{\A\trans \frac{\partial
                               q_{\lambda}(\wh \bb_{\mM^c})}{\partial
                               \bb_{\mM^c}}\right\}_{\mM \cup S}
                           -\lambda(\A\trans \wh{\z})_{\mM\cup
                             S}\right]\bigg\|_\infty\n\\
&& +\| \{\wh{\Q}_{\mM \cup S, \mM\cup S}(\bb^*)\}^{-1} \A_2 \{\wh{\Q}_{\mM
  \cup S, \mM\cup S}(\bb^*)\}^{-1}\nonumber\\
&&\times \left[- \left\{\frac{\partial
                               \mL(\wc \bb)}{\partial \bb}\right\}_{\mM
                             \cup S} + \left\{\A\trans \frac{\partial
                               q_{\lambda}(\wh\bb_{\mM^c})}{\partial
                               \bb_{\mM^c}}\right\}_{\mM \cup S}
                           -\lambda(\A\trans \wh{\z})_{\mM\cup S}
                         \right]\|_\infty\n\\
&\leq  &o_p(\lambda) + \lambda \|\{\wh{\Q}_{\mM \cup S,
                             \mM\cup S}(\bb^*)\}^{-1}\|_\infty + \lambda \|\{\wh{\Q}_{\mM \cup S, \mM\cup S}(\bb^*)\}^{-1} \A_2 \{\wh{\Q}_{\mM
  \cup S, \mM\cup S}(\bb^*)\}^{-1}\|_\infty \n\\
&=& o_p(\lambda) + 2\lambda c_\infty +  2\lambda c_{\infty},
\ee
where $o_p(\lambda)$ in  the second inequality is obtained by using the same
argument as those lead to
(\ref{eq:QQL}). The other two terms are obtained by
                           using Lemma \ref{lem:fromlemm5support} such
                           that $$\bigg\|\left[\left\{\A\trans \frac{\partial
                               q_{\lambda}(\wh \bb_{\mM^c})}{\partial
                               \bb_{\mM^c}}\right\}_{\mM \cup S}
                           -\lambda(\A\trans \wh{\z})_{\mM\cup
                             S}\right]\bigg\|_\infty \leq \lambda.$$
The last equality holds because
\bse
&& \|\{\wh{\Q}_{\mM \cup S,
                             \mM\cup S}(\bb^*)\}^{-1}\|_\infty\\
&\leq & \|\{\wh{\Q}_{\mM \cup S,
                             \mM\cup S}(\bb^*)\}^{-1} - \{\Q_{(\mM \cup S), \mM \cup S}(\bb_t)\}^{-1}\|_\infty+ \| \{\Q_{(\mM \cup S), \mM \cup S}(\bb_t)\}^{-1}\|_\infty \\
&\leq&2 c_\infty
\ese
by (\ref{eq:qdiffsup}).
And similarly, we have
\bse
&& \|\{\wh{\Q}_{\mM \cup S, \mM\cup S}(\bb^*)\}^{-1} \A_2 \{\wh{\Q}_{\mM
  \cup S, \mM\cup S}(\bb^*)\}^{-1}\|_\infty\\
&\leq &\|\{\wh{\Q}_{\mM \cup S, \mM\cup S}(\bb^*)\}^{-1} \A_2 \{\wh{\Q}_{\mM
  \cup S, \mM\cup S}(\bb^*)\}^{-1} \\
&&-\{ {\Q}_{\mM \cup S, \mM\cup S}(\bb_t)\}^{-1} \A_2 \{{\Q}_{\mM
  \cup S, \mM\cup S}(\bb_t)\}^{-1}\|_\infty \\
&&+\| \{{\Q}_{\mM \cup S, \mM\cup S}(\bb_t)\}^{-1} \A_2 \{{\Q}_{\mM
  \cup S, \mM\cup S}(\bb_t)\}^{-1}\|_\infty\\
&=&  2 c_{\infty}.
\ese
Now because $\min(|\beta_{j}|)\geq \lambda (\gamma + 5c_\infty )$ for $j
\in S$, we have
\bse
\min |\wh{\beta}_j| \geq \lambda \gamma
\ese
and $q_\lambda'(\wh{\beta})_j = \lambda {\rm sign}(\wh{\beta}_j) =
\lambda (\A\trans\wh{\z})_j$.
Hence by Lemma \ref{lem:fromlemma4loh2017}  in the supplementary
material, we have
\be\label{eq:zeros}
\left\{\A\trans \frac{\partial
                               q_{\lambda}(\wh \bb_{\mM^c})}{\partial
                               \bb_{\mM^c}}\right\}_{\mM \cup S}
                           -\lambda(\A\trans \wh{\z})_{\mM\cup
                             S} = \bf 0.
\ee
Inserting (\ref{eq:zeros}) into (\ref{eq:AZ}), the first two terms of
(\ref{eq:AZ}) are zero hence
$
\|(\A\trans\z)_{(\mM \cup S)^c}\| =o_p(1)<1.
$
This implies that $\wh{\bb}$ has support in $\mM \cup
S$. Further, because of (\ref{eq:difbb1}),
$\wh{\bb}$ is unique in the feasible set.
Therefore, using (\ref{eq:diffbb}) together with (\ref{eq:zeros}) and (\ref{eq:a3}),
we have
\bse
&&\wh{\bb}_{\mM \cup S}- \wc \bb_{ \mM \cup S}\\
&=& (\{\wh{\Q}_{\mM \cup S,
  \mM\cup S}(\bb^*)\}^{-1} - \{\wh{\Q}_{\mM \cup S, \mM\cup S}(\bb^*)\}^{-1} \A_2 \{\wh{\Q}_{\mM
  \cup S, \mM\cup S}(\bb^*)\}^{-1} ) \\
&&\times \left[- \left\{\frac{\partial
                               \mL(\wc \bb)}{\partial \bb}\right\}_{\mM
                             \cup S} \right] \\
&=& - (\{{\Q}_{\mM \cup S,
  \mM\cup S}(\bb)\}^{-1} -\{ {\Q}_{\mM \cup S, \mM\cup S}(\bb_t)\}^{-1} \A_2^*\{ {\Q}_{\mM
  \cup S, \mM\cup S}(\bb_t)\}^{-1} ) \\
&&\times \left[\left\{\frac{\partial
                               \mL(\wc \bb)}{\partial \bb}\right\}_{\mM
                             \cup S} \right] \{1 + o_p(1)\}\\
\ese
and $\wh{\bb}_{(\mM \cup S)^c} = \bf 0$.
Now because
\bse
&& \frac{\partial
                               \mL(\wc \bb)}{\partial \bb} - \frac{\partial
                               \mL(\bb_t)}{\partial \bb} \\
&=& \frac{\partial^2
                               \mL(\bb^* )}{\partial \bb \partial
                               \bb\trans} (\wc{\bb} - \bb_t)\\
&=& \Q(\bb) (\wc{\bb} - \bb_t) +
\left\{\frac{\partial^2
                               \mL(\bb_t )}{\partial \bb \partial
                               \bb\trans} - \Q(\bb_t)\right\}(\wc{\bb} - \bb_t)
                             + \left\{\frac{\partial^2
                               \mL(\bb^* )}{\partial \bb \partial
                               \bb\trans} - \frac{\partial^2
                               \mL(\bb_t )}{\partial \bb \partial
                               \bb\trans} \right\}  (\wc{\bb} -
                           \bb_t) \\
&=&  \Q(\bb_t)(\wc{\bb}  - \bb_t)  \{1 + o_p(1)\}.
\ese
The last equality holds because
\bse
\|\left\{\frac{\partial^2
                               \mL(\bb_t )}{\partial \bb \partial
                               \bb\trans} - \Q(\bb_t)\right\}\|_2 = o_p(1),
\ese
by (\ref{eq:2normbound2}),
\bse
\|\left\{\frac{\partial^2
                               \mL(\wt \bb )}{\partial \bb \partial
                               \bb\trans} - \frac{\partial^2
                               \mL(\bb_t )}{\partial \bb \partial
                               \bb\trans} \right\}\|_2 = o_p(1)
\ese by (\ref{eq:2normbound1}),  and $\|\Q(\bb_t)\|_2 = O(1)$.  Therefore,
\bse
\wh{\bb}_{\mM \cup S}- \wc \bb_{ \mM \cup S}
&=& - (\{{\Q}_{\mM \cup S,
  \mM\cup S}(\bb_t)\}^{-1} -\{ {\Q}_{\mM \cup S, \mM\cup S}(\bb_t)\}^{-1} \A_2^* \{{\Q}_{\mM
  \cup S, \mM\cup S}(\bb_t)\}^{-1} ) \nonumber\\
&&\times \left[\left\{\frac{\partial
                               \mL( \bb_t)}{\partial \bb} + \Q (\bb_t) (\wc{\bb} - \bb_t) \right\}_{\mM
                             \cup S} \right] \{1 + o_p(1)\}.
\ese
Further, we have
\bse
&&\wh{\bb}_{\mM \cup S}- \bb_{ t\mM \cup S}\\
&=& \wh{\bb}_{\mM \cup S}- \wc \bb_{ \mM \cup S} -
[(\C\C\trans)^{-1}\C, {\bf 0}_{ r \times k}]\trans \h_n\\
&=& - (\{{\Q}_{\mM \cup S,
  \mM\cup S}(\bb_t)^{-1} \}- \{{\Q}_{\mM \cup S, \mM\cup S}(\bb_t)\}^{-1} \A_2^* \{{\Q}_{\mM
  \cup S, \mM\cup S}(\bb_t)\}^{-1} ) \left[\left\{\frac{\partial
                               \mL( \bb_t)}{\partial \bb}\right\}_{\mM
                             \cup S} \right.\\
&&\left.+ {\Q}_{\mM
  \cup S, \mM\cup S}(\bb_t)  [(\C\C\trans)^{-1}\C, {\bf 0}_{ r \times
  k}]\trans \h_n  \mystrut \right] \{1 + o_p(1)\} -
[(\C\C\trans)^{-1}\C, {\bf 0}_{ r \times k}]\trans \h_n\\
&=& - \{({\Q}_{\mM \cup S,
  \mM\cup S}(\bb_t)\}^{-1} - \{{\Q}_{\mM \cup S, \mM\cup S}(\bb_t)\}^{-1} \A_2^* \{{\Q}_{\mM
  \cup S, \mM\cup S}(\bb_t)\}^{-1} ) \left\{\frac{\partial
                               \mL( \bb_t)}{\partial \bb}\right\}_{\mM
                             \cup S}\{1 + o_p(1)\}\\
&& -\{{\Q}_{\mM \cup S, \mM\cup S}(\bb_t)\}^{-1} \A_2^* [\{(\C\C\trans)^{-1}\C, {\bf 0}_{ r \times
  k}\}\trans \h_n ]\{1 + o_p(1)\}.
\ese
\qed

\subsubsection{Proof of Theorem \ref{from:th:2}}
The proof is very similar to that of Theorem \ref{th:2} hence is omitted.

\subsection{Proof of Results in Section \ref{sec:asytest}}

\subsubsection{Proof of Lemma \ref{th:chisq}}
 First note that
\bse
&&\bPsi^{-1/2}(\bSig, \Q, \bb_t)\bomega_n \\
&=&- \sqrt{n}\bPsi^{-1/2}(\bSig, \Q, \bb_t) \C [\I_{m\times m}, {\bf 0}_{m \times k}] {\Q}_{ \mM \cup
  S,  \mM \cup S}^{-1} (\bb_t)\left\{\frac{\partial \mL(\bb_t)}{\partial
    \bb}\right\}_{\mM \cup S}\\
&=&\bPsi^{-1/2}(\bSig, \Q, \bb_t)\C [\I_{m\times m}, {\bf 0}_{m \times k}] {\Q}_{ \mM \cup
  S,  \mM \cup S}^{-1} (\bb_t)\\
&&\times n^{-1/2 }\sumi \{Y_i \W_i - \exp(\bb_t\trans\W_i - \bb_t\trans
\bOmega \bb_t/2)  (\W_i -
\bOmega \bb_t)\}_{\mM \cup S}.
\ese
It is easy to see that
\be\label{eq:var1}
&&\sumi \cov\left(\bPsi^{-1/2}(\bSig, \Q, \bb_t)\C [\I_{m\times m}, {\bf 0}_{m \times k}] {\Q}_{ \mM \cup
  S,  \mM \cup S}^{-1} (\bb_t)\right.\nonumber\\
&&\left.\times  n^{-1/2} \{Y_i \W_i - \exp(\bb_t\trans\W_i - \bb_t\trans
\bOmega \bb_t/2)  (\W_i -
\bOmega  \bb_t)\}_{\mM \cup S}\right) = \I_{r\times r}.
\ee
Further,
\be\label{eq:3ordop1}
&& r^{1/4} \sumi E \|n^{-1/2 } \bPsi^{-1/2}(\bSig, \Q, \bb_t)\C [\I_{m\times m}, {\bf 0}_{m \times k}] {\Q}_{ \mM \cup
  S,  \mM \cup S}^{-1} (\bb_t) \nonumber\\
&&\times \{Y_i \W_i - \exp(\bb_t\trans\W_i - \bb_t\trans
\bOmega \bb_t/2)  (\W_i -
\bOmega \bb_t)\}_{\mM \cup S}\|_2^3 \nonumber\\
&\leq&  r^{1/4} \sumi E \|n^{-1/2 } \{Y_i \bPsi^{-1/2} (\bSig, \Q, \bb_t)\C [\I_{m\times m}, {\bf 0}_{m \times k}]
{\Q}_{ \mM \cup
  S,  \mM \cup S}^{-1} (\bb_t)  \W_{i\mM\cup S} \nonumber\\
&& - \exp(\bb_t\trans\W_i - \bb_t\trans
\bOmega \bb_t/2)  \bPsi^{-1/2}(\bSig, \Q, \bb_t)\C [\I_{m\times m}, {\bf 0}_{m \times k}] {\Q}_{ \mM \cup
  S,  \mM \cup S}^{-1} (\bb_t)   (\W_i -
\bOmega \bb_t) _{\mM \cup S}\}\|_2^3 \nonumber\\
&\leq&  r^{1/4} \sumi E \|n^{-1/2}\{Y_i \bPsi^{-1/2} (\bSig, \Q, \bb_t)\C [\I_{m\times m}, {\bf 0}_{m \times k}]
{\Q}_{ \mM \cup
  S,  \mM \cup S}^{-1} (\bb_t)  \W_{i\mM\cup S} \nonumber\\
&& - \exp(\bb_t\trans\W_i - \bb_t\trans
\bOmega \bb_t/2)  \bPsi^{-1/2}(\bSig, \Q, \bb_t)\C [\I_{m\times m}, {\bf 0}_{m \times k}] {\Q}_{ \mM \cup
  S,  \mM \cup S}^{-1} (\bb_t)   (\W_i -
\bOmega \bb_t) _{\mM \cup S}\}\|_2^3 \nonumber\\
&\leq&  r^{1/4} D \sumi   n^{-3/2 } (\|\bPsi^{-1/2} (\bSig, \Q, \bb_t)\C [\I_{m\times m}, {\bf 0}_{m \times k}]{\Q}_{\mM \cup
  S,  \mM \cup S}^{-1}(\bb_t)  \W_{i \mM\cup S} \|_2^3 \nonumber\\
&&+ \|\bPsi^{-1/2} (\bSig, \Q, \bb_t)\C [\I_{m\times m}, {\bf 0}_{m \times k}]{\Q}_{\mM \cup
  S,  \mM \cup S}^{-1}(\bb_t) (\bOmega  \bb_t)_{\mM\cup S}\|_2^3) \nonumber\\
&=& o_p(1),
\ee
where $D$ is a positive constant. The second to the last inequality
holds by  Condition \ref{con:3ord}.   Also because
\bse
&& \|\bPsi^{-1/2} (\bSig, \Q, \bb_t)\C [\I_{m\times m}, {\bf 0}_{m \times k}]{\Q}_{\mM \cup
  S,  \mM \cup S}^{-1}(\bb_t)  \W_{i \mM\cup S} \|_2\\
&\leq& \|\bPsi^{-1/2} (\bSig, \Q, \bb_t)\C [\I_{m\times m}, {\bf 0}_{m \times k}]{\Q}_{\mM \cup
  S,  \mM \cup S}^{-1}(\bb_t) \|_2 \|\W_{i \mM\cup S} \|_2\\
&\leq& \|\bPsi^{-1/2} (\bSig, \Q, \bb_t)\C [\I_{m\times m}, {\bf 0}_{m \times k}]{\Q}_{\mM \cup
  S,  \mM \cup S}^{-1}(\bb_t) \|_2 |\W_i \trans \v |/\|\v\|_2\\
&\leq& O(1) M_W \sqrt{m +k}= O(\sqrt{m+k}),
\ese
where $\v$ is a $p$ dimensional vector with $\|\v\|_0 = m+k$, the last
inequality holds by Condition \ref{con:boundX} (a). Further since $\|\bPsi^{-1/2} (\bSig, \Q, \bb_t)\C [\I_{m\times m}, {\bf 0}_{m \times k}]{\Q}_{\mM \cup
  S,  \mM \cup S}^{-1}(\bb_t) (\bOmega  \bb_t)_{\mM\cup S}\|_2 =
O(1)$,
\bse
 && r^{1/4} D \sumi   n^{-3/2 } (\|\bPsi^{-1/2} (\bSig, \Q, \bb_t)\C [\I_{m\times m}, {\bf 0}_{m \times k}]{\Q}_{\mM \cup
  S,  \mM \cup S}^{-1}(\bb_t)  \W_{i \mM\cup S} \|_2^3 \nonumber\\
&&+ \|\bPsi^{-1/2} (\bSig, \Q, \bb_t)\C [\I_{m\times m}, {\bf 0}_{m \times k}]{\Q}_{\mM \cup
  S,  \mM \cup S}^{-1}(\bb_t) (\bOmega  \bb_t)_{\mM\cup S}\|_2^3)\\
&=& O\{r^{1/4} (m+k)^{3/2} n^{- 1/2}\}= o_p(1),
\ese
by the Condition that $n\geq c_\infty (m + k)^4 \log(p)$.
Combine (\ref{eq:var1})  and (\ref{eq:3ordop1}), and using Lemma
\ref{lem:frombentkus} in Section \ref{sec:restlemma}, let $\Z$ is a
standard Gaussian random variable we have
\be\label{eq:asynormal}
\lim_{n\to \infty}\sup_{\mathcal{C}} |\Pr(\bPsi^{-1/2}(\bSig, \Q,
\bb_t)\bomega_n \in {\mathcal{C}} ) - \Pr(\Z \in {\mathcal C} )| = 0.
\ee
Now we choose
\bse
{\mathcal{C}_x} = \{\z: \|\z - \sqrt{n}\bPsi^{-1/2}\h_n \|_2 \leq x\},
\ese
then
\bse
\lim_{n\to \infty}  |\Pr(\bPsi^{-1/2}(\bSig, \Q,
\bb_t)\bomega_n \in {\mathcal{C}_x} ) - \Pr(\Z \in {\mathcal C}_x )| =
0,
\ese
which implies
\bse
&& \lim_{n\to \infty} |\Pr(T_0 \leq x) -
\Pr\{\chi^2(r, n \h_n\trans \bPsi^{-1}(\bSig, \Q,
\bb_t) \h_n ) \leq x\}| = 0,
\ese
where $\chi^2(r, \gamma)$ is a non-central chi-square
distribution, with non-centrality parameter $\gamma$. \qed

\subsubsection{Proof of Theorem \ref{th:TW}}
 From Theorem \ref{from:th:2}, we have
\bse
\wh{\bb}_{a \mM \cup S}- \bb_{t\mM \cup S}
&=&  - \{{\Q}_{\mM \cup S,
  \mM\cup S}(\bb_t)\}^{-1} \left\{\frac{\partial
                               \mL( \bb_t)}{\partial \bb}\right\}_{\mM
                             \cup S}\{1 + o_p(1)\}
\ese
Hence,
\bse
\sqrt{n}\C(\wh{\bb}_{a \mM }- \bb_{t\mM}) &=& - \sqrt{n} \C
[\I_{m\times m}, {\bf 0}_{m \times k}] \{ {\Q}_{\mM \cup S,
  \mM\cup S}(\bb_t)\}^{-1} \left\{\frac{\partial
                               \mL( \bb_t)}{\partial \bb}\right\}_{\mM
                             \cup S}\{1 + o_p(1)\}\\
&=& \bomega_n\{1 + o_p(1)\}.
\ese
Further, because $\C \bb_{t\mM} = \t + \h_n$, we have
\bse
&&\sqrt{n}\bPsi^{-1/2} (\wh\bSig,\wh\bQ,\wh \bb_a)(\C\wh{\bb}_{a \mM }- \t) = \bPsi^{-1/2} (\wh\bSig,
\wh\bQ, \wh\bb_a)\bomega_n\{1 + o_p(1 )\} +
\sqrt{n}\bPsi^{-1/2} (\wh\bSig,
\wh\bQ, \wh\bb_a)  \h_n.
\ese
Further because $\h_n = O\{\sqrt{\max (m + k - r,
    r)/n}\}$,
\be\label{eq:TWT0}
&&n (\C\wh{\bb}_{a \mM }- \t)  \trans \bPsi^{-1} (\wh\bSig,\wh\bQ,\wh
\bb_a)(\C\wh{\bb}_{a \mM }- \t) \nonumber\\
&=&  (\bomega_n + \sqrt{n}\h_n)\trans \bPsi^{-1} (\bSig,
\bQ, \bb_t) (\bomega_n + \sqrt{n}\h_n) \{1 + o_p(1)\}\nonumber\\
&&+ (\bomega_n + \sqrt{n}\h_n)\trans \{\bPsi^{-1} (\wh\bSig,
\wh\bQ, \wh \bb_a) -  \bPsi^{-1} (\bSig,
\bQ, \bb_t)\} (\bomega_n + \sqrt{n}\h_n) \nonumber\\
&\leq&(\bomega_n + \sqrt{n}\h_n)\trans \bPsi^{-1} (\bSig,
\bQ, \bb_t) (\bomega_n + \sqrt{n}\h_n) \{1 + o_p(1)\}\nonumber\\
&&+ \|\bomega_n + \sqrt{n}\h_n\|_2^2 \|\bPsi^{-1} (\wh\bSig,
\wh\bQ, \wh \bb_a) -  \bPsi^{-1} (\bSig,
\bQ, \bb_t)\|_2 \nonumber \\
&=& T_0 + o_p(r).
\ee
The last equality holds because $T_0$ converge in distribution to a
non-central chi-square distribution with degree freedom $r$ as
shown in Lemma \ref{th:chisq}  and
hence $T_0$ is  of the order $O_p(r)$. Also because \bse
\|\bomega_n + \sqrt{n}\h_n\|_2^2&\leq& (\bomega_n + \sqrt{n}\h_n)\trans \bPsi^{-1} (\bSig,
\bQ, \bb_t) (\bomega_n + \sqrt{n}\h_n)  /\alpha_{\min} \{\bPsi^{-1} (\bSig,
\bQ, \bb_t)\}\\
&\leq&  (\bomega_n + \sqrt{n}\h_n)\trans \bPsi^{-1} (\bSig,
\bQ, \bb_t) (\bomega_n + \sqrt{n}\h_n)  /c_{\bPsi}
\ese
by Condition \ref{con:psieigen},
we have
\bse
&& \|\bPsi^{-1} (\wh\bSig,
\wh\bQ, \wh \bb_a) -  \bPsi^{-1} (\bSig,
\bQ, \bb_t)\|_2  \|\bomega_n + \sqrt{n}\h_n\|_2^2 \\
&\leq&\left [\left\{\frac{\log(p)}{n}\right\}^{1/4}  (m+ k)\right]O_p(T_0).
\ese
Therefore,  $T_W - T_0 = o_p(r)$.
\qed

\subsubsection{Proof of Theorem \ref{th:TS}}
From Theorem \ref{th:2}, we have
\bse
&&\wh{\bb}_{\mM \cup S}- \bb_{t\mM \cup S}\\
&=&  - (\{{\Q}_{\mM \cup S,
  \mM\cup S}(\bb_t)\}^{-1} - \{{\Q}_{\mM \cup S, \mM\cup S}(\bb_t)\}^{-1} \A_2^* \{{\Q}_{\mM
  \cup S, \mM\cup S}(\bb_t)\}^{-1} ) \left\{\frac{\partial
                               \mL( \bb_t)}{\partial \bb}\right\}_{\mM
                             \cup S}\{1 + o_p(1)\}\\
&& +\{{\Q}_{\mM \cup S, \mM\cup S}(\bb_t)\}^{-1} \A_2^*  [\{(\C\C\trans)^{-1}\C, {\bf 0}_{ r \times
  k}\}\trans \h_n ]\{1 + o_p(1)\}
\ese
and $\wh{\bb}_{(\mM \cup S)^c} = \bf 0$.
 By Taylor expansion we have
\bse
&& \left\{\frac{\partial
    \mL(\wh{\bb})}{\partial \bb }\right\}_{\mM \cup S} \\
& =& \left\{\frac{\partial
    \mL({\bb_t})}{\partial \bb }\right\}_{\mM \cup S} + \left\{ \frac{\partial^2
    \mL({\bb}^*)}{\partial \bb \partial \bb\trans } \right\}_{\mM \cup
S}(\wh{\bb} - \bb_t)_{\mM \cup S} \\
&=&\left\{\frac{\partial
    \mL({\bb_t})}{\partial \bb }\right\}_{\mM \cup S} + \Q _{\mM \cup S,\mM \cup S } (\bb_t)(\wh{\bb} - \bb_t)_{\mM \cup S}
\{1 + o_p(1)\}\\
&=& \A_2^* \{{\Q}_{\mM
  \cup S, \mM\cup S}(\bb_t)\}^{-1}  \left\{\frac{\partial
                               \mL( \bb_t)}{\partial \bb}\right\}_{\mM
                             \cup S}\{1 + o_p(1)\}\\
&& +\{\A_2^*  [\{(\C\C\trans)^{-1}\C, {\bf 0}_{ r \times
  k}\}\trans \h_n ]\{1 + o_p(1)\}
\ese
where $\bb^*$ is a point in between $\bb_t$ and $\wh{\bb}$. The second
equality holds by (\ref{eq:2normbound1}) and
(\ref{eq:2normbound2}) in Appendix.
Therefore, we have
\bse
&&\sqrt{n}\bPsi(\wh{\bSig}, \wh{\Q}, \wh\bb)^{-1/2}\C [\I_{m\times m}, {\bf 0}_{m \times k}]\{\wh {\Q}_{\mM \cup S,
  \mM\cup S}(\wh \bb)\}^{-1 }\A_2^* \{{\Q}_{\mM
  \cup S, \mM\cup S}(\bb_t)\}^{-1}   \left\{\frac{\partial
    \mL(\bb_t)}{\partial \bb }\right\}_{\mM \cup S} \\
&=& \sqrt{n}\bPsi(\wh{\bSig}, \wh{\Q}, \wh\bb)^{-1/2}\C [\I_{m\times m}, {\bf 0}_{m \times k}]\{ {\Q}_{\mM \cup S,
  \mM\cup S}(\bb_t)\}^{-1 } \left\{\frac{\partial
    \mL({\bb_t})}{\partial \bb }\right\}_{\mM \cup S} \{1 + o_p(1)\}\\
&=& \bPsi(\wh{\bSig}, \wh{\Q}, \wh\bb)^{-1/2}\bomega_n\{1 +
o_p(1)\},
\ese
The first
equality holds by (\ref{eq:2normbound1}) and
(\ref{eq:2normbound2}) in the supplementary material.
And similarly
\bse
&&\sqrt{n}\bPsi(\wh{\bSig}, \wh{\Q}, \wh\bb)^{-1/2}\C [\I_{m\times m}, {\bf 0}_{m \times k}]\{\wh {\Q}_{\mM \cup S,
  \mM\cup S}(\wh \bb)\}^{-1 }  \{\A_2^*  [\{(\C\C\trans)^{-1}\C, {\bf 0}_{ r \times
  k}\}\trans \h_n ]\\
&=& \sqrt{n}\bPsi(\wh{\bSig}, \wh{\Q}, \wh\bb)^{-1/2} \h_n\{1 + o_p(1)\}
\ese
and hence
\bse
T_S =  (\bomega_n + \sqrt{n}\h_n)\trans \bPsi(\wh{\bSig}, \wh{\Q},
\wh\bb)^{-1}  (\bomega_n + \sqrt{n}\h_n)\{1 + o_p(1)\}
\ese
Now  the same steps in (\ref{eq:TWT0}) lead to $T_S - T_0 =
o_p(r)$.
\qed

\setcounter{section}{0}
\setcounter{equation}{0}
\def\theequation{S\arabic{section}.\arabic{equation}}
\def\thesection{S\arabic{section}}

\fontsize{12}{14pt plus.8pt minus .6pt}\selectfont

\renewcommand{\theequation}{B.\arabic{equation}}
\renewcommand{\thetable}{B.\arabic{table}}
\renewcommand{\thefigure}{B.\arabic{figure}}
\renewcommand{\thesubsection}{B.\arabic{section}.\arabic{subsection}}
\renewcommand{\thesection}{B.\arabic{section}}

\setcounter{Lem}{0}\renewcommand{\theLem}{B.\arabic{Lem}}
\setcounter{Th}{0}\renewcommand{\theTh}{B.\arabic{Th}}
\setcounter{Cor}{0}\renewcommand{\theCor}{B.\arabic{Cor}}

\section*{Appendix B}
\section{Some Lemmas on the Penalty Function}

\begin{Lem}\label{lem:fromlemma4}
Conditions (A1)--(A4)
imply that $\rho_{\lambda}$ is $\lambda $-Lipschitz and  all
sub gradients and derivatives of $\rho_{\lambda}$ are bounded by
$\lambda $ in magnitude. 
Conditions (A1)--(A5) imply
\bse
\lambda \|\bb_{\mM^c}\|_1\leq \rho_\lambda(\bb_{\mM^c}) + \mu/2 \|\bb_{\mM^c}\|_2^2,
\forall \bb_{\mM^c}\in \mathbb {\R}^{p-m}
\ese
\end{Lem}
\noindent Proof: This lemma is a direct consequence of Lemma 4 in
\cite{loh2015}.

\begin{Lem}\label{lem:fromlemma5}
Suppose $\rho_\lambda$ satisfies Conditions (A1)--(A5). Let
$\v\in\mathbb{R}^{p-m}$,  let $\mathcal{A}$ be the index set of $k$
largest elements of $\v$ in magnitude, and 
 let $\mathcal{A}^c$ be the index set of the remaining $p-m-k$
elements of $\v$. 
Suppose $\xi>0$ and satisfies 
\bse
\xi \rho_{\lambda}(\v_{\mathcal{A}}) -
\rho_{\lambda}(\v_{\mathcal{A}^c})\geq 0. 
\ese
Then 
\bse
\xi \rho_{\lambda}(\v_{\mathcal{A}}) -
\rho_{\lambda}(\v_{\mathcal{A}^c}) \leq \lambda (\xi\|\v_{\mathcal{A}}\|_1 - \|\v_{\mathcal{A}^c}\|_1).
\ese
Moreover, if $\bb_{t\mM^c} \in \mathbb{R}^{p-m}$ is $k$-sparse, then for a vector
$\bb_{\mM^c}\in \mathbb{R}^{p-m}$ such that $\xi \rho_{\lambda}(\bb_{t\mM^c} ) -
\rho_{\lambda}(\bb_{\mM^c}) >0$ and $\xi \geq 1$, we have 
\bse
\xi \rho_{\lambda}(\bb_{t\mM^c}) - \rho_{\lambda}(\bb_{\mM^c}) \leq \lambda (\xi
\|\v_{\mathcal{A}}\|_1 - \|\v_{\mathcal{A}^c}\|_1), 
\ese
where $\v= \bb_{\mM^c}- \bb_{t\mM^c}$, $\mathcal{A}$ is the index set of the $k$
largest elements of $\v$ in magnitude,  and
$\mathcal{A}^c$ is the index set of the remaining $p-m-k$ elements of $\v$.
\end{Lem}
\noindent Proof: This lemma is a direct consequence of Lemma 5 in
\cite{loh2015}. \qed

\begin{Lem}\label{lem:fromlemma4loh2017}
Suppose $\rho_{\lambda}$ is $(\mu, \gamma)$-amenable,
$|\wh\beta_j|\geq 
\lambda \gamma $ for $j \in \mM \cup S$, then
$q_{\lambda}'(\wh\beta_j) = \lambda {\rm sign}(\wh\beta_j)$. 
\end{Lem}
\noindent Proof: Because $\rho_{\lambda}$ is $(\mu, \gamma)$-amenable,
$\rho'(\wh{\beta}_j) = 0$ by Condition (A6) and (A1). Hence $q_{\lambda}'(\wh\beta_j)  =\partial
\lambda |\wh \beta_j |/\partial \wh \beta_j = \lambda {\rm
  sign}(\wh\beta_j)$. This proves the result. \qed

\begin{Lem}\label{lem:fromlemm5support}
Consider a $\mu$-amenable regularizer $\rho_{\lambda}$. Then \\
(a) $|\rho_{\lambda}'(t)| \leq \lambda$ for all $t\neq 0$. \\
(b) The function $q_{\lambda}(t) - \mu/2 t^2$ is concave and
everywhere differentiable, where $q_{\lambda}(t) = \lambda
|t| - \rho_{\lambda} (t)$. 
\end{Lem}
\noindent Proof: This lemma is a direct consequence of Lemma 5 in
\cite{loh2017}. \qed

\section{Some Lemmas on the Criterion Function}

\begin{Lem}\label{lem:fromCor1}
Assume Conditions \ref{con:boundX} -- \ref{con:pn}
hold.  There exists a constant $c_1>0$ so that
\bse
\pr \left[n^{-1}\|\sumi \{Y_i \W_i - \exp(\bb_t\trans\W_i - \bb_t\trans
\bOmega \bb_t/2)  (\W_i - \bOmega\bb_t)\} \|_\infty 
 > c_1\sqrt{\log(p) /n}\mystrut\right]\leq 6p^{-1}. 
\ese
\end{Lem}
\noindent Proof: The lemma is the direct consequence of Corollary 1 in
\cite{jiang2021}. We omit the proof here.

\begin{Lem}\label{lem:RE1}
Assume that Conditions \ref{con:boundX}, \ref{con:pn} and
\ref{con:ww2} hold,  then 
for any $\bb$ with  $\|\bb\|_2 \leq R_2$, and 
for sufficiently large $n$ and $p$, 
with probability $1  - O[\exp\{-\sqrt{n\log p}\}]$,
$n^{-1} \sumi \exp(\bb\trans \W_i - \bb\trans
\bOmega\bb/2)  \{(\W_i -
\bOmega\bb) ^{\otimes2} - \bOmega\}$ satisfies the lower and upper-RE conditions
with   $$\alpha_1 = \min_{\|\bb\|_1\le R_1, \|\bb\|_2\le R_2}
\alpha_{\min}
  [E\{\exp(\bb\trans\X_i)\X_i\X_i\trans\}]  /2,$$ 
$$\alpha_2 =\max_{\|\bb\|_1\le R_1, \|\bb\|_2\le R_2}
3 \alpha_{\max}
  [E\{\exp(\bb\trans\X_i)\X_i\X_i\trans\}]/ 2 $$
and 
$\tau(n, p) 
=\tau_1\sqrt{\log(p)/n}
$
for a bounded positive constant $\tau_1$.
\end{Lem}
\noindent Proof: The lemma is a direct consequence of Lemma 12 in
\cite{jiang2021} with $ c=1$. \qed

\begin{Cor}\label{cor:modLemma1}
Assume Conditions \ref{con:boundX}--\ref{con:pn} to hold, $\|\h_n\| _2
= O\{\sqrt{\max (m + k - r, 
    r)/n}\}$ and $m + k = o(n^{1/3})$, then there exists
    a positive constant $c_{10}$ so that
\bse
\|\frac{\partial\mL(\wc{\bb})}{\partial \bb}\|_\infty
\leq c_{10}\max[\sqrt{(m + k)/n}, \sqrt{\log(p)/n}], 
\ese
with probability $1 - O(p^{-1}) -O[\exp\{-
\sqrt{n\log (p)}\}]$.
\end{Cor}
\noindent Proof:
\bse
&& \|\frac{\partial\mL(\wc{\bb})}{\partial \bb}\|_\infty\\
& =& n^{-1}\|\sumi \{Y_i \W_i - \exp(\wc{\bb}\trans\W_i - \wc{\bb}\trans
\bOmega \wc{\bb}/2)  (\W_i - 
\bOmega \wc{\bb})\} \|_\infty\\
&\leq & n^{-1}\|\sumi \{Y_i \W_i - \exp({\bb_t}\trans\W_i - {\bb_t}\trans
\bOmega {\bb_t}/2)  (\W_i -
\bOmega {\bb_t})\} \|_\infty \\
&&+ n^{-1}\|\sumi \exp(\wc{\bb}\trans\W_i - \wc{\bb}\trans
\bOmega \wc{\bb}/2)  (\W_i - \wc{\bb}\trans
\bOmega)\}  - \exp({\bb_t}\trans\W_i - {\bb_t}\trans
\bOmega {\bb_t}/2)  (\W_i -
\bOmega  {\bb_t})\}  \|_\infty. 
\ese
Let $\mathbb{K}_0=\{\v\in{\mathbb R}^p:\v_{\mM^c}=0, \|\v\|_2 = 1\}.$ Then,
\bse
&& n^{-1}\|\sumi \exp(\wc{\bb}\trans\W_i - \wc{\bb}\trans
\bOmega \wc{\bb}/2)  (\W_i - \wc{\bb}\trans
\bOmega)  - \exp({\bb_t}\trans\W_i - {\bb_t}\trans
\bOmega {\bb_t}/2)  (\W_i - 
\bOmega {\bb_t}) \|_\infty\\
&\leq  & n^{-1}\| \sumi \exp(\bb^{*\rm T} \W_i - \bb^{*\rm T} 
\bOmega\bb^*/2)  \{(\W_i -
\bOmega\bb^*) ^{\otimes2} - \bOmega\} \{\wc{\bb} - {\bb_t}\}\|_2\\
&\leq& \sup_{\v\in \mathbb{K}}  n^{-1}\sumi \v\trans \exp(\bb^{*\rm T} \W_i - \bb^{*\rm T} 
\bOmega\bb^*/2)  \{(\W_i -
\bOmega\bb^*) ^{\otimes2} - \bOmega\} \v  \|\wc{\bb} -
{\bb_t}\|_2\\
&\leq & 3 \alpha_{\max}
  [E\{\exp(\bb\trans\X_i)\X_i\X_i\trans\}]/ 2 \|\C\trans(\C\C\trans)^{-1}\|_2\|\h_n\|_2+ \tau_1 m\sqrt{\log(p)/n}\|\C\trans(\C\C\trans)^{-1}\|_2\|\h_n\|_2\\
&\leq& c_{10}\max[\sqrt{(m + k)/n}, \sqrt{\log(p)/n}]
\ese
for some constants $c_{10}$. The first inequality holds by the Taylor
expansion with $\bb^*$ be  the point on the line connecting $\wc{\bb}$ and
$\bb$. The second inequality holds because $\wc{\bb} - \bb$ only
has $m$ nonzero elements supported in $\mM$, and hence only the 
corresponding $m\times m$ sub-matrix  in $n^{-1}\sumi \exp(\bb^{*\rm T} \W_i -
\bb^{*\rm T}  
\bOmega\bb^*/2)  \{(\W_i -
\bOmega\bb^*) ^{\otimes2} - \bOmega\} $ contributes to the $L_2$
norm. The third inequality holds by  Lemma \ref{lem:RE1} and the fact that
$\v$ only as $m$ non-zero elements. The last equality holds because
$\|\h_n\| _2 = O\{\sqrt{\max (m + k - r, r)/n} \}\leq O(\sqrt{m+ k/n} )$, and
$m \sqrt{m+ k/n} \to 0$ because $m + k = o(n^{1/3})$. This proves the
results. \qed

\begin{Lem}\label{lem:fromlemma9jiang}
Assume that Conditions \ref{con:boundX} and \ref{con:ww2} hold.  If
$\X_i, \U_i \in {\mathbb R}^p$, define ${\mathbb K}
\equiv \{\v\in {\mathbb R}^{\mM\cup S}: \|\v\|_2\leq 1\}$, ${\mathbb K}_1
\equiv \{\v\in {\mathbb R}^{(\mM\cup S)^{c}}: \|\v\|_2= 1, \|\v\|_0
= 1\}$  \bse 
&&\Pr\left(\sup_{\v, \w \in {\mathbb K}} |\sumi [A(\bb\trans\W_i)
  g(\W_i, \bb, \v, \w) - \v\trans E\{\exp(\bb\trans 
  \X_i) \X_i\X_i\trans\} \w ]| > nt \right)\\
& \leq& 2 \exp\left(-  \min \left[ \frac{ n  t^2}{324 e^2 M_4 },\frac{n
    t}{36  e M_5 \log (n) } \right] + 2 (m + k)\log(9)\right), 
\ese
 \bse 
&&\Pr\left(\sup_{\v \in {\mathbb K}_1, \w \in {\mathbb K}} |\sumi [A(\bb\trans\W_i)
  g(\W_i, \bb, \v, \w) - \v\trans E\{\exp(\bb\trans 
  \X_i) \X_i\X_i\trans\} \w ]| > nt \right)\\
& \leq& 2 \exp\left(-  \min \left[ \frac{ n  t^2}{36 e^2 M_4 },\frac{n
    t}{12  e M_5 \log (n) } \right] +  (m + k)\log (9) +\log(p)
\right), 
\ese
and 
\bse
&& \Pr\left(\sup_{\v\in \mathbb{K}_1, \w\in \mathbb{K}_1} |\sumi [A(\bb\trans\W_i)
  g(\W_i, \bb, \v, \w) - \v\trans E\{\exp(\bb\trans
  \X_i) \X_i\X_i\trans\} \w] | > nt \right)\\
&\leq&    2 \exp\left(-  \min \left[ \frac{ n  t^2}{16 e^2 M_4 },
\frac{nt}{8  e M_5 \log (n) } \right] +  2\log(p) \right).
\ese
\end{Lem}
\noindent Proof: 
By Lemma 1 statement 3 and Lemma 3 statement 3 in \cite{jiang2021}, we can see 
that the square of a conditional sub-Gaussian
variable is  sub-exponential. Now because $\v \trans(\W_i -
\bb\trans \bOmega) $ given $\X_i$ and $\bb\trans\W_i$ is normal, and
hence 
$
 \v \trans (\W_i -
\bb\trans \bOmega) ^{\otimes2} \v
$
is conditional sub-exponential. Now by the Cauchy--Schwarz inequality, and
without loss of generality we assume $ \v \trans (\W_i -
\bb\trans \bOmega) ^{\otimes2} \v \geq  \w \trans (\W_i -
\bb\trans \bOmega) ^{\otimes2} \w$, we 
have 
\bse
 |\v \trans (\W_i -
\bb\trans \bOmega) ^{\otimes2} \w| \leq \v \trans (\W_i -
\bb\trans \bOmega) ^{\otimes2} \v. 
\ese
Hence, by Lemma 3 statement 3 in \cite{jiang2021}, we have for some
bounded positive $K_3(\bb\trans \W_i, \X_i)$, 
\bse
&& E  [\exp\{|\v \trans (\W_i -
\bb\trans \bOmega) ^{\otimes2} \w|/K_3(\bb\trans \W_i,
\X_i)\}|\bb\trans \W_i, \X_i] \\
&\leq& E  [\exp\{|\v \trans (\W_i -
\bb\trans \bOmega) ^{\otimes2} \v|/K_3(\bb\trans
\W_i, \X_i)\}|\bb\trans\W_i, \X_i] \leq e
\ese 
Hence, $ \v \trans (\W_i -
\bb\trans \bOmega) ^{\otimes2} \w$ is also conditional sub-exponential
variable. Therefore, 
we have that
\bse
g(\W_i, \bb, \v, \w) - E\{g(\W_i, \bb, \v, \w)
|\bb\trans \W_i, \X_i\}
\ese
is centered sub-exponential. 
Then using the same argument as those that lead to Lemma 8
in \cite{jiang2021} and Condition \ref{con:ww2},  for any  unit
vectors $\v , \w$, 
we can show that
\bse
&& \Pr\left(|\sumi [A(\bb\trans\W_i)
  g(\W_i, \bb, \v, \w) - \v\trans E\{\exp(\bb\trans 
  \X_i) \X_i\X_i\trans\} \w ]| > nt \right) \\
&\leq& 2 \exp\left(-  \min \left[ \frac{ n  t^2}{16 e^2 M_4 },\frac{n
    t}{8  e M_5 \log (n) } \right] \right).
\ese 
Now we define 
$\mathcal{B} = \{\u_1, \ldots, \u_r\}\subset {\mathbb K}$   to be  a 1/3-cover of ${\mathbb K}$, if for every $\v \in {\mathbb K}$, there is some $\u_i \in
\mathcal{B} $ such that $\|\v - \u_i\|_2 \leq 1/3$. Define
$\Delta\v=\v-\u_j$ where $\u_j = \arg\min_{\u_i}\|\v -\u_i\|_2$. We have $\|\Delta\v\|_2  \leq
1/3$. Similarly define $\u_k = \arg\min_{\u_i}\|\w -\u_i\|_2$ for $\w
\in \mathbb{K}$. 
By
\cite{ledoux2013}, we can construct $\mathcal{B}$ with $|\mathcal{B}|
< 9^{(m+ k)}$.
Now for $\v_1, \v_2 \in \mathbb{K}$, define 
\bse
\Phi (\v_1, \v_2) = \v_1 \trans \left[\sumi  \frac{A(\bb\trans\W_i)\{ (\W_i -
\bOmega\bb) ^{\otimes2} - \bOmega \} -  E\{\exp(\bb\trans
  \X_i) \X_i\X_i\trans\}}{n} \right]\v_2. 
\ese 
We have 
\bse
|\Phi (\v, \w)|
&= & |\Phi (\Delta \v + \u_j , \Delta \w
+ \u_k )|  \\
&\le&  \max_{j, k}  |\Phi (\u_j,
\u_k)| + \max_i  |\Phi (\Delta
\v,\u_i)  | + \max_i | \Phi (\u_i,\Delta
\w)|  +   | \Phi (\Delta
\v,\Delta
\w)  |.
\ese 
Since $\|3\Delta\v\|_2 \leq 1$ and $\supp(3 \Delta\v)\subseteq  \mathbb{K}$, 
$3\Delta\v  \in \mathbb{K}$. 
It follows that 
\bse
&&\sup_{\v, \w \in \mathbb{K}}|\Phi (\v, \w)| \\
&\leq&
\max_{j, k} |\Phi (\u_j,
\u_k)| + 1/3 \sup_{\v \in \mathbb{K}} \max_i |\Phi (3\Delta
\v,\u_i)  | + 1/3 \sup_{\w \in \mathbb{K}} \max_i |\Phi (\u_k, 3\Delta \w)  | +  1/9\sup_{\v, \w \in \mathbb{K} } | \Phi (3\Delta
\v,3\Delta
\w)  |\\
&\le&
\max_{j, k} |\Phi (\u_j,
\u_k)| + 2/3 \{\sup_{\v \in \mathbb{K}} |\Phi (3\Delta
\v, 3\Delta\v)  |\}^{1/2}\{ \max_i|\Phi (
\u_i,\u_i)  | \}^{1/2}
+ 1/9\sup_{\v, \w\in \mathbb{K}} | \Phi (\v,
\w)  |\\
&\le&
 \max_{j, k} |\Phi (\u_j,
\u_k)| + \sup_{\v, \w \in \mathbb{K} } \{2/3 |\Phi (\v, \w)| + 1/9
|\Phi (\v, \w)|\}. 
\ese
Hence,  
$\sup_{\v, \w \in \mathbb{K} }|\Phi (\v, \w)|\leq 9/2 \max_{j, k} |\Phi (\u_j,
\u_k)|$. 
By a union bound while considering $|\mathcal{B}|
< 9^{2(m+ k)}$, we have 
\bse
&& \Pr\left(\sup_{\v, \w \in \mathbb{K}} |\sumi [A(\bb\trans\W_i)
  g(\W_i, \bb, \v, \w) - \v\trans E\{\exp(\bb\trans
  \X_i) \X_i\X_i\trans\} \w] | > 9/2 nt \right)\\
&\le&\Pr\left(\max_{j, k}|\sumi [A(\bb\trans\W_i)
  g(\W_i, \bb, \v, \w) - \v\trans E\{\exp(\bb\trans
  \X_i) \X_i\X_i\trans\} \w] | >  nt \right)\\
&\le&\Pr\left(\max_{j, k}|\sumi [A(\bb\trans\W_i)
  g(\W_i, \bb, \v/\|\v\|_2, \w/\|\w\|_2) - \v/\|\v\|_2\trans E\{\exp(\bb\trans
  \X_i) \X_i\X_i\trans\} \w/\|\w\|_2] | >  nt \right)\\
&\leq&  2 \exp\left(-  \min \left[ \frac{ n  t^2}{16 e^2 M_4 },\frac{n
    t}{8  e M_5 \log (n) } \right] + 2 (m+k)\log(9)  \right).
\ese
It also follows that 
\bse
&&\sup_{\v\in \mathbb{K}_1,\w\in\mathbb{K}}|\Phi (\v, \w)| \\
&\leq&
\max_{\v\in\mathbb{K}_1,j} |\Phi (\v,\u_j)| 
+ 1/3 \sup_{\v\in \mathbb{K}_1, \w \in \mathbb{K}} |\Phi (
\v, 3\Delta\w)  | \\
&\le&
 \max_{\v \in\mathbb{K}_1,j} |\Phi (\v,\u_j)| + 1/3\sup_{\v \in \mathbb{K}_1,\w\in\mathbb{K} } \{ |\Phi (\v, \w)|,
\ese
so
$
\sup_{\v\in \mathbb{K}_1,\w\in\mathbb{K}}|\Phi (\v, \w)| 
\le 3/2\max_{\v \in\mathbb{K}_1,j} |\Phi (\v,\u_j)|$.
Because $\v \in \mathbb{K}_1 $ is a vector with a single
nonzero entry  $1$,  there are only $p-m-k$ elements in
$\mathbb{K}_1$. We thus have
\bse
&& \Pr\left(\sup_{\v\in \mathbb{K}_1, \w\in \mathbb{K}} |\sumi [A(\bb\trans\W_i)
  g(\W_i, \bb, \v, \w) - \v\trans E\{\exp(\bb\trans
  \X_i) \X_i\X_i\trans\} \w] | > 3/2 nt \right)\\
&\le&\Pr\left(\max_{\v\in \mathbb{K}_1,\u_j} |\sumi [A(\bb\trans\W_i)
  g(\W_i, \bb, \v, \u_j) - \v\trans E\{\exp(\bb\trans
  \X_i) \X_i\X_i\trans\} \u_j] | >  nt \right)\\
&\leq&    2 \exp\left(-  \min \left[ \frac{ n  t^2}{16 e^2 M_4 },\frac{n
    t}{8  e M_5 \log (n) } \right] +  (m+k)\log(9) + \log(p) \right).
\ese
This proves the result. 
Further, since ${\mathbb K}_1$ contains only unit
  vectors, we have 
\bse
&& \Pr\left(\sup_{\v\in \mathbb{K}_1, \w\in \mathbb{K}_1} |\sumi [A(\bb\trans\W_i)
  g(\W_i, \bb, \v, \w) - \v\trans E\{\exp(\bb\trans
  \X_i) \X_i\X_i\trans\} \w] | > nt \right)\\
&\leq&    2 \exp\left(-  \min \left[ \frac{ n  t^2}{16 e^2 M_4 },\frac{n
    t}{8  e M_5 \log (n) } \right] +  2\log(p) \right).
\ese
\qed

\begin{Cor}\label{cor:fromlemma9jiang}
Assume that Conditions \ref{con:boundX}, \ref{con:pn} and \ref{con:ww2} hold and $ m +
      k  = o\{n /\log^2(n) \} $.  If
$\X_i, \U_i \in {\mathbb R}^p$, define ${\mathbb K}
\equiv \{\v\in {\mathbb R}^{\mM\cup S}: \|\v\|_2\leq 1\}$, ${\mathbb K}_1
\equiv \{\v\in {\mathbb R}^{(\mM\cup S)^{c}}: \|\v\|_2= 1, \|\v\|_0
= 1\}$, there are $a_0>0, a_1>0, b_1 >0$  such that \bse 
&&\Pr\left(\sup_{\v, \w \in {\mathbb K}} |\sumi [A(\bb\trans\W_i)
  g(\W_i, \bb, \v, \w) \right.\\
&&\left.- \v\trans E\{\exp(\bb\trans 
  \X_i) \X_i\X_i\trans\} \w ]| > n a_0
  \sqrt{\frac{ m + k}{n}} \right)\\
& \leq& 2 \exp\left\{-  (m + k)\right\}, 
\ese
 \bse 
&&\Pr\left(\sup_{\v \in {\mathbb K}_1, \w \in {\mathbb K}} |\sumi [A(\bb\trans\W_i)
  g(\W_i, \bb, \v, \w) \right.\\
&&\left.- \v\trans E\{\exp(\bb\trans 
  \X_i) \X_i\X_i\trans\} \w ]| > n a_1
  \sqrt{\frac{ \max[m + k, \log(p)]}{n}}  \right)\\
& \leq& 2 \exp[-\max \{\log(p), (m + k)\}], 
\ese
and 
 \bse 
&&\Pr\left(\sup_{\v \in {\mathbb K}_1, \w \in {\mathbb K}_1} |\sumi [A(\bb\trans\W_i)
  g(\W_i, \bb, \v, \w) \right.\\
&&\left.- \v\trans E\{\exp(\bb\trans 
  \X_i) \X_i\X_i\trans\} \w ]| > n b_1
  \sqrt{\frac{ \log(p)}{n}}  \right)\\
& \leq& 2 \exp[-\log(p)].
\ese
\end{Cor}
\noindent Proof: By Lemma \ref{lem:fromlemma9jiang}, take $t = a_0
  \sqrt{ (m + k)/n}$, we have
\bse 
&&\Pr\left(\sup_{\v, \w \in {\mathbb K}} |\sumi [A(\bb\trans\W_i)
  g(\W_i, \bb, \v, \w) - \v\trans E\{\exp(\bb\trans 
  \X_i) \X_i\X_i\trans\} \w ]| > n t \right)\\
& \leq& 2 \exp\left(-  \min \left[ \frac{a _0^2  (m +
      k) }{324 e^2 M_4 },\frac{a_0\sqrt{ n  (m +
        k)}}{36  e M_5 \log (n) } \right] + 2 (m +
  k)\log( 9 )\right)\\
&=& 2 \exp\left(-  \frac{   a _0^2  (m +
      k) }{324 e^2 M_4 }+ 2 (m +
  k)\right)\\
&=& 2 \exp\left\{-  (m +
  k)\right\}
\ese 
The second to the last equality holds because $ m +
      k = o\{n /\log^2(n) \} $. The last equality holds
      by choosing $a _0 = \sqrt{972 e^2 M_4}$. 

Further, by the second relation in  Lemma \ref{lem:fromlemma9jiang},
take $t = a_1 \sqrt{\max\{\log(p), (m+k)\}/n}$, we have 
\bse 
&&\Pr\left(\sup_{\v \in {\mathbb K}_1, \w \in {\mathbb K}} |\sumi [A(\bb\trans\W_i)
  g(\W_i, \bb, \v, \w) - \v\trans E\{\exp(\bb\trans 
  \X_i) \X_i\X_i\trans\} \w ]| > nt \right)\\
& \leq& 2 \exp\left(-  \min \left[ \frac{ a_1^2 \max\{\log(p), (m+k)\} }{36 e^2 M_4 }, \right.\right.\\
&&\left.\left.\frac{a_1\sqrt{  n \max\{\log(p), (m+k)\}}}{12  e M_5
      \log (n) } \right] +   (m + k)\log (9)
  +\log(p) \right)\\
&=&  2 \exp\left\{-  \min \left( \frac{ a_1^2 \max\{\log(p), (m+k)\} }{36 e^2 M_4 }, \right.\right.\\
&&\left. a_1 \max[ \log(p) \sqrt{n }/\{12 e M_5\sqrt{\log (p)} \log(n)
  \}, \sqrt{ n (m+k)} /\{12 e M_5 \log(n)\}  ]\right)\\
&&\left.+  2(m + k)
  +\log(p) \mystrut \right)\\
&\leq& 2 \exp\left(-  \min \left(\frac{ a_1^2 \max\{\log(p), (m+k)\} }{36 e^2 M_4 }, \right.\right.\\
&&\left.a_1  \max[ \log(p) /(12  e M_5 C), \sqrt{ n (m+k)} /\{12 e M_5 \log(n)\}  ]\right) \\
&&\left.+  2(m + k)
  +\log(p) \mystrut \right\}\\
&\leq& 2 \exp\left\{-  \min \left[ \frac{ a_1^2 \max\{\log(p), (m+k)\} }{36 e^2 M_4 }, \right.\right.\\
&&\left.\left. a_1\max\{  \log(p) /(12 e M_5 C), (m+k)/(12  e M_5) \} \right] +  3\max\{(m + k), \log(p)
    \} \mystrut \right\}\\
&\leq& 2 \exp\left(-  a_2  \max\{\log(p), (m+k)\}  +  3\max\{\log(p), (m + k)\}\right)\\
&=& 2 \exp[- \max \{\log(p), (m + k)\}], 
\ese
where $a_2 = \min\{a_1^2/(36 e^2 M_4),  a_1/(12 e M_5 C),
a_1/(12 e M_5)\}$, and we select $a_1$ such that $a_2 \geq 4$. 
The third equality holds because $\log(n) \leq
C\sqrt{n/\log(p)}$ by Condition \ref{con:pn}. The fourth equality holds $ (m +
      k)  = o\{n /\log^2(n) \} $. 

In addition,  by the third relation in  Lemma \ref{lem:fromlemma9jiang},
take $t = b_1 \sqrt{\log(p)/n}$, we have 
\bse 
&&\Pr\left(\sup_{\v \in {\mathbb K}_1, \w \in {\mathbb K}_1} |\sumi [A(\bb\trans\W_i)
  g(\W_i, \bb, \v, \w) - \v\trans E\{\exp(\bb\trans 
  \X_i) \X_i\X_i\trans\} \w ]| > nt \right)\\
& \leq& 2 \exp\left(-  \min \left[ \frac{ b_1^2  \log(p) }{16 e^2 M_4 }, \frac{b_1\sqrt{  n \log(p)}}{8  e M_5
      \log (n) } \right] 
  +2 \log(p) \right)\\
&=&  2 \exp\left\{-  \min \left( \frac{ b_1^2 \log(p)}{16 e^2 M_4 },
 b_1  \log(p) \sqrt{n }/\{8 e M_5\sqrt{\log (p)} \log(n)
  \}\right) 
  +2 \log(p) \mystrut \right)\\
&\leq& 2 \exp\left(-  \min \left(\frac{ b_1^2 \log(p) }{16 e^2 M_4 },
    b_1  \log(p) /(8  e M_5 C)\right) 
  +2 \log(p) \mystrut \right\}\\
&\leq& 2 \exp\left(-  b_{12}\log(p)  +  2 \log(p)\right)\\
&=& 2 \exp[- \log(p)], 
\ese
where $b_{12}  = \min\{b_1^2/(16 e^2 M_4),  b_1/(8 e M_5 C)\}$, and we select $b_1$ such that $b_{12} \geq 3$. 
The second inequality holds because $\log(n) \leq
C\sqrt{n/\log(p)}$ by Condition \ref{con:pn}. 
\qed

\begin{Lem}\label{lem:fromlemma9jiangext}
Assume that Conditions \ref{con:boundX} and \ref{con:ww3} hold.  If
$\X_i, \U_i \in {\mathbb R}^p$, define ${\mathbb K}
\equiv \{\v\in {\mathbb R}^{\mM\cup S}: \|\v\|_2\leq 1\}$, ${\mathbb K}_1
\equiv \{\v\in {\mathbb R}^{(\mM\cup S)^{c}}: \|\v\|_2= 1, \|\v\|_0
= 1\}$. Then \bse 
&&\Pr\left(\sup_{\v, \w \in {\mathbb K}} |\sumi [A^2(\bb\trans\W_i)
  g_1(\W_i, \bb, \v, \w) \right.\\
&&\left.- \v\trans E\{\exp(2 \bb\trans \W_i - \bb\trans
  \bOmega \bb) (\W_i - \bOmega\bb)^{\otimes2}\} \w ] > nt \right)\\
& \leq& 2 \exp\left(-  \min \left[ \frac{ n  t^2}{324 e^2 M_6 },\frac{n
    t}{36  e M_7 \log (n) } \right] +  2(m + k)\log(9)\right), 
\ese
\bse 
&&\Pr\left(\sup_{\v \in {\mathbb K}_1, \w \in {\mathbb K}} |\sumi [A^2(\bb\trans\W_i)
  g_1(\W_i, \bb, \v, \w) \right.\\
&&\left.- \v\trans E\{\exp(2 \bb\trans \W_i - \bb\trans
  \bOmega \bb) (\W_i - \bOmega\bb)^{\otimes2}\} \w ]| > nt \mystrut \right)\\
& \leq& 2 \exp\left\{-  \min \left[ \frac{ n  t^2}{36 e^2 M_6 },\frac{n
    t}{12  e M_7 \log (n) } \right] +  (m + k)\log(9) +
\log(p)\right\}, 
\ese
and 
\bse 
&&\Pr\left(\sup_{\v \in {\mathbb K}_1, \w \in {\mathbb K}_1} |\sumi [A^2(\bb\trans\W_i)
  g_1(\W_i, \bb, \v, \w) \right.\\
&&\left.- \v\trans E\{\exp(2 \bb\trans \W_i - \bb\trans
  \bOmega \bb) (\W_i - \bOmega\bb)^{\otimes2}\} \w ]| > nt \mystrut \right)\\
& \leq& 2 \exp\left\{-  \min \left[ \frac{ n  t^2}{16 e^2 M_6 },\frac{n
    t}{8  e M_7 \log (n) } \right] +  2 \log(p)\right\}.
\ese
\end{Lem}
\noindent Proof: The lemma holds by using the same arguments as those
lead to Lemma \ref{lem:fromlemma9jiang}. \qed

\begin{Cor}\label{cor:fromlemma9jiangext}
Assume that Conditions \ref{con:boundX}, \ref{con:pn} \ref{con:ww3}
hold, $ m +
      k   = o\{n /\log^2(n) \} $.  If
$\X_i, \U_i \in {\mathbb R}^p$, define ${\mathbb K}
\equiv \{\v\in {\mathbb R}^{\mM\cup S}: \|\v\|_2\leq 1\}$, ${\mathbb K}_1
\equiv \{\v\in {\mathbb R}^{(\mM\cup S)^{c}}: \|\v\|_2= 1, \|\v\|_0
= 1\}$. There are $a_2 >0, a_3 >0, b_3 >0 $ such that \bse 
&&\Pr\left(\sup_{\v, \w \in {\mathbb K}} |\sumi [A^2(\bb\trans\W_i)
  g_1(\W_i, \bb, \v, \w) \right.\\
&&\left.- \v\trans E\{\exp(2 \bb\trans \W_i - \bb\trans
  \bOmega \bb) (\W_i - \bOmega\bb)^{\otimes2}\} \w ]| > n a_2
  \sqrt{\frac{ (m + k)}{n}}  \right)\\
& \leq& 2 \exp\left\{-  (m + k)\right\}, 
\ese
\bse 
&&\Pr\left(\sup_{\v \in {\mathbb K}_1, \w \in {\mathbb K}} |\sumi [A^2(\bb\trans\W_i)
  g_1(\W_i, \bb, \v, \w) \right.\\
&&- \v\trans E\{\exp(2 \bb\trans \W_i - \bb\trans
  \bOmega \bb) (\W_i - \bOmega\bb)^{\otimes2}\} \w ]|\\
&&\left. > n  a_3
  \sqrt{\frac{ \max[(m + k), \log(p)]}{n}}  \right)\\
& \leq& 2 \exp\left[- \max \{\log(p), (m + k)\}\right].
\ese
and
\bse 
&&\Pr\left(\sup_{\v \in {\mathbb K}_1, \w \in {\mathbb K}_1} |\sumi [A^2(\bb\trans\W_i)
  g_1(\W_i, \bb, \v, \w) \right.\\
&&- \v\trans E\{\exp(2 \bb\trans \W_i - \bb\trans
  \bOmega \bb) (\W_i - \bOmega\bb)^{\otimes2}\} \w ]|\\
&&\left. > n  b_3
  \sqrt{\frac{ \log(p)}{n}}  \right)
\leq 2 \exp\{- \log(p)\}.
\ese
\end{Cor}
\noindent Proof: The corollary follows the same arguments as those
lead to Corollary \ref{cor:fromlemma9jiang}.\qed

\begin{Cor}\label{cor:fromlemma9jiangext1}
Assume that Conditions \ref{con:boundX}, \ref{con:pn} \ref{con:ww3}
hold, $ m +
      k   = o\{n /\log^2(n) \} $.  If
$\X_i, \U_i \in {\mathbb R}^p$, define ${\mathbb K}
\equiv \{\v\in {\mathbb R}^{\mM\cup S}: \|\v\|_2\leq 1\}$, ${\mathbb K}_1
\equiv \{\v\in {\mathbb R}^{(\mM\cup S)^{c}}: \|\v\|_2= 1, \|\v\|_0
= 1\}$. There are $a_{21} >0, a_{31} >0, b_{31} >0 $ such that \bse 
&&\Pr\left(\sup_{\v, \w \in {\mathbb K}} |\sumi [A(\bb\trans\W_i)
  g_1(\W_i, \bb, \v, \w) \right.\\
&&\left.- \v\trans E\{\exp(\bb\trans \W_i - \bb\trans
  \bOmega \bb/2) (\W_i - \bOmega\bb)^{\otimes2}\} \w ]| > n a_{21}
  \sqrt{\frac{ (m + k)}{n}}  \right)\\
& \leq& 2 \exp\left\{-  (m + k)\right\}, 
\ese
\bse 
&&\Pr\left(\sup_{\v \in {\mathbb K}_1, \w \in {\mathbb K}} |\sumi [A(\bb\trans\W_i)
  g_1(\W_i, \bb, \v, \w) \right.\\
&&- \v\trans E\{\exp(\bb\trans \W_i - \bb\trans
  \bOmega \bb/2) (\W_i - \bOmega\bb)^{\otimes2}\} \w ]|\\
&&\left. > n  a_{31}
  \sqrt{\frac{ \max[(m + k), \log(p)]}{n}}  \right)\\
& \leq& 2 \exp\left[- \max \{\log(p), (m + k)\}\right], 
\ese
and 
\bse 
&&\Pr\left(\sup_{\v \in {\mathbb K}_1, \w \in {\mathbb K}_1} |\sumi [A(\bb\trans\W_i)
  g_1(\W_i, \bb, \v, \w) \right.\\
&&- \v\trans E\{\exp(\bb\trans \W_i - \bb\trans
  \bOmega \bb/2) (\W_i - \bOmega\bb)^{\otimes2}\} \w ]|\\
&&\left. > n  b_{31}
  \sqrt{\frac{ \log(p)}{n}}  \right) \leq 2 \exp\{-\log(p)\}.
\ese
\end{Cor}
\noindent Proof: The corollary follows the same arguments as those that
lead to Corollary \ref{cor:fromlemma9jiangext}.\qed

\begin{Lem}\label{lem:fromlemma15loh2012}
Assume that Conditions \ref{con:boundX} holds, $ m +
      k   = o\{n /\log^2(n) \} $.  If
$\X_i, \U_i \in {\mathbb R}^p$, ${\mathbb K}
\equiv \{\v\in {\mathbb R}^{\mM\cup S}: \|\v\|_2= 1\}$, ${\mathbb K}_1
\equiv \{\v\in {\mathbb R}^{(\mM\cup S)^{c}}: \|\v\|_2= 1, \|\v\|_0
= 1\}$. Then there exists constants $c_2, c_3, c_4, c_5, c_6, c_7,
c_8, c_9, c_{10} >0 $\bse 
&&\Pr\left(\sup_{\v, \w\in {\mathbb K}} |\sumi  \{\v\trans(\W_i -
  \bOmega\bb) ^{\otimes2} \w- E( \v\trans(\W_i -
  \bOmega\bb) ^{\otimes2} \w )\} | > nt \right)\\
& \leq& 2 \exp\left(-  c_2n \min (c_3 t^2, c_4 t)  + 2 (m + k)\log(9)\right), 
\ese
\bse 
&&\Pr\left(\sup_{\v \in {\mathbb K}_1, \w\in {\mathbb K}} |\sumi \{ \v\trans(\W_i -
  \bOmega\bb) ^{\otimes2} \w- E( \v\trans(\W_i -
  \bOmega\bb) ^{\otimes2} \w )\}|> nt \right)\\
& \leq& 2 \exp\left(-  c_7n \min (c_5 t^2, c_6 t)  + (m + k)\log(9)+
  \log(p)\right), 
\ese
and 
\bse 
&&\Pr\left(\sup_{\v \in {\mathbb K}_1, \w\in {\mathbb K}_1} |\sumi \{ \v\trans(\W_i -
  \bOmega\bb) ^{\otimes2} \w- E( \v\trans(\W_i -
  \bOmega\bb) ^{\otimes2} \w )\}|> nt \right)\\
& \leq& 2 \exp\left(-  c_{10}n \min (c_8 t^2, c_9 t)  + 2\log(p)\right).
\ese
\end{Lem}
\noindent Proof: The lemma is a  consequence of Lemma 15 in \cite{loh2012} and using the same arguments as those
lead to Lemma \ref{lem:fromlemma9jiang}.\qed

\begin{Cor}\label{cor:fromlemma15loh2012}
Assume that Conditions \ref{con:boundX} holds.  If
$\X_i, \U_i \in {\mathbb R}^p$, ${\mathbb K}
\equiv \{\v\in {\mathbb R}^{\mM\cup S}: \|\v\|_2\leq 1\}$, ${\mathbb K}_1
\equiv \{\v\in {\mathbb R}^{(\mM\cup S)^{c}}: \|\v\|_2= 1, \|\v\|_0
= 1\}$. There are constants $a_4>0, a_5>0, b_5 >0 $ such that \bse 
&&\Pr\left(\sup_{\v, \w\in {\mathbb K}} |\sumi  \{\v\trans(\W_i -
  \bOmega\bb) ^{\otimes2} \w- E( \v\trans(\W_i -
  \bOmega\bb) ^{\otimes2} \w )\} |\right.\\
&&\left. > n  a_4
  \sqrt{\frac{ (m + k)}{n}} \right) \leq 2 \exp\left\{-  (m + k)\right\}, 
\ese
\bse 
&&\Pr\left(\sup_{\v \in {\mathbb K}_1, \w\in {\mathbb K}} |\sumi \{ \v\trans(\W_i -
  \bOmega\bb) ^{\otimes2} \w- E( \v\trans(\W_i -
  \bOmega\bb) ^{\otimes2} \w )\} |\right.\\
&&\left.> n a_5
  \sqrt{\frac{ \max[(m + k), \log(p)]}{n}}  \right)\leq 2 \exp[- \max
\{\log(p), (m + k)\}], 
\ese
and 
\bse 
&&\Pr\left(\sup_{\v \in {\mathbb K}_1, \w\in {\mathbb K}_1} |\sumi \{ \v\trans(\W_i -
  \bOmega\bb) ^{\otimes2} \w- E( \v\trans(\W_i -
  \bOmega\bb) ^{\otimes2} \w )\} |\right.\\
&&\left.> n b_5
  \sqrt{\frac{ \log(p)}{n}}  \right)\leq 2 \exp\{-\log(p)\}.
\ese
\end{Cor}
\noindent Proof: The corollary follows the same arguments as those
lead to Corollary \ref{cor:fromlemma9jiang}.\qed

\section{Lemmas on Criterion Function and Penalty Function}
\begin{Lem}\label{lem:fromlemma1}
Consider a $\mu$-amenable regularizer $\rho_{\lambda}$, with $\mu <
\alpha_1$ and  $n>\log (p) \tau_1^2 (m + k)^2/(\alpha_1-\mu)^2$,
where $\alpha_1, 
\tau_1$ are defined in Lemma \ref{lem:RE1}. Then the function $\mL(\bb) -
\mu \|\bb_{\mM^c}\|_2^2/2$ and $\mL(\bb) + 
\rho_\lambda(\bb_{\mM^c})$ are strictly convex on
$\bb\in \mathbb{R}^{\mM \cup S}$, and hence the restricted program
(\ref{eq:lapHTsub}) is also strictly convex. 
\end{Lem}
\noindent Proof: 
First define a vector $\v\in{\mathbb R}^p$ with the $j$th element
$|v_j| >0$ if $j \in \mM \cup S$, and $\|v_j\| =0$, otherwise.   By Lemma \ref{lem:RE1}, we have for $\bb$ in the feasible set of
program (\ref{eq:lapHTsub}), 
\bse
\v\trans \frac{\partial ^2 \mL(\bb)}{\partial \bb\partial \bb\trans }\v
\geq \alpha_1 \|\v\|_2^2 - \tau_1 \sqrt{\frac{\log(p)}{n}} \|\v\|_1^2, 
\ese
and   $\|\v\|_1 \leq \sqrt{m+k} \|\v\|_2$, and hence we
have 
\bse
\v\trans \frac{\partial ^2 \mL(\bb)}{\partial \bb\partial \bb\trans }\v
\geq \left\{\alpha_1 - \tau_1 (m + k)\sqrt{\frac{\log(p)}{n}} \right\}\|\v\|_2^2.
\ese
Therefore, 
\bse
&&\v_{\mM \cup S}\trans \left\{\frac{\partial ^2 \mL(\bb)}{\partial \bb\partial \bb\trans }
\right\}_{(\mM + S) (\mM + S)}\v_{\mM \cup S} -\mu \v_S\trans\v_S\\
&\geq& \left\{\alpha_1 - \mu - \tau_1 (m + k)\sqrt{\frac{\log(p)}{n}} \right\}\|\v\|_2^2, 
\ese
where $\left\{ \partial ^2 \mL(\bb)/(\partial \bb\partial \bb\trans ) 
\right\}_{(\mM \cup S), (\mM\cup S)}$ is the $(m+k)
\times (m+ k)$ block of $\left\{ \partial ^2 \mL(\bb)/\partial \bb\partial \bb\trans 
\right\}$ corresponding to $\mM\cup S$.
Hence $\mL(\bb) - \mu \|\bb_{\mM^c}\|_2^2/2$ is strictly convex on
$\mathbb{R}^{\mM\cup S}$. 

Finally, because
\bse
\mL(\bb) - q_{\lambda}(\bb_{\mM^c}) = (\mL(\bb) - \mu\|\bb_{\mM^c}\|_2^2/2)+
\{\mu\|\bb_{\mM^c}\|_2^2/2 - q_{\lambda}(\bb_{\mM^c}) \},
\ese
where the second part is convex over $\mathbb{R}^{\mM\cup S}$ by Lemma \ref{lem:fromlemm5support}, 
hence $\mL(\bb) - q_{\lambda}(\bb_{\mM^c})$ restricted to
$\mathbb{R}^{\mM \cup S}$ is strictly
convex. Because $\mL(\bb) + 
\rho_{\lambda} (\bb_{\mM^c}) = \mL(\bb) - q_{\lambda}(\bb_{\mM^c}) +
\lambda \|\bb_{\mM^c}\|_1$,  the strict convexity of $\mL(\bb) + 
\rho_{\lambda} (\bb_{\mM^c})$ over $\mathbb{R}^{\mM\cup S}$ follows. This proves the
result. \qed

Now as we know $\C$ contains $r$ independent columns, without loss of
generality, we write $\C = (\C_r, \C_{m -r})$ where $\C_r$ is a full
rank square matrix.

\begin{Lem}\label{lem:fromlemma11loh2017}
Let $\A, \B \in \mathbb{R}^{p\times p}$ be invertible. For any matrix
norm $\|\|$, we have 
\bse
\|\A^{-1} - \B^{-1}\| \leq \frac{\|\A^{-1}\|^2 \|\A - \B\|}{1
  -\|\A^{-1}\| \|\A - \B\| }. 
\ese
In particular, if $\|\A^{-1}\|\|\A - \B\| \leq 1/2$, then $\|\A^{-1} -
\B^{-1}\| = O(\|\A^{-1}\|^2\|\A - \B\|)$. 
\end{Lem}
\noindent Proof: This lemma is  Lemma 11 in
\cite{loh2017}. \qed
\begin{Lem}\label{lem:fromlemma10}
Suppose $\x^*$ is feasible for the program
\be\label{eq:G1}
\min_{\x}\{f(\x) - g(\x_{\mM^c}) + \lambda \|\x_{\mM^c}\|_1\}, \text{  such that }
\|\x\|_1\leq R_1, \|\x\|_2\leq R_2, \text{ and } \C\x_{\mM} = \t,
\ee
where $f \in C^2$, $g \in C^1$ and $g(\x_{\mM^c}) -
\kappa/2\|\x_{\mM^c}\|_2^2$ is concave and $\C$ is an $r\times m$
matrix. Define $\A = \{\0_{(p-m)\times m}, \I_{(p-m)\times (p-m)}\}$ and
$\A_1 = \{\I_{m\times m}, \0_{m\times (p-m)}\}$. Assume there are
$\v^* \in \partial \|\x^*_{\mM^c}\|_1$, $\w_1^*
\in \partial \|\x^*\|_1$, $\w_2^* \in \partial \|\x^*\|_2$,
$\mu_1^*\geq 0$, $\mu_2^*\geq 0$, $\bmu_3^*$
such that 
\be
\mu_1^* (R_1 - \|\x^*\|_1) = 0\label{eq:x1}\\
\mu_2^* (R_2 - \|\x^*\|_2) = 0\label{eq:x2}\\
\frac{\partial f(\x^*)}{\partial \x^*} - \A\trans \frac{\partial
  g(\x^*_{\mM^c})}{\partial \x^*_{\mM^c}} + \lambda \A\trans \v^* + \mu_1^* \w_1^*
+
\mu_2^*\w_2^* + \A_1\trans  \C \trans \bmu_3^*= 0 \label{eq:x1x2}\\
\s\trans \frac{\partial^2 f(\x)}{\partial \x\partial \x\trans}  \s >\kappa, 
\forall \s \in G^*, \label{eq:G2c}
\ee
where 
\bse
G^* &:=& \left\{\s\in{\mathbb R}^p: \|\s\|_2 = 1; \sup_{\w_1 \in \partial \|\x^*\|_1} \s\trans
\w_1 \leq 0 \text{ if } \|\x^*\|_1  = R_1;\sup_{\w_2 \in \partial \|\x^*\|_2} \s\trans
\w_2 \leq 0 \text{ if } \|\x^*\|_2 = R_2; \right.\\
&& \sup_{\v \in \partial
  \|\x^*_{\mM^c}\|_1}  \s\trans \left(\frac{\partial f(\x^*)}{\partial \x^*} -\A\trans \frac{\partial
g(\x_{\mM^c}^*)}{\partial \x_{\mM^c}^* } \right)  + \lambda \s\trans\A\trans\v = {
0 };\\
&&\left.
\mu_1^* \sup_{\w_1 \in \partial \|\x^*\|_1} \s\trans \w_1 = 0, \mu_2^*
\sup_{\w_2 \in \partial \|\x^*\|_2} \s\trans \w_2 = 0; \C \s_{\mM} =
\bf 0\right\}.
\ese
Then $\x^*$ is an isolated local minimum of the program (\ref{eq:G1}). 
\end{Lem}
\noindent Proof: The proof of this lemma is similar to the proof
of Theorem 3 in \cite{fletcher1980} and that of Lemma 10 in \cite{loh2017},
except that we allow additional constraints $\|\x\|_2\leq R_2$ and
$\C\x_{\mM} = \t$. 

Suppose $\x^*$ is not an isolated local minimum. Then there is a
sequence $(\x^{k})$, so that $\x^k\to \x^*$ and
\bse
\phi(\x^k) \leq \phi(\x^*), 
\ese
where $\phi(\x) = f(\x) - g(\x_{\mM^c}) + \lambda
\|\x_{\mM^c}\|_1$. Let $\s^k := (\x^k - \x^*)/\|\x^k - \x^*\|_2$, so
$(\s^k)$ is a set of feasible directions. Since $(\s^k) \subset
\mathbb{B}_2(1)$, where $\mathbb{B}_2(1)$ is the ball with radius
1, the set must possess a point of accumulation $\s \in
\mathbb{B}_2(1)$, and we can extract a convergence subsequence such
that $(\s^k)
\to \s$. With a slight abuse of notation,  we still use $(\s^k)$ to denote
the subsequence. We will show that $\s \in  G^*$. 

First of all by the construction, $\x^k$'s are all feasible, and hence
$\C\x^k_{\mM} = \t$. Therefore, $\C\s^k_{\mM} = 0$, take the limits on
the left and right of the equation  implies
\be\label{eq:s0}
\C\s _{\mM} = \bf 0.
\ee

Further, because the feasible region is closed, $\s$ is also in the
feasible direction at $\x^*$. If $\|\x^*\|_1 = R_1$, by the
sub-gradient of convexity function $\|\x^*\|_1$ we have 
\bse
0 \geq \|\x^k\|_1 - \|\x^*\|_1 = \|\x^* +\|\x^k - \x^*\|_2 \s^k \|_1 -
\|\x^*\|_1 
\geq \|\x^k - \x^*\|_2  {\s^k}\trans \w_1,
\ese
for any $\w_1 \in \partial \|\x^*\|_1$. When $k\to\infty$, 
this leads to 
\be\label{eq:G31}
\sup_{\w_1 \in \partial \|\x^*\|_1} \s\trans\w_1\leq 0. 
\ee 
Further (\ref{eq:x1}) also implies that if $\|\x^*\|_1\neq R_1 $, then
$\mu_1^*  =
0$. Since $\mu_1^*\ge 0$, hence we have 
\be\label{eq:mu1}
\mu_1^* \sup_{\w_1 \in \partial \|\x^*\|_1} \s\trans\w_1 \leq 0. 
\ee 

Similarly, by the
sub-gradient of convexity function $\|\x^*\|_2$ we have 
\bse
0 \geq \|\x^k\|_2 - \|\x^*\|_2= \|\x^* +\|\x^k - \x^*\|_2 \s^k \|_2 -
\|\x^*\|_2 
\geq \|\x^k - \x^*\|_2  {\s^k}\trans \w_2,
\ese
for any $\w_2 \in \partial \|\x^*\|_2$. When $k\to\infty$, 
this leads to 
\be\label{eq:G32}
\sup_{\w_2 \in \partial \|\x^*\|_2} \s\trans\w_2\leq 0. 
\ee 
Further (\ref{eq:x2}) also implies that if $\|\x^*\|_2\neq R_2 $, then
$\mu_2^*  =
0$. Since $\mu_2^*\ge 0$, hence we have 
\be\label{eq:mu2}
\mu_2^* \sup_{\w_2 \in \partial \|\x^*\|_2} \s\trans\w_2 \leq 0. 
\ee 

Further by (\ref{eq:x1x2}), we have 
\bse
\s\trans \frac{\partial f(\x^*)}{\partial \x^*} - \s\trans \A\trans \frac{\partial
  g(\x^*_{\mM^c})}{\partial \x^*_{\mM^c}} + \s\trans \lambda \A\trans  \v^* +
\mu_1^* \s\trans \w_1^*
+
\mu_2^*\s\trans \w_2^* + \s\trans \A_1\trans  \C \trans \bmu_3^{*\trans}= 0
\ese
which, together with (\ref{eq:s0}),  implies 
\be\label{eq:G5}
\s\trans \frac{\partial f(\x^*)}{\partial \x^*} - \s\trans \A\trans \frac{\partial
  g(\x^*_{\mM^c})}{\partial \x^*_{\mM^c}} + \s\trans \lambda \A\trans
\v^*
=-
\mu_1^* \s\trans \w_1^*-
\mu_2^*\s\trans \w_2^*\geq 0. 
\ee 
By the definition of sub-gradient, we have 
\bse
\|\x_{\mM^c}^* +\|\x^k -\x^*\|_2 \A\s^k  \|_1 -
\|\x_{\mM^c}^*\|_1 \geq \|\x^k -\x^*\|_2 \s^{k\rm T}
\A\trans \v
\ese
for all $\v \in \partial \|\x_{\mM^c}^*\|_1$ and all $k$. Further
because $\|\x_{\mM^c}^* +\|\x^k -\x^*\|_2 \A\s^k  \|_1 -
\|\x_{\mM^c}^*\|_1= \|\x_{\mM^c}^k\|_1 -\|\x_{\mM^c}^*\|_1 $, we have 
\be\label{eq:G6}
\s\trans \A\trans \v = \lim_{k \to \infty }\s^{k\rm T}\A \trans\v \leq  \lim_{k \to \infty }\frac{\|\x_{\mM^c}^k\|_1 -\|\x_{\mM^c}^*\|_1 }{\|\x^k -\x^*\|_2}
\ee 
for all $\v \in \partial \|\x_{\mM^c}^*\|_1$.
Furthermore, 
\be\label{eq:G7}
&&\s\trans \frac{\partial f(\x^*)}{\partial \x^*} - \s\trans \A\trans \frac{\partial
  g(\x^*_{\mM^c})}{\partial \x^*_{\mM^c}}\nonumber\\
& =& \lim_{k\to \infty}\s^{k\rm T} \frac{\partial f(\x^*)}{\partial \x^*} - \s^{k\rm T} \A\trans \frac{\partial
  g(\x^*_{\mM^c})}{\partial \x^*_{\mM^c}}\nonumber\\
&=&\lim_{k\to \infty}\frac{\langle\x^k - \x^*,\frac{\partial
    f(\x^*)}{\partial \x^*} - \A\trans \frac{\partial
  g(\x^*_{\mM^c})}{\partial \x^*_{\mM^c}}\rangle}{\|\x^k -\x^*\|_2}\nonumber\\
&=& \lim_{k \to \infty }\frac{f(\x^{k}) - g(\x_{\mM^c}^k) - f(\x^{*})
  +  g(\x_{\mM^c}^*)}{\|\x^k -\x^*\|_2}, 
\ee
since $\x^k \to \x^*$ and $f(\x^{k}) - g(\x_{\mM^c}^k) \in
C^1$. Combining (\ref{eq:G6}) and (\ref{eq:G7}), we conclude that 
\bse
&& \sup_{\v \in \partial \|\x_{\mM^c}^*\|_1} \s\trans \frac{\partial f(\x^*)}{\partial \x^*} - \s\trans \A\trans \frac{\partial
  g(\x^*_{\mM^c})}{\partial \x^*_{\mM^c}} + \lambda \s\trans \A\trans
\v\\
&\leq&\lim_{k\to \infty}\frac{\phi(\x^k) - \phi(\x^*)}{\|\x^k - \x*\|}
\leq 0.
\ese 
Combining with (\ref{eq:G5}), we have 
\bse
0&\le& \s\trans \frac{\partial f(\x^*)}{\partial \x^*} - \s\trans \A\trans \frac{\partial
  g(\x^*_{\mM^c})}{\partial \x^*_{\mM^c}}  + \lambda \s\trans \A\trans
\v^*\nonumber\\
&\le& \sup_{\v \in \partial \|\x^*_{\mM^c}\|_1} \s\trans \frac{\partial f(\x^*)}{\partial \x^*} - \s\trans \A\trans \frac{\partial
  g(\x^*_{\mM^c})}{\partial \x^*_{\mM^c}}  + \lambda \s\trans \A\trans
\v \nonumber \\
&\le& 0,
\ese
hence
\be\label{eq:G8}
\sup_{\v \in \partial \|\x^*_{\mM^c}\|_1} \s\trans \frac{\partial f(\x^*)}{\partial \x^*} - \s\trans \A\trans \frac{\partial
  g(\x^*_{\mM^c})}{\partial \x^*_{\mM^c}}  + \lambda \s\trans \A\trans
\v=0
\ee
Now together with (\ref{eq:x1x2}), we have 
\bse
\mu_1^*\s\trans \w_1^* + \mu_2^*\s\trans \w_2^*  = 0. 
\ese 
Further by (\ref{eq:mu1}) and (\ref{eq:mu2}), we have 
\be
\mu_1^*\sup_{\w_1 \in \partial \|\x^*\|_1}\s\trans \w_1 =
\mu_2^*\sup_{\w_2 \in \partial \|\x^*\|_2}\s\trans \w_2 =
0. \label{eq:G9}
\ee 
Combine (\ref{eq:s0}), (\ref{eq:G31}), (\ref{eq:G32}), (\ref{eq:G8})
and (\ref{eq:G9}), we conclude that $\s \in G^*$. 

By the convexity of the $L_1$ norm, we have 
\bse
\|\x^k\|_1 - \|\x^*\|_1 =  \|\x^* + (\x^k - \x^*) \|_1 - \|\x^*\|_1
\geq (\x^k - \x^*)\trans \w_1 =  \x^{k\rm T} \w_1- \|\x^*\|_1
\ese
for all $\w_1 \in \partial \|\x^* \|_1$. Therefore, $\x^{k\rm T} \w_1^* \leq
\|\x^k\|_1\leq R_1$. Similarly,
by the convexity of the $L_2$ norm, we have 
\bse
\|\x^k\|_2 - \|\x^*\|_2 =  \|\x^* + (\x^k - \x^*) \|_2 - \|\x^*\|_2
\geq (\x^k - \x^*)\trans \w_2 =  \x^{k\rm T} \w_2- \|\x^*\|_2
\ese
for all $\w_2 \in \partial \|\x^* \|_2$. Therefore, $\x^{k\rm T} \w_2^* \leq
\|\x^k\|_2\leq R_2$. 
Further, $\x^{k \rm T}  \A_1\trans  \C \trans
\bmu_3^* - \t\trans \bmu_3^* = 0$. 
Hence 
\bse
\phi(\x^k) &=&  f(\x^k) - g(\x_{\mM^c}^k) + \lambda\|\x_{\mM^c}^k\|_1 \\
&\geq& f(\x^k) - g(\x_{\mM^c}^k) + \lambda \x_{\mM^c}^{k\rm T} \v^* +
\mu_1^* (\x^{k\rm T} \w_1^* - R_1) + \mu_2^* (\x^{k\rm T} \w_2^* -
R_2) \\
&&+ \x^{k\rm T} \A_1\trans  \C \trans \bmu_3^* - \t\trans \bmu_3^* 
\ese
for all $\v^*\in\partial\|\x_{\mM^c}^*\|_1$,
and 
\bse
\phi(\x^*) &=&  f(\x^*) - g(\x_{\mM^c}^* ) + \lambda
\x_{\mM^c}^{* \rm T} \v^* + \mu_1^* (\x^{*\rm T} \w_1^* - R_1) +
\mu_2^* (\x^{* \rm T} \w_2^* - R_2) \\
&&+ \x^{* \rm T}\A_1\trans  \C \trans \bmu_3^* - \t\trans \bmu_3^*. 
\ese
The equality holds because $\mu_1^* (\x^{*\rm T} \w_1^* - R_1)  =0$
and $\mu_2^* (\x^{* \rm T} \w_2^* - R_2) = 0$. 
Hence 
\be\label{eq:G10}
&&\lim_{k\to\infty}\{\phi(\x^k)  - \phi(\x^* )\}/\|\x^k-\x^*\|_2^2 \nonumber\\
&\geq&\lim_{k\to\infty} \{f(\x^k)  - f(\x^* )- g(\x_{\mM^c}^k) + g(\x_{\mM^c}^*) + \langle\lambda\A\trans
\v^* + \mu_1^* \w_1^* + \mu_2^*\w_2^* \nonumber\\
&&+ \A_1\trans \C\trans \bmu_3^*, \x^k - \x^*\rangle \}/\|\x^k-\x^*\|_2^2\nonumber\\
&=&\lim_{k\to\infty} \{ f(\x^k)  - f(\x^* )- g(\x_{\mM^c}^k) +
g(\x_{\mM^c}^*) \}/\|\x^k-\x^*\|_2^2 \n\\
&&- \left\langle \frac{\partial f(\x^*)}{\partial \x^*} - \A\trans \frac{\partial
  g(\x^*_{\mM^c})}{\partial \x^*_{\mM^c}}, \x^k - \x^*\right\rangle /\|\x^k-\x^*\|_2^2 \nonumber\\
&=& \lim_{k\to\infty} \left\{f(\x^k)  - f(\x^* ) - \left \langle \frac{\partial f(\x^*)}{\partial
  \x^*}, \x^k - \x^*\right\rangle\right\}/\|\x^k-\x^*\|_2^2\nonumber\\
&& - \left\{g(\x_{\mM^c}^k)  - g(\x_{\mM^c}^*)  -\left\langle \A\trans \frac{\partial
  g(\x^*_{\mM^c})}{\partial \x^*_{\mM^c}}, \x^k - \x^*\right\rangle
\right\}/\|\x^k-\x^*\|_2^2.
\ee
By the concavity of $g(\x_{\mM^c}) - \kappa/2 \|\x_{\mM^c}\|^2_2$, we have 
\bse
\left\{g(\x_{\mM^c}^k)  - g(\x_{\mM^c}^*)  -\left\langle \A\trans \frac{\partial
  g(\x^*_{\mM^c})}{\partial \x^*_{\mM^c}}, \x^k - \x^*\right\rangle
\right\} \leq  \kappa/2 \|\x^k_{\mM^c} - \x^*_{\mM^c}  \|_2^2.
\ese
Further note that $\phi(\x^k) - \phi(\x^*) \leq 0$. Combine with
(\ref{eq:G10}), we have 
\bse
 &&\lim_{k\to \infty }\left\{f(\x^k)  - f(\x^* ) - \left \langle \frac{\partial f(\x^*)}{\partial
  \x^*}, \x^k - \x^*\right\rangle\right\} /\|\x^k-\x^*\|_2^2\\
&&- \kappa/2 \|\x^k_{\mM^c} -
\x^*_{\mM^c}  \|_2^2 /\|\x^k-\x^*\|_2^2 \leq 0, 
\ese
which by Taylor expansion implies 
\bse
\s\trans  \frac{\partial^2 f(\x^*)}{\partial \x\partial \x\trans }\s 
&=&\lim_{k \to \infty }(\x^k -
\x^* ) \trans \frac{\partial^2 f(\x^{*})}{\partial\x\partial\x\trans }(\x^k -
\x^*) /\|\x^k-\x^*\|_2^2\\
&\leq& \lim_{k \to \infty} \{\kappa \|\x^k_{\mM^c} -
\x^*_{\mM^c}  \|_2^2  + o(\|(\x^k_{\mM^c} -
\x^*_{\mM^c} )\|_2^2)\}/\|\x^k-\x^*\|_2^2\\
&\leq&\kappa,
\ese
which contradicts with (\ref{eq:G2c}). Hence $\x^*$ must be an isolated
local minimum. \qed

\begin{Cor}\label{cor:lem:fromlemma10}
Suppose $\x^*$ is feasible for the program
\be\label{eq:corG1}
\min_{\x}\{f(\x) - g(\x_{\mM^c}) + \lambda \|\x_{\mM^c}\|_1\}, \text{  such that }
\|\x\|_1\leq R_1, \|\x\|_2\leq R_2, 
\ee
where $f \in C^2$, $g \in C^1$ and $g(\x_{\mM^c}) -
\kappa/2\|\x_{\mM^c}\|_2^2$ is concave and $\C$ is an $r\times m$
matrix. Define $\A = \{\0_{(p-m)\times m}, \I_{(p-m)\times (p-m)}\}$ and
$\A_1 = \{\I_{m\times m}, \0_{m\times (p-m)}\}$. Assume there are
$\v^* \in \partial \|\x^*_{\mM^c}\|_1$, $\w_1^*
\in \partial \|\x^*\|_1$, $\w_2^* \in \partial \|\x^*\|_2$,  $\mu_1^*\geq 0$, $\mu_2^*\geq 0$
such that 
\bse
\mu_1^* (R_1 - \|\x^*\|_1) = 0\label{eq:corx1}\\
\mu_2^* (R_2 - \|\x^*\|_2) = 0\label{eq:corx2}\\
\frac{\partial f(\x^*)}{\partial \x^*} - \A\trans \frac{\partial
  g(\x^*_{\mM^c})}{\partial \x^*_{\mM^c}} + \lambda \A\trans \v^* + \mu_1^* \w_1^*
+
\mu_2^*\w_2^* = 0 \label{eq:corx1x2}\\
\s\trans \frac{\partial^2 f(\x)}{\partial \x\partial \x\trans}  \s >\kappa, 
\forall \s \in G^*, \label{eq:corG2c}
\ese
where 
\bse
G^* &:=& \left\{\s\in{\mathbb R}^p: \|\s\|_2 = 1; \sup_{\w_1 \in \partial \|\x^*\|_1} \s\trans
\w_1 \leq 0 \text{ if } \|\x^*\|_1  = R_1;\sup_{\w_2 \in \partial \|\x^*\|_2} \s\trans
\w_2 \leq 0 \text{ if } \|\x^*\|_2 = R_2; \right.\\
&& \sup_{\v \in \partial
  \|\x^*_{\mM^c}\|_1}  \s\trans \left(\frac{\partial f(\x^*)}{\partial \x^*} -\A\trans \frac{\partial
g(\x_{\mM^c}^*)}{\partial \x_{\mM^c}^* } \right)  + \lambda \s\trans\A\trans\v = {
0 };\\
&&\left.
\mu_1^* \sup_{\w_1 \in \partial \|\x^*\|_1} \s\trans \w_1 = 0, \mu_2^*
\sup_{\w_2 \in \partial \|\x^*\|_2} \s\trans \w_2 = 0\right\}.
\ese
Then $\x^*$ is an isolated local minimum of the program (\ref{eq:corG1}). 
\end{Cor}
\noindent Proof: The Corollary holds by using the same argument as
those lead to Lemma \ref{lem:fromlemma10}, while ignoring the equality
constraint.\qed

\begin{Lem}\label{lem:fromlem3}
Assume  Conditions \ref{con:boundX}--\ref{con:ww2} hold,   $ \lambda \leq \alpha_1/ (8 R_1)$, $R_1 \leq [n \alpha_1^2/\{64 \tau_1^2 \log
(p)\}]^{1/4}$, $\delta \in [4 R_1 \tau_1 \sqrt{\log (p)/n} /\lambda,
1]$, $n \geq 4 \log (p) \tau_1^2 \{\sqrt{r}\|\C_r^{-1}\C_{m-r}\|_2 + \sqrt{m
  -r} + (2/ \delta + 1/2) \sqrt{k} \}^4/(\alpha_1-\mu)^2$, and
$\alpha_1 >\mu$. Define $\A = \{\0_{(p-m)\times m}, \I_{(p-m)\times (p-m)}\}$ and
$\A_1 = \{\I_{m\times m}, \0_{m\times (p-m)}\}$.  Suppose 
$\wt{\bb}$ is a stationary point of program (\ref{eq:lapHT}) and
$\wh{\bb}$ is the interior local minimizer of (\ref{eq:lapHT}) and satisfies
$\supp (\wh{\bb})  \subseteq \mM \cup S$. Then $\supp
(\wt{\bb}) \subseteq \mM \cup S$. 
\end{Lem}
\noindent Proof: Let $\wt{\v} := \wt{\bb} - \wh{\bb}$, by the Taylor
expansion of the first order derivative
\bse
\{\partial \mL(\wt \bb)/\partial \bb\trans - \partial \mL(\wh {\bb})/\partial \bb\trans\}\wt{\v}
=  \wt \v \trans \partial ^2 \mL(\bb^*)/\partial \bb\partial \bb\trans
 \wt \v, 
\ese
where $\bb^* $ is a point on the line connecting $\wh{\bb}$ and
$\wt{\bb}$ and hence in the feasible set. Therefore, by Lemma
\ref{lem:RE1} we have 
\bse
\{\partial \mL(\wt \bb)/\partial \bb\trans - \partial
\mL(\wh{\bb})/\partial \bb\trans\}\wt{\v} \geq \alpha_1\|\wt \v\|_2^2 -
\tau_1\sqrt{\log(p)/n} \|\wt \v\|_1^2.
\ese  
We first show
that $\|\wt {\v}\|_2\leq 1$. Suppose that $\|\wt{\v}\|_2 >1$, we have 
\be\label{eq:RElem9}
\{\partial \mL(\wt \bb)/\partial \bb\trans - \partial
\mL(\wh{\bb})/\partial \bb\trans\}\wt{\v}  \geq \alpha_1\|\wt \v\|_2 -2 
\tau_1 R_1 \sqrt{\log(p)/n} \|\wt \v\|_1. 
\ee
The first order condition implies 
\be\label{eq:lem9first}
\left\{\partial \mL(\wt \bb)/\partial \bb+ \A\trans \partial
\rho_\lambda(\wt \bb_{\mM^c})/\partial \bb_{\mM^c} \right\}\trans
(\wh{\bb} - \wt{\bb})\geq 0
\ee
and hence 
\bse
\partial \mL(\wt \bb)/\partial \bb\trans \wt{\v} \leq - \A\trans \partial
\rho_\lambda(\wt \bb_{\mM^c})/\partial \bb_{\mM^c} \trans\wt{\v}.
\ese 
Combine with (\ref{eq:RElem9}), we have 
\be \label{eq:RElem91}
\{- \A\trans \partial
\rho_\lambda(\wt \bb_{\mM^c})/\partial \bb_{\mM^c}   - \partial
\mL(\wh{\bb})/\partial \bb\}\trans\wt{\v}  \geq \alpha_1\|\wt \v\|_2 -2 
\tau_1R_1 \sqrt{\log(p)/n} \|\wt \v\|_1. 
\ee 
Further because $\wh{\bb}$ is an interior local minimum, 
we have 
\bse
\partial \mL(\wh \bb)/\partial \bb + \A\trans \partial
\rho_\lambda(\wh \bb_{\mM^c})/\partial \bb_{\mM^c} + \A_1\trans  \C\trans \bmu_4= \bf 0. 
\ese
for some  Lagrange  multiplier $\bmu_4 >\0$. 
Note that $(\A_1\trans  \C\trans \bmu_4)\trans\wt\v
=\bmu_t\trans\C(\wt\bb_{\mM}-\wh\bb_{mM})=\0$.
Therefore, plug in (\ref{eq:RElem91}), 
we have 
\bse
&&\alpha_1\|\wt \v\|_2 -2 
\tau_1R_1 \sqrt{\log(p)/n} \|\wt \v\|_1 \\
&\leq&  \{ \A\trans \partial
\rho_\lambda(\wh \bb_{\mM^c})/\partial \bb_{\mM^c}
+\A_1\trans  \C\trans \bmu_4
- \A\trans \partial
\rho_\lambda(\wt \bb_{\mM^c})/\partial \bb_{\mM^c} \}\trans\wt{\v} \\
&=&  \{ \A\trans \partial
\rho_\lambda(\wh \bb_{\mM^c})/\partial \bb_{\mM^c}- \A\trans \partial
\rho_\lambda(\wt \bb_{\mM^c})/\partial \bb_{\mM^c} \}\trans\wt{\v}  \\
&\leq & \{ \|\A\trans \partial
\rho_\lambda(\wh \bb_{\mM^c})/\partial \bb_{\mM^c}\|_\infty+ \|\A\trans \partial
\rho_\lambda(\wt \bb_{\mM^c})/\partial \bb_{\mM^c}
\|_\infty\}\|\wt{\v}\|_1 \\
&\leq& 2\lambda \|\wt{\v}\|_1, 
\ese
where $ \|\A\trans \partial
\rho_\lambda(\bb_{\mM^c})/\partial \bb_{\mM^c}\|_\infty \leq
\lambda$ holds by Lemma \ref{lem:fromlemm5support}. Hence we have 
\bse
\|\wt{\v}\|_2 \leq  \{2\lambda +  2 
\tau_1R_1 \sqrt{\log(p)/n}\} \|\wt{\v}\|_1/\alpha_1 \leq 2R_1 \{2 \lambda + 2 
\tau_1R_1 \sqrt{\log(p)/n}\}/\alpha_1. 
\ese
The right hand side is at most 1 because $\lambda\leq \alpha_1/(8 R_1)$,
and  $n \geq \log (p) 64 \tau_1^2 R_1^4/\alpha_1^2$, which
contradicts  with $\|\wt\v\|_2>1$. Therefore $\|\wt \v\|_2\leq 1$. 

Now note that 
\bse
\mL(\bb) - q_{\lambda}(\bb_{\mM^c}) = (\mL(\bb) - \mu\|\bb_{\mM^c}\|_2^2/2)+
\{\mu\|\bb_{\mM^c}\|_2^2/2 - q_{\lambda}(\bb_{\mM^c}) \}
\ese
and $\{\mu\|\bb_{\mM^c}\|_2^2/2 - q_{\lambda}(\bb_{\mM^c}) \}$ is convex by Lemma
\ref{lem:fromlemm5support}, and hence for any $\bb$ in the feasible
set, we have 
\be\label{eq:lem9second}
\wt{\v}\trans \frac{\partial^2 \{\mL(\bb) -
  q_{\lambda}(\bb_{\mM^c})\}}{\partial \bb \partial
  \bb\trans} \wt{\v} 
&\geq& \wt{\v}\trans \frac{\partial^2 \mL(\bb) }{\partial \bb\partial
  \bb\trans} \wt{\v} - \mu \wt{\v}\trans \A\trans\A \wt{\v}\nonumber\\
&\geq&( \alpha_1 - \mu) \|\wt {\v}\|_2^2- \tau_1 \sqrt{\frac{\log(p)}{n}}\|\wt {\v}\|_1^2.
\ee 
Further by (\ref{eq:lem9first}), we have 
\bse
0\leq \left\{\frac{\partial\mL(\wt\bb)}{\partial \bb\trans} - \frac{\partial q_{\lambda}(\wt
  \bb_{\mM^c})}{\partial \bb_{\mM^c}\trans }\A\right\} (\wh{\bb} - \wt{\bb}) + \lambda\wt{\z}\trans \A(\wh{\bb} -
\wt{\bb}), 
\ese
where $\wt{\z} \in \partial \|\wt\bb_{\mM^c}\|_1$. Further, because $\wh{\bb}$
is an interior local minimum,   for $\wh{\z} \in \partial
\|\wh\bb_{\mM^c}\|_1$,  we have 
\bse
\0&=&\partial \mL(\wh \bb)/\partial \bb + \A\trans \partial
\rho_\lambda(\wh \bb_{\mM^c})/\partial \bb_{\mM^c} 
+ \A_1\trans \C\trans \bmu_4\\
& =& \left\{\frac{\partial\mL(\wh\bb)}{\partial \bb} - \A\trans\frac{\partial q_{\lambda}(\wh
  \bb_{\mM^c})}{\partial \bb_{\mM^c} }\right\}   +
\lambda \A\trans\wh{\z} + \A_1\trans \C\trans \bmu_4,
\ese
which leads to
\bse
&&\left\{\frac{\partial\mL(\wh\bb)}{\partial \bb\trans} - \frac{\partial q_{\lambda}(\wh
  \bb_{\mM^c})}{\partial \bb_{\mM^c}\trans }\A\right\}   \wt{\v} +
\lambda (\wh{\z}\trans \A) \wt{\v}  = \bf 0. 
\ese
Hence
\bse\label{eq:firstorderless}
0 &\leq&  \left\{\frac{\partial\mL(\wh\bb)}{\partial \bb\trans} - \frac{\partial q_{\lambda}(\wh
  \bb_{\mM^c})}{\partial \bb_{\mM^c}\trans }\A\right\} \wt{\v}  -  \left\{\frac{\partial\mL(\wt\bb)}{\partial \bb\trans} - \frac{\partial q_{\lambda}(\wt
  \bb_{\mM^c})}{\partial \bb_{\mM^c}\trans }\A\right\} \wt{\v } \\
&&+ \lambda (\wh{\z}\trans \A) (\wt{\bb} - \wh{\bb}) -
\lambda (\wt{\z}\trans \A) (\wt{\bb} - \wh{\bb}) \\
&=&  \left\{\frac{\partial\mL(\wh\bb)}{\partial \bb\trans} - \frac{\partial q_{\lambda}(\wh
  \bb_{\mM^c})}{\partial \bb_{\mM^c}\trans }\A\right\}  \wt{\v}  -  \left\{\frac{\partial\mL(\wt\bb)}{\partial \bb\trans} - \frac{\partial q_{\lambda}(\wt
  \bb_{\mM^c})}{\partial \bb_{\mM^c}\trans }\A\right\} \wt{\v } \\
&&+ \lambda (\wh{\z}\trans \A) \wt{\bb} - \lambda
\|\wh\bb_{\mM^c}\|_1- \lambda
\|\wt\bb_{\mM^c}\|_1 + \lambda (\wt{\z}\trans \A) \wh{\bb}, 
\ese
which implies 
\be\label{eq:lem939}
&&\lambda
\|\wt\bb_{\mM^c}\|_1  - \lambda (\wh{\z}\trans \A) \wt{\bb} \nonumber\\
&\leq&   \left\{\frac{\partial\mL(\wh\bb)}{\partial \bb\trans} - \frac{\partial q_{\lambda}(\wh
  \bb_{\mM^c})}{\partial \bb_{\mM^c}\trans }\A\right\}   \wt{\v}  -  \left\{\frac{\partial\mL(\wt\bb)}{\partial \bb\trans} - \frac{\partial q_{\lambda}(\wt
  \bb_{\mM^c})}{\partial \bb_{\mM^c}\trans }\A\right\}  \wt{\v } - \lambda
\|\wh\bb_{\mM^c}\|_1+ \lambda (\wt{\z}\trans \A) \wh{\bb}\nonumber\\
&\leq& \left\{\frac{\partial\mL(\wh\bb)}{\partial \bb\trans} - \frac{\partial q_{\lambda}(\wh
  \bb_{\mM^c})}{\partial \bb_{\mM^c}\trans }\A\right\}  \wt{\v}  -  \left\{\frac{\partial\mL(\wt\bb)}{\partial \bb\trans} - \frac{\partial q_{\lambda}(\wt
  \bb_{\mM^c})}{\partial \bb_{\mM^c}\trans }\A\right\}  \wt{\v } \nonumber\\
&\leq& \tau_1 \sqrt{\frac{\log (p)}{n}}\|\wt{\v}\|_1^2 - (\alpha_1 -
\mu)\|\wt{\v}\|_2^2. 
\ee
The second inequality holds because $|(\wt{\z}\trans \A) \wh{\bb}| 
=|\wt{\z}\trans  \wh\bb_{\mM^c}| 
\leq
\|\wt{\z}\|_\infty \|\wh \bb_{\mM^c}\|_1\leq \|\wh \bb_{\mM^c}\|_1$, and
last inequality holds by the Taylor expansion and 
(\ref{eq:lem9second}). 

Now we first assume that the following statement holds: If  $\|(\A\trans \wh{\z})_{\mM\cup S}\|_\infty \leq 1
- \delta$ and  $\lambda \geq 4 R_1 \tau_1 \sqrt{\log (p)/(n)}/\delta$, then 
\bse
\|\wt \v \|_1 \leq \{\sqrt{r}\|\C_r^{-1}\C_{m-r}\|_2 + \sqrt{m
  -r}+  (2/ \delta + 1/2) \sqrt{k} \}\|\wt \v \|_2. 
\ese
We then will have 
\bse
&&\lambda
\|\wt\bb_{\mM^c}\|_1  - \lambda (\wh{\z}\trans \A) \wt{\bb} \\
&\leq& \left[\tau_1\sqrt{\frac{\log(p)}{n}}\{\sqrt{r}\|\C_r^{-1}\C_{m-r}\|_2 + \sqrt{m
  -r}+  (2/ \delta + 1/2) \sqrt{k} \}^2  - (\alpha_1 - \mu)\right]\|\wt \v\|_2^2.
\ese
Now because  $n \geq 4 \log (p) \tau_1^2 \{\sqrt{r}\|\C_r^{-1}\C_{m-r}\|_2 + \sqrt{m
  -r} +  (2/ \delta + 1/2)  \sqrt{k} \}^4/(\alpha_1-\mu)^2$, we have 
\bse
0 = \lambda \|\wt\bb_{\mM^c}\|_1 -\lambda \|\wt\bb_{\mM^c}\|_1 \leq \lambda
\|\wt\bb_{\mM^c}\|_1  - \lambda (\wh{\z}\trans \A) \wt{\bb} \leq -
(\alpha_1 - \mu)/2 \|\wt \v\|_2^2 \leq 0. 
\ese
The first inequality holds by the fact that $ (\wh{\z}\trans
\A) \wt{\bb} 
=\wh{\z}\trans\wt\bb_{\mM^c}
 \leq \|\wh{\z}\|_\infty \|\wt\bb_{\mM^c}\|_1 \leq \|\wt\bb_{\mM^c}\|_1$. Hence we have 
\bse
\lambda
\|\wt\bb_{\mM^c}\|_1  = \lambda (\wh{\z}\trans \A) \wt{\bb}.
\ese
Now because $\|(\wh{\z}\trans \A)_{(\mM \cup S)^c}\|_\infty < 1$,
we conclude that $\wt{\bb}_{(\mM \cup S)^c} = \bf 0$. Hence, $\supp
(\wt{\bb})\subset \mM \cup S$. This would prove the claim in the statement. 

Thus, 
 we only need to show that  if $\|(\A\trans \wh{\z})_{\mM\cup S}\|_\infty \leq 1
- \delta$ and $\lambda \geq 4 R_1 \tau_1 \sqrt{\log (p)/(n)}/\delta
$, then
\bse
\|\wt \v \|_1 \leq \{(\sqrt{r}\|\C_r^{-1}\C_{m-r}\|_2 + \sqrt{m
  -r}) +  (2/ \delta + 1/2) \sqrt{k} \}\|\wt \v \|_2. 
\ese
First from (\ref{eq:lem939}), we have 
\be\label{eq:lem940}
(\alpha_1 - \mu)\|\wt{\v}\|_2^2 - \tau_1 \sqrt{\frac{\log(p)}{n}} \|\wt{\v}\|_1^2
&\leq&  \left\{\frac{\partial\mL(\wt\bb)}{\partial \bb\trans} - \frac{\partial q_{\lambda}(\wt
  \bb_{\mM^c})}{\partial \bb_{\mM^c}\trans }\A\right\}  \wt{\v}  -  \left\{\frac{\partial\mL(\wh\bb)}{\partial \bb\trans} - \frac{\partial q_{\lambda}(\wh
  \bb_{\mM^c})}{\partial \bb_{\mM^c}\trans }\A\right\}  \wt{\v } \nonumber\\
&\leq& \lambda (\wh{\z}\trans \A) \wt{\bb} - \lambda
\|\wh\bb_{\mM^c}\|_1- \lambda
\|\wt\bb_{\mM^c}\|_1 + \lambda (\wt{\z}\trans \A) \wh{\bb}\nonumber\\
&=& \lambda (\wt{\z}\trans \A) \wh{\bb} - \lambda
\|\wt\bb_{\mM^c}\|_1  + \lambda  (\wh{\z}\trans \A) \wt{\v}. 
\ee
The second inequality holds by (\ref{eq:firstorderless}). Now since
$\supp(\wh {\bb}) \subseteq \mM \cup S$, we have 
\be\label{eq:lem941}
&&\lambda (\wt{\z}\trans \A) \wh{\bb} - \lambda
\|\wt\bb_{\mM^c}\|_1 \n\\
 &\leq&  \lambda \|\wh\bb_{\mM^c}\|_1 -  \lambda
\|\wt\bb_{\mM^c}\|_1\nonumber\\
&=&  \lambda \|\wh{\bb}_{ S}\|_1 -  \lambda
\|\wt\bb_{S}\|_1 - \lambda \|\wt\bb_{(\mM\cup S)^c}\|_1\nonumber\\
&=& \lambda \|\wh{\bb}_{ S}\|_1 -  \lambda
\|\wt\bb_{S}\|_1 - (\lambda \|\wt\bb_{(\mM\cup S)^c}\|_1 - \lambda
\|\wh\bb_{(\mM\cup S)^c}\|_1)\nonumber\\
&=&  \lambda \|\wt {\bb}_{ S} + \wh{\bb}_S - \wt{\bb}_S\|_1 -  \lambda
\|\wt\bb_{S}\|_1 
- (\lambda \|\wt\bb_{(\mM\cup S)^c}\|_1 - \lambda
\|\wt \bb_{(\mM\cup S)^c} + \wh \bb_{(\mM\cup S)^c} - \wt
\bb_{(\mM\cup S)^c}  \|_1)\nonumber\\
&\leq& \lambda \|\wt{\v}_{S}\|_1 -  \lambda \|\wt{\v}_{(\mM \cup
  S)^c}\|_1. 
\ee
In addition, 
\be\label{eq:lem942}
 \lambda  (\wh{\z}\trans \A) \wt{\v} &=& \lambda \{(\A\trans
 \wh{\z})_{S} \}\trans \wt{\v}_S + \lambda \{(\A\trans
 \wh{\z})_{(\mM \cup S) ^c } \}\trans \wt{\v}_{(\mM \cup S) ^c }\nonumber\\
&\leq& \lambda \|(\A\trans
 \wh{\z})_{S}  \|_\infty \|\wt{\v}_S\|_1 + \lambda \|(\A\trans
 \wh{\z})_{(\mM \cup S)^c}  \|_\infty \|\wt{\v}_{(\mM \cup S) ^c
 }\|_1\nonumber\\
&\leq & \lambda \|\wt{\v}_S\|_1 + \lambda (1 - \delta ) \|\wt{\v}_{(\mM \cup S) ^c}\|_1.
\ee
Combine (\ref{eq:lem940}), (\ref{eq:lem941}), and (\ref{eq:lem942}),
we have 
\bse
- \tau_1 \sqrt{\frac{\log(p)}{n}}\|\wt{\v}\|_1^2 &\leq& (\alpha_1-
\mu)\|\wt{\v}\|_2^2 - \tau_1 \sqrt{\frac{\log(p)}{n}}\|\wt{\v}\|_1^2
\\
&\leq& \lambda \|\wt{\v}_{S}\|_1 -  \lambda \|\wt{\v}_{(\mM \cup
  S)^c}\|_1  + \lambda \|\wt{\v}_S\|_1 + \lambda (1 - \delta )
\|\wt{\v}_{(\mM \cup S) ^c}\|_1\\
&=& \lambda 2 \|\wt{\v}_{S}\|_1  - \lambda \delta \|\wt{\v}_{(\mM \cup S)
  ^c}\|_1. 
\ese
Now because $\delta \lambda/2 \geq 2 \tau_1 R_1 \sqrt{\log(p)/n} \geq
\tau_1 \sqrt{\log(p)/n}\|\wt{\v}\|_1$,  the above display implies 
\bse
- 2^{-1}\delta \lambda  \|\wt{\v}_{S}\|_1  \leq  \lambda 2 \|\wt{\v}_{S}\|_1  - \lambda \delta \|\wt{\v}_{(\mM \cup S)
  ^c}\|_1
\ese
which leads to 
\bse
\delta \|\wt{\v}_{(\mM \cup S)^c}\|_1 \leq (2 +
\delta/2)\|\wt{\v}_{S}\|_1. 
\ese
Then 
\bse
\|\wt \v\|_1 &=& \|\wt \v_{\mM} \|_1 + \|\wt \v_{S } \|_1 + \|\wt
\v_{(\mM \cup S)^c } \|_1 \\
&\leq& \|\wt \v_{\mM} \|_1 +  \|\wt \v_{S } \|_1 + (4/ \delta + 1)
\|\wt \v_{S } \|_1 \\
&\leq& \|\wt \v_{\mM} \|_1 + (2/ \delta + 1/2) \sqrt{k}
\|\wt \v\|_2 \\
&\leq& \{\sqrt{r}\|\C_r^{-1}\C_{m-r}\|_2 + \sqrt{m
  -r} +  (2/ \delta + 1/2)  \sqrt{k} \}\|\wt \v \|_2.
\ese
The last equality holds by using the same argument as those lead to
(\ref{eq:wtv1}). 
This completes the proof. \qed

\begin{Cor}\label{cor:lem:fromlem3}
Assume  Conditions \ref{con:boundX}--\ref{con:ww2} hold,   $ \lambda \leq \alpha_1/ (8 R_1)$, $R_1 \leq [n \alpha_1^2/\{64 \tau_1^2 \log
(p)\}]^{1/4}$, $\delta \in [4 R_1 \tau_1 \sqrt{\log (p)/n} /\lambda,
1]$, $n \geq 4 \log (p) \tau_1^2 \{ \sqrt{m}+ (2/ \delta + 1/2)
\sqrt{k} \}^4/(\alpha_1-\mu)^2$, and  $\alpha_1>\mu$. Suppose
$\wt{\bb}_a$ is a stationary point of program (\ref{eq:lapHT1}) and
$\wh{\bb}_a$ is the interior local minimize of (\ref{eq:lapHT1}) and satisfies
$\supp  (\wh{\bb}_a)  \subseteq \mM \cup S$. Then $\supp
(\wt{\bb}_a) \subseteq \mM \cup S$. 
\end{Cor}
\noindent Proof: The Corollary holds by using the same arguments as
those lead to Lemma \ref{lem:fromlem3}, while using the consistency
result in Theorem \ref{from:lem:th1}. \qed

\section{Supporting Results related to Test Statistics}
\subsection{Some Definitions}
We define 
\bse
\bSig(\bb)&\equiv& E [\{Y_i \W_i - \exp(\bb\trans\W_i - \bb\trans
\bOmega \bb/2)  (\W_i -
\bOmega \bb)\}^{\otimes2}]
\ese
By using the relation (\ref{eq:merr1})--(\ref{eq:merr3}) and the
additional relations 
\bse
E\{\exp(2\bb_t\trans \W_i - \bb_t\trans \bOmega\bb_t )\mid\X_i\} &=&
\exp(2\bb_t\trans\X_i),\\
E\{\exp(2\bb_t\trans \W_i - \bb_t\trans \bOmega\bb_t )(\W_i-\bOmega\bb_t)\mid\X_i\} &=&
\exp(2\bb_t\trans\X_i)\X_i,\\
E[\exp(2\bb_t\trans \W_i - \bb_t\trans \bOmega\bb_t
) \{(\W_i-\bOmega\bb_t)^{\otimes2} - \bOmega/2\} \mid\X_i] &=&
\exp(2\bb_t\trans\X_i)\X_i^{\otimes2}, 
\ese
we have
\bse
\bSig(\bb_t)&=& E[\{\exp(\bb_t\trans\X_i) - \exp(2
\bb_t\trans\X_i)\}(\X_i\X_i\trans + \bOmega) \\
&&- \exp(2\bb_t\trans
\X_i) \bOmega \bb_t\X_i\trans - \exp(2\bb_t\trans
\X_i) \X_i\bb_t\trans \bOmega + \exp(2\bb_t\trans \W_i - \bb_t\trans
\bOmega \bb_t) (\W_i - \bOmega \bb_t)^{\otimes2} ]
\ese
and the sample version
\bse
\wh \bSig(\bb)&=& n^{-1}\sumi \exp(\bb\trans \W_i -
\bb\trans\bOmega\bb/2)\{ (\W_i - \bOmega \bb)^{\otimes2}  - \bOmega\}\\
&& - \exp(2\bb\trans \W_i -
\bb\trans\bOmega\bb)\{ (\W_i - \bOmega \bb)^{\otimes2} - \bOmega/2\} \\
&&+ \exp(\bb\trans \W_i -
\bb\trans\bOmega\bb/2) \bOmega - \exp(2\bb\trans \W_i -
\bb\trans\bOmega\bb)\bOmega\\
&& - \exp(2\bb\trans\W_i - \bb\trans\bOmega\bb) \bOmega \bb (\W_i - \bOmega\bb)\trans 
- \exp(2\bb\trans\W_i - \bb\trans\bOmega\bb) (\W_i - \bOmega\bb)\bb\trans \bOmega\\
&&+  \exp(2 \bb\trans \W_i - \bb\trans \bOmega\bb ) (\W_i -
\bOmega \bb)^{\otimes2}\\
&=& n^{-1}\sumi \exp(\bb\trans \W_i -
\bb\trans\bOmega\bb/2)\{ (\W_i - \bOmega \bb)^{\otimes2} \}\\
&& - \exp(2\bb\trans \W_i -
\bb\trans\bOmega\bb)\{ (\W_i - \bOmega \bb)^{\otimes2} + \bOmega/2\} \\
&& - \exp(2\bb\trans\W_i - \bb\trans \bOmega \bb) \bOmega \bb (\W_i - \bOmega\bb)\trans -
\exp(2\bb\trans\W_i - \bb\trans \bOmega \bb) (\W_i - \bOmega\bb)\bb\trans \bOmega\\
&&+ \exp(2 \bb\trans \W_i - \bb\trans \bOmega\bb ) (\W_i -
\bOmega \bb)^{\otimes2}.
\ese
Further, let 
\bse
{\bPsi} (\bSig, \Q, \bb)= (\C [\I_{m\times m}, {\bf 0}_{m \times k}]  {\Q}_{\mM \cup S,
  \mM\cup S}^{-1}({\bb}) \bSig_{\mM \cup S,
  \mM\cup S} ({\bb}) {\Q}_{\mM \cup S,
  \mM\cup S}^{-1} ({\bb}) [\I_{m\times m}, {\bf 0}_{m \times k}]\trans
\C\trans), 
\ese
and 
\bse
T_0 = (\bomega_n + \sqrt{n}\h_n) \trans  \bPsi^{-1}(\bSig, \Q, \bb_t)  (\bomega_n +
\sqrt{n}\h_n), 
\ese
where 
\bse
\bomega_n = - \sqrt{n} \C [\I_{m\times m}, {\bf 0}_{m \times k}] {\Q}_{ \mM \cup
  S,  \mM \cup S}^{-1} (\bb_t)\left\{\frac{\partial \mL(\bb_t)}{\partial
    \bb_t}\right\}_{\mM \cup S}
\ese
Further define Wald statistics as 
\bse
T_W =n (\C \wh{\bb}_{a\mM} - \t)\trans  {\bPsi} (\wh\bSig, \wh\Q, \wh\bb_a) ^{-1}(\C
\wh{\bb}_{a\mM} - \t); 
\ese
the score statistics 
\bse
T_S &=& n \left\{\frac{\partial \mL(\wh{\bb})}{\partial \bb\trans
  }\right\}_{\mM \cup S} (\C [\I_{m\times m}, {\bf 0}_{m \times k}] {\Q}_{ \mM \cup
  S,  \mM \cup S}^{-1} (\wh\bb))\trans  \\
&&\times \bPsi^{-1}(\wh\bSig, \wh\Q, \wh\bb) \C [\I_{m\times m}, {\bf 0}_{m \times k}] {\Q}_{ \mM \cup
  S,  \mM \cup S}^{-1} (\wh\bb) \left\{\frac{\partial
    \mL(\wh{\bb})}{\partial \bb }\right\}_{\mM \cup S}.  
\ese
Because without knowing the distribution of $\X$, the full likelihood
is unknown and hence we do not discuss the likelihood ratio test
here. 

\subsection{Assumptions}{\label{sec:assD}}
\begin{enumerate}[label=(D\arabic*)]
\item \label{con:3ord}
Assume 
\bse
&&\max_{1\leq i\leq n}(\|\bPsi^{-1/2} (\bSig, \Q, \bb_t)\C [\I_{m\times m}, {\bf 0}_{m \times k}]{\Q}_{\mM \cup
  S,  \mM \cup S}^{-1}(\bb_t)  \W_{i \mM\cup S} \|_2^3\\
&&+ \|\bPsi^{-1/2} (\bSig, \Q, \bb_t)\C [\I_{m\times m}, {\bf 0}_{m \times k}]{\Q}_{\mM \cup
  S,  \mM \cup S}^{-1}(\bb_t)  (\bOmega \bb_t)_{\mM \cup S}\|_2^3) ^{-1} \\
&&\times E [\|Y_i \bPsi^{-1/2} (\bSig, \Q, \bb_t)\C [\I_{m\times m}, {\bf 0}_{m \times k}]{\Q}_{\mM \cup
  S,  \mM \cup S}^{-1}(\bb_t)   \W_i- \exp(\bb_t\trans\W_i - \bb_t\trans
\bOmega \bb_t/2) \\
&&\times \bPsi^{-1/2} (\bSig, \Q, \bb_t)\C [\I_{m\times m}, {\bf 0}_{m \times k}]{\Q}_{\mM \cup
  S,  \mM \cup S}^{-1}(\bb_t)   (\W_i - \bb_t\trans
\bOmega)_{\mM \cup S}\|_2^3|\W_i ]= O(1)
\ese
\item \label{con:psieigen}$c_{\bPsi}\leq \alpha_{\min}\{{\bPsi} (\bSig, \Q, \bb)\} \leq
  \alpha_{\max}\{{\bPsi} (\bSig, \Q, \bb)\} \leq C_{\bPsi}.$
\end{enumerate}

\subsection{Some Lemmas}\label{sec:restlemma}
\begin{Lem}\label{lem:frombentkus}
Suppose $\X_1, \ldots, \X_n$ are independent $m$ dimensional random
vector which satisfies, $E(X_i)= 0$ and $\sumi\cov(\X_i) = \I_m$. Let
$\Z$ be a m-dimensional standard  multivariate normal vector, then 
\bse
\sup_{\mathcal C} |\Pr(\sumi \X_i \in {\mathcal C}) - \Pr(\Z\in
\mathcal{C})|= O(m^{1/2}\sumi E\|\X_i\|_2^3). 
\ese
\end{Lem}
\noindent Proof: The lemma  follows Theorem 1.1 in
\cite{bentkus2005}.

\begin{Lem}\label{lem:psi}
Assume Conditions \ref{con:boundX} -- \ref{con:ww3} hold, then 
\bse
\|{\bPsi}^{-1} (\wh \bSig, \wh\Q, \wh\bb_a) - {\bPsi}^{-1} (\bSig, \Q,
\bb_t)\|_2 = \left\{\frac{\log(p)}{n}\right\}^{1/4}  (m+ k), 
\ese
and 
\bse
\|{\bPsi}^{-1} (\wh \bSig, \wh\Q, \wh\bb) - {\bPsi}^{-1} (\bSig, \Q,
\bb_t)\|_2 = \left\{\frac{\log(p)}{n}\right\}^{1/4}  (m+ k). 
\ese
\end{Lem}
\noindent Proof: First 
\be\label{eq:bPsi}
&&\|{\bPsi}(\wh\bSig, \wh\Q, \wh{\bb}_a) - {\bPsi}(\bSig, \Q,
\bb_t)\|_2\nonumber\\
&=& O_p\{\max(\|\wh\Q_{\mM \cup S, \mM \cup S}(\wh{\bb}_a) ^{-1} -
\Q_{\mM \cup S, \mM \cup S}({\bb_t}) ^{-1}\|_2,\nonumber\\
&&
\|\wh \bSig_{\mM \cup S, \mM \cup S}(\wh{\bb}_a) - \bSig_{\mM \cup S, \mM \cup S}({\bb_t})\|_2\}. 
\ee
Using the same arguments as those lead to (\ref{eq:2normdiff}),  we
have $\|\wh\Q(\wh{\bb}_a) ^{-1} -  \Q({\bb_t}) ^{-1}\|_2 =
\left\{\frac{\log(p)}{n}\right\}^{1/4}  (m+ k) $. 
Now  recall that ${\mathbb K} \equiv \{\v\in {\mathbb R}^{\mM\cup S}: \|\v\|_2\leq 1\}$, for $\v \in
{\mathbb K}$, 
by the Taylor expansion
\bse
&&  \v\trans\{ \wh \bSig(\wh{\bb}_a) - \wh{\bSig}({\bb_t})\}\v\\
&=&n^{-1}\sumi \exp(\bb_a\trans \W_i -
\bb_a\trans\bOmega\bb_a/2)\v \trans \{ (\W_i - \bOmega \bb_a)^{\otimes2}
\}\v (\W_i - \bOmega \bb_a )\trans (\wh{\bb}_a - \bb_t)\\
&& + n^{-1}\sumi   2 \exp(\bb_a\trans \W_i -
\bb_a\trans\bOmega\bb_a/2)  \v\trans  (\W_i - \bOmega \bb_a )\v\trans \bOmega
(\wh{\bb}_a - \bb_t)\\
&& - n^{-1}\sumi 2\exp(2\bb_a\trans \W_i -
\bb_a\trans\bOmega\bb_a)\v \trans \{ (\W_i - \bOmega \bb_a)^{\otimes2} + \bOmega/2
\}\v (\W_i - \bOmega \bb_a)\trans (\wh{\bb}_a - \bb_t) \\
&& - n^{-1}\sumi   2 \exp(2\bb_a\trans \W_i -
\bb_a\trans\bOmega\bb_a)  \v\trans  (\W_i - \bOmega \bb_a )\v\trans \bOmega
(\wh{\bb}_a - \bb_t)\\
&& - n^{-1}\sumi 2 \exp(2\bb_a\trans\W_i - \bb_a\trans \bOmega
\bb_a)\v\trans\bOmega \bb_a (\W_i - \bOmega\bb_a)\trans \v (\W_i -
\bOmega\bb_a)\trans (\wh{\bb}_a - \bb_t)\\
&& - n^{-1}\sumi  2\exp(2\bb_a\trans\W_i - \bb_a\trans \bOmega
\bb_a)\{\v\trans (\W_i - \bOmega\bb_a)\v\trans \bOmega -
\v\trans\bOmega\bb_a \v\trans
\bOmega \} (\wh{\bb}_a - \bb_t)\\
&& 
- n^{-1}\sumi 2 \exp(2\bb_a\trans\W_i - \bb_a\trans \bOmega
\bb_a)\v\trans (\W_i - \bOmega\bb_a)\bb_a\trans \bOmega \v (\W_i -
\bOmega\bb_a)\trans (\wh{\bb}_a - \bb_t)\\
&& - n^{-1}\sumi  \exp(2\bb_a\trans\W_i - \bb_a\trans \bOmega
\bb_a)\{\v\trans (\W_i - \bOmega\bb_a)\v\trans \bOmega -
\v\trans\bOmega\bb_a \v\trans
\bOmega \} (\wh{\bb}_a - \bb_t)\\
&&+ n^{-1}\sumi \exp(2\bb_a\trans \W_i -
\bb_a\trans\bOmega\bb_a)\v \trans \{ (\W_i - \bOmega \bb_a)^{\otimes2}
\}\v (\W_i - \bOmega \bb_a)\trans (\wh{\bb}_a - \bb_t) \\
&& +  n^{-1}\sumi   2 \exp(2\bb_a\trans \W_i -
\bb_a\trans\bOmega\bb_a)  \v\trans  (\W_i - \bOmega \bb_a)\v\trans \bOmega
(\wh{\bb}_a - \bb_t)\\
&=& n^{-1}\sumi \exp(\bb_a\trans \W_i -
\bb_a\trans\bOmega\bb_a/2)\v \trans \{ (\W_i - \bOmega \bb_a)^{\otimes2}
\}\v (\W_i - \bOmega \bb_a)\trans (\wh{\bb}_a - \bb_t)\\
&& + n^{-1}\sumi   2 \exp(\bb_a\trans \W_i -
\bb_a\trans\bOmega\bb_a/2)  \v\trans  (\W_i - \bOmega \bb_a )\v\trans \bOmega
(\wh{\bb}_a - \bb_t)\\
&& - n^{-1}\sumi 2 \exp(2\bb_a\trans \W_i -
\bb_a\trans\bOmega\bb_a)\v \trans \bOmega \v/2  (\W_i - \bOmega \bb_a
)\trans (\wh{\bb}_a - \bb_t) \\
&& - n^{-1}\sumi 4 \exp(2\bb_a\trans\W_i - \bb_a\trans \bOmega
\bb_a)\v\trans\bOmega \bb_a (\W_i - \bOmega\bb_a)\trans \v (\W_i -
\bOmega\bb_a)\trans (\wh{\bb}_a - \bb_t)\\
&& - n^{-1}\sumi 2 \exp(2\bb_a\trans\W_i - \bb_a\trans \bOmega
\bb_a)\{\v\trans (\W_i - \bOmega\bb_a)\v\trans \bOmega -
\v\trans\bOmega\bb_a  \v\trans
\bOmega \} (\wh{\bb}_a - \bb_t), 
\ese
where $\bb_a$ is a point in between $\wh{\bb}_a$ and $\bb_t$. Using
the same arguments as those lead to (\ref{eq:2normbound1}) and (\ref{eq:2normbound2}), we have \bse
 \| \wh \bSig(\wh{\bb}_a) - \wh{\bSig}({\bb_t})\|_2=
 O_p\left[\left\{\frac{\log(p)}{n}\right\}^{1/4} (m + k) \right], 
\ese
and 
\bse
 \| \wh \bSig(\bb) - {\bSig}({\bb_t})\|_2=
o_p\left[\left\{\frac{\log(p)}{n}\right\}^{1/4} (m + k) \right]. 
\ese
Therefore, 
\bse
\| \wh \bSig(\wh{\bb}_a) - \wh{\bSig}({\bb_t})\|_2 =
O_p\left[\left\{\frac{\log(p)}{n}\right\}^{1/4} (m + k)
\right]. 
\ese
Combine with (\ref{eq:bPsi}) and the fact that $\|\wh\Q(\wh{\bb}_a) ^{-1} -  \Q({\bb_t}) ^{-1}\|_2 =
\left\{\frac{\log(p)}{n}\right\}^{1/4}  (m+ k) $, by Lemma
\ref{lem:fromlemma11loh2017} we have 
\bse
&& \|{\bPsi}^{-1} (\wh \bSig, \wh\Q, \wh\bb_a) - {\bPsi}^{-1} (\bSig, \Q,
\bb_t)\|_2\\
&=& O_p\{\|{\bPsi} (\wh \bSig, \wh\Q, \wh\bb_a) - {\bPsi} (\bSig, \Q,
\bb_t)\|_2\}\\
&=& \left\{\frac{\log(p)}{n}\right\}^{1/4}  (m+ k).
\ese
The second relation in the statement holds by using the same arguments as those lead
the above results. This proves the result.\qed

\section{The composite ROIs}\label{sec:composite}
Based on \cite{braak1991, Landau2016, scholl2016},  we define the
composite regions as follow, where letter $L$ and $R$ represent the
left and right hemispheres, respectively. 
\begin{itemize}

\item {\bf Braak 1 and 2 composite region (Braak12):}
  L\_entorhinal, R\_entorhinal
  \item {\bf Braak 3 and 4 composite region (Braak34):
    }\\
    L\_parahippocampal, L\_fusiform, 
L\_lingual, L\_amygdala, R\_parahippocampal, R\_fusiform, R\_lingual,
R\_amygdala, L\_middletemporal, L\_caudantcing, L\_rostantcing,
L\_postcing, L\_isthmuscing,  L\_insula, L\_inferiortemporal, L\_temppole,
R\_middletemporal, R\_caudantcing, R\_rostantcing, R\_postcing,
R\_isthmuscing, R\_insula, R\_inferiortemporal, R\_temppole.

\item {\bf Braak 5 and 6 composite region (Braak56):}\\
L\_superior\_frontal, L\_lateral\_orbitofrontal, L\_medial\_orbitofrontal,\\
L\_frontal\_pole, L\_caudal\_middle\_frontal, L\_rostral\_middle\_frontal,
L\_pars\_opercularis, L\_pars\_orbitalis, L\_pars\_triangularis,
L\_lateraloccipital, L\_parietalsupramarginal, L\_parietalinferior,
L\_superiortemporal, L\_parietalsuperior, L\_precuneus,
L\_bankSuperiorTemporalSulcus, L\_tranvtemp, R\_superior\_frontal,
R\_lateral\_orbitofrontal, R\_medial\_orbitofrontal, R\_frontal\_pole,
R\_caudal\_middle\_frontal, R\_rostral\_middle\_frontal, R\_pars\_opercularis,
R\_pars\_orbitalis, R\_pars\_triangularis, R\_lateraloccipital,
R\_parietalsupramarginal, R\_parietalinferior, \\ R\_superiortemporal,
R\_parietalsuperior, R\_precuneus, R\_bankSuperiorTemporalSulcus,
R\_tranvtemp, L\_pericalcarine, L\_postcentral, L\_cuneus, L\_precentral,
L\_paracentral, R\_pericalcarine, R\_postcentral, R\_cuneus,
R\_precentral, R\_paracentral

\end{itemize}

\section{Simulations when measurement error covariance is
  unknown}\label{sec:addsimulation}

{\color{cyan} To evaluate the proposed adjustment estimator in Section \ref{sec:con}, we added the
  following three
  simulation settings where the covariance matrices of the measurement
  errors contain different numbers of unknown parameters. 
	\begin{enumerate}
		\item $\bOmega$ is a matrix with $p/4$ unknown
                  parameters, corresponding to $p/4$ nonzero diagonal entries.
		\item  
            $\bOmega$ is a matrix with $p/2$ unknown parameters,
                corresponding to $p/2$ nonzero entries.
		\item $\bOmega$ is a matrix with $6p-15$ unknown
                  parameters.  Specifically, 
$\bOmega=(\sigma_{ij})_{p\times p}$,  where
                  $\sigma_{ij}=0.05(1-|i-j|/5)$ for $|i-j| \leq 5$ and
                  $\sigma_{ij}=0$ for  $|i-j| > 5$.
                \end{enumerate}
                In these settings, the number of unknown
                  parameters in $\bOmega$
increases, while all 
other settings are the same  as in Section
\ref{sec:simulation}. We evaluate the empirical sizes and powers of
the Wald test for  $n = 300$, $p=50$ in Table \ref{S1S} and $p=350$ in Table \ref{S2S}.
In the first two settings, where $\bOmega$ contains relatively few
parameters hence their convergence is
sufficiently fast, the 
empirical sizes are well controlled around the
nominal  level 0.05 in all hypothesis tests. In the last setting,
the parameter number is too large to achieve sufficiently
  fast convergence, hence
the Wald test cannot control the Type I errors. These results
imply that the proposed adjustment can control Type I error rate
when the error covariance matrix contains small or
  moderate number of unknown parameters.

\begin{table}
	\centering
	\caption{ The empirical size and power of Wald tests for
          linear hypothesis testing with $n=300$. NUP standards for
          the number of unknown parameter in $\bOmega$.  }\label{S1S}
	\begin{tabular}{c|c|c|c|c|cccc}\hline
		&&NUP=$p/4$& NUP=$p/2$& NUP $=6p-15$\\\hline
		\multirow{35}{*}{$p=50$}&	$\beta_{2}$&\multicolumn{7}{c}{$H_{0,1}:\beta_{2}=-0.75,\ v.s. \ H_{a,1}:  \beta_{2} \neq-0.75$}\\\cline{2-9}
		&-0.75 &0.054&0.044&0.262\\
		&-0.65 &0.334&0.384&0.348\\
		&-0.55 &0.832&0.852&0.542\\
		&-0.35 &1.000&1.000&0.910\\\cline{2-9}	                                                        			
		&$\beta_{3}$&\multicolumn{7}{c}{$H_{0,2}:\beta_{3}=0, \ v.s. \ H_{a,2}:  \beta_{3} \neq0$ } \\\cline{2-9}		
		&0.0&0.064&0.044&0.558\\
		&0.1&0.262&0.244&0.330\\
		&0.2&0.714&0.768&0.320\\
		&0.4&1.000&0.994&0.644\\\cline{2-9}
		&$\beta_{p}$&\multicolumn{7}{c}{$H_{0,3}:\beta_{p}=0, \ v.s.  \  H_{a,3}:\beta_{p} \neq0$ }    \\\cline{2-9}	            			
		&0.0&0.074&0.066&0.070\\
		&0.1&0.302&0.316&0.228\\
		&0.2&0.778&0.804&0.676\\
		&0.4&1.000&1.000&0.992\\\cline{2-9}   
		&	$\beta_{1} + \beta_{2}$&\multicolumn{7}{c}{$H_{0,4}:\beta_{1}+\beta_{2}=0, \ v.s. \ H_{a,4}:  \beta_{1}+\beta_{2}\neq0$}\\\cline{2-9}	
		&0.0&0.048&0.052&0.056\\
		&0.1&0.164&0.152&0.220\\
		&0.2&0.488&0.484&0.650\\
		&0.4&0.966&0.972&0.998\\\cline{2-9}
		&$\beta_{3} +\beta_{4}$&\multicolumn{7}{c}{$H_{0,5}:\beta_{3}+\beta_{4}=0,  \ v.s. \ H_{a,5}: \beta_{3}+\beta_{4}\neq 0.$}\\ \cline{2-9}	   
		&0.0&0.056&0.056&0.542\\
		&0.1&0.182&0.156&0.326\\
		&0.2&0.456&0.518&0.278\\
		&0.4&0.972&0.952&0.588\\\cline{2-9}
		&$\beta_{1} + \beta_{p}$&\multicolumn{7}{c}{$H_{0,6}:\beta_{1}+\beta_{p}=0.75, \ v.s.  \ H_{a,6}: \beta_{1}+\beta_{p}\neq0.75$}\\\cline{2-9}	   
		&0.75 &0.050&0.048&0.328\\
		&0.85 &0.154&0.168&0.338\\
		&0.95 &0.472&0.486&0.456\\
		&1.15 &0.974&0.974&0.668\\\cline{2-9}
		&$\beta_{2} + \beta_{3}$&\multicolumn{7}{c}{$H_{0,7}: \beta_{2}+ \beta_{3}=-0.75,  \ v.s.  \ H_{a,7}: \beta_{2} + \beta_{3} \neq -0.75$}\\ \cline{2-9}	   
		&-0.75 &0.048&0.040&0.246\\
		&-0.65 &0.170&0.166&0.292\\
		&-0.55 &0.480&0.490&0.414\\
		&-0.35 &0.976&0.968&0.856\\\hline		
	\end{tabular}	
\end{table}

\begin{table}
	\centering
	\caption{ The empirical size and power of Wald tests for linear hypothesis testing with $n=300$. NUP standards for
          the number of unknown parameter in $\bOmega$.}\label{S2S}
	\begin{tabular}{c|c|c|c|c|cccc}\hline
	  &&NUP=$p/4$&NUP=$p/2$& NUP $=6p-15$ \\\hline
	  \multirow{35}{*}{$p=350$}&
                                     $\beta_{2}$&\multicolumn{7}{c}
{$H_{0,1}:\beta_{2}=-0.75,\ v.s. \ H_{a,1}:  \beta_{2} \neq-0.75$}\\\cline{2-9}	
	  &-0.75 &0.060&0.052&0.672\\
	  &-0.65 &0.326&0.314&0.792\\
	  &-0.55 &0.846&0.818&0.910\\
	  &-0.35 &0.898&0.886&0.738\\\cline{2-9}	                                                          			
	  &$\beta_{3}$&\multicolumn{7}{c}{$H_{0,2}:\beta_{3}=0, \ v.s. \ H_{a,2}:  \beta_{3} \neq0$ } \\\cline{2-9}		
	  &0.0&0.046&0.062&0.466\\
	  &0.1&0.290&0.294&0.150\\
	  &0.2&0.774&0.732&0.150\\
	  &0.4&0.996&0.998&0.762\\\cline{2-9}                                                      			
	  &$\beta_{p}$&\multicolumn{7}{c}{$H_{0,3}:\beta_{p}=0, \ v.s.  \  H_{a,3}:\beta_{p} \neq0$ }    \\\cline{2-9}	            			
	  &0.0&0.054&0.068&0.044\\
	  &0.1&0.298&0.286&0.206\\
	  &0.2&0.790&0.778&0.662\\
	  &0.4&0.986&0.986&1.000\\\cline{2-9}   
	&	$\beta_{1} + \beta_{2}$&\multicolumn{7}{c}{$H_{0,4}:\beta_{1}+\beta_{2}=0, \ v.s. \ H_{a,4}:  \beta_{1}+\beta_{2}\neq0$}\\\cline{2-9}	
		&0.0&0.044&0.052&0.060\\
		&0.1&0.190&0.184&0.226\\
		&0.2&0.506&0.498&0.648\\
		&0.4&0.872&0.860&0.998\\\cline{2-9}
		&$\beta_{3} +\beta_{4}$&\multicolumn{7}{c}{$H_{0,5}:\beta_{3}+\beta_{4}=0,  \ v.s. \ H_{a,5}: \beta_{3}+\beta_{4}\neq 0.$}\\ \cline{2-9}	   
		&0.0&0.064&0.066&0.508\\
		&0.1&0.168&0.180&0.164\\
		&0.2&0.472&0.478&0.116\\
		&0.4&0.950&0.944&0.684\\\cline{2-9}
		&$\beta_{1} + \beta_{p}$&\multicolumn{7}{c}{$H_{0,6}:\beta_{1}+\beta_{p}=0.75, \ v.s.  \ H_{a,6}: \beta_{1}+\beta_{p}\neq0.75$}\\\cline{2-9}	   
		&0.75 &0.056&0.054&0.680\\
		&0.85 &0.182&0.152&0.738\\
		&0.95 &0.476&0.464&0.778\\
		&1.15 &0.966&0.960&0.512\\\cline{2-9}
		&$\beta_{2} + \beta_{3}$&\multicolumn{7}{c}{$H_{0,7}: \beta_{2}+ \beta_{3}=-0.75,  \ v.s.  \ H_{a,7}: \beta_{2} + \beta_{3} \neq -0.75$}\\ \cline{2-9}	   
		&-0.75 &0.068&0.056&0.678\\
		&-0.65 &0.172&0.182&0.680\\
		&-0.55 &0.492&0.446&0.696\\
		&-0.35 &0.950&0.942&0.900\\\hline		
	\end{tabular}                                                                                                    
\end{table}

\bibliographystyle{plainnat}
\bibliography{hdmerrorinfer}

%% file: hdinferrev3.bbl
\begin{thebibliography}{xx}

\harvarditem[Baker et~al.]{Baker, Lockhart, Price, He, Huesman, Schonhaut,
  Faria, Rabinovici \harvardand\ Jagust}{2017}{baker2017}
Baker, S.~L., Lockhart, S.~N., Price, J.~C., He, M., Huesman, R.~H., Schonhaut,
  D., Faria, J., Rabinovici, G. \harvardand\ Jagust, W.~J.  \harvardyearleft
  2017\harvardyearright , `Reference tissue--based kinetic evaluation of
  18f-av-1451 for tau imaging', {\em Journal of Nuclear Medicine} {\bf
  58}(2),~332--338.

\harvarditem{Belloni, Chernozhukov \harvardand\
  Kaul}{2017}{belloni2017inference}
Belloni, A., Chernozhukov, V. \harvardand\ Kaul, A.  \harvardyearleft
  2017\harvardyearright , `Confidence bands for coefficients in high
  dimensional linear models with error-in-variables', {\em arXiv preprint
  arXiv:1703.00469} .

\harvarditem{Belloni \harvardand\ Rosenbaum}{2016}{belloni2016}
Belloni, A. \harvardand\ Rosenbaum, M.  \harvardyearleft 2016\harvardyearright
  , `An $\{l_1$, $l_2$, $l_\infty\}$-regularization approach to
  high-dimensional errors-in-variables models', {\em Electronic Journal of
  Statistics} {\bf 10},~1935--7524.

\harvarditem{Belloni, Rosenbaum \harvardand\ Tsybakov}{2017}{belloni2017linear}
Belloni, A., Rosenbaum, M. \harvardand\ Tsybakov, A.~B.  \harvardyearleft
  2017\harvardyearright , `Linear and conic programming estimators in high
  dimensional errors-in-variables models', {\em Journal of the Royal
  Statistical Society: Series B (Statistical Methodology)} {\bf
  79}(3),~939--956.

\harvarditem{Benjamini \harvardand\ Hochberg}{1995}{benjamini1995}
Benjamini, Y. \harvardand\ Hochberg, Y.  \harvardyearleft 1995\harvardyearright
  , `Controlling the false discovery rate: a practical and powerful approach to
  multiple testing', {\em Journal of the Royal statistical society: series B
  (Methodological)} {\bf 57}(1),~289--300.

\harvarditem{Bentkus}{2005}{bentkus2005}
Bentkus, V.  \harvardyearleft 2005\harvardyearright , `A lyapunov-type bound in
  $r^d$', {\em Theory of Probability \& Its Applications} {\bf
  49}(2),~311--323.

\harvarditem{Bickel \harvardand\ Levina}{2008}{bickel2008}
Bickel, P.~J. \harvardand\ Levina, E.  \harvardyearleft 2008\harvardyearright ,
  `Covariance regularization by thresholding', {\em The Annals of Statistics}
  {\bf 36}(6),~2577--2604.

\harvarditem{Braak \harvardand\ Braak}{1991}{braak1991}
Braak, H. \harvardand\ Braak, E.  \harvardyearleft 1991\harvardyearright ,
  `Neuropathological stageing of alzheimer-related changes', {\em Acta
  neuropathologica} {\bf 82}(4),~239--259.

\harvarditem[Brown et~al.]{Brown, Weaver \harvardand\
  Wolfson}{2019}{brown2019meboost}
Brown, B., Weaver, T. \harvardand\ Wolfson, J.  \harvardyearleft
  2019\harvardyearright , `Meboost: Variable selection in the presence of
  measurement error', {\em Statistics in medicine} {\bf 38}(15),~2705--2718.

\harvarditem{Cai \harvardand\ Liu}{2011}{cai2011}
Cai, T. \harvardand\ Liu, W.  \harvardyearleft 2011\harvardyearright ,
  `Adaptive thresholding for sparse covariance matrix estimation', {\em Journal
  of the American Statistical Association} {\bf 106}(494),~672--684.

\harvarditem[Carroll et~al.]{Carroll, Ruppert, Stefanski \harvardand\
  Crainiceanu}{2006}{carroll2006}
Carroll, R.~J., Ruppert, D., Stefanski, L.~A. \harvardand\ Crainiceanu, C.~M.
  \harvardyearleft 2006\harvardyearright , {\em Measurement error in nonlinear
  models: a modern perspective}, Chapman and Hall/CRC, New Tork.

\harvarditem{Chen \harvardand\ Gu}{2014}{chen2014}
Chen, L. \harvardand\ Gu, Y.  \harvardyearleft 2014\harvardyearright , `The
  convergence guarantees of a non-convex approach for sparse recovery', {\em
  IEEE Transactions on Signal Processing} {\bf 62}(15),~3754--3767.

\harvarditem{Cook \harvardand\ Stefanski}{1995}{cook1995}
Cook, J. \harvardand\ Stefanski, L.~A.  \harvardyearleft 1995\harvardyearright
  , `A simulation extrapolation method for parametric measurement error
  models', {\em Journal of the American Statistical Association} {\bf
  89},~1314–1328.

\harvarditem{Datta \harvardand\ Zou}{2017}{datta2017cocolasso}
Datta, A. \harvardand\ Zou, H.  \harvardyearleft 2017\harvardyearright ,
  `Cocolasso for high-dimensional error-in-variables regression', {\em The
  Annals of Statistics} {\bf 45}(6),~2400--2426.

\harvarditem[Duchi et~al.]{Duchi, Shalev-Shwartz, Singer \harvardand\
  Chandra}{2008}{duchi2008}
Duchi, J., Shalev-Shwartz, S., Singer, Y. \harvardand\ Chandra, T.
  \harvardyearleft 2008\harvardyearright , Efficient projections onto the l
  1-ball for learning in high dimensions, {\em in} `Proceedings of the 25th
  international conference on Machine learning', ACM, pp.~272--279.

\harvarditem[Fallah et~al.]{Fallah, Mitnitski \harvardand\
  Rockwood}{2011}{fallah2011applying}
Fallah, N., Mitnitski, A. \harvardand\ Rockwood, K.  \harvardyearleft
  2011\harvardyearright , `Applying neural network poisson regression to
  predict cognitive score changes', {\em Journal of Applied Statistics} {\bf
  38}(9),~2051--2062.

\harvarditem{Fan \harvardand\ Li}{2001}{fan2001}
Fan, J. \harvardand\ Li, R.  \harvardyearleft 2001\harvardyearright , `Variable
  selection via nonconcave penalized likelihood and its oracle properties',
  {\em Journal of the American statistical Association} {\bf
  96}(456),~1348--1360.

\harvarditem[Fan et~al.]{Fan, Liao \harvardand\ Mincheva}{2011}{fan2011}
Fan, J., Liao, Y. \harvardand\ Mincheva, M.  \harvardyearleft
  2011\harvardyearright , `High dimensional covariance matrix estimation in
  approximate factor models', {\em Annals of statistics} {\bf 39}(6),~3320.

\harvarditem{Fletcher \harvardand\ Watson}{1980}{fletcher1980}
Fletcher, R. \harvardand\ Watson, G.~A.  \harvardyearleft 1980\harvardyearright
  , `First and second order conditions for a class of nondifferentiable
  optimization problems', {\em Mathematical Programming} {\bf 18}(1),~291--307.

\harvarditem{Jiang \harvardand\ Ma}{2021}{jiang2021}
Jiang, F. \harvardand\ Ma, Y.  \harvardyearleft 2021\harvardyearright ,
  `Poisson regression with error corrupted high dimensional features', {\em
  Statistica Sinica} .
\newblock In press.

\harvarditem[Jiang et~al.]{Jiang, Zhou \harvardand\ Ma}{2021}{github}
Jiang, F., Zhou, Y. \harvardand\ Ma, Y.  \harvardyearleft 2021\harvardyearright
  , `Code and data for high dimensional poisson models with noisy data:
  hypothesis testing for nonlinear nonconvex optimization'.
\newblock
  https://github.com/feigroup/high-dimenisional-inference-for-count-data-with-errors.

\harvarditem[Katz et~al.]{Katz, Wang, Nester, Derby, Zimmerman, Lipton,
  Sliwinski \harvardand\ Rabin}{2021}{katz2021t}
Katz, M.~J., Wang, C., Nester, C.~O., Derby, C.~A., Zimmerman, M.~E., Lipton,
  R.~B., Sliwinski, M.~J. \harvardand\ Rabin, L.~A.  \harvardyearleft
  2021\harvardyearright , `T-moca: A valid phone screen for cognitive
  impairment in diverse community samples', {\em Alzheimer's \& Dementia:
  Diagnosis, Assessment \& Disease Monitoring} {\bf 13}(1),~e12144.

\harvarditem[Landau et~al.]{Landau, Ward, Murphy \harvardand\
  Jagust}{2016}{Landau2016}
Landau, S., Ward, T.~J., Murphy, A. \harvardand\ Jagust, W.  \harvardyearleft
  2016\harvardyearright , `Flortaucipir (av-1451) processing methods'.
\newline\harvardurl{\url{https://adni.bitbucket.io/reference/docs/UCBERKELEYAV1451/UCBERKELEY_AV1451_Methods_Aug2018.pdf}}

\harvarditem[Leisman et~al.]{Leisman, Moustafa \harvardand\
  Shafir}{2016}{leisman2016}
Leisman, G., Moustafa, A.~A. \harvardand\ Shafir, T.  \harvardyearleft
  2016\harvardyearright , `Thinking, walking, talking: integratory motor and
  cognitive brain function', {\em Frontiers in public health} {\bf 4},~94.

\harvarditem[Li et~al.]{Li, Li \harvardand\ Ma}{2021}{lilima2021}
Li, M., Li, R. \harvardand\ Ma, Y.  \harvardyearleft 2021\harvardyearright ,
  `Inference in high dimensional linear measurement error models', {\em Journal
  of Multivariate Analysis} {\bf in press}.

\harvarditem{Loh \harvardand\ Wainwright}{2012}{loh2012}
Loh, P.-L. \harvardand\ Wainwright, M.~J.  \harvardyearleft
  2012\harvardyearright , `High-dimensional regression with noisy and missing
  data: Provable guarantees with nonconvexity', {\em The Annals of Statistics}
  pp.~1637--1664.

\harvarditem{Loh \harvardand\ Wainwright}{2015}{loh2015}
Loh, P.-L. \harvardand\ Wainwright, M.~J.  \harvardyearleft
  2015\harvardyearright , `Regularized m-estimators with nonconvexity:
  Statistical and algorithmic theory for local optima', {\em Journal of Machine
  Learning Research} {\bf 16},~559--616.

\harvarditem{Loh \harvardand\ Wainwright}{2017}{loh2017}
Loh, P.-L. \harvardand\ Wainwright, M.~J.  \harvardyearleft
  2017\harvardyearright , `Support recovery without incoherence: A case for
  nonconvex regularization', {\em The Annals of Statistics} {\bf
  45}(6),~2455--2482.

\harvarditem{McCullagh \harvardand\ Nelder}{2019}{mccullagh2019}
McCullagh, P. \harvardand\ Nelder, J.~A.  \harvardyearleft
  2019\harvardyearright , {\em Generalized linear models}, Routledge.

\harvarditem[Mitnitski et~al.]{Mitnitski, Fallah, Dean \harvardand\
  Rockwood}{2014}{mitnitski2014multi}
Mitnitski, A.~B., Fallah, N., Dean, C.~B. \harvardand\ Rockwood, K.
  \harvardyearleft 2014\harvardyearright , `A multi-state model for the
  analysis of changes in cognitive scores over a fixed time interval', {\em
  Statistical methods in medical research} {\bf 23}(3),~244--256.

\harvarditem{Ning \harvardand\ Liu}{2017}{ning2017}
Ning, Y. \harvardand\ Liu, H.  \harvardyearleft 2017\harvardyearright , `A
  general theory of hypothesis tests and confidence regions for sparse high
  dimensional models', {\em The Annals of Statistics} {\bf 45},~158--195.

\harvarditem[Sch{\"o}ll et~al.]{Sch{\"o}ll, Lockhart, Schonhaut, O’Neil,
  Janabi, Ossenkoppele, Baker, Vogel, Faria, Schwimmer
  et~al.}{2016}{scholl2016}
Sch{\"o}ll, M., Lockhart, S.~N., Schonhaut, D.~R., O’Neil, J.~P., Janabi, M.,
  Ossenkoppele, R., Baker, S.~L., Vogel, J.~W., Faria, J., Schwimmer, H.~D.
  et~al.  \harvardyearleft 2016\harvardyearright , `Pet imaging of tau
  deposition in the aging human brain', {\em Neuron} {\bf 89}(5),~971--982.

\harvarditem{{{S}}ent{\"u}rk \harvardand\ M{\"u}ller}{2005}{csenturk2005}
{{S}}ent{\"u}rk, D. \harvardand\ M{\"u}ller, H.-G.  \harvardyearleft
  2005\harvardyearright , `Covariate-adjusted regression', {\em Biometrika}
  {\bf 92}(1),~75--89.

\harvarditem[Shi et~al.]{Shi, Song, Chen, Li et~al.}{2019}{shi2019linear}
Shi, C., Song, R., Chen, Z., Li, R. et~al.  \harvardyearleft
  2019\harvardyearright , `Linear hypothesis testing for high dimensional
  generalized linear models', {\em The Annals of statistics} {\bf
  47}(5),~2671--2703.

\harvarditem{Smith \harvardand\ Kosslyn}{2008}{smith2008}
Smith, E. \harvardand\ Kosslyn, S.  \harvardyearleft 2008\harvardyearright ,
  {\em Cognitive Psychology: Mind and Brain}, Pearson Education, Pearson
  Prentice Hall.
\newline\harvardurl{https://books.google.com/books?id=YIPSNwAACAAJ}

\harvarditem[S{\o}rensen et~al.]{S{\o}rensen, Frigessi \harvardand\
  Thoresen}{2015}{sorensen2015}
S{\o}rensen, {\O}., Frigessi, A. \harvardand\ Thoresen, M.  \harvardyearleft
  2015\harvardyearright , `Measurement error in lasso: Impact and likelihood
  bias correction', {\em Statistica sinica} pp.~809--829.

\harvarditem[S{\o}rensen et~al.]{S{\o}rensen, Hellton, Frigessi \harvardand\
  Thoresen}{2018}{sorensen2018covariate}
S{\o}rensen, {\O}., Hellton, K.~H., Frigessi, A. \harvardand\ Thoresen, M.
  \harvardyearleft 2018\harvardyearright , `Covariate selection in
  high-dimensional generalized linear models with measurement error', {\em
  Journal of Computational and Graphical Statistics} {\bf 27}(4),~739--749.

\harvarditem{Stefanski}{1989}{stefanski1989}
Stefanski, L.~A.  \harvardyearleft 1989\harvardyearright , `Unbiased estimation
  of a nonlinear function a normal mean with application to measurement err
  oorf models', {\em Communications in Statistics-Theory and Methods} {\bf
  18}(12),~4335--4358.

\harvarditem[Van~de Geer et~al.]{Van~de Geer, B{\"u}hlmann, Ritov \harvardand\
  Dezeure}{2014}{van2014asymptotically}
Van~de Geer, S., B{\"u}hlmann, P., Ritov, Y. \harvardand\ Dezeure, R.
  \harvardyearleft 2014\harvardyearright , `On asymptotically optimal
  confidence regions and tests for high-dimensional models', {\em Annals of
  Statistics} {\bf 42}(3),~1166--1202.

\harvarditem{Vial}{1982}{vial1982}
Vial, J.-P.  \harvardyearleft 1982\harvardyearright , `Strong convexity of sets
  and functions', {\em Journal of Mathematical Economics} {\bf
  9}(1-2),~187--205.

\harvarditem{Wainwright}{2009}{wainwright2009}
Wainwright, M.~J.  \harvardyearleft 2009\harvardyearright , `Sharp thresholds
  for high-dimensional and noisy sparsity recovery using $l_1$ - constrained
  quadratic programming (lasso)', {\em IEEE transactions on information theory}
  {\bf 55}(5),~2183--2202.

\harvarditem[Weiner et~al.]{Weiner, Veitch, Aisen, Beckett, Cairns, Green,
  Harvey, Jack~Jr, Jagust, Morris et~al.}{2017}{weiner2017}
Weiner, M.~W., Veitch, D.~P., Aisen, P.~S., Beckett, L.~A., Cairns, N.~J.,
  Green, R.~C., Harvey, D., Jack~Jr, C.~R., Jagust, W., Morris, J.~C. et~al.
  \harvardyearleft 2017\harvardyearright , `The alzheimer's disease
  neuroimaging initiative 3: Continued innovation for clinical trial
  improvement', {\em Alzheimer's \& Dementia} {\bf 13}(5),~561--571.

\harvarditem[Zhang et~al.]{Zhang, Zhang et~al.}{2012}{zhang2012}
Zhang, C.-H., Zhang, T. et~al.  \harvardyearleft 2012\harvardyearright , `A
  general theory of concave regularization for high-dimensional sparse
  estimation problems', {\em Statistical Science} {\bf 27}(4),~576--593.

\harvarditem{Zhang \harvardand\ Cheng}{2017}{zhang2017simultaneous}
Zhang, X. \harvardand\ Cheng, G.  \harvardyearleft 2017\harvardyearright ,
  `Simultaneous inference for high-dimensional linear models', {\em Journal of
  the American Statistical Association} {\bf 112}(518),~757--768.

\end{thebibliography}


\begin{thebibliography}{11}
\providecommand{\natexlab}[1]{#1}
\providecommand{\url}[1]{\texttt{#1}}
\expandafter\ifx\csname urlstyle\endcsname\relax
  \providecommand{\doi}[1]{doi: #1}\else
  \providecommand{\doi}{doi: \begingroup \urlstyle{rm}\Url}\fi

\bibitem[Bentkus(2005)]{bentkus2005}
Vidmantas Bentkus.
\newblock A lyapunov-type bound in $r^d$.
\newblock \emph{Theory of Probability \& Its Applications}, 49\penalty0
  (2):\penalty0 311--323, 2005.

\bibitem[Braak and Braak(1991)]{braak1991}
Heiko Braak and Eva Braak.
\newblock Neuropathological stageing of alzheimer-related changes.
\newblock \emph{Acta neuropathologica}, 82\penalty0 (4):\penalty0 239--259,
  1991.

\bibitem[Fletcher and Watson(1980)]{fletcher1980}
Roger Fletcher and G~Alistair Watson.
\newblock First and second order conditions for a class of nondifferentiable
  optimization problems.
\newblock \emph{Mathematical Programming}, 18\penalty0 (1):\penalty0 291--307,
  1980.

\bibitem[Jiang and Ma(2021)]{jiang2021}
Fei Jiang and Yanyuan Ma.
\newblock Poisson regression with error corrupted high dimensional features.
\newblock \emph{Statistica Sinica}, 2021.
\newblock In press.

\bibitem[Landau et~al.(2016)Landau, Ward, Murphy, and Jagust]{Landau2016}
Susan Landau, Tyler~J. Ward, Alice Murphy, and William Jagust.
\newblock Flortaucipir (av-1451) processing methods, 2016.
\newblock URL
  \url{\url{https://adni.bitbucket.io/reference/docs/UCBERKELEYAV1451/UCBERKELEY_AV1451_Methods_Aug2018.pdf}}.

\bibitem[Ledoux and Talagrand(2013)]{ledoux2013}
Michel Ledoux and Michel Talagrand.
\newblock \emph{Probability in Banach Spaces: isoperimetry and processes}.
\newblock Springer Science \& Business Media, 2013.

\bibitem[Loh and Wainwright(2012)]{loh2012}
Po-Ling Loh and Martin~J Wainwright.
\newblock High-dimensional regression with noisy and missing data: Provable
  guarantees with nonconvexity.
\newblock \emph{The Annals of Statistics}, pages 1637--1664, 2012.

\bibitem[Loh and Wainwright(2015)]{loh2015}
Po-Ling Loh and Martin~J Wainwright.
\newblock Regularized m-estimators with nonconvexity: Statistical and
  algorithmic theory for local optima.
\newblock \emph{Journal of Machine Learning Research}, 16:\penalty0 559--616,
  2015.

\bibitem[Loh and Wainwright(2017)]{loh2017}
Po-Ling Loh and Martin~J. Wainwright.
\newblock Support recovery without incoherence: A case for nonconvex
  regularization.
\newblock \emph{The Annals of Statistics}, 45\penalty0 (6):\penalty0
  2455--2482, 2017.

\bibitem[Sch{\"o}ll et~al.(2016)Sch{\"o}ll, Lockhart, Schonhaut, O’Neil,
  Janabi, Ossenkoppele, Baker, Vogel, Faria, Schwimmer, et~al.]{scholl2016}
Michael Sch{\"o}ll, Samuel~N Lockhart, Daniel~R Schonhaut, James~P O’Neil,
  Mustafa Janabi, Rik Ossenkoppele, Suzanne~L Baker, Jacob~W Vogel, Jamie
  Faria, Henry~D Schwimmer, et~al.
\newblock Pet imaging of tau deposition in the aging human brain.
\newblock \emph{Neuron}, 89\penalty0 (5):\penalty0 971--982, 2016.

\bibitem[Wainwright(2009)]{wainwright2009}
Martin~J Wainwright.
\newblock Sharp thresholds for high-dimensional and noisy sparsity recovery
  using $l_1$ - constrained quadratic programming (lasso).
\newblock \emph{IEEE transactions on information theory}, 55\penalty0
  (5):\penalty0 2183--2202, 2009.

\end{thebibliography}
